\documentclass[12pt,reqno,psamsfonts]{amsart}
\usepackage[alphabetic]{amsrefs}
\usepackage{feynmp}
\usepackage{latexsym}
\usepackage{amsfonts}
\usepackage{amssymb}
\usepackage{amscd}
\usepackage{amsbsy}
\usepackage{amsmath}
\usepackage{bbm}
\usepackage{euscript}
\usepackage{mathrsfs}
\usepackage{yhmath}

\usepackage{dashbox}
\usepackage{environ}

\usepackage[dvips]{graphicx}
\usepackage{epsfig}    
\usepackage{xcolor}
\usepackage{tikz}
\usetikzlibrary{shapes}
\usetikzlibrary{decorations.markings}
\usetikzlibrary{decorations.pathreplacing}
\usetikzlibrary{patterns}
\usetikzlibrary{matrix,arrows}
\usepgflibrary{decorations.pathmorphing}

\usepackage{tikz-cd}

\usepackage[all]{xy}
\xyoption{knot}
\xyoption{graph}
\xyoption{tile}
\xyoption{arc}
\xyoption{poly}
\xyoption{pdf}

\usepackage{setspace}

\def\ie{{\it i.e.\ }}

\def\cf{{\it cf.\ }}
\def\rhs{{\it r.h.s.\ }}
\def\lhs{{\it l.h.s.\ }}

\def\Hom{\mathop{{\rm Hom}}\nolimits}

\def\deg{ \mathop{{\rm deg}}\nolimits }
\def\p{^{\prime}}

\def\promod{\mathop{{\rm mod}}\nolimits}
\def\mod#1{\;(\promod #1)}

\def\pr#1#2{ \noindent{\em Proof of #1~\ref{#2}.} }

\def\qed{ \hfill $\Box$ }

\def\lrbc#1{ \left( #1 \right) }

\parskip=\medskipamount                 
\thicklines                     
\def\inbar{\vrule height1.5ex width.4pt depth0pt}
\def\IC{\relax\,\hbox{$\inbar\kern-.3em{\rm C}$}}
\def\IN{\relax{\rm I\kern-.18em N}}
\def\IQ{\relax\,\hbox{$\inbar\kern-.3em{\rm Q}$}}
\def\IR{\relax{\rm I\kern-.18em R}}
\def\ZZ{\relax{\sf Z\kern-.4em Z}}

\newtheorem{theorem}{Theorem}[section]
\newtheorem{proposition}[theorem]{Proposition}
\newtheorem{corollary}[theorem]{Corollary}
\newtheorem{conjecture}[theorem]{Conjecture}
\newtheorem{lemma}[theorem]{Lemma}
\newtheorem{definition}[theorem]{Definition}
\newtheorem{remark}[theorem]{Remark}

\catcode`\@=11

\newif\if@fewtab\@fewtabtrue

\catcode`\@=11

\newif\if@fewtab\@fewtabtrue

{\count255=\time\divide\count255 by 60
\xdef\hourmin{\number\count255} \multiply\count255
by-60\advance\count255 by\time
\xdef\hourmin{\hourmin:\ifnum\count255<10 0\fi\the\count255}}
\def\ps@draft{\let\@mkboth\@gobbletwo
    \def\@oddhead{}
    \def\@oddfoot
      {\hbox to 7 cm{\footnotesize {\em Draft of \jobname:} \draftdate
       \hfil}\hskip -7cm\hfil\rm\thepage \hfil}
    \def\@evenhead{}\let\@evenfoot\@oddfoot}

\def\ceqno{\global\@fewtabfalse
    \ifcase\@eqcnt \def\@tempa{& & &}\or \def\@tempa{& &}
      \or \def\@tempa{&}
      \or\def\@tempa{}\fi\@tempa
{\rm(\theequation)}}

\def\aeqno#1{\global\@fewtabfalse
    \ifcase\@eqcnt \def\@tempa{& & &}\or \def\@tempa{& &}
      \or \def\@tempa{&}
      \or\def\@tempa{}\fi\@tempa
{\rm(\theequation,#1)}}

\def\label#1{\ifnum\draftcontrol=1
 \global\def\draftnote{$\scriptstyle #1$}\fi
 \@bsphack\if@filesw {\let\thepage\relax
   \def\protect{\noexpand\noexpand\noexpand}%
\xdef\@gtempa{\write\@auxout{\string
      \newlabel{#1}{{\@currentlabel}{\thepage}}}}}\@gtempa
   \if@nobreak \ifvmode\nobreak\fi\fi\fi
  \@esphack}

\def\alabel#1#2{\label{#1}\global\@fewtabfalse
    \ifcase\@eqcnt \def\@tempa{& & &}\or \def\@tempa{& &}
      \or \def\@tempa{&}
      \or\def\@tempa{}\fi\@tempa
{\hbox to 3cm{\phantom{\rm(\theequation,#2)} \draftnote
\hfil}\hskip -3cm {\rm(\theequation,#2)}}}

\def\clabel#1{\label{#1}\global\@fewtabfalse
    \ifcase\@eqcnt \def\@tempa{& & &}\or \def\@tempa{& &}
      \or \def\@tempa{&}
      \or\def\@tempa{}\fi\@tempa
{\hbox to 3cm{\phantom{\rm(\theequation)} \draftnote \hfil}\hskip
-3cm{\rm(\theequation)}}}

\def\eqnarray{\def\draftnote{{}}\global\@fewtabtrue
\stepcounter{equation}\let\@currentlabel=\theequation
\global\@eqnswtrue
\global\@eqcnt\z@\tabskip\@centering\let\\=\@eqncr
$$\halign to \displaywidth\bgroup\@eqnsel\hskip\@centering\@eqcnt\z@
  $\displaystyle\tabskip\z@{##}$&\global\@eqcnt\@ne
  \hskip 1\arraycolsep \hfil$\displaystyle{##}$\hfil
  &\global\@eqcnt\tw@ \hskip 1\arraycolsep
$\displaystyle\tabskip\z@{##}$ \hfil
\tabskip\@centering&\global\@eqcnt\thr@@\llap{##}\tabskip\z@ \cr}

\def\endeqnarray{\@@eqncr\egroup
      \global\advance\c@equation\m@ne$$\global\@ignoretrue}

\def\@eqnnum{\hbox to 3cm{\phantom{\rm(\theequation)} \draftnote
                         \hfil}\hskip -3cm {\rm(\theequation)}}

\def\@@eqncr{\let\@tempa\relax
    \ifcase\@eqcnt \def\@tempa{& & &}\or \def\@tempa{& &}
      \or \def\@tempa{&}
      \or\def\@tempa{}
\fi\@tempa \if@eqnsw \if@fewtab\@eqnnum\fi
\stepcounter{equation}\fi\global
\@eqnswtrue\global\@eqcnt\z@\global\@fewtabtrue\cr}

\def\draftcite#1{\ifnum\draftcontrol=1#1\else{}\fi}

\def\@lbibitem[#1]#2{\item{}\hskip -3cm \hbox to 2cm
{\hfil$\scriptstyle\draftcite{#2}$}\hskip
1cm[\@biblabel{#1}]\if@filesw
     {\def\protect##1{\string ##1\space}\immediate
      \write\@auxout{\string\bibcite{#2}{#1}}}\fi\ignorespaces}

\def\@bibitem#1{\item\hskip -3cm \hbox to 2cm
{\hfil $\scriptstyle\draftcite{#1}$}\hskip 1cm \if@filesw
\immediate\write\@auxout
       {\string\bibcite{#1}{\the\value{\@listctr}}}\fi\ignorespaces}

\def\draftdate{\number\month/\number\day/\number\year\ \ \ \hourmin }
\def\draft{\pagestyle{draft}\thispagestyle{draft}
\global\def\draftcontrol{1}} \global\def\draftcontrol{0}
\catcode`\@=12

\def\theequation{{\thesection.\arabic{equation}}}

\def\qq{\begin{eqnarray}}
\def\qqq{\end{eqnarray}}
\def\ee{\begin{eqnarray}}
\def\eee{\end{eqnarray}}
\def\rx#1{~(\ref{#1})}
\def\rxw#1{(\ref{#1})}
\def\ex#1{eq.\hspace*{-3pt}\rx{#1}}
\def\eex#1{eqs.\hspace*{-3pt}\rx{#1}}
\def\cx#1{~\cite{#1}}
\def\rw#1{~\ref{#1}}

\def\xlee#1{ \begin{eqnarray} \label{#1} }
\def\xeee{ \end{eqnarray} }
\def\ylee#1{ \begin{eqnarray}\nonumber }
\def\yeee{ \end{eqnarray} }
\def\zlee#1{ \begin{displaymath} }
\def\zleee{ \end{displaymath} }
\def\wlee#1{ $ }
\def\weee{ $ }

\def\fg#1{Fig.~\ref{#1}}

\newlength{\shiftwidth}
\def\shift#1{&&\hbox to \shiftwidth{\hfill $\displaystyle#1$}}
\newlength{\sshiftwidth}
\def\sshift#1{\lefteqn{\hbox to
\sshiftwidth{\hfill$\displaystyle#1$}}}

\def\qbezier{\bezier{120}}
\def\DottedCircle{
\bezier{4}(0.966,-0.259)(1.04,0)(0.966,0.259)
\bezier{4}(0.966,0.259)(0.897,0.518)(0.707,0.707)
\bezier{4}(0.707,0.707)(0.518,0.897)(0.259,0.966)
\bezier{4}(0.259,0.966)(0,1.04)(-0.259,0.966)
\bezier{4}(-0.259,0.966)(-0.518,0.897)(-0.707,0.707)
\bezier{4}(-0.707,0.707)(-0.897,0.518)(-0.966,0.259)
\bezier{4}(-0.966,0.259)(-1.04,0)(-0.966,-0.259)
\bezier{4}(-0.966,-0.259)(-0.897,-0.518)(-0.707,-0.707)
\bezier{4}(-0.707,-0.707)(-0.518,-0.897)(-0.259,-0.966)
\bezier{4}(-0.259,-0.966)(0,-1.04)(0.259,-0.966)
\bezier{4}(0.259,-0.966)(0.518,-0.897)(0.707,-0.707)
\bezier{4}(0.707,-0.707)(0.897,-0.518)(0.966,-0.259) }
\def\Endpoint[#1]{
\ifcase#1 \put(1,0){\circle*{0.15}}
\or\put(0.866,0.5){\circle*{0.15}}
\or\put(0.5,0.866){\circle*{0.15}} \or\put(0,1){\circle*{0.15}}
\or\put(-0.5,0.866){\circle*{0.15}}
\or\put(-0.866,0.5){\circle*{0.15}} \or\put(-1,0){\circle*{0.15}}
\or\put(-0.866,-0.5){\circle*{0.15}}
\or\put(-0.5,-0.866){\circle*{0.15}} \or\put(0,-1){\circle*{0.15}}
\or\put(0.5,-0.866){\circle*{0.15}}
\or\put(0.866,-0.5){\circle*{0.15}} \fi}
\def\Arc[#1]{
\thicklines         
\ifcase#1 \bezier{25}(0.966,-0.259)(1.04,0)(0.966,0.259) \or
\bezier{25}(0.966,0.259)(0.897,0.518)(0.707,0.707) \or
\bezier{25}(0.707,0.707)(0.518,0.897)(0.259,0.966) \or
\bezier{25}(0.259,0.966)(0,1.04)(-0.259,0.966) \or
\bezier{25}(-0.259,0.966)(-0.518,0.897)(-0.707,0.707) \or
\bezier{25}(-0.707,0.707)(-0.897,0.518)(-0.966,0.259) \or
\bezier{25}(-0.966,0.259)(-1.04,0)(-0.966,-0.259) \or
\bezier{25}(-0.966,-0.259)(-0.897,-0.518)(-0.707,-0.707) \or
\bezier{25}(-0.707,-0.707)(-0.518,-0.897)(-0.259,-0.966) \or
\bezier{25}(-0.259,-0.966)(0,-1.04)(0.259,-0.966) \or
\bezier{25}(0.259,-0.966)(0.518,-0.897)(0.707,-0.707) \or
\bezier{25}(0.707,-0.707)(0.897,-0.518)(0.966,-0.259) \fi}
\def\DottedArc[#1]{
\ifcase#1 \bezier{4}(0.966,-0.259)(1.04,0)(0.966,0.259) \or
\bezier{4}(0.966,0.259)(0.897,0.518)(0.707,0.707) \or
\bezier{4}(0.707,0.707)(0.518,0.897)(0.259,0.966) \or
\bezier{4}(0.259,0.966)(0,1.04)(-0.259,0.966) \or
\bezier{4}(-0.259,0.966)(-0.518,0.897)(-0.707,0.707) \or
\bezier{4}(-0.707,0.707)(-0.897,0.518)(-0.966,0.259) \or
\bezier{4}(-0.966,0.259)(-1.04,0)(-0.966,-0.259) \or
\bezier{4}(-0.966,-0.259)(-0.897,-0.518)(-0.707,-0.707) \or
\bezier{4}(-0.707,-0.707)(-0.518,-0.897)(-0.259,-0.966) \or
\bezier{4}(-0.259,-0.966)(0,-1.04)(0.259,-0.966) \or
\bezier{4}(0.259,-0.966)(0.518,-0.897)(0.707,-0.707) \or
\bezier{4}(0.707,-0.707)(0.897,-0.518)(0.966,-0.259) \fi}
\def\Chord[#1,#2]{
\thinlines \ifnum#1>#2\Chord[#2,#1] \else\ifnum#1<#2 \ifcase#1
\ifcase#2 \or\qbezier(1,0)(0.516,0.138)(0.866,0.5)
\or\qbezier(1,0)(0.45,0.26)(0.5,0.866)
\or\qbezier(1,0)(0.327,0.327)(0,1)
\or\qbezier(1,0)(0.179,0.311)(-0.5,0.866)
\or\qbezier(1,0)(0.0536,0.2)(-0.866,0.5) \or\put(1, 0){\line(-2,
0){2}} \or\qbezier(1,0)(0.0536,-0.2)(-0.866,-0.5)
\or\qbezier(1,0)(0.179,-0.311)(-0.5,-0.866)
\or\qbezier(1,0)(0.327,-0.327)(0,-1)
\or\qbezier(1,0)(0.45,-0.26)(0.5,-0.866)
\or\qbezier(1,0)(0.516,-0.138)(0.866,-0.5) \fi \or\ifcase#2\or
\or\qbezier(0.866,0.5)(0.378,0.378)(0.5,0.866)
\or\qbezier(0.866,0.5)(0.26,0.45)(0,1)
\or\qbezier(0.866,0.5)(0.12,0.446)(-0.5,0.866)
\or\qbezier(0.866,0.5)(0,0.359)(-0.866,0.5)
\or\qbezier(0.866,0.5)(-0.0536,0.2)(-1,0) \or\put(0.866,
0.5){\line(-5, -3){1.73}}
\or\qbezier(0.866,0.5)(0.146,-0.146)(-0.5,-0.866)
\or\qbezier(0.866,0.5)(0.311,-0.179)(0,-1)
\or\qbezier(0.866,0.5)(0.446,-0.12)(0.5,-0.866)
\or\qbezier(0.866,0.5)(0.52,0)(0.866,-0.5) \fi \or\ifcase#2\or\or
\or\qbezier(0.5,0.866)(0.138,0.516)(0,1)
\or\qbezier(0.5,0.866)(0,0.52)(-0.5,0.866)
\or\qbezier(0.5,0.866)(-0.12,0.446)(-0.866,0.5)
\or\qbezier(0.5,0.866)(-0.179,0.311)(-1,0)
\or\qbezier(0.5,0.866)(-0.146,0.146)(-0.866,-0.5) \or\put(0.5,
0.866){\line(-3, -5){1}} \or\qbezier(0.5,0.866)(0.2,-0.0536)(0,-1)
\or\qbezier(0.5,0.866)(0.359,0)(0.5,-0.866)
\or\qbezier(0.5,0.866)(0.446,0.12)(0.866,-0.5) \fi
\or\ifcase#2\or\or\or \or\qbezier(0,1.)(-0.138,0.516)(-0.5,0.866)
\or\qbezier(0,1.)(-0.26,0.45)(-0.866,0.5)
\or\qbezier(0,1.)(-0.327,0.327)(-1,0)
\or\qbezier(0,1.)(-0.311,0.179)(-0.866,-0.5)
\or\qbezier(0,1.)(-0.2,0.0536)(-0.5,-0.866) \or\put(0, 1){\line(0,
-2){2}} \or\qbezier(0,1.)(0.2,0.0536)(0.5,-0.866)
\or\qbezier(0,1.)(0.311,0.179)(0.866,-0.5) \fi
\or\ifcase#2\or\or\or\or
\or\qbezier(-0.5,0.866)(-0.378,0.378)(-0.866,0.5)
\or\qbezier(-0.5,0.866)(-0.45,0.26)(-1,0)
\or\qbezier(-0.5,0.866)(-0.446,0.12)(-0.866,-0.5)
\or\qbezier(-0.5,0.866)(-0.359,0)(-0.5,-0.866)
\or\qbezier(-0.5,0.866)(-0.2,-0.0536)(0,-1) \or\put(-0.5,
0.866){\line(3, -5){1}}
\or\qbezier(-0.5,0.866)(0.146,0.146)(0.866,-0.5) \fi
\or\ifcase#2\or\or\or\or\or
\or\qbezier(-0.866,0.5)(-0.516,0.138)(-1,0)
\or\qbezier(-0.866,0.5)(-0.52,0)(-0.866,-0.5)
\or\qbezier(-0.866,0.5)(-0.446,-0.12)(-0.5,-0.866)
\or\qbezier(-0.866,0.5)(-0.311,-0.179)(0,-1)
\or\qbezier(-0.866,0.5)(-0.146,-0.146)(0.5,-0.866) \or\put(-0.866,
0.5){\line(5, -3){1.73}} \fi \or\ifcase#2\or\or\or\or\or\or
\or\qbezier(-1,0)(-0.516,-0.138)(-0.866,-0.5)
\or\qbezier(-1,0)(-0.45,-0.26)(-0.5,-0.866)
\or\qbezier(-1,0)(-0.327,-0.327)(0,-1)
\or\qbezier(-1,0)(-0.179,-0.311)(0.5,-0.866)
\or\qbezier(-1,0)(-0.0536,-0.2)(0.866,-0.5) \fi
\or\ifcase#2\or\or\or\or\or\or\or
\or\qbezier(-0.866,-0.5)(-0.378,-0.378)(-0.5,-0.866)
\or\qbezier(-0.866,-0.5)(-0.26,-0.45)(0,-1)
\or\qbezier(-0.866,-0.5)(-0.12,-0.446)(0.5,-0.866)
\or\qbezier(-0.866,-0.5)(0,-0.359)(0.866,-0.5) \fi
\or\ifcase#2\or\or\or\or\or\or\or\or
\or\qbezier(-0.5,-0.866)(-0.138,-0.516)(0,-1)
\or\qbezier(-0.5,-0.866)(0,-0.52)(0.5,-0.866)
\or\qbezier(-0.5,-0.866)(0.12,-0.446)(0.866,-0.5) \fi
\or\ifcase#2\or\or\or\or\or\or\or\or\or
\or\qbezier(0,-1.)(0.138,-0.516)(0.5,-0.866)
\or\qbezier(0,-1.)(0.26,-0.45)(0.866,-0.5) \fi
\or\ifcase#2\or\or\or\or\or\or\or\or\or\or
\or\qbezier(0.5,-0.866)(0.378,-0.378)(0.866,-0.5) \fi\fi\fi\fi}
\def\FullChord[#1,#2]{
\Endpoint[#1] \Endpoint[#2] \Arc[#1] \Arc[#2] \Chord[#1,#2] }
\def\EndChord[#1,#2]{
\Endpoint[#1] \Endpoint[#2] \Chord[#1,#2] }
\def\Picture#1{
\begin{picture}(2,1)(-1,-0.167)
#1
\end{picture}
}
\def\DottedChordDiagram[#1,#2]{
\Picture{\DottedCircle \FullChord[#1,#2]} }
\def\ZZ{ \mathbb{Z} }
\def\IQ{ \mathbb{Q} }
\def\IC{ \mathbb{C} }
\def\IR{ \mathbb{R} }

\def\hlf{ {1\over 2} }

\def\xS{ \mathbb{S} }
\def\xSv#1{ \xS^{#1} }
\def\xSo {\xSv{1} }
\def\xSt{\xSv{2} }
\def\xStSo{\xSt\times\xSo}

\def\xD{ \mathbb{D} }
\def\xDp{ \xD' }
\def\xB{ \mathbb{B} }
\def\xT{ \mathbb{T} }
\def\xTv#1{ \xT^{#1} }
\def\xTt{ \xTv{2} }

\def\yT{ \mathrm{T} }
\def\yTvv#1#2{ \yT_{#1,#2} }
\def\yTmn{ \yTvv{m}{n} }
\def\yTmmn{ \yTvv{m}{-n} }

\def\Hom{ \mathop{\mathrm{Hom}}\nolimits }

\def\xId{ \mathbbm{1} }

\def\hlf{ \frac{1}{2} }
\def\thlf{ \frac{3}{2} }
\def\frth{ \frac{1}{4} }
\def\shlf{ \tfrac{1}{2} }
\def\sfrth{ \tfrac{1}{4} }
\def\sthlf{ \tfrac{3}{2} }

\def\hem{\bullet}

\def\TQFT{TQFT}

\def\tKbr{Kauffman bracket}

\def\tJpol{Jones polynomial}

\def\tKhom{Khovanov homology}
\def\tKhoms{Khovanov homologies}
\def\tKhbr{Khovanov bracket}

\def\tJW{Jones-Wenzl}
\def\tJWp{\tJW\ projector}
\def\cJWp{categorified \tJWp}

\def\tTL{Temperley-Lieb}
\def\taTL{TL}
\def\tTLc{\tTL\ category}

\def\tTLt{\tTL\ tangle}
\def\taTLt{\taTL\ tangle}
\def\tTLa{\tTL\ algebra}

\def\trbr{torus braid}

\def\tmcn{multi-cone}

\def\thead{tail}

\def\Asplng{A-splicing}
\def\Bsplng{B-splicing}

\def\Arpl{A-replacement}
\def\Arpld{A-replaced}
\def\Brpl{B-replacement}
\def\Brpld{B-replaced}

\def\tBcr{B-circle}

\def\tadq{adequate}
\def\tBadq{B-\tadq}
\def\tBiadq{B-inadequate}
\def\tBdg{B-diagram}
\def\trBdg{reduced \tBdg}

\def\tinadq{inadequate}
\def\tnBadq{not \tBadq}

\def\tBrdc{B-reduction}

\def\uclrd{unicolored}

\def\splng{splicing}

\def\nsplcd{negatively spliced}

\def\twdfc{\twd\ deficit}
\def\twd{width}

\def\tbdgr{bi-degree}

\def\ZZt{ \ZZ_2 }

\def\ZZZtt{\ZZ\times\ZZ\times\ZZt}

\def\ZZtqqi{ \ZZ[q^{\pm 2} ] }

\def\QQ{ \mathbb{Q} }

\def\QQqqi{ \QQ[q^{\pm 1} ] }
\def\QQqqip{ \QQ[[q,q^{-1}] }

\def\ctfont{ \mathsf }
\def\chfont{ \mathbf }
\def\stfont{ \mathfrak }

\def\tcat#1{ \ddot{#1} }

\def\cTng{ \ctfont{Tng} }
\def\ctTng{ \tcat{\cTng} }

\def\cTL{ \ctfont{TL} }
\def\cTLp{ \cTL^+ }

\def\ctTL{ \tcat{\cTL} }

\def\caC{ \ctfont{C} }
\def\caCt{ \tilde{\caC} }
\def\ctaC{ \tcat{\caC} }

\def\caA{ \ctfont{A} }

\def\cKom{ \mathop{\mathbf{Kom}} }
\def\cKomm{ \cKom^+ }

\def\cKommA{ \cKomm(\caA) }

\def\chA{ \chfont{A} }
\def\chB{ \chfont{B} }

\def\stA{ \stfont{A} }

\def\oba{ \alpha }

\def\obO{ O }

\def\hmord#1{ |#1|_{\mathrm h} }
\def\hlmord#1{ \left| #1 \right|_{\mathrm h} }

\def\drsys{direct system}
\def\drlim{direct limit}

\def\SUv#1{ \SU({#1}) }
\def\SUt{ \SUv{2} }

\def\aTL{ \mathrm{TL} }
\def\aTLv#1{ \aTL_{#1} }
\def\aTLmn{ \aTLv{m,n} }

\def\pJ{ \mathrm{J} }
\def\pJv#1{ \pJ_{#1} }
\def\pJqv#1{ \pJv{#1}(q) }
\def\pJvv#1#2{ \pJ_{#1,#2} }
\def\pJqvv#1#2{ \pJvv{#1}{#2}(q) }
\def\pJqL{ \pJqv{\xL} }
\def\pJqLi{ \pJqvv{\xL}{\infty} }

\def\pJqaL{ \pJqvv{\xca}{\xL} }
\def\pJqaD{ \pJqvv{\xca}{\xD} }

\def\zpol{ \mathrm{T} }
\def\zpolnvv#1#2{ \zpol^{#1}_{#2} }
\def\zpolnqvv#1#2{ \zpolnvv{#1}{#2}(q) }
\def\zpolnqnL{ \zpolnqvv{\xnu}{\xL} }

\def\ypolvv#1#2{ \zpol_{#1,#2} }
\def\ypolqvv#1#2{ \ypolvv{#1}{#2}(q) }
\def\ypolqnL{ \ypolqvv{\xn}{\xL} }
\def\ypolqnD{ \ypolqvv{\xn}{\xD} }

\def\ypolv#1{ \zpol_{#1} }
\def\ypolqtv#1{ \ypolv{#1}(\xt,q) }
\def\ypolqtL{ \ypolqtv{\xL} }
\def\ypolqtD{ \ypolqtv{\xD} }

\def\xdmm{ - }

\def\xKbrv#1{ \langle #1 \rangle }
\def\xKbrBv#1{ \Big \langle\, #1 \,\Big \rangle }

\def\xKbrd{ \xKbrv{\xdmm} }

\def\xKhv#1{ \langle\!\langle #1 \rangle\!\rangle }

\def\xvKhv#1{ \left \langle\!\!\!\left\langle #1 \right\rangle\!\!\!\right\rangle }

\def\Cnv#1{ \boxed{#1} }
\def\sdff{ \circlearrowleft }
\def\Pcnv#1{ \boxed{#1}_{\,\displaystyle \sdff} }

\def\Conv#1{ \mathrm{Cone} (#1) }

\def\hteqv{ \sim }

\def\dlm{ \lim\limits_{\rightarrow}}

\def\xbul{ \bullet }

\def\Kh{\scriptscriptstyle{\mathrm{Kh}} }
\def\Hm{ \mathrm{H} }
\def\KHm{ \Hm^{\Kh} }

\def\KHmvv#1#2{ \KHm_{#1,#2} }
\def\KHmvb#1{ \KHmvv{#1}{\hem} }

\def\KHmib{\KHmvv{i}{\xbul} }

\def\tKHm{ \tilde{\Hm}^{\Kh} }
\def\tKHmvv#1#2{ \tKHm_{#1,#2} }
\def\tKHmvb#1{ \tKHmvv{#1}{\hem} }

\def\tlH{ \Hm^{\infty} }
\def\tlHvv#1#2{ \tlH_{#1,#2} }

\def\ttlH{ \Hm^{\boldsymbol{\infty}} }
\def\ttlH{ \Hm^{\thicksim} }
\def\ttlHvvv#1#2#3{ \ttlH_{#1,#2,#3} }

\def\hgrshv#1#2#3{ [#1,#2,#3] }

\def\dgh{ \deg_{\mathrm{h}} }
\def\dgq{ \deg_{\mathrm{q}} }
\def\dgb{ \deg_{\mathrm{b}} }

\def\tqdgr{$q$-degree}
\def\thdgr{$h$-degree}
\def\tbdgr{$b$-degree}

\def\tbgrd{$b$-grading}

\def\smxnzi{ \sum_{\xn=0}^\infty }

\def\wdv#1{ |#1|_{\mathrm{wd}} }
\def\swdv#1{ \left| #1 \right|_{\mathrm wd} }
\def\wdfcv#1{ |#1|_{\mathrm{df}} }

\def\sTL{ \stfont{T} }
\def\sTLv#1{ \sTL_{#1} }

\def\sTLa{ \sTLv{a} }
\def\sTLab{ \sTLv{a,b} }
\def\sTLabc{ \sTLv{a,b}^{\supset} }

\def\shcr{ \mathsf{h} }
\def\shfr{ \mathsf{q} }
\def\shcrh{\shcr^{\hlf} }
\def\shcrmh{ \shcr^{-\hlf} }

\def\clN{ N }
\def\clNo{ \clN + 1}

\def\xDclv#1{ \xD_{#1} }
\def\xDclN{ \xDclv{\clN} }
\def\xDclNo{ \xDclv{\clNo} }

\def\xDclvv#1#2{ \xD_{#1,#2} }

\def\xDNoa{ \xDclvv{\xca+1}{\incra} }
\def\xDNoamo{ \xDclvv{\xca+1}{\incra-1} }

\def\xDpclvv#1#2{ \xD'_{#1,#2} }

\def\xDpNoa{ \xDpclvv{\xca+1}{\incra} }
\def\xDpNof{ \xDpclvv{\xca+1}{\ncrD} }

\def\xtD{ \tilde{\xD} }
\def\xtDv#1{ \xtD_{#1} }
\def\xtDNo{ \xtDv{\xca+1} }
\def\xtDNo{ \xtDv{\xca} }
\def\xtDN{ \xtDv{\xca} }

\def\xtDNov#1{ \xtD_{\xca+1;#1} }
\def\xtDNob{ \xtDNov{\crcb} }
\def\xtDNobo{ \xtDNov{\crcb+1} }
\def\xtDNobv#1{ \xtDNov{\crcb,#1} }
\def\xtDNobg{ \xtDNobv{\prjg} }

\def\xtDNobgo{ \xtDNobv{\prjg+1} }
\def\xtDNobf{ \xtDNobv{\pncrb-1} }

\def\xhD{ \hat{\xD} }

\def\xhDNov#1{ \xhD_{\xca+1,#1} }
\def\xhDNob{ \xhDNov{\crcb} }

\def\xDs{ \xD_{\spmp} }
\def\xDscr{ \xD_{\spmp,\scs} }

\def\xDcir{ \xD_{\circ} }

\def\xDov#1{ \xD_{#1} }
\def\xDon{ \xDov{n} }
\def\xDom{ \xDov{m} }

\def\xBov#1{ \xB_{#1} }
\def\xBmpn{ \xBov{m+n} }

\def\eqdiam{equatorial diameter}

\def\ttngl{tangle}
\def\ttnglv#1#2{$(#1,#2)$-\ttngl}
\def\TLttnglv#1#2{\taTL\ \ttnglv{#1}{#2}}
\def\ttnglmn{\ttnglv{m}{n}}
\def\ttnglnn{\ttnglv{n} {n}}
\def\TLttnglmn{\TLttnglv{m}{n}}
\def\TLttnglnn{\TLttnglv{n}{n}}

\def\qi{ q^{-1} }
\def\qpqi{ q + \qi }
\def\mqpqi{ - (\qpqi) }

\def\tHom{ \widetilde{\Hom} }

\def\xlam{ \lambda }
\def\xKhl{ \xKhv{\xlam} }

\def\shm{ \mu }
\def\hgrshklm{ \hgrshv{k}{l}{\shm} }

\def\elcf{ f }
\def\elcfv#1{ f_{#1} }

\def\Zcat{ \mathrm{Z} }
\def\Zcatv#1{ \Zcat(#1) }

\def\gAB{ g_{AB} }

\def\xIdv#1{ \xId_{#1} }
\def\xIdn{ \xIdv{n} }
\def\xIdvv#1#2{ \xId_{#1,#2} }
\def\xIdnL{ \xIdvv{n}{\xL} }

\def\xL{ L }
\def\xLp{ \xL' }
\def\xLb{ \bar{\xL} }
\def\xLcv#1{ \xL_{#1} }
\def\xLcN{ \xLcv{\xca} }
\def\xLcNo{ \xLcv{\xca+1} }

\def\xLpcv#1{ \xLp_{#1} }
\def\xLpcN{ \xLpcv{\xca} }

\def\xD{ D }

\def\xczt{ \zeta }
\def\xcztv#1{ \xczt(#1) }
\def\xcztL{ \xcztv{\xL} }

\def\cmfont{ \mathbf }
\def\cmA{ \cmfont{A} }

\def\xca{ N }
\def\xnu{ \nu }

\def\sB{ \mathrm{B} }

\def\nBv#1{ \gvv{#1} }
\def\nBD{ \nBv{\xD} }

\def\ncr{ \mathrm{cr} }
\def\ncr{ \times }
\def\ncr{ \chi }
\def\ncrv#1{ \ncr_{\!#1} }
\def\ncrD{ \ncrv{\xD} }
\def\ncrL{ \ncrv{\xL} }

\def\ncrm{ \chi^{!} }
\def\ncrmv#1{ \ncrm_{\!#1} }
\def\ncrmL{ \ncrmv{\xL} }

\def\ncriv#1{ \ncr_{\!#1}^{\mathrm in} }
\def\ncriD{ \ncriv{\xD} }

\def\ncrt{ \tilde{\ncr} }
\def\ncrtv#1{ \ncrt_{#1} }
\def\ncrtD{ \ncrtv{\xD} }

\def\sBv#1{ \sB(#1) }
\def\sBD{ \sBv{\xD} }
\def\sBDp{ \sBv{\xDp} }

\def\xfs{ s }
\def\xfsv#1{ \xfs_{#1} }

\def\xfsvv#1#2{ \xfsv{#1}(#2) }
\def\xfsLN{ \xfsvv{\xL}{\xca} }

\def\xn{ n }
\def\xt{ t }

\def\degq{ \deg_q }

\def\qhlf{ q^{ \frac{1}{2} } }
\def\qmhlf{ q^{-\frac{1}{2} } }

\def\gl{\mathrm{ l} }
\def\glv#1{ \gl_{#1} }
\def\glD{ \glv{\xD} }
\def\glL{ \glv{\xL} }
\def\ge{ \mathrm{e} }
\def\gev#1{ \ge_{#1} }
\def\geD{ \gev{\xD} }
\def\geL{ \gev{\xL} }
\def\gv{\kappa}
\def\gvv#1{ \gv_{#1} }
\def\gvD{ \gvv{\xD} }
\def\gvL{ \gvv{\xL} }

\def\xfrm{ \phi }
\def\xfrmv#1{ \xfrm_{#1} }
\def\xfrmD{ \xfrmv{\xD} }

\def\pncr{ \pi }
\def\pncrv#1{ \pncr_{#1} }
\def\pncrb{ \pncrv{\crcb} }

\def\xM{ M }

\def\bdA{ A }

\def\Nsgn{ (-1)^{\xca} }

\def\brbet{ \beta }
\def\xnp{ n_+ }
\def\xnm{n_- }

\def\stlcl#1{ {\scriptstyle #1 } }

\def\drlmA{ \dlm \chA_i }

\def\xalg{ \IQ }

\def\mnf{ f }
\def\mnfv#1{ \mnf_{#1} }
\def\mnfN{ \mnfv{\xca} }
\def\mnfi{ \mnfv{\mathrm i} }
\def\mnff{ \mnfv{\mathrm f} }

\def\mntf{ \tilde{\mnf} }
\def\mntfv#1{ \mntf_{#1} }
\def\mntfN{ \mntfv{\xca} }

\def\incra{ \alpha }
\def\incrb{ \alpha\p }
\def\xlbv#1{ v_{#1} }
\def\xlba{ \xlbv{\incra} }
\def\xlbao{ \xlbv{\incra+1} }
\def\xlbb{ \xlbv{\incrb} }

\def\crcb{ \beta }
\def\crcbp{ \beta\p }
\def\ylbv#1{ c_{#1} }
\def\ylbb{ \ylbv{\crcb} }
\def\ylbbp{ \ylbv{\crcbp} }

\def\prjg{ \gamma }
\def\prjgp{ \gamma\p }

\def\zlbvv#1#2{ p_{#1,#2} }
\def\zlbbv#1{ \zlbvv{\crcb}{#1} }
\def\zlbbg{ \zlbbv{\prjg} }
\def\zlbbgp{ \zlbbv{\prjgp} }

\def\xmg{ g }
\def\xmgNv#1{ \xmg_{\xca,#1} }
\def\xmgNa{ \xmgNv{\incra} }
\def\xmgN{ \xmg_{\xca} }

\def\xtmg{ \tilde{\xmg} }
\def\xtmgv#1{ \xtmg_{#1} }
\def\xtmgbg{ \xtmgv{\crcb,\prjg} }

\def\ltrf{local transformation}
\def\Ltrf{Local transformation}

\def\lrpl{local replacement}
\def\Lrpl{Local replacment}

\def\ytng{ \tau }

\def\ytngsv#1{ \ytng_{#1} }
\def\ytngsi{ \ytngsv{\mathrm i} }
\def\ytngsf{ \ytngsv{\mathrm f} }
\def\ytngse{ \ytngsv{\mathrm e} }
\def\ytngsfp{ \ytngsf\p }
\def\ytngsc{ \ytngsv{\mathrm c} }
\def\ytngsci{ \ytngsv{\mathrm{c},\xki} }

\def\ytngki{ \xKhv{\ytngsi} }
\def\ytngkf{ \xKhv{\ytngsf} }
\def\ytngkc{ \xKhv{\ytngsc} }

\def\ytngkfp{ \xKhv{\ytngsf'}}

\def\xDsv#1{ \xD_{#1} }
\def\xDsi{ \xDsv{\mathrm i} }
\def\xDsf{ \xDsv{\mathrm f} }

\def\xDsc{ \xDsv{\mathrm c} }
\def\xDscv#1{ \xDsv{\mathrm{c},#1} }
\def\xDsci{ \xDscv{\xki} }
\def\xDscip{ \xDsci' }
\def\xDscipp{ \xDsci''}
\def\xDse{\xDsv{\mathrm e} }
\def\xDsep{ \xDse'}
\def\xDscp{ \xDsc\p }

\def\hbnd{ M_{\mathrm h} }

\def\sgmm{ s }

\def\yncr{ n_{\times} }
\def\yncrv#1{ \yncr(#1) }
\def\yncrpv#1{ \yncr'(#1) }
\def\yncrD{ \yncrv{\xD} }

\def\yncc{ n_{\circ} }
\def\ynccv#1{ \yncc(#1) }

\def\tdgpr{degree preserving}

\def\xhsh{ m }
\def\xhshv#1{ \xhsh_{#1} }
\def\xhshf{ \xhshv{\mathrm f}}

\def\xqsh{ n }
\def\xqshv#1{ \xqsh_{#1} }
\def\xqshf{ \xqshv{\mathrm f}}

\def\ztau{ \tau }

\def\svrt{ \mathfrak{V} }
\def\svrta{ \svrt_{\mathrm{ad}} }
\def\svrti{ \svrt_{\mathrm{in}} }
\def\xvrt{ v }
\def\spmp{ s }
\def\sipvr{ \spmp_{\xvrt} }
\def\spvr{ \sipvr }
\def\xabms{ |\!|\spmp|\!| }

\def\tstt{state}

\def\tstrt{strut}
\def\tstrtl{\tstrt\ line}
\def\tjmp{jumping}
\def\tjmpc{\tjmp\ circle}
\def\trlx{relaxed}
\def\trlxc{\trlx\ circle}

\def\crb{ c }

\def\nvcr#1{ n_{\mathrm #1} }
\def\njcr{ \nvcr{j} }
\def\nscr{ \nvcr{s} }
\def\nscrb{ n_{\mathrm{s},\crb} }
\def\nwcr{ \nvcr{w} }
\def\nwcrb{ n_{\mathrm{w},\crb} }

\def\njrc{ n_{\mathrm{r},c} }

\def\tstrght{straight}
\def\txstr{straight}
\def\txstrs{\txstr\ segment}
\def\tstrghtc{\tstrght\ circle}

\def\twndg{winding}
\def\twndgc{\twndg\ circle}

\def\ntwd{not widening}
\def\scs{ c }

\def\xtusc{ \tau_{s,c} }

\def\xE{ E }
\def\xEo{\xE^1}
\def\xEovv#1#2{ \xEo_{#1,#2} }
\def\xEoij{ \xEovv{i}{j} }

\def\xki{ k }
\def\yki{ k }
\def\yvki{ \yki(i) }
\def\ysvki{ \yki^2(i) }

\def\zmi{ i }
\def\zmj{ j }

\def\xfd{ d }
\def\xfdv#1{ \xfd_{#1} }
\def\xfdN{ \xfdv{\xca} }
\def\xti{ \tilde{i}}

\def\ptJqLi { \mathrm{J}_{\xL,\thicksim}(b,q)}

\def\Stosic{Sto\v si\'c}

\def\SUv#1{ \mathrm{SU}(#1) }

\def\xltone{I}
\def\xlttwo{II}
\def\xltthree{III}

\def\prmlt{ \mu }
\def\prmltv#1{ \prmlt_{#1} }
\def\prmltijg{ \prmltv{ij,\gamma} }
\def\prmltij{ \prmltv{ij} }
\def\prmltijk{ \prmltv{ij,\xki} }
\def\tprmlt{ \tilde{\prmlt} }
\def\tprmltv#1{ \tprmlt_{#1} }
\def\tprmltij{ \tprmltv{ij} }

\def\lumps{lump sum}

\def\mtot{ m^{\mathrm{tot}} }
\def\mtotv#1{ \mtot_{#1} }
\def\mtotja{ \mtotv{j,a} }

\def\betbr{ \beta }
\def\xLv#1{ \xL_{#1} }
\def\xLb{ \xLv{\betbr} }

\def\xcrsp{
\xygraph{
!{0;/r1.5pc/:}
[u(0.5)]
!{\xoverv}
}
}

\def\xpver{
\xygraph{
!{0;/r1.5pc/:}
[u(0.5)]
!{\xunoverv}
}
}

\def\xphor{
\xygraph{
!{0;/r1.5pc/:}
[u(0.5)]
!{\xunoverh}
}
}

\def\zoverv{
\begin{tikzpicture} \draw + (1,0) -- ++(0,1);
\draw  [line width=6pt, draw=black] (0,0) -- ++(1,1);
\draw (0,0) -- ++(1,1);
\end{tikzpicture}
 }

 \def\zoverv#1#2{
 \begin{tikzpicture}[scale=#1,baseline=11*#1-#2*#1]
 \path[use as bounding box] (-0.1,-0.1) rectangle (1.1,1.1);
 \draw (0,1) -- ++(1,-1);
 \draw[line width=6pt, draw=white] (0,0) -- ++(1,1);
 \draw (0,0) -- ++(1,1);
\end{tikzpicture}
 }

\def\zunoverv#1#2{
 \begin{tikzpicture}[scale=#1,baseline=11*#1-#2*#1]
 \path[use as bounding box] (-0.1,-0.1) rectangle (1.1,1.1);
 \draw (0,0) to [out=45,in=-45] (0,1);
 \draw (1,0) to [out=135,in=-135] (1,1);
\end{tikzpicture}
 }

 \def\zunoverh#1#2{
 \begin{tikzpicture}[scale=#1,baseline=11*#1-#2*#1]
 \path[use as bounding box] (-0.1,-0.1) rectangle (1.1,1.1);
 \draw (0,0) to [out=45,in=135] (1,0);
 \draw (0,1) to [out=-45,in=-135] (1,1);
\end{tikzpicture}
 }

\def\zcirc#1#2{
 \begin{tikzpicture}[scale=#1,baseline=11*#1-#2*#1]
 \path[use as bounding box] (-0.1,-0.1) rectangle (1.1,1.1);
 \draw (0.5,0.5) circle (0.5);
\end{tikzpicture}
 }

\def\zposfr#1#2{
 \begin{tikzpicture}[scale=#1,baseline=11*#1-#2*#1]
 \path[use as bounding box] (-0.1,-0.1) rectangle (0.9,1.1);
 \draw (0.6,0.25) to [out=180,in=-90] (0,1);
 \draw [line width=6pt,draw=white] (0,0) to [out=90,in=180] (0.6,0.75);
 \draw (0,0) to [out=90,in=180] (0.6,0.75);
 \draw (0.6,0.75) .. controls +(0:0.4) and +(0:0.4) .. (0.6,0.25);
\end{tikzpicture}
 }

 \def\zstvr#1#2{
 \begin{tikzpicture}[scale=#1,baseline=11*#1-#2*#1]
 \path[use as bounding box] (-0.1,-0.1) rectangle (0.6,1.1);
 \draw (0.25,0) -- (0.25,1);
\end{tikzpicture}
 }

 \def\zcposfr#1#2{
 \begin{tikzpicture}[scale=#1,baseline=11*#1-#2*#1-0.3cm]
 \path[use as bounding box] (-1.2,-0.1) rectangle ++(2.4,0.8);
 \begin{scope}[rotate=90]
 \draw [thkln] (0.6,0.25) to [out=180,in=-90] (0,1);
 \draw [line width=6pt,draw=white] (0,0) to [out=90,in=180] (0.6,0.75);
 \draw [thkln] (0,0) to [out=90,in=180] (0.6,0.75);
 \draw [thkln] (0.6,0.75) .. controls +(0:0.4) and +(0:0.4) .. (0.6,0.25);
 \draw [ptzer] (-0.6,-0.3) rectangle ++ (1.2,0.3);
 \draw [thkln] (0,-0.3) -- (0,-0.9) node [near end, above] {$\scriptstyle a$};
 \end{scope}
\end{tikzpicture}
 }

  \def\zcstvr#1#2{
 \begin{tikzpicture}[scale=#1,baseline=11*#1-#2*#1-0.3cm]
 \path[use as bounding box] (-0.9,-0.5) rectangle (0.9,0.5);
\draw [ptzer] (-0.15,-0.6) rectangle ++(0.3,1.2);
\draw [thkln] (-0.75,0) -- (-0.15,0) (0.15,0) -- (0.75,0) node [near end,above] {$\scriptstyle a$};
\end{tikzpicture}
 }

\def\cblth{1.3pt}
\def\ljwp{1pt}
\def\ocrw{6pt}

\tikzset{ptzer/.style={line width=\ljwp} }
\tikzset{ptzert/.style={line width=\ljwp,fill=white}}
\tikzset{ptone/.style={line width=\ljwp,fill=gray!30}}
\tikzset{pttwo/.style={line width=\ljwp,fill=white,pattern=horizontal lines} }
\tikzset{ptthr/.style={line width=\ljwp,pattern=north east lines} }
\tikzset{ptfour/.style={line width=\ljwp,pattern=crosshatch} }
\tikzset{vthln/.style={line width=\ljwp} }

\tikzset{menvone/.style={scale=0.5,baseline=-1.5} }
\tikzset{menvtwo/.style={scale=0.70,baseline=-1.5} }
\tikzset{menvthree/.style={scale=0.35, baseline=-1.5} }
\tikzset{menvfour/.style={scale=0.25, baseline=-1.5} }
\tikzset{menvzer/.style={scale=1,baseline=-1.5} }

\tikzset{menvrtone/.style={scale=0.45,rotate=-45,baseline=-3.5} }
\tikzset{menvrtsone/.style={scale=0.65,rotate=-45,baseline=-3.5} }
\tikzset{menvrtonep/.style={scale=0.5,baseline=-3.5} }

\tikzset{lnovr/.style={line width=6pt,color=white}}
\tikzset{lnovrtw/.style={line width=4pt,color=white}}
\tikzset{thkln/.style={line width=\cblth} }
\tikzset{thnln/.style={} }
\tikzset{thkc/.style={line width=0.6pt} }
\tikzset{bcrc/.style={line width=0.2cm,color=gray!40,opacity=0.5} }
\tikzset{bcrct/.style={line width=0.25cm,color=gray!20} } 

\pgfdeclaredecoration{ignore}{final}
{
\state{final}{}
}

\pgfdeclaremetadecoration{middle}{initial}{
    \state{initial}[
        width={(\pgfmetadecoratedpathlength - \the\pgfdecorationsegmentlength)/2},
        next state=middle
    ]
    {\decoration{moveto}}

    \state{middle}[
        width={\the\pgfdecorationsegmentlength},
        next state=final
    ]
    {\decoration{curveto}}

    \state{final}
    {\decoration{ignore}}
}

\tikzset{middle segment/.style={decoration={middle},decorate, segment length=#1}}

\numberwithin{equation}{section}

\title[A tail of Khovanov homologyl]
{Khovanov homology of a unicolored B-adequate link has a tail}

\author[L.~Rozansky]{Lev Rozansky}
\address{
L.~Rozansky\\
Department of Mathematics\\
University of North Carolina at Chapel Hill\\
CB \# 3250, Phillips Hall\\
Chapel Hill, NC 27599
}
\email{rozansky@math.unc.edu}
\thanks{This work was supported in part by the NSF grants DMS-0808974 and DMS-1108727}

\input xy.sty
\input xyarrow.tex

\begin{document}
\draft
\maketitle
\begin{abstract}

C.~Armond~\cite{CK10} and S.~Garoufalidis and T.~Le~\cite{GL11} have shown that a \uclrd\ Jones polynomial of a \tBadq\ link has a stable tail at large colors. We categorify this tail by showing that Khovanov homology of a \uclrd\ link also has a stable tail, whose graded Euler characteristic coincides with the tail of the Jones polynomial.

\end{abstract}

\begin{spacing}{0.65}
\tableofcontents
\end{spacing}

\section{Notations and basic facts}



\subsection{Links, \tadq\ links and their diagrams}
All diagrams of tangles and links  in this paper are framed, we assume blackboard framing. Links are presumed unframed, unless specified otherwise.

Each crossing of a link diagram $\xD$ can be `spliced' in two ways, we call them \Asplng\ and \Bsplng:
\def\zpsh{1}
%
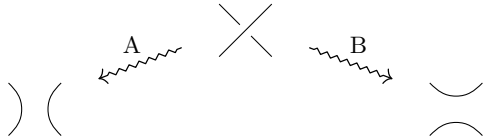
\[
\begin{tikzpicture}[menvthree,commutative diagrams/every diagram]
\node (ul) at ( -2,-0.5){};
\node (dl) at (-6,-2){};
\path[commutative diagrams/.cd, every arrow, every label]
(ul) edge[commutative diagrams/squiggly] node[pos=0.4,swap] {A} (dl);
\node (ur) at ( 2,-0.5){};
\node (dr) at (6,-2){};
\path[commutative diagrams/.cd, every arrow, every label]
(ur) edge[commutative diagrams/squiggly] node[pos=0.4] {B} (dr);
\begin{scope}
\draw (-\zpsh,\zpsh) -- (\zpsh,-\zpsh);
\draw [lnovr] (-\zpsh,-\zpsh) -- (\zpsh,\zpsh);
\draw (-\zpsh,-\zpsh) -- (\zpsh,\zpsh);
\end{scope}
\begin{scope}[xshift=-8cm,yshift=-3cm,rotate=90]
\draw(-\zpsh,-\zpsh) to [out=45,in=-180] (0,-\zpsh*0.5) to [out=0,in=135] (\zpsh,-\zpsh);
\draw(-\zpsh,\zpsh) to [out=-45,in=-180] (0,\zpsh*0.5) to [out=0,in=-135] (\zpsh,\zpsh);
\end{scope}
\begin{scope}[xshift=8cm,yshift=-3cm]
\draw(-\zpsh,-\zpsh) to [out=45,in=-180] (0,-\zpsh*0.5) to [out=0,in=135] (\zpsh,-\zpsh);
\draw(-\zpsh,\zpsh) to [out=-45,in=-180] (0,\zpsh*0.5) to [out=0,in=-135] (\zpsh,\zpsh);
\end{scope}
\end{tikzpicture}
\]

Let $\sBD$ denote the diagram which consists of two parts: the circles resulting from \Bsplng\ of all crossings of $\xD$ (\tBcr s) and segments connecting those circles at places where crossings were in $\xD$ (\tstrt s). Schematically, one passes from $\xD$ to $\sBD$ in the following way:
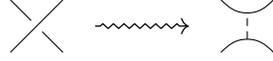
\def\zpsh{1}
\[
\begin{tikzpicture}[menvthree]
\node (l) at (-2.2,0) {};
\node (r) at (2.2,0) {};
\path[commutative diagrams/.cd, every arrow, every label]
(l) edge[commutative diagrams/squiggly] (r);
\begin{scope}[xshift=-4cm]
\draw (-\zpsh,\zpsh) -- (\zpsh,-\zpsh);
\draw [lnovr] (-\zpsh,-\zpsh) -- (\zpsh,\zpsh);
\draw (-\zpsh,-\zpsh) -- (\zpsh,\zpsh);
\end{scope}
\begin{scope}[xshift=4cm]
\draw(-\zpsh,-\zpsh) to [out=45,in=-180] (0,-\zpsh*0.5) to [out=0,in=135] (\zpsh,-\zpsh);
\draw(-\zpsh,\zpsh) to [out=-45,in=-180] (0,\zpsh*0.5) to [out=0,in=-135] (\zpsh,\zpsh);
\draw [densely dashed] (0,-\zpsh*0.5) -- (0,\zpsh*0.5);
\end{scope}
\end{tikzpicture}
\]
where the arcs in the right diagram are parts of \tBcr s and the dashed segment is a \tstrt.

A crossing of $\xD$ and its \tstrt\ in $\sBD$ are called \tBadq, if the \tstrt\ connects two different \tBcr s. A framed diagram $\xD$ is \tBadq, if all of its crossings are \tBadq. Finally, a framed link is \tBadq, if it can be represented by a \tBadq\ framed diagram. An unframed link $\xL$ is called \tBadq, if there is at least one framed \tBadq\ diagram that represents it. Note that if an unframed link is \tBadq, then, generally, it can not be represented by a \tBadq\ framed diagram for all framings.

All alternating links are \tBadq, but not all \tBadq\ links are alternating: an example of this is a torus knot $\yTmmn$, $n \geq m \geq 3$. More generally, a link constructed by closing a totally negative braid is \tBadq. Torus knots $\yTmn$, $n\geq m \geq 3$ provide examples of links which are not \tBadq.

Here are some notations associated with a link diagram $\xD$ throughout the paper:
\begin{center}
\begin{tabular}{c c l}
$\svrt$ & -- &a set of crossings (\tstrt s) in $\xD$ or in $\sBD$
\\
$\ncrD$ & -- &the number of crossings in $\xD$
\\
$\ncriD$ & -- & the number of B-\tinadq\ crossings in $\xD$
\\
$\gvD$ & -- &  the number of \tBcr s in $\sBD$
\\
$\xfrmD$ & -- & the framing number of $\xD$
\end{tabular}
\end{center}

The following is an easy corollary of the results of section 7.7 of~\cite{Kh99}:
\begin{theorem}
\label{thm:frm}
The numbers $\ncrv{\xL}$ and $\gvv{\xL}$ are topological invariants of a \tBadq\ framed link $\xL$, because they do not depend on the choice of representative \tBadq\ diagram $\xD$ for $\xL$. Moreover, if
\tBadq\ framed links $\xL$ and $\xLp$ differ only by framing, then
\begin{equation}
\label{eq:afrm}
\ncrv{\xLp} - \ncrv{\xL} = \gvv{\xLp} - \gvv{\xL} = - (\xfrmv{\xLp} - \xfrmv{\xL}).
\end{equation}
\end{theorem}

For an unframed \tBadq\ link $\xL$ we define the minimal crossing number $\ncrmL$ as the minimum among the numbers $\ncrD$ for \tBadq\ framed diagrams $\xD$ representing $\xL$.


\subsection{The \tKbr\ and the Jones polynomial}
The \tKbr\ of a framed tangle diagram is defined by the \splng\ relation and the unknot normalization condition:
\begin{equation}
\label{eq:dkbr}
\xKbrBv{\zoverv{0.75}{1}}
\;=\;
\qhlf\xKbrBv{\zunoverv{0.75}{1}} \;+\; \qmhlf \xKbrBv{\zunoverh{0.75}{1}},
\qquad
\xKbrBv{\zcirc{0.75}{1}} \; = \; -(\qpqi).
\end{equation}
Thus defined, the bracket is framing-dependent:
\[
\xKbrBv{
\zposfr{0.75}{1}
}\;= \; -q^{\frac{3}{2}}\xKbrBv{
\zstvr{0.75}{1}
}.
\]

The Jones polynomial of a framed link $\xL$ is the \tKbr\ of its diagram: $\pJqL=\xKbrv{\xL}$.

\subsection{Cables and coloring}

We introduce coloring of tangle and link components through cabling and \tJWp s. A cable of a strand is depicted by using a thicker line with the label indicating the number of strands, and the \tJWp\ is depicted by a box:
\[
\begin{tikzpicture}
\draw [thkln] (0,0) -- (1,0); \node [above] at (0.5,0) {$\stlcl{a}$};
\end{tikzpicture}
\;=\;
\begin{tikzpicture}[baseline=-4]
\draw (0,0.4) -- (1,0.4);
\draw (0,-0.4) -- (1,-0.4);
\node at (0.5,0.1) {$\vdots$};
\draw[decorate,decoration={brace, amplitude=4pt}]
        (1.2,0.4) -- (1.2,-0.4) node[midway, right=2pt]{$a$};
\end{tikzpicture},
\qquad
\begin{tikzpicture}[scale=0.75,baseline=-2.5]
\draw [ptzer] (-0.15,-0.6) rectangle ++(0.3,1.2);
\draw [thkln] (-.65,0) -- (-0.15,0) node [near start, above] {$\scriptstyle a$}
(0.15,0) -- (.65,0);
\end{tikzpicture}
\]

For a positive integer $\xca$ let $\pJqaL$ denote the \uclrd\ Jones polynomial of $\xL$, that is, all components of $\xL$ are colored by the same color $\xca$. A coloring of a link component by $\xca$ means that we assign the $(\xca+1)$-dimensional irreducible representation of $\SUt$ to it. Equivalently, the color $\xca$ means that the link component is $\xca$-cabled and we place the \tJWp\ on this cable.

In this paper we consider \uclrd\ links, that is, links, all of whose components are colored by the same number $\xca$.
Their colored Jones polynomial $\pJqaL$ is a Laurent polynomial of $q^2$ up to an overall factor: if $\xL$ is presented by a diagram $\xD$, then
\[
q^{\hlf\ncrD\xca^2 + \gvD \xca } \pJqaL \in\ZZtqqi.
\]

%

\subsection{Homological notations}
Let $\caA$ be a finitely generated additive category: objects of $\caA$ are finite sums of elements of a finite set $\stA$.
Let $\cKommA$ denote the homotopy category of its complexes bounded from below: an object of $\cKommA$ is a chain
\begin{equation}
\label{eq:complex}
\chA = ( \cdots \rightarrow A_{i+1}\rightarrow A_i \rightarrow\cdots\rightarrow A_m),
\end{equation}
where $A_i = \bigoplus_{\oba\in\stA} m_{i,\oba}\,\oba$ and $ m_{i,\oba}\in\ZZ_{\geq 0}$ are the multiplicities of generators. The notation $m$ for multiplicity is treated in this paper as an arbitrary constant, so the appearance of $m$ in different expressions does not imply that there is a relation between the multiplicities, unless it is stated specifically. The special multiplicities appearing in a presentation of the \cJWp\ are denoted by $\prmlt$.

We use a non-standard notation for the translation functor: $\shcr \chA = \chA[1]$, which allows us to define a functor $p(\shcr)$ for any polynomial $p(x)$ with integer non-negative coefficients. In particular, we use a functor ${i \brace j}_{\shcr}$ based on a combinatorial polynomial
\begin{equation}
\label{eq:combp}
{i \brace j}_{x} =
\frac
{
(1-x^{2i})(1-x^{2i-2})\cdots(1-x^{2i -2j + 2})
}
{
(1-x^2)(1-x^4)\cdots(1-x^{2j})
}.
\end{equation}

We also use a non-standard notation for the cone of two complexes:
\begin{equation}
\label{eq:cnnt}
\Cnv{\shcr \chA \xrightarrow{f} \chB} = \Conv{\chA\xrightarrow{f} \chB}.
\end{equation}
in order to emphasize the fact that the cone $\Conv{\chA\rightarrow\chB}$ can be presented as a sum $\shcr \chA \oplus \chB$ deformed by an extra differential $\chA\xrightarrow{f} \chB$. Moreover, when we work with bi-graded Khovanov complexes, there may be some confusion about which of two gradings is homological, but our non-standard notation\rx{eq:cnnt} specifies all degree shifts explicitly.

The homological order $\hmord{\obO}$ of an object $\obO\in\cKommA$ is the minimum number $m$, for which $\obO$ can be presented by a complex\rx{eq:complex}.

Consider a direct system of complexes of $\cKommA$:
$\chA_0\rightarrow\chA_1\rightarrow\cdots.$ If this system is `Cauchy', that is,
if for the cones $\chB_i =  \Conv{\chA_{i-1} \rightarrow \chA_{i}}$ there is a limit
$\lim_{i\rightarrow\infty}\hmord{\chB_i} = \infty$,
then there exists a direct limit $\drlmA$.

Since $\chA_{i}\hteqv \Conv{\shcr^{-1} \chB_i\rightarrow\chA_{i-1}}$, 
the direct limit $\drlmA$ can be viewed as a result of attaching the complexes $\chB_i$ one after another to the initial complex $\chB_0=\chA_0$, hence we use the following notation for the complex $\drlmA$:
\begin{equation}
\label{eq:mtcn}
\drlmA \hteqv \boxed{
\cdots\rightarrow \chB_i \rightarrow \cdots
}_{\;i=0}^{\;\infty}
\end{equation}
In fact, if all $\chB_i$ are `homologially minimal' representatives of their equivalence classes, then the sum $\bigoplus_{i=0}^{\infty} \chB_i$ is well-defined (every chain object is finitely generated) and $\drlmA$ is homotopy equivalent to $\bigoplus_{i=0}^{\infty} \chB_i$ defomed by adding extra differentials $\chB_i\xrightarrow{f_{ij}}\chB_j$ for all pairs $i>j$.

We refer to the \rhs of \ex{eq:mtcn} as a \emph{\tmcn}, and we also use a similar notation for the complex\rx{eq:complex}: $\chA = \boxed{\cdots\rightarrow \shcr^iA_i\rightarrow\cdots}_{\;i=m}^{\;\infty}$. Note the use of the functor $\shcr$ to set explicitly the correct homological degree of the chain object $A_i$ in the \tmcn.

If a \tmcn\ $\chA$ is generated by complexes $\chB_a$:
\begin{equation}
\label{eq:lsfmc}
\chA = \Cnv{
\cdots\rightarrow \bigoplus_{j,a} m_{ij,a}\shcr^j\, \chB_a\rightarrow\cdots
}\;,
\end{equation}
(where $m_{ij,a}$ are multiplicities) but we do not care how those complexes are arranged within the \tmcn, then we use a `\lumps' notation
\[
\chA = \Pcnv{ \bigoplus_{j,a} \mtotja\,\shcr^j\,\chB_a
}, \qquad
\mtotja = \sum_i m_{ij,a},
\]
because, as a complex, $\chA$ is a sum of $\chB_a$ with total multiplicities $\sum_i m_{ij,a}$ deformed by an extra differential depicted as $\sdff$.

If the category $\caA$ is abelian, then we can compute the homology of the \tmcn\rx{eq:mtcn} with the help of the filtered complex spectral sequence. The $\xEo$ term of this spectral sequence is the sum of homologies of $\chB_i$: $\xEo = \bigoplus_{i=0}^\infty \Hm(\chB_i)$ and it is determined by the \lumps\ form of the \tmcn.
\begin{remark}
\label{rmk:bndss}
Since subsequent terms in the spectral sequence get only smaller, there is a bound on the homological order of the homology of\rx{eq:mtcn} in terms of homological orders of its constituent complexes:
\begin{equation*}
\hlmord{
\Hm\left( \;\boxed{
\cdots\rightarrow \chB_i \rightarrow \cdots
}_{\;i=0}^{\;\infty}
\;
\right)
}
\geq
\min_{i} \hmord{\Hm(\chB_i)}.
\end{equation*}
In particular, for the \lumps\ \tmcn\rx{eq:lsfmc}
\[
\hmord{ \Hm(\chA) } \geq \min\{\hmord{\Hm(\chB_a)} + j\colon \mtotja\neq 0\}.
\]
\end{remark}

\subsection{Khovanov homology}
\hyphenation{ca-no-po-ly}
In defining Khovanov complexes~\cite{Kh99} for tangles we follow the cobordism based approach of D.~Bar-Natan~\cite{BN05}, albeit with a different grading convention. We still have two degrees: \thdgr\ $\dgh$ and \tqdgr\ $\dgq$, and we use the notations $\shcr$ and $\shfr$ for their translation functors (these functors increase the corresponding degrees by 1). The \tqdgr\ is the genuine homological degree: it takes values in $\ZZ$ and its parity determines the sign factors. The \thdgr\ is `pseudo-homological', it takes values in $\hlf\ZZ$ and it has no impact on signs, however it is the \thdgr\ shift functor $\shcr$ which is present explicitly in the \tKhbr.

In our notations, \tKhbr\ of a crossing and of the unknot are
\begin{equation}
\label{eq:dkhbr}
\xvKhv{
\zoverv{0.75}{1}
}
\;=\;
\Cnv{
\shcr^{\hlf}
\xvKhv{\zunoverv{0.75}{1}
}
\xrightarrow{\;\sgmm \; }
\shcr^{-\hlf}
\xvKhv{\zunoverh{0.75}{1}
}
}\;,\qquad
\xvKhv{
\zcirc{0.75}{1}
}
\;=\;
(\shfr + \shfr^{-1})\xalg,
\end{equation}
where $\shcr$ and $\shfr$ are degree shift functors, while
$\sgmm$ is the morphism corresponding to the saddle cobordism. Note that $\dgh\sgmm = -1$, while $\dgq\sgmm = 1$, so $\sgmm$ is odd.

Thus defined, \tKhbr\ is invariant under the first Reidemeister move only up to a degree shift:
\begin{equation}
\label{eq:frsh}
\xvKhv{
\zposfr{0.75}{1}
}\;= \; \shcr^{\hlf}\shfr \xvKhv{
\zstvr{0.75}{1}
}.
\end{equation}

Relations~\eqref{eq:dkhbr} transform into the relations~\eqref{eq:dkbr} after the substitution
$\shcr\mapsto q$, $\shfr\mapsto - q$, hence in our notations the graded Euler characteristic of \tKhom\ of a framed link equals its \tJpol:
\[
\pJqL = \sum_{i,j}(-1)^j q^{i+j} \dim\KHmvv{i}{j}(\xL).
\]

We will use the \tKhbr\ notation $\xKhv{-}$ very sparingly, because it clutters the pictures, especially when the diagrams are big. Nevertheless, we hope that the distinction between diagrams and their Khovanov complexes will be clear. Actually, we blur this distinction further by allowing the presence of \cJWp s within diagrams, since, strictly speaking, projectors are not diagrams but rather complexes within Bar-Natan's universal category.

\subsection{A \cJWp}
%

An $(a,b)$-tangle is an embedding of circles and segments, the segment endpoints coinciding with initial $a$ points or final $b$ points. Imagine that the tangle goes from the bottom up. Depending on the position of its endpoints, the segment is either \tstrght, or a cap, or a cup. If one of its endpoints is initial and the other is final (so the segment goes straight through the tangle), then the segment is \emph{\tstrght}, if both endpoints are initial, then the segment is a \emph{cap}, and if both of its  the segments are final, then the segment is a \emph{cup}.

The \emph{\twd} $\wdv{\tau}$ of a tangle $\tau$ is the number of its \txstrs s.
An $(a,a)$ tangle has an equal number of cups and caps, we call this number a \emph{\twdfc} and denote it as $\wdfcv{\tau}$. Obviously, $\wdfcv{\tau} = \hlf(a-\wdv{\tau})$.

A \emph{\tTL} (\taTL) tangle is a flat tangle which contains no circles. Let $\sTLa$ be the set of all $(a,a)$ \taTLt s. The \cJWp\
\begin{tikzpicture}[scale=0.5,baseline=-2.5]
\draw [ptzer] (0,-0.5) rectangle (0.25,0.5);
\draw [thkln] (0.25,0) -- (0.75,0);
\draw [thkln] (-0.75,0) -- (0,0);
\node [above] at (-0.5,0) {$\scriptstyle a$};
\end{tikzpicture}
was constructed independently by Frenkel, Stroppel and Sussan~\cite{FSS11}, Cooper and Krushkal~\cite{CK10} and by the author~\cite{Ro11}. It satisfies three essential properties:
it is a projector:
\begin{equation}
\label{eq:cmppr}
\begin{tikzpicture}[menvtwo]
\draw [ptzer] (-0.15,-0.6) rectangle ++(0.3,1.2);
\draw [ptzer] (0.65,-0.6) rectangle ++(0.3,1.2);
\draw [thkln] (-.65,0) -- (-0.15,0) node [near start, above] {$\scriptstyle a$}
(0.15,0) -- (.65,0) (0.65+0.3,0) -- (0.65+0.3+0.5,0);
\end{tikzpicture}
\;\hteqv\;
\begin{tikzpicture}[menvtwo]
\draw [ptzer] (-0.15,-0.6) rectangle ++(0.3,1.2);
\draw [thkln] (-.65,0) -- (-0.15,0) node [near start, above] {$\scriptstyle a$}
(0.15,0) -- (.65,0);
\end{tikzpicture}
\;,
\end{equation}
it annihilates cups and caps:
\begin{equation}
\label{eq:ccann}
\def\zext{0.4}
\begin{tikzpicture}[menvtwo]
\draw [ptzer] (-0.15,-0.6) rectangle ++(0.3,1.2);
\draw [thkln] (-.65-\zext,0) -- (-0.15,0) node [near start, above] {$\scriptstyle a$};
\draw [thkln] (0.15,0.45) -- (.65+\zext,0.45) node [near end,above] {$\scriptstyle b$};
\draw [thkln] (0.15,-0.45) -- (.65+\zext,-0.45) node [near end,below] {$\scriptstyle a-b-2$};
\draw (0.15,0.25) to [out=0,in=90] (0.5,0) to [out=-90,in=0] (0.15,-0.25);
\end{tikzpicture}
\;\hteqv\;
\begin{tikzpicture}[menvtwo,xscale=-1]
\draw [ptzer] (-0.15,-0.6) rectangle ++(0.3,1.2);
\draw [thkln] (-.65-\zext,0) -- (-0.15,0) node [near start, above] {$\scriptstyle a$};
\draw [thkln] (0.15,0.45) -- (.65+\zext,0.45) node [near end,above] {$\scriptstyle b$};
\draw [thkln] (0.15,-0.45) -- (.65+\zext,-0.45) node [near end,below] {$\scriptstyle a-b-2$};
\draw (0.15,0.25) to [out=0,in=90] (0.5,0) to [out=-90,in=0] (0.15,-0.25);
\end{tikzpicture}
\;\hteqv\; 0,
\end{equation}
and it
has a presentation as a cone of an identity braid and a complex
$
\begin{tikzpicture}[scale=0.4,baseline=-2.5]
\draw[ptone] (-0.15,-0.6) rectangle ++(0.3,1.2);
\draw [thkln] (-.65,0) -- (-0.15,0) node [near start, above] {$\scriptstyle a$}
(0.15,0) -- (.65,0);
\end{tikzpicture}
$
generated by \taTLt s with positive \twdfc\ and with non-negative \thdgr\ and \tqdgr\ shifts
\begin{equation}
\label{eq:projcn}
\begin{tikzpicture}[scale=0.75,baseline=-2.5]
\draw [ptzer] (-0.15,-0.6) rectangle ++(0.3,1.2);
\draw [thkln] (-.65,0) -- (-0.15,0) node [near start, above] {$\scriptstyle a$}
(0.15,0) -- (.65,0);
\end{tikzpicture}
\;\hteqv\;
\boxed{
\shcr
\begin{tikzpicture}[scale=0.75,baseline=-2.5]
\draw[ptone] (-0.15,-0.6) rectangle ++(0.3,1.2);
\draw [thkln] (-.65,0) -- (-0.15,0) node [near start, above] {$\scriptstyle a$}
(0.15,0) -- (.65,0);
\end{tikzpicture}
\longrightarrow
\begin{tikzpicture}
\draw [thkln] (-0.5,0) -- (0.5,0) node [midway,above] {$\scriptstyle a$};
\end{tikzpicture}
}\;,
\end{equation}
where
\begin{equation}
\label{eq:grproj}
\begin{tikzpicture}[scale=0.75,baseline=-2.5]
\draw[ptone] (-0.15,-0.6) rectangle ++(0.3,1.2);
\draw [thkln] (-.65,0) -- (-0.15,0) node [near start, above] {$\scriptstyle a$}
(0.15,0) -- (.65,0);
\end{tikzpicture}
\;=\;
\boxed{
\cdots\longrightarrow\shcr^i\bigoplus_{\substack{0\le j \le i \\ \gamma\in\sTLa, \wdfcv{\gamma}>0}} \prmltijg\,\shfr^j\,\xKhv{\gamma}
\longrightarrow\cdots
}_{\;i=0}^{\;\infty}
\end{equation}
and $\prmltijg$ are multiplicities.

\subsection{\tKhbr\ of colored tangles}
We define \tKhbr\ of colored tangles by cabling tangle components and adding at least one  \cJWp\ to each tangle component. This means that we allow semi-infinite complexes which may extend infinitely far into positive homological degree.

The colored \tKhbr\ is independent of the framing up to a degree shift:
\begin{equation}
\label{eq:cfrsh}
\zcposfr{0.75}{1}
\;= \; \shcr^{\hlf a^2}\shfr^a
\zcstvr{0.75}{1}
.
\end{equation}

\section{Results}
\subsection{Overview}
\subsubsection{Bounds on colored \tKhom}
Let $\xLcN$  denote the $\xca$-\uclrd\ version of the link $\xL$. For a \tBadq\ framed link $\xL$, a \emph{shifted} \tKhom\ of $\xLcN$ is defined by the formula
\begin{equation}
\label{eq:shdgt}
\tKHm(\xLcN) = \shcr^{\hlf\clN^2\ncrL} \shfr^{\gvL\xca} \KHm(\xLcN).
\end{equation}
In view of \ex{eq:afrm}, if links $\xL$ and $\xLp$ differ only by framing, then their shifted Khovanov homologies are isomorphic: $\tKHm(\xLpcN) = \tKHm(\xLcN)$.
\begin{theorem}
\label{thm:bndd}
There are bounds on degrees of shifted homology of a \tBadq\ link $\xL$:
$\tKHmvv{i}{j}(\xLcN) = 0$, if one of the following conditions is met:
\begin{align}
i&<0
\label{eq:bd1}
\\
\label{eq:bd2}
j&< -\shlf i -\shlf \ncrmL
\\
\label{eq:bd3}
j& < -i
\\
\label{eq:bd4}
j& =-i \neq 0,
\end{align}
where $\ncrmL$ is the minimum crossing number of a diagram representing $\xL$.
Moreover,
\begin{equation}
\label{eq:endi}
\dim\tKHmvv{0}{0} = 1.
\end{equation}
\end{theorem}

\subsubsection{Tail of Khovanov homology}
The main result of this paper is a definition of special degree-preserving maps between shifted Khovanov homologies of \uclrd\ \tBadq\ links, such that these maps are isomorphisms at low \thdgr s.
\begin{theorem}
\label{thm:llviso}
A \tBadq\ diagram of a link $\xL$ determines a sequence of \tdgpr\ maps
\begin{equation}
\label{eq:splmps}
\tKHm(\xLcN)\xrightarrow{\mnfN}
\tKHm(\xLcNo)
\end{equation}
which are isomorphisms on $\tKHmvv{i}{\hem}(\xLcN)$  for $i\leq \xca-1$.
\end{theorem}
Let $\tlH(\xL)$ be the \drlim\ of the \drsys\ determined by the sequence of maps $\mnfN$, $\xca\in\ZZ_+$ : 
\begin{equation}
\label{eq:drlm}
\tlH(\xL) = \dlm \tKHm(\xLcN).
\end{equation}

%
\begin{corollary}
\label{cor:mpis}
The maps
\begin{equation}
\label{eq:shti}
\tKHm(\xLcN)\rightarrow\tlH(\xL)
\end{equation}
 associated with the \drlim\rx{eq:drlm} are isomorphisms on $\tKHmvv{i}{\hem}(\xLcN)$ for $i\leq\xca-1$, hence
the \drlim\ $\tlH(\xL)$ is finite-dimensional in every bi-degree, it satisfies the bounds $\tlHvv{i}{j}(\xL)=0$ at the conditions\rx{eq:bd1}--\rxw{eq:bd4} and
\begin{equation}
\label{eq:hdtl}
\dim\tlHvv{0}{0}(\xLcN) = 1.
\end{equation}
\end{corollary}
%
\begin{remark}
To define the maps\rx{eq:splmps} we have to choose a diagram $\xD$ representing $\xL$, however we expect that the maps can be defined canonically, that is, independently of that choice.
\end{remark}

\subsubsection{Relation to the tail of the Jones polynomial}
The bounds\rx{eq:bd1} and\rx{eq:bd2} on $\tlH(\xL)$ mean that the graded Euler characteristic of the tail homology is well-defined, because in its presentation as an alternating sum of homology dimensions
\[
\pJqLi = \sum_{i,j} (-1)^j\,q^{i+j}\,\dim\tlHvv{i}{j}(\xD)
\]
there is only a finite number of non-trivial terms for any given value of $i+j$.

The bound\rx{eq:bd2} indicates that $\tKHmvv{i}{j}(\xLcN)$ and $\tlHvv{i}{j}(\xL)$ are trivial when $i+j< \hlf i - \hlf\ncrmL$, hence their high \thdgr s  contribute only to coefficients at high powers of $q$ in the graded Euler characteristic. Since the map\rx{eq:shti} is an isomorphism at low \thdgr s, we come to the following:
\begin{theorem}
The graded Euler characteristic of the tail homology determines the lower powers of $q$ in the \uclrd\ Jones polynomial of a \tBadq\ link:
\[
\pJqaL=
(-1)^{\gvL\xca}
 q^{-\hlf\clN^2\ncrL-\gvL\xca} \left(\pJqLi + O\left(q^{\hlf\xca - \hlf\ncrmL}\right)\right).
\]
\end{theorem}
This means that the tail homology categorifies the tail of the \uclrd\ Jones polynomial of \tBadq\ links studied by C.~Armond~\cite{Arm11} and by S.~Garoufalidis and T.~Le~\cite{GL11}.

\subsubsection{Tail homology is determined by a \trBdg\ of a link}


A link diagram $\xDp$ is a \tBrdc\ of a link diagram $\xD$, if the diagram $\sBDp$ is constructed from $\sBD$ in two stages: at first stage for each pair of distinct \tBcr s of $\sBD$ connected by more that one \tstrt\ we remove all connecting \tstrt s but one; at the second stage we remove all \tBcr s which have only one \tstrt\ attached to them. Obviously, if $\xD$ is \tBadq, then so is $\xDp$.

A link $\xLp$ is a \tBrdc\ of a \tBadq\ link $\xL$, if $\xLp$ can be presented by a diagram which is a \tBrdc\ of a \tBadq\ diagram presenting $\xL$.

\begin{theorem}
\label{thm:rddg}
If a link $\xLp$ is a \tBrdc\ of a \tBadq\ link $\xL$, then their tail homologies are isomorphic: $\tlH(\xLp)\cong\tlH(\xL)$.
\end{theorem}
C.~Armond and O.~Dasbach~\cite{AD12} are working on a similar result for the tail of the \uclrd\ Jones polynomial.
\begin{corollary}
If $\xLb$ is a circular closure of a connected negative braid $\betbr$, then the tail homology of $\xLb$ is that of an unknot: $\tlH(\xLb)\cong\tlH(\bigcirc )$.
\end{corollary}
\begin{proof}
It is easy to see that $\xLb$ is \tBadq\ and a \trBdg\ of $\xLb$ consists of a single circle without \tstrt s.
\end{proof}
\begin{corollary}[Invariance under strut doubling]
If $\xL$ is a \tBadq\ link and $\xLp$ is constructed by performing a replacement
\begin{equation}
\label{eq:crdcr}
\def\mxsh{2}
\def\mysh{1.5}
\begin{tikzpicture}[menvfour]
\node (l) at (-3.2,0) {};
\node (r) at (3.2,0) {};
\path[commutative diagrams/.cd, every arrow, every label]
(l) edge[commutative diagrams/squiggly] (r);
\begin{scope}[xshift=-6cm]
\draw [thnln] (-\mxsh+0.15,\mysh) to [out=0,in=180] (\mxsh-0.15,-\mysh);
\draw [lnovr] (-\mxsh+0.15,-\mysh) to [out=0,in=180] (\mxsh-0.15,\mysh);
\draw [thnln] (-\mxsh+0.15,-\mysh) to [out=0,in=180] (\mxsh-0.15,\mysh);
\end{scope}
\begin{scope}[xshift=6cm]
\draw [thnln] (-\mxsh+0.15,\mysh) to [out=0,in=180] (\mxsh-0.15,-\mysh);
\draw [lnovr] (-\mxsh+0.15,-\mysh) to [out=0,in=180] (\mxsh-0.15,\mysh);
\draw [thnln] (-\mxsh+0.15,-\mysh) to [out=0,in=180] (\mxsh-0.15,\mysh);
\begin{scope}[xshift=4cm-0.3cm]
\draw [thnln] (-\mxsh+0.15,\mysh) to [out=0,in=180] (\mxsh-0.15,-\mysh);
\draw [lnovr] (-\mxsh+0.15,-\mysh) to [out=0,in=180] (\mxsh-0.15,\mysh);
\draw [thnln] (-\mxsh+0.15,-\mysh) to [out=0,in=180] (\mxsh-0.15,\mysh);
\end{scope}
\end{scope}
\end{tikzpicture}
\end{equation}
in a \tBadq\ diagram of $\xL$, then tail homologies of $\xL$ and $\xLp$ are isomorphic: $\tlH(\xL)\cong\tlH(\xLp)$.
\end{corollary}
\begin{proof}
Obviously, $\xL$ and $\xLp$ have the same \tBrdc.
\end{proof}
%
The latter corollary is a categorification of a similar property of the tail of the \uclrd\ Jones polynomial observed by C.~Armond and O.~Dasbach~\cites{AD11,Arm11}. This property suggests that a single crossing plays the role of a \cJWp\ in the tail homology. In order to make this statement precise, tail homology has to be defined for knotted graphs, which may include both finite and infinite colors, so that the essential property of contracting cups/caps can be formulated. We hope to address this issue in a subsequent paper. Meanwhile, we prove in Appendix that for large $\xca$ a crossing of two $\xca$-cables is, indeed, homologically close to the \tJWp\ placed on two parallel $\xca$-cables.

\subsection{Technicalities}

We prove most statements of the previous subsection not just for \tBadq\ links, but for any link diagram. However, we conjecture that the results are trivial for \tBiadq\ diagrams, because their tail homology is trivial, if defined as a \drlim\ of the system that we construct. We expect that \tBiadq\ links also have tail homology, but the proof that the tail of their Khovanov homology stabilizes in the limit of large color requires new ideas.

\subsubsection{Shifted Khovanov homology}
Let $\xD$ be a diagram of a tangle which may include single lines, cables and \tJWp s. We define $\yncrv{\xD}$ to be the total number of single line crossings in $\xD$  (that is, a crossing between an $a$-cable and a $b$-cable contributes $ab$ to $\yncrv{\xD}$). The following is an obvious corollary of \ex{eq:dkhbr} and the fact that, according to\rx{eq:projcn} and\rx{eq:grproj},
$\hmord{
\begin{tikzpicture}[scale=0.4,baseline=-2.5]
\draw[line width=\ljwp] (-0.15,-0.6) rectangle ++(0.3,1.2);
\draw [line width=\cblth] (-.65,0) -- (-0.15,0) 
(0.15,0) -- (.65,0);
\end{tikzpicture}
}=0$:
\begin{theorem}
\label{thm:smfr}
The complex $\xKhv{\xD}$ has a lower homological bound:
$\hmord{\xKhv{\xD}}\geq - \hlf\yncrv{\xD}.$
\end{theorem}

Let $\xD$ be a diagram of a link which may include single lines, cables and \tJWp s. Define $\ynccv{\xD}$ to be the total number of circles in the diagram constructed from $\xD$ by replacing the \tJWp s with identity braids and performing \Bsplng s on all crossings. Now we define the shifted Khovanov homology of $\xD$:
\[
\tKHm(\xD) = \shcr^{\hlf\yncrv{\xD}} \shfr^{\ynccv{\xD}} \KHm(\xD).
\]
The following is a particular case of Theorem\rw{thm:smfr}:
\begin{theorem}
\label{thm:eaest}
If $i<0$, then $\tKHmvv{i}{\hem}(\xD)=0$.
\end{theorem}

For a link diagram $\xD$ let $\xDclN$ denote the corresponding \uclrd\ diagram (that is, every link components is $\xca$-cabled and contains at least one \tJWp).
Then, obviously, $\yncrv{\xDclN} = \clN^2\ncrD$ and
$\ynccv{\xDclN} =\xca \gvD$, so
\begin{equation}
\label{eq:shdgd}
\tKHm(\xDclN) = \shcr^{\hlf\clN^2\ncrD} \shfr^{\gvD\xca} \KHm(\xDclN).
\end{equation}
Hence, if $\xD$ is \tBadq\ and represents a link $\xL$, then $\tKHm(\xDclN)$ coincides with the shifted homology defined by \ex{eq:shdgt}.

Now Theorem\rw{thm:bndd} is a corollary of Theorem\rw{thm:eaest} and the following:
\begin{theorem}
\label{thm:kqbnd}
A shifted homology of a \uclrd\ diagram $\xDclN$ has a bound:
$\tKHmvv{i}{j}(\xDclN)=0$ if one of the following conditions is satisfied:
\begin{align}
\label{eq:bd2a}
j&<-\shlf i - \shlf\ncrD - \sthlf\ncriD
\\
\label{eq:bd3a}
j& < -i - \ncriD
\end{align}
Moreover, if $\xD$ is \tBadq, then
\begin{equation}
\label{eq:cndbd}
\dim \tKHmvv{i}{-i} =
\begin{cases}
0,&\text{if $i>0$,}
\\
1,
&\text{if $i=0$.}
\end{cases}
\end{equation}
\end{theorem}
Theorem\rw{thm:llviso} is a special case of the following:
\begin{theorem}
\label{thm:leviso}
For any link diagram $\xD$ there is a sequence of \tdgpr\ maps
\begin{equation}
\label{eq:spmaps}
\tKHm(\xDclN)\xrightarrow{\mnfN}
\tKHm(\xDclNo),
\end{equation}
which are isomorphisms on $\tKHmvv{i}{\hem}$ for $i\leq \xca-1$.
\end{theorem}
This theorem implies that the \drsys\ formed by maps\rx{eq:spmaps} has a limit
\begin{equation}
\label{eq:dglim}
\tlH(\xD) = \dlm \tKHm(\xDclN).
\end{equation}
which is finite-dimensional in every bi-degree.
\begin{conjecture}
If the diagram $\xD$ is not \tBadq, then the \drlim\rx{eq:dglim} is trivial: $\tlH(\xD)=0$.
\end{conjecture}

\subsection{Discussion}
We conjecture that \tBadq\ links have a tri-graded homology $\ttlH(\xL)$, which has an additional `\tbgrd', whose zero-degree part coincides with the tail homology $\tlH(\xL)$:
\[
\ttlH(\xL) = \bigoplus_{\substack{i,j\\k\geq 0}}\ttlHvvv{i}{j}{k}(\xL),\qquad
\ttlHvvv{i}{j}{0}(\xL) = \tlHvv{i}{j}(\xL).
\]
This homology should have a family of mutually anti-commuting differentials $\xfdN$, $\xca\in\ZZ_+$, with degrees $\dgq\xfdN=-1$, $\dgh\xfdN=1$ and $\dgb=-1$
such that homology of $\ttlH(\xL)$ with respect to $\xfdN$ matches the shifted Khovanov homology $\tKHm$ up to a level proportional to $\xca^2$, after the \tbdgr\ is converted into \thdgr:
\begin{equation}
\label{eq:trconj}
\tKHmvv{\xti}{\hem}(\xL) = \bigoplus_{i+\xca k = \xti}\Hm_{i,\hem,k}^{\xfdN}\big(
\ttlH(\xL)\big),\quad\text{if $\xti\geq a_{\xL}\xca^2$},
\end{equation}
where $a_{\xL}$ is a constant determined by $\xL$.

There are three reasons to formulate this conjecture. The first reason is that the proof of Theorem\rw{thm:llviso} is based on numerous long exact sequences\rx{eq:les}, in which the `correction homology' starts at homological degree proportional to $\xca$.

The second reason is that, according to Garoufalidis and Le~\cite{GL11}, the tail of the Jones polynomial of \tBadq\ links has a `telescopic' structure. They show that if $\xL$ is alternating, then there exists a family of Laurent series $\Phi_n(q) = \sum_m a_m q^m$ such that for any $k>0$ the combined series $F_k(q) = \sum_{n=0}^k \Phi_n(q)$ is a better approximation for the tail of the colored Jones polynomial that just the first term $\Phi_0(q)= \pJqLi $. With the help of the colored Kauffman bracket (\cf\rx{eq:colKhbr}) we can prove a similar result for all \tBadq\ links and we expect that the 2-variable series
\begin{equation}
\label{eq:dsrs}
\ptJqLi = \sum_{n=0}^\infty b^n \Phi_n(q)
\end{equation}
is the bi-graded Euler characteristic of the tri-graded homology $\ttlH(\xL)$.

The third reason for our conjecture comes from the paper by Gukov and \Stosic~\cite{GS11}. Based on QFT models of Khovanov homology, they suggest that its dependence on color should be similar to the dependence of the $\SUv{n}$ homology on $n$: this homology  may be presented as a homology of a special differential acting on $\SUv{N}$ homology if $N>n$. We suggest to go half step further. Ultimately, the $\SUv{N}$ homology may be presented, at least, conjecturally, from the tri-graded HOMFLY-PT homology with the help of special differentials $\xfdv{n}$ and we expect that a similar process may work for the tail homology.

We expect that the formation of a stable tail of a \uclrd\ \tBadq\ link is a general feature which originates in the tri-graded homology when the Young diagram describing the color has a very large value of one of the differences between the lengths of rows or columns. In particular, it could be easy to follow the tail formation in case when the diagram consists of a single very large column.

Witten suggested~\cite{Wi11} that a series of the form\rx{eq:dsrs} should represent the graded Euler characteristic of Khovanov homology in the background of a flat $\mathrm{U}(1)$-reducible $\SUt$ connection in the link complement. We conjecture that if a link can be presented as a circular closure of a totally negative braid, then the tail homology coincides with the one related to the flat $\mathrm{U}(1)$-reducible $\SUt$ connection.

\subsection{Acknowledgements}
The author thanks Eugene Gorsky for sharing the results of his unfinished research and, in particular, the conjecture about the structure of Khovanov homology of a colored unknot and its stabilization in the limit of high color.

This work was supported in part by the NSF grants DMS-0808974 and DMS-1108727.

\section{Five tools}

The proof of Theorems\rw{thm:kqbnd} and\rw{thm:leviso} requires five tools: \ltrf s, purging, braid straightening, colored \tKhbr\ and recurrence relations between \cJWp s. \Ltrf s relate homologies of similar diagrams. Purging gets rid of redundant \taTLt s in complexes, which contain \tJWp s, thus improving estimates of homological order. Straightening a braid is a simple observation that braiding within a cable attached to a projector results only in degree shifts. The colored \tKhbr\ is a special presentation of a crossing of two cables attached to \tJWp s. Finally, recurrence relations relate \tJWp s on $\xca$ and $\xca+1$ strands.

\subsection{\Lrpl s and \ltrf s}
A \emph{\lrpl} is a pair of tangles $\ytngsi$ and $\ytngsf$, which may contain single and cabled lines, as well as \tJWp s. Both tangles should have the same sets of incoming legs and the same sets of outgoing legs. Hence if an initial diagram $\xDsi$ contains the tangle $\ytngsi$ attached by its legs to the rest of the diagram, then we can construct a final diagram $\xDsf$ by replacing $\ytngsi$ with $\ytngsf$. If $\ytngsi$ or $\ytngsf$ is not an actual diagram, but rather a complex of diagrams within the universal category, the \lrpl\ still makes sense as a construction of $\xKhv{\xDsf}$ from $\xKhv{\xDsi}$.

A \emph{\ltrf} is a \lrpl\ together with a specified \tdgpr\ morphism $\ytngkfp\xrightarrow{\xmg}\ytngki$, where we use a shortcut $\ytngkfp= \shfr^{\xhshf}\shcr^{\xqshf}\ytngkf$\footnote{The unnatural direction of morphism is chosen for future convenience.} The morphism $\xmg$ determines the presentation of $\ytngki$ as a cone
\begin{equation}
\label{eq:cnrel}
\ytngki \hteqv
\boxed{
\ytngkc
\rightarrow\ytngkfp
}\;,
\qquad\text{where}\quad
\ytngkc = \boxed{
\shcr\ytngkfp\xrightarrow{\xmg}\ytngki
}\;.
\end{equation}
Up to a degree shift, the `correction' complex $\ytngkc$ may be the categorification complex of an actual tangle $\ytngsc$, or it may be just a convenient shortcut.

Let $\xDsi$ be a diagram of a link which contains $\ytngsi$ and let $\xDsf$ and $\xDsc$ be the diagrams constructed by replacing $\ytngsi$ with $\ytngsf$ and $\ytngsc$. The relations\rx{eq:cnrel} imply a long exact sequence
\begin{equation}
\label{eq:les}
\shcr^{-1}\KHm(\xDsc) \longrightarrow  \shfr^{\xhshf}\shcr^{\xqshf}\KHm(\xDsf) \xrightarrow{\;\;\xmg\;\;} \KHm(\xDsi)
\longrightarrow \KHm(\xDsc).
\end{equation}
For all \ltrf s considered in this paper, there are relations
\[
\xqshf = \shlf\yncrv{\xDsf} - \shlf\yncrv{\xDsi},\qquad
\xhshf = \ynccv{\xDsf} - \ynccv{\xDsi},
\]
hence \ex{eq:les} turns into the following sequence of \tdgpr\ maps between shifted homologies:
\[
\shcr^{-1} \KHm(\xDscp) \longrightarrow
\tKHm(\xDsf) \xrightarrow{\;\;\xmg\;\;} \tKHm(\xDsi)
\longrightarrow \KHm(\xDscp),
\]
where $
\KHm(\xDscp) = \shcr^{\shlf\yncrv{\xDsi}}\shfr^{\ynccv{\xDsi}}\KHm(\xDsc)$. This exact sequence implies the following:
\begin{proposition}
\label{prp:gestdg}
If $\KHmib(\xDsc)=0$ for $i \leq \hbnd-1$, then
the \tdgpr\ map
\begin{equation}
\label{eq:dgprmg}
\tKHm(\xDsf) \xrightarrow{\;\xmg\;} \tKHm(\xDsi),
\end{equation}
is an isomorphism on $\tKHmvv{i}{\hem}$ for
\begin{equation}
\label{eq:ineqo}
i\leq\hbnd+
\shlf\yncrv{\xDsi}-2.
\end{equation}
\end{proposition}

%
%
%
%

\subsection{Purging}

Purging is a process of using \ex{eq:ccann} to remove constituent \taTLt s of a complex, whose cups or caps are connected directly to a \tJWp. Let $\sTLab$ be the set of all (a,b) \taTLt s and let
$\sTLabc =\{ \gamma\in\sTLabc\colon \swdv{\gamma} = b\}$. In other words, $\sTLabc$ is a subset of $(a,b)$ \taTLt s, which contain no cups, but only caps and \txstrs s.
\begin{proposition}
\label{prp:purge}
There is a homotopy equivalence
\begin{multline}
\label{eq:bprg}
\boxed{
\cdots\longrightarrow\shcr^i\bigoplus_{\substack{j \\ \gamma\in\sTLab}} m_{ij,\gamma}\,\shfr^j\,
\begin{tikzpicture}[menvtwo]
\draw [ptzer] (-0.6,-0.6) rectangle (0.6,0.6); \node at (0,0) {$\gamma$};
\draw [thkln] (-1.2,0) -- (-0.6,0) node [near start,below] {$\scriptstyle a$};
\draw [thkln] (1.2,0) -- (0.6,0) node [near start,below] {$\scriptstyle b$};
\end{tikzpicture}
\longrightarrow\cdots
}
\begin{tikzpicture}[menvtwo]
\path [use as bounding box] (-1,-1) rectangle (2,1);
\draw [ptzer] (-0.15,-0.6) rectangle (.15,.6);
\draw [thkln] (-0.75-0.22,0) -- (-0.15,0)
node [midway,below] {$\scriptstyle b$}
(0.15,0) -- (.75,0)
node [near end,below] {$\scriptstyle b$};
\end{tikzpicture}
\\
\hteqv\;
\boxed{
\cdots\longrightarrow\shcr^i\bigoplus_{\substack{j \\ \gamma\in\sTLabc}} m_{ij,\gamma}\,\shfr^j\,
\begin{tikzpicture}[menvtwo]
\draw [ptzer] (-0.6,-0.6) rectangle (0.6,0.6); \node at (0,0) {$\gamma$};
\draw [thkln] (-1.2,0) -- (-0.6,0) node [near start,below] {$\scriptstyle a$};
\draw [thkln] (1.8,0) -- (0.6,0) node [midway,below] {$\scriptstyle b$}
(2.1,0) -- (2.7,0) node [near end,below] {$\scriptstyle b$};
\draw[ptzer] (1.8,-0.6) rectangle (2.1,0.6);
\end{tikzpicture}
\longrightarrow\cdots
}
\end{multline}
\end{proposition}


\subsection{Straightening a braid attached to a \tJWp}

\begin{theorem}
The categorification complex of the tangle composition of the \tJWp\ with a braid $\brbet$ is homotopy equivalent to the shifted complex of the \tJWp:
\[
\begin{tikzpicture}[menvtwo]
\draw [ptzer] (-0.6,-0.6) rectangle (0.6,0.6); \node at (0,0) {$\brbet$};
\draw [thkln] (-2.7,0) -- (-2.1,0) node [near start,below] {$\scriptstyle a$} (-1.8,0) -- (-0.6,0)
 node [midway, below] {$\scriptstyle a$} (0.6,0) -- (1.2,0) node [near end, below] {$\scriptstyle a$};
\draw [ptzer] (-2.1,-0.6) rectangle ++(0.3,1.2);
\end{tikzpicture}
\;\hteqv\;
\shcr^{\hlf(\xnp-\xnm)}
\begin{tikzpicture}[menvtwo]
\draw [ptzer] (-0.15,-0.6) rectangle ++(0.3,1.2);
\draw [thkln] (-0.75,0) -- (-0.15,0) node [near start,below] {$\scriptstyle a$}
(0.15,0) -- (0.75,0);
\end{tikzpicture}
\]
where $\xnp$ ($\xnm$) is the number of positive (negative) elementary crossings in a presentation of $\brbet$.
\end{theorem}

\begin{proof}
It is sufficient to prove this equivalence in the case of $\xnp=1$ and $\xnm=0$ (the case of $\xnp=0$ and $\xnm=1$ is similar and other cases can be proved by consequent composition of elementary crossings).
Thus we replace the positive crossing by its Khovanov complex
\[
\begin{tikzpicture}[xscale=0.6,yscale=-0.6,baseline=-2.5]
\draw [thkln] (-0.5,0) -- (0.1,0);
\draw (0.5,0.3) to [out=0,in=180] (2,-0.3);
\draw [line width=6pt, draw=white] (0.5,-0.3) to [out=0,in=180] (2,0.3);
\draw (0.5,-0.3) to [out=0,in=180] (2,0.3);
\draw [thkln] (0.5,0.6) -- (2,0.6);
\draw [thkln] (0.5,-0.6) -- (2,-0.6);
\draw [ptzer] (0.1,-0.8) rectangle (0.5,0.8);
\end{tikzpicture}
\; = \;
\Cnv{
\shcrh\;
\begin{tikzpicture}[scale=0.6,baseline=-2.5]
\draw [thkln] (-0.5,0) -- (0.1,0);
\draw (0.5,0.3) -- (0.75,0.3) to [out=0,in=180] (1.25,0.15) to [out=0,in=180] (1.75,0.3) -- (1.75,0.3) -- (2,0.3);
\draw (0.5,-0.3) -- (0.75,-0.3) to [out=0,in=180] (1.25,-0.15) to [out=0,in=180] (1.75,-0.3) -- (1.75,-0.3) -- (2,-0.3);
\draw [thkln] (0.5,0.6) -- (2,0.6);
\draw [thkln] (0.5,-0.6) -- (2,-0.6);
\draw [ptzer] (0.1,-0.8) rectangle (0.5,0.8);
\end{tikzpicture}
\xrightarrow{\hspace{0.8cm}}
\shcrmh\;
\begin{tikzpicture}[scale=0.6,baseline=-2.5]
\draw [thkln] (-0.5,0) -- (0.1,0);
\draw (0.5,0.3) -- (0.75,0.3) .. controls +(0:0.4) and +(0:0.4) .. (0.75,-0.3) -- (0.5,-0.3);
\draw (2,0.3) -- (1.75,0.3) .. controls +(180:0.4) and +(180:0.4) .. (1.75,-0.3) -- (2,-0.3);
\draw [thkln] (0.5,0.6) -- (2,0.6);
\draw [thkln] (0.5,-0.6) -- (2,-0.6);
\draw [ptzer] (0.1,-0.8) rectangle (0.5,0.8);
\end{tikzpicture}
}
\]
and observe that the second term in the resulting cone is contractible.
\end{proof}

This proposition has important special cases:
\begin{equation}
\label{eq:projtw}
\begin{tikzpicture}[scale=0.75,baseline=-3]
\draw [thkln] (-0.55,0)
-- (0.2,0)
node [near start,above] {$\scriptstyle a+1$};
\draw  [thkln] (0.5,0.3) to [out=0,in=180] (1.25,-0.3)
node[right] {$\scriptstyle a$};
\draw [lnovr] (0.5,-0.3) to [out=0,in=180] (1.25,0.3);
\draw (0.5,-0.3) to [out=0,in=180] (1.25,0.3);
\draw [ptzer] (0.2,-0.6) rectangle (0.5,0.6);
\end{tikzpicture}
\;\hteqv\;
\shcr^{\frac{a}{2}}
\begin{tikzpicture}[scale=0.75,baseline=-3]
\draw [line width=\cblth] (-0.55,0)
-- (0.2,0)
node [midway,above] {$\scriptstyle a+1$\;};
\draw [line width=\cblth] (0.5,0) -- (1.25,0);
\draw [line width=1pt] (0.2,-0.6) rectangle (0.5,0.6);
\end{tikzpicture}
,\qquad
\begin{tikzpicture}[scale=0.75,baseline=-3]
\draw [line width=\cblth] (-0.55,0)
-- (0.2,0)
node [near start,above] {$\scriptstyle a+1$};
\begin{scope}[yscale=-1]
\draw  [line width=\cblth] (0.5,0.3) to [out=0,in=180] (1.25,-0.3)
node [right] {$\scriptstyle a$};
\draw [line width=6pt, draw=white] (0.5,-0.3) to [out=0,in=180] (1.25,0.3);
\draw (0.5,-0.3) to [out=0,in=180] (1.25,0.3);
\end{scope}
\draw [line width=1pt] (0.2,-0.6) rectangle (0.5,0.6);
\end{tikzpicture}
\;\hteqv\;
\shcr^{-\frac{a}{2}}
\begin{tikzpicture}[scale=0.75,baseline=-3]
\draw [line width=\cblth] (-0.55,0)
-- (0.2,0)
node [midway,above] {$\scriptstyle a+1$\;};
\draw [line width=\cblth] (0.5,0) -- (1.25,0);
\draw [line width=1pt] (0.2,-0.6) rectangle (0.5,0.6);
\end{tikzpicture}
\end{equation}
and two similar cases with opposite powers of $\shcr$ when the cable runs over the single line.


\subsection{Colored \tKhbr}
\begin{theorem}[Single strand splicing]
\label{lm:sss}
A \tKhbr\ of the crossing of two equally colored strands can be presented as a cone:
\begin{equation}
\label{eq:colssKhbr}
\begin{tikzpicture}[scale=0.5,rotate=-45,baseline=-3.5]
\draw [thkln] (-1.5,0) -- (1.5,0);
\draw [lnovr] (0,-1.5) -- (0,1.5);
\draw [thkln] (0,-1.5) -- (0,1.5);
%
%
\draw [thkln] (-2.2,0)  -- (-1.8,0);
\node at (-2.6,0) {$\scriptstyle a+1$};
\draw [thkln] (1.8,0) -- (2.2,0);
\node at (2.6,0) {$\scriptstyle a+1$};
\draw [thkln] (0,1.8) -- (0,2.2);
\node at (0,2.6) {$\scriptstyle a+1$};
\draw [thkln] (0,-1.8) -- (0,-2.2);
\node at (0,-2.6) {$\scriptstyle a+1$};
\draw [thkln] (-1.8,-0.6) rectangle ++(0.3,1.2);
\draw [thkln] (1.5,-0.6) rectangle ++(0.3,1.2);
\draw [thkln] (-0.6,1.5) rectangle ++(1.2,0.3);
\draw [thkln] (-0.6,-1.5) rectangle ++(1.2,-0.3);
\end{tikzpicture}
\hteqv
\;\;
\Cnv{
\shcr^{a+\hlf}
\begin{tikzpicture}[scale=0.5,rotate=-45,baseline=-3.5]
\draw  (-1.5,-0.3) to [out=0, in=90] (-0.3,-1.5);
\draw  (1.5,0.3) to [out=180, in=-90] (0.3, 1.5);
\draw [thkln] (-1.5,0) -- (1.5,0);
\draw [line width=\ocrw,draw=white] (0,-1.5) -- (0,1.5);
\draw [thkln] (0,-1.5) -- (0,1.5);
\draw [thkln] (-2.2,0)  -- (-1.8,0);
\node at (-2.6,0) {$\scriptstyle a+1$};
\draw [thkln] (1.8,0) -- (2.2,0);
\node at (2.6,0) {$\scriptstyle a+1$};
\draw [thkln] (0,1.8) -- (0,2.2);
\node at (0,2.6) {$\scriptstyle a+1$};
\draw [thkln] (0,-1.8) -- (0,-2.2);
\node at (0,-2.6) {$\scriptstyle a+1$};
\draw [ptzer] (-1.8,-0.6) rectangle ++(0.3,1.2);
\draw [ptzer] (1.5,-0.6) rectangle ++(0.3,1.2);
\draw [ptzer] (-0.6,1.5) rectangle ++(1.2,0.3);
\draw [ptzer] (-0.6,-1.5) rectangle ++(1.2,-0.3);
\end{tikzpicture}
\longrightarrow
\shcr^{-(a+\hlf)}
\begin{tikzpicture}[scale=0.5,rotate=-45,baseline=-3.5]
\draw (-1.5,0.3) to [out=0, in=-90] (-0.3,1.5);
\draw (1.5,-0.3) to [out=180, in=90] (0.3,-1.5);
\draw [thkln] (-1.5,0) -- (1.5,0);
\draw [line width=\ocrw,draw=white] (0,-1.5) -- (0,1.5);
\draw [thkln] (0,-1.5) -- (0,1.5);
\draw [thkln] (-2.2,0)  -- (-1.8,0);
\node at (-2.6,0) {$\scriptstyle a+1$};
\draw [thkln] (1.8,0) -- (2.2,0);
\node at (2.6,0) {$\scriptstyle a+1$};
\draw [thkln] (0,1.8) -- (0,2.2);
\node at (0,2.6) {$\scriptstyle a+1$};
\draw [thkln] (0,-1.8) -- (0,-2.2);
\node at (0,-2.6) {$\scriptstyle a+1$};
\draw [ptzer] (-1.8,-0.6) rectangle ++(0.3,1.2);
\draw [ptzer] (1.5,-0.6) rectangle ++(0.3,1.2);
\draw [ptzer] (-0.6,1.5) rectangle ++(1.2,0.3);
\draw [ptzer] (-0.6,-1.5) rectangle ++(1.2,-0.3);
\end{tikzpicture}
}
\end{equation}
\end{theorem}
\begin{proof}
Split off a single strand from each crossing cable, apply the \tKhbr\ relation\rx{eq:dkhbr} to their crossing:
\[
\begin{tikzpicture}[scale=0.5,rotate=-45,baseline=-3.5]
\draw [thkln] (-1.5,0.4) -- (1.5,0.4);
\draw (-1.5,-0.4) -- (1.5,-0.4);
\draw [line width=\ocrw,draw=white] (-0.4,-1.5) -- (-0.4,1.5);
\draw [thkln] (-0.4,-1.5) -- (-0.4,1.5);
\draw [line width=\ocrw,draw=white] (0.4,-1.5) -- (0.4,1.5);
\draw  (0.4,-1.5) -- (0.4,1.5);
%
%
\draw [thkln] (-2.2,0)  -- (-1.8,0);
\node at (-2.6,0) {$\scriptstyle a+1$};
\draw [thkln] (1.8,0) -- (2.2,0);
\node at (2.6,0) {$\scriptstyle a+1$};
\draw [thkln] (0,1.8) -- (0,2.2);
\node at (0,2.6) {$\scriptstyle a+1$};
\draw [thkln] (0,-1.8) -- (0,-2.2);
\node at (0,-2.6) {$\scriptstyle a+1$};
\draw [ptzer] (-1.8,-0.6) rectangle ++(0.3,1.2);
\draw [ptzer] (1.5,-0.6) rectangle ++(0.3,1.2);
\draw [ptzer] (-0.6,1.5) rectangle ++(1.2,0.3);
\draw [ptzer] (-0.6,-1.5) rectangle ++(1.2,-0.3);
\end{tikzpicture}
\hteqv\;\;
\Cnv{
\shcrh
\begin{tikzpicture}[scale=0.5,rotate=-45,baseline=-3.5]
\draw [thkln] (-1.5,0.4) -- (1.5,0.4);
\draw (-1.5,-0.4) -- (-0.7,-0.4);
\draw (-0.7,-0.4) to [out=0,in=90] (0.4,-1.5);
\draw [line width=\ocrw,draw=white] (-0.4,-1.5) -- (-0.4,1.5);
\draw [thkln] (-0.4,-1.5) -- (-0.4,1.5);
\draw  (0.4,0.7) -- (0.4,1.5);
\draw [line width=\ocrw,draw=white] (1.5,-0.4) to [out=180,in=-90] (0.4,0.7);
\draw (1.5,-0.4) to [out=180,in=-90] (0.4,0.7);
%
%
\draw [thkln] (-2.2,0)  -- (-1.8,0);
\node at (-2.6,0) {$\scriptstyle a+1$};
\draw [thkln] (1.8,0) -- (2.2,0);
\node at (2.6,0) {$\scriptstyle a+1$};
\draw [thkln] (0,1.8) -- (0,2.2);
\node at (0,2.6) {$\scriptstyle a+1$};
\draw [thkln] (0,-1.8) -- (0,-2.2);
\node at (0,-2.6) {$\scriptstyle a+1$};
\draw [ptzer] (-1.8,-0.6) rectangle ++(0.3,1.2);
\draw [ptzer] (1.5,-0.6) rectangle ++(0.3,1.2);
\draw [ptzer] (-0.6,1.5) rectangle ++(1.2,0.3);
\draw [ptzer] (-0.6,-1.5) rectangle ++(1.2,-0.3);
\end{tikzpicture}
\longrightarrow
\shcrmh
\begin{tikzpicture}[scale=0.5,rotate=-45,baseline=-3.5]
\draw [thkln] (-1.5,0.4) -- (1.5,0.4);
\draw (-1.5,-0.4) -- (-0.7,-0.4);
%
%
\draw  (0.4,0.7) -- (0.4,1.5);
%
%
\draw (0.4,-1.5) to [out=90,in=180] (1.5,-0.4);
\draw [line width=\ocrw,draw=white] (-0.7,-0.4) to [out=0,in=-90] (0.4,0.7);
\draw (-0.7,-0.4) to [out=0,in=-90] (0.4,0.7);
\draw [line width=\ocrw,draw=white] (-0.4,-1.5) -- (-0.4,1.5);
\draw [thkln] (-0.4,-1.5) -- (-0.4,1.5);
%
%
\draw [thkln] (-2.2,0)  -- (-1.8,0);
\node at (-2.6,0) {$\scriptstyle a+1$};
\draw [thkln] (1.8,0) -- (2.2,0);
\node at (2.6,0) {$\scriptstyle a+1$};
\draw [thkln] (0,1.8) -- (0,2.2);
\node at (0,2.6) {$\scriptstyle a+1$};
\draw [thkln] (0,-1.8) -- (0,-2.2);
\node at (0,-2.6) {$\scriptstyle a+1$};
\draw [ptzer] (-1.8,-0.6) rectangle ++(0.3,1.2);
\draw [ptzer] (1.5,-0.6) rectangle ++(0.3,1.2);
\draw [ptzer] (-0.6,1.5) rectangle ++(1.2,0.3);
\draw [ptzer] (-0.6,-1.5) rectangle ++(1.2,-0.3);
\end{tikzpicture}
}
\]
and then use the relations\rx{eq:projtw} to bring both tangles to the form of \ex{eq:colKhbr}.
\end{proof}
%

%
\begin{theorem}[Colored \tKhbr]
\label{thm:colkhovbr}
A \tKhbr\ of the crossing of two equally colored strands can be presented as a \tmcn\ of crossingless colored tangles:
\begin{equation}
\label{eq:lmtcn}
\begin{tikzpicture}[menvrtone]
\draw [thkln] (-1.5,0) -- (1.5,0);
\draw [line width=\ocrw,draw=white] (0,-1.5) -- (0,1.5);
\draw [thkln] (0,-1.5) -- (0,1.5);
%
%
\draw [thkln] (-2.2,0)  -- (-1.8,0);
\node at (-2.6,0) {$\scriptstyle a$};
\draw [thkln] (1.8,0) -- (2.2,0);
\node at (2.6,0) {$\scriptstyle a$};
\draw [thkln] (0,1.8) -- (0,2.2);
\node at (0,2.6) {$\scriptstyle a$};
\draw [thkln] (0,-1.8) -- (0,-2.2);
\node at (0,-2.6) {$\scriptstyle a$};
\draw [ptzer] (-1.8,-0.6) rectangle ++(0.3,1.2);
\draw [ptzer] (1.5,-0.6) rectangle ++(0.3,1.2);
\draw [ptzer] (-0.6,1.5) rectangle ++(1.2,0.3);
\draw [ptzer] (-0.6,-1.5) rectangle ++(1.2,-0.3);
\end{tikzpicture}
\;\hteqv\;
\shcr^{-\hlf a^2}
\;
\Cnv{
\cdots\longrightarrow
\shcr^{\ysvki + n_i}
\begin{tikzpicture}[menvrtsone]
\node at (135:0.55) {$\scriptstyle a-\yvki$};
\node at (135:-0.55) {$\scriptstyle a-\yvki$};
\node at (45:1.5) {$\scriptstyle \yvki$};
\node at (45:-1.5) {$\scriptstyle \yvki$};
\draw [thkln] (-1.5,0.3) to [out=0, in=-90] (-0.3,1.5);
\draw [thkln] (-1.5,-0.3) to [out=0, in=90] (-0.3,-1.5);
\draw [thkln] (1.5,0.3) to [out=180, in=-90] (0.3, 1.5);
\draw [thkln] (1.5,-0.3) to [out=180, in=90] (0.3,-1.5);
%
\draw [thkln] (-2.2,0)  -- (-1.8,0);
\node at (-2.6,0) {$\scriptstyle a$};
\draw [thkln] (1.8,0) -- (2.2,0);
\node at (2.6,0) {$\scriptstyle a$};
\draw [thkln] (0,1.8) -- (0,2.2);
\node at (0,2.6) {$\scriptstyle a$};
\draw [thkln] (0,-1.8) -- (0,-2.2);
\node at (0,-2.6) {$\scriptstyle a$};
\draw [ptzer] (-1.8,-0.6) rectangle ++(0.3,1.2);
\draw [ptzer] (1.5,-0.6) rectangle ++(0.3,1.2);
\draw [ptzer] (-0.6,1.5) rectangle ++(1.2,0.3);
\draw [ptzer] (-0.6,-1.5) rectangle ++(1.2,-0.3);
\end{tikzpicture}
\longrightarrow\cdots
}_{\;i=0}^{\;\infty}\;,
\end{equation}
such that
\begin{equation}
\label{eq:twineq}
n_i\geq 0,\qquad  i\geq 2^{\yki}-1
\end{equation}
and the \lumps\ form of this \tmcn\ is
\begin{equation}
\label{eq:colKhbr}
\begin{tikzpicture}[menvrtone]
\draw [thkln] (-1.5,0) -- (1.5,0);
\draw [line width=\ocrw,draw=white] (0,-1.5) -- (0,1.5);
\draw [thkln] (0,-1.5) -- (0,1.5);
%
%
\draw [thkln] (-2.2,0)  -- (-1.8,0);
\node at (-2.6,0) {$\scriptstyle a$};
\draw [thkln] (1.8,0) -- (2.2,0);
\node at (2.6,0) {$\scriptstyle a$};
\draw [thkln] (0,1.8) -- (0,2.2);
\node at (0,2.6) {$\scriptstyle a$};
\draw [thkln] (0,-1.8) -- (0,-2.2);
\node at (0,-2.6) {$\scriptstyle a$};
\draw [ptzer] (-1.8,-0.6) rectangle ++(0.3,1.2);
\draw [ptzer] (1.5,-0.6) rectangle ++(0.3,1.2);
\draw [ptzer] (-0.6,1.5) rectangle ++(1.2,0.3);
\draw [ptzer] (-0.6,-1.5) rectangle ++(1.2,-0.3);
\end{tikzpicture}
\;\hteqv\;
\shcr^{-\hlf a^2}\;
\Pcnv{
\bigoplus_{\xki=0}^a
\shcr^{\xki^2}
{a \brace \yki}_{\shcr}
\begin{tikzpicture}[menvrtone]
\node at (135:0.65) {$\scriptstyle a-\yki$};
\node at (135:-0.65) {$\scriptstyle a-\yki$};
\node at (45:1.5) {$\scriptstyle \yki$};
\node at (45:-1.5) {$\scriptstyle \yki$};
\draw [thkln] (-1.5,0.3) to [out=0, in=-90] (-0.3,1.5);
\draw [thkln] (-1.5,-0.3) to [out=0, in=90] (-0.3,-1.5);
\draw [thkln] (1.5,0.3) to [out=180, in=-90] (0.3, 1.5);
\draw [thkln] (1.5,-0.3) to [out=180, in=90] (0.3,-1.5);
%
\draw [thkln] (-2.2,0)  -- (-1.8,0);
\node at (-2.6,0) {$\scriptstyle a$};
\draw [thkln] (1.8,0) -- (2.2,0);
\node at (2.6,0) {$\scriptstyle a$};
\draw [thkln] (0,1.8) -- (0,2.2);
\node at (0,2.6) {$\scriptstyle a$};
\draw [thkln] (0,-1.8) -- (0,-2.2);
\node at (0,-2.6) {$\scriptstyle a$};
\draw [ptzer] (-1.8,-0.6) rectangle ++(0.3,1.2);
\draw [ptzer] (1.5,-0.6) rectangle ++(0.3,1.2);
\draw [ptzer] (-0.6,1.5) rectangle ++(1.2,0.3);
\draw [ptzer] (-0.6,-1.5) rectangle ++(1.2,-0.3);
\end{tikzpicture}
}
\end{equation}
\end{theorem}
\begin{proof}
We prove this theorem by induction over $a$. At $a=1$ it amounts to \tKhbr\rx{eq:dkhbr}. Suppose that it holds for some $a$ and consider the crossing of two $(a+1)$-cables. We split each cable into an $a$-cable and a single line and apply \eex{eq:lmtcn} and\rx{eq:colKhbr} to the crossing of $a$-cables:
\def\xscl{0.6}
\begin{equation}
\label{eq:spcrmt}
\begin{split}
\begin{tikzpicture}[scale=\xscl,rotate=-45,baseline=-3.5]
\draw [thkln] (-1.5,0.4) -- (1.5,0.4);
\draw (-1.5,-0.4) -- (1.5,-0.4);
\draw [line width=\ocrw,draw=white] (-0.4,-1.5) -- (-0.4,1.5);
\draw [thkln] (-0.4,-1.5) -- (-0.4,1.5);
\draw [line width=\ocrw,draw=white] (0.4,-1.5) -- (0.4,1.5);
\draw  (0.4,-1.5) -- (0.4,1.5);
\draw [thkln] (-2.2,0)  -- (-1.8,0);
\node at (-2.6,0) {$\scriptstyle a+1$};
\draw [thkln] (1.8,0) -- (2.2,0);
\node at (2.6,0) {$\scriptstyle a+1$};
\draw [thkln] (0,1.8) -- (0,2.2);
\node at (0,2.6) {$\scriptstyle a+1$};
\draw [thkln] (0,-1.8) -- (0,-2.2);
\node at (0,-2.6) {$\scriptstyle a+1$};
\draw [ptzer] (-1.8,-0.6) rectangle ++(0.3,1.2);
\draw [ptzer] (1.5,-0.6) rectangle ++(0.3,1.2);
\draw [ptzer] (-0.6,1.5) rectangle ++(1.2,0.3);
\draw [ptzer] (-0.6,-1.5) rectangle ++(1.2,-0.3);
\end{tikzpicture}
&\;\hteqv\;
\shcr^{-\hlf a^2}\;
\Pcnv{
\bigoplus_{\yki=0}^{a}
\shcr^{\yki^2}
{a \brace \yki}_{\shcr}
\begin{tikzpicture}[scale=\xscl,
rotate=-45,
baseline=-3.5]
\draw [thkln] (-1.5,0.4) to [out=0,in=-90] (-0.4,1.5);
\node at (-1,1) {$\scriptstyle a-\yki$};
\draw (-1.5,-0.4) -- (1.5,-0.4);
\draw [lnovrtw] (-1.5,0) to [out=0,in=90] (-0.4,-1.5);
\draw [thkln] (-1.5,0) to [out=0,in=90] (-0.4,-1.5) ;
\draw [lnovrtw] (0,-1.5) to [out=90,in=-135] (0,0);
\draw [thkln] (0,-1.5) to [out=90,in=-135] (0,0) to [out=45,in=180] (1.5,0);
\draw [thkln]  (1.5,0.4) to [out = 180, in=-90] (0,1.5);
\draw [lnovrtw] (0.4,-1.5) -- (0.4,1.5);
\draw  (0.4,-1.5) -- (0.4,1.5);
%
%
\draw [thkln] (-2.2,0)  -- (-1.8,0);
\node at (-2.6,0) {$\scriptstyle a+1$};
\draw [thkln] (1.8,0) -- (2.2,0);
\node at (2.6,0) {$\scriptstyle a+1$};
\draw [thkln] (0,1.8) -- (0,2.2);
\node at (0,2.6) {$\scriptstyle a+1$};
\draw [thkln] (0,-1.8) -- (0,-2.2);
\node at (0,-2.6) {$\scriptstyle a+1$};
\draw [ptzer] (-1.8,-0.6) rectangle ++(0.3,1.2);
\draw [ptzer] (1.5,-0.6) rectangle ++(0.3,1.2);
\draw [ptzer] (-0.6,1.5) rectangle ++(1.2,0.3);
\draw [ptzer] (-0.6,-1.5) rectangle ++(1.2,-0.3);
\end{tikzpicture}
}
\\
&\;\hteqv\;
\shcr^{-\hlf a^2}\;
\Cnv{
\cdots\longrightarrow
\shcr^{\ysvki + n_i}
\begin{tikzpicture}[scale=\xscl,
rotate=-45,
baseline=-3.5]
\draw [thkln] (-1.5,0.4) to [out=0,in=-90] (-0.4,1.5);
\node at (-1+0.45,1-0.45) {$\scriptstyle a-\yvki$};
\draw (-1.5,-0.4) -- (1.5,-0.4);
\draw [lnovrtw] (-1.5,0) to [out=0,in=90] (-0.4,-1.5);
\draw [thkln] (-1.5,0) to [out=0,in=90] (-0.4,-1.5) ;
\draw [lnovrtw] (0,-1.5) to [out=90,in=-135] (0,0);
\draw [thkln] (0,-1.5) to [out=90,in=-135] (0,0) to [out=45,in=180] (1.5,0);
\draw [thkln]  (1.5,0.4) to [out = 180, in=-90] (0,1.5);
\draw [lnovrtw] (0.4,-1.5) -- (0.4,1.5);
\draw  (0.4,-1.5) -- (0.4,1.5);
%
%
\draw [thkln] (-2.2,0)  -- (-1.8,0);
\node at (-2.6,0) {$\scriptstyle a+1$};
\draw [thkln] (1.8,0) -- (2.2,0);
\node at (2.6,0) {$\scriptstyle a+1$};
\draw [thkln] (0,1.8) -- (0,2.2);
\node at (0,2.6) {$\scriptstyle a+1$};
\draw [thkln] (0,-1.8) -- (0,-2.2);
\node at (0,-2.6) {$\scriptstyle a+1$};
\draw [ptzer] (-1.8,-0.6) rectangle ++(0.3,1.2);
\draw [ptzer] (1.5,-0.6) rectangle ++(0.3,1.2);
\draw [ptzer] (-0.6,1.5) rectangle ++(1.2,0.3);
\draw [ptzer] (-0.6,-1.5) rectangle ++(1.2,-0.3);
\end{tikzpicture}
\longrightarrow\cdots
}_{\;i=0}^{\;\infty}\;
\end{split}
\end{equation}
The categorification complex of a constituent tangle of the resulting \tmcn\ can be simplified:
\begin{equation}
\label{eq:sepbr}
\begin{split}
\begin{tikzpicture}[scale=\xscl,
rotate=-45,
baseline=-3.5]
\draw [thkln] (-1.5,0.4) to [out=0,in=-90] (-0.4,1.5);
\node at (-1,1) {$\scriptstyle a-\yki$};
\draw (-1.5,-0.4) -- (1.5,-0.4);
\draw [lnovrtw] (-1.5,0) to [out=0,in=90] (-0.4,-1.5);
\draw [thkln] (-1.5,0) to [out=0,in=90] (-0.4,-1.5) ;
\draw [lnovrtw] (0,-1.5) to [out=90,in=-135] (0,0);
\draw [thkln] (0,-1.5) to [out=90,in=-135] (0,0) to [out=45,in=180] (1.5,0);
\draw [thkln]  (1.5,0.4) to [out = 180, in=-90] (0,1.5);
\draw [lnovrtw] (0.4,-1.5) -- (0.4,1.5);
\draw  (0.4,-1.5) -- (0.4,1.5);
%
%
\draw [thkln] (-2.2,0)  -- (-1.8,0);
\node at (-2.6,0) {$\scriptstyle a+1$};
\draw [thkln] (1.8,0) -- (2.2,0);
\node at (2.6,0) {$\scriptstyle a+1$};
\draw [thkln] (0,1.8) -- (0,2.2);
\node at (0,2.6) {$\scriptstyle a+1$};
\draw [thkln] (0,-1.8) -- (0,-2.2);
\node at (0,-2.6) {$\scriptstyle a+1$};
\draw [ptzer] (-1.8,-0.6) rectangle ++(0.3,1.2);
\draw [ptzer] (1.5,-0.6) rectangle ++(0.3,1.2);
\draw [ptzer] (-0.6,1.5) rectangle ++(1.2,0.3);
\draw [ptzer] (-0.6,-1.5) rectangle ++(1.2,-0.3);
\end{tikzpicture}
&
\;\hteqv\;
\shcr^{2\xki-a}
\begin{tikzpicture}[scale=\xscl,
rotate=-45,
baseline=-3.5]
\draw [thkln] (-1.5,0.4) to [out=0,in=-90] (-0.4,1.5);
\node at (-1.05,1.05) {$\scriptstyle a-\yki$};
\node at (1,-1) {$\scriptstyle a-\yki$};
\node at (1.,1.) {$\scriptstyle \yki$};
\node at (-1.,-1.) {$\scriptstyle \yki$};
\draw (-1.5,0) -- (1.5,0);
\draw [thkln] (-1.5,-0.4) to [out=0,in=90] (-0.4,-1.5) ;
\draw [thkln] (0.4,-1.5) to [out=90,in=180]  (1.5,-0.4);
\draw [thkln]  (1.5,0.4) to [out = 180, in=-90] (0.4,1.5);
\draw [lnovrtw] (0,-1.5) -- (0,1.5);
\draw  (0,-1.5) -- (0,1.5);
%
%
\draw [thkln] (-2.2,0)  -- (-1.8,0);
\node at (-2.6,0) {$\scriptstyle a+1$};
\draw [thkln] (1.8,0) -- (2.2,0);
\node at (2.6,0) {$\scriptstyle a+1$};
\draw [thkln] (0,1.8) -- (0,2.2);
\node at (0,2.6) {$\scriptstyle a+1$};
\draw [thkln] (0,-1.8) -- (0,-2.2);
\node at (0,-2.6) {$\scriptstyle a+1$};
\draw [ptzer] (-1.8,-0.6) rectangle ++(0.3,1.2);
\draw [ptzer] (1.5,-0.6) rectangle ++(0.3,1.2);
\draw [ptzer] (-0.6,1.5) rectangle ++(1.2,0.3);
\draw [ptzer] (-0.6,-1.5) rectangle ++(1.2,-0.3);
\end{tikzpicture}
\\
&
\;\hteqv\;
\shcr^{2\xki-a}\;
\boxed{
\,\shcr^{\hlf}
\begin{tikzpicture}[scale=\xscl,
rotate=-45,
baseline=-3.5]
\draw [thkln] (-1.5,0.3) to [out=0,in=-90] (-0.3,1.5);
\node at (-0.4,0.4) {$\scriptstyle a-\yki$};
\node at (0.4,-0.4) {$\scriptstyle a-\yki$};
\node at (1.15,1.15) {$\scriptstyle \yki+1$};
\node at (-1.15,-1.15) {$\scriptstyle \yki+1$};
\draw [thkln] (-1.5,-0.3) to [out=0,in=90] (-0.3,-1.5) ;
\draw [thkln] (0.3,-1.5) to [out=90,in=180]  (1.5,-0.3);
\draw [thkln]  (1.5,0.3) to [out = 180, in=-90] (0.3,1.5);
%
\draw [thkln] (-2.2,0)  -- (-1.8,0);
\node at (-2.6,0) {$\scriptstyle a+1$};
\draw [thkln] (1.8,0) -- (2.2,0);
\node at (2.6,0) {$\scriptstyle a+1$};
\draw [thkln] (0,1.8) -- (0,2.2);
\node at (0,2.6) {$\scriptstyle a+1$};
\draw [thkln] (0,-1.8) -- (0,-2.2);
\node at (0,-2.6) {$\scriptstyle a+1$};
\draw [ptzer] (-1.8,-0.6) rectangle ++(0.3,1.2);
\draw [ptzer] (1.5,-0.6) rectangle ++(0.3,1.2);
\draw [ptzer] (-0.6,1.5) rectangle ++(1.2,0.3);
\draw [ptzer] (-0.6,-1.5) rectangle ++(1.2,-0.3);
\end{tikzpicture}
\longrightarrow
\shcr^{-\hlf}
\begin{tikzpicture}[scale=\xscl,
rotate=-45,
baseline=-3.5]
\draw [thkln] (-1.5,0.3) to [out=0,in=-90] (-0.3,1.5);
\node at (-0.4,0.4) {$\scriptstyle a-\yki+1$};
\node at (0.4,-0.4) {$\scriptstyle a-\yki+1$};
\node at (0.9,0.9) {$\scriptstyle \yki$};
\node at (-0.9,-0.9) {$\scriptstyle \yki$};
\draw [thkln] (-1.5,-0.3) to [out=0,in=90] (-0.3,-1.5) ;
\draw [thkln] (0.3,-1.5) to [out=90,in=180]  (1.5,-0.3);
\draw [thkln]  (1.5,0.3) to [out = 180, in=-90] (0.3,1.5);
%
\draw [thkln] (-2.2,0)  -- (-1.8,0);
\node at (-2.6,0) {$\scriptstyle a+1$};
\draw [thkln] (1.8,0) -- (2.2,0);
\node at (2.6,0) {$\scriptstyle a+1$};
\draw [thkln] (0,1.8) -- (0,2.2);
\node at (0,2.6) {$\scriptstyle a+1$};
\draw [thkln] (0,-1.8) -- (0,-2.2);
\node at (0,-2.6) {$\scriptstyle a+1$};
\draw [ptzer] (-1.8,-0.6) rectangle ++(0.3,1.2);
\draw [ptzer] (1.5,-0.6) rectangle ++(0.3,1.2);
\draw [ptzer] (-0.6,1.5) rectangle ++(1.2,0.3);
\draw [ptzer] (-0.6,-1.5) rectangle ++(1.2,-0.3);
\end{tikzpicture}
}
\end{split}
\end{equation}
Here the first homotopy equivalence follows from \ex{eq:projtw} and the second one is the application of Khovanov bracket\rx{eq:dkhbr} to the crossing of two single lines. We substitute \ex{eq:sepbr} for every constituent tangle in both \tmcn s of \ex{eq:spcrmt}. The \lumps\ \tmcn\ transforms into the \rhs of \ex{eq:colKhbr} for the intersection of two $(a+1)$-cables with the help of a simple identity
%
\[
{a \brace \yki-1}_{\shcr} + \shcr^{2\yki} {a\brace \yki}_{\shcr} =
{a+1\brace \yki}_{\shcr}.
\]
Associativity of the cone operation implies that the second \tmcn\ of \ex{eq:spcrmt} can be brought to the linear form of the \rhs of \ex{eq:lmtcn}, so it remains to verify inequalities\rx{eq:twineq}. The first inequality follows from the \lumps\ \tmcn\ formula\rx{eq:colKhbr} which we have just proved. Let us verify the second inequality for two tangles of the cone\rx{eq:sepbr} after they appear through the substitution in the second \tmcn\ of \ex{eq:spcrmt}. Since every constituent tangle of\rx{eq:spcrmt} is replaced by a cone of two tangles, the second tangle of the cone\rx{eq:sepbr} will appear at the \tmcn\ position $i'=2i$ and the second inequality of\rx{eq:twineq} for it obviously holds. The first tangle of \rx{eq:spcrmt} appears at the position $i'=2i+1$ and it carries $\yki'=\yki+1$. The inequality $2i'+1\geq 2^{\yki'}-1$ follows easily from the assumed inequality $i\leq 2^{\yki}-1$.
\end{proof}

\subsection{Recurrence relations for \cJWp s}

\begin{proposition}
A larger \tJWp\ absorbs a smaller one:
\begin{equation}
\label{eq:absorb}
\begin{tikzpicture}[scale=0.5,baseline=-1.5]
\draw [thkln] 
(-0.,0) -- (1-0.15,0)
node [near start,below] {$\scriptstyle \xca$}
 (1+.15,0) -- (1.75,0)
(2.05,0) -- (2.9,0)
node [near end, below] {$\scriptstyle \xca$};
\draw (0,0.9) -- (1.75,0.9) (2.05,0.9) -- (2.9,0.9);
\draw [ptzer] (1 - 0.15,-0.6) rectangle ++(0.3,1.2);
\draw [ptzer] (1.75,-1.2) rectangle ++(0.3,2.4);
\end{tikzpicture}
\;\hteqv\;
\begin{tikzpicture}[scale=0.5,baseline=-1.5]
\draw [ptzer] (1.75,-0.6) rectangle ++(0.3,1.2);
\draw [thkln] (1.75-1.3,0) -- (1.75,0)
node [near start, below] {$\scriptstyle \xca+1$}
(2.05,0) -- (2.05+1.3,0)
node [near end, below] {$\scriptstyle \xca+1$};
\end{tikzpicture}\;.
\end{equation}
\end{proposition}
\begin{proof}
In view of \eex{eq:projcn} and\rx{eq:grproj} for $a=\xca$, this equivalence is a result of purging the smaller projector with the larger one.
\end{proof}

Let us introduce a shortcut notation:
\begin{equation}
\label{eq:xthrpr}
\begin{tikzpicture}[menvone]
\draw[pttwo] (-.15,-0.6) rectangle ++(0.3,1.2);
\draw[thkln] (-1.35,0) -- (-.15,0)
node[near start,below] {$\scriptstyle \xca+1$ }
(.15,0) -- (1.35,0)
node [near end, below] {$\scriptstyle \xca+1$};
\end{tikzpicture}
=
\begin{tikzpicture}[scale=0.5,baseline=-1.5]
\draw[ptone] (-.15,-1.2) rectangle ++(.3,2.4);
\draw[ptzer] (-1-.15,-0.6) rectangle ++(.3,1.2);
\draw[ptzer] (1-.15,-0.6) rectangle ++(.3,1.2);
\draw[vthln] (-2,0) -- (-1-.15,0)
 node[near start, below] {$\scriptstyle \xca$}
 (-1+.15,0) -- (-.15,0) (.15,0) -- (1-.15,0)  (1.15,0) -- (2,0)
 node[near end,below] {$\scriptstyle \xca$};
 \draw (-2,0.9) -- (-.15,0.9) (.15,0.9) -- (2,0.9);
\end{tikzpicture},
\end{equation}
where the complex
$
\begin{tikzpicture}[scale=0.4,baseline=-3]
\draw[ptone] (-.15,-0.6) rectangle ++(0.3,1.2);
\draw[vthln] (-0.75,0) -- (-.15,0) (.15,0) -- (.75,0);
\end{tikzpicture}
$
is defined by \ex{eq:projcn}.
\begin{proposition}
\label{pr:hteqvgr}
Thus defined, the complex
$
\begin{tikzpicture}[scale=0.4,baseline=-3]
\draw[pttwo] (-.15,-0.6) rectangle ++(0.3,1.2);
\draw[vthln] (-0.75,0) -- (-.15,0) (.15,0) -- (.75,0);
\end{tikzpicture}
$
has a \tmcn\ presentation
\begin{equation}
\label{eq:swgrpr}
\begin{tikzpicture}[menvone]
\draw[pttwo] (-.15,-0.6) rectangle ++(0.3,1.2);
\draw[thkln] (-1.35,0) -- (-.15,0)
node[near start,below] {$\scriptstyle \xca+1$ }
(.15,0) -- (1.35,0)
node [near end, below] {$\scriptstyle \xca+1$};
\end{tikzpicture}
\;
\hteqv
\;
\boxed{
\cdots\longrightarrow
\shcr^i \bigoplus_{j=0}^i
\prmltij\shfr^j
\begin{tikzpicture}[menvone]
\draw [line width=6pt, color=white] (-1+0.15,0) -- (0,0);
\draw [thkln] (-1.75,0) -- (-1 - .15,0)
node [near start, below] {$\scriptstyle \xca$}
(-1+0.15,0) -- (0,0) -- (1-0.15,0) (1+.15,0) -- (1.75,0)
node [near end, below] {$\scriptstyle \xca$};
\draw (1, 0.9) to [out=180,in=90] (0.4,0.6) to [out=-90,in=180] (1-.15,0.3);
\draw [xscale=-1] (1, 0.9) to [out=180,in=90] (0.4,0.6) to [out=-90,in=180] (1-.15,0.3);
\draw (-1.75,0.9) -- (-1,0.9) (1,0.9) -- (1.75,0.9);
\draw [ptzer] (-1.15,-0.6) rectangle ++(0.3,1.2);
\draw [ptzer] (1 - 0.15,-0.6) rectangle ++(0.3,1.2);
\draw [dotted] (-0.6,-0.7) rectangle (0.6,1.2);
\end{tikzpicture}
\longrightarrow\cdots
}_{\;i=0}^{\infty},
\end{equation}
where $\prmltij$ are the multiplicities of the \tTLt\ inside the dotted box, with which it appears in the \rhs of \ex{eq:grproj}.
\end{proposition}
\begin{proof}
We purge the complex $
\begin{tikzpicture}[scale=0.4,baseline=-3]
\draw[pttwo] (-.15,-0.6) rectangle ++(0.3,1.2);
\draw[vthln] (-0.75,0) -- (-.15,0) (.15,0) -- (.75,0);
\end{tikzpicture}
$ with the help of two $\xca$-strand \tJWp s.
The tangle in the dotted box is the only
$(\xca+1,\xca+1)$ \taTLt\
which is not contracted when sandwiched between them.
\end{proof}

\begin{theorem}
\label{thm:projpr}
The $(\xca+1)$-strand \cJWp\ is homotopy equivalent to a cone
\begin{equation}
\label{eq:jwpar}
\begin{tikzpicture}[menvone]
\draw[ptzer] (-.15,-0.6) rectangle ++(0.3,1.2);
\draw[thkln] (-1.35,0) -- (-.15,0)
node[near start,below] {$\scriptstyle \xca+1$ }
(.15,0) -- (1.35,0)
node [near end, below] {$\scriptstyle \xca+1$};
\end{tikzpicture}
\hteqv
\boxed{
\shcr
\begin{tikzpicture}[menvone]
\draw[pttwo] (-.15,-0.6) rectangle ++(0.3,1.2);
\draw[thkln] (-1.35,0) -- (-.15,0)
node[near start,below] {$\scriptstyle \xca+1$ }
(.15,0) -- (1.35,0)
node [near end, below] {$\scriptstyle \xca+1$};
\end{tikzpicture}
\longrightarrow
\begin{tikzpicture}[menvone]
\draw[color=white] (-0.15,-1.2) rectangle (0.15,1.2);
\draw[ptzer] (-0.15,-0.6) rectangle (0.15,0.6);
\draw [thkln] (-0.75,0) -- (-0.15,0)
node [near start, below] {$\scriptstyle \xca$} (0.15,0) -- (0.75,0)
node [near end, below] {$\scriptstyle \xca$};
\draw (-0.75,0.9) -- (0.75, 0.9);
\end{tikzpicture}
}
\end{equation}
\end{theorem}
\begin{proof}
Consider a sequence of homotopy equivalences
\begin{multline}
\label{eq:thrprj}
\begin{tikzpicture}[menvone]
\draw[ptzer] (-.15,-0.6) rectangle ++(0.3,1.2);
\draw[thkln] (-1.35,0) -- (-.15,0)
node[near start,below] {$\scriptstyle \xca+1$ }
(.15,0) -- (1.35,0)
node [near end, below] {$\scriptstyle \xca+1$};
\end{tikzpicture}
\;\hteqv\;
\begin{tikzpicture}[scale=0.5,baseline=-1.5]
\draw[ptzer] (-.15,-1.2) rectangle ++(.3,2.4);
\draw[ptzer] (-1-.15,-0.6) rectangle ++(.3,1.2);
\draw[ptzer] (1-.15,-0.6) rectangle ++(.3,1.2);
\draw[vthln] (-2,0) -- (-1-.15,0)
 node[near start, below] {$\scriptstyle \xca$}
 (-1+.15,0) -- (-.15,0) (.15,0) -- (1-.15,0)  (1.15,0) -- (2,0)
 node[near end,below] {$\scriptstyle \xca$};
 \draw (-2,0.9) -- (-.15,0.9) (.15,0.9) -- (2,0.9);
\end{tikzpicture}
\;
\hteqv
\;
\boxed{
\shcr
\begin{tikzpicture}[scale=0.5,baseline=-1.5]
\draw[ptone] (-.15,-1.2) rectangle ++(.3,2.4);
\draw[ptzer] (-1-.15,-0.6) rectangle ++(.3,1.2);
\draw[ptzer] (1-.15,-0.6) rectangle ++(.3,1.2);
\draw[vthln] (-2,0) -- (-1-.15,0)
 node[near start, below] {$\scriptstyle \xca$}
 (-1+.15,0) -- (-.15,0) (.15,0) -- (1-.15,0)  (1.15,0) -- (2,0)
 node[near end,below] {$\scriptstyle \xca$};
 \draw (-2,0.9) -- (-.15,0.9) (.15,0.9) -- (2,0.9);
\end{tikzpicture}
\longrightarrow
\begin{tikzpicture}[menvone]
\draw[ptzer] (-.15,-0.6) rectangle ++(0.3,1.2);
\draw[ptzer] (1-.15,-0.6) rectangle ++(0.3,1.2);
\draw[vthln] (-1,0) -- (-.15,0)
 node[near start, below] {$\scriptstyle \xca$}
(.15,0) -- (1-.15,0) (1+.15,0) --(2,0)
 node[near end, below] {$\scriptstyle \xca$};
 \draw (-1,0.9) -- (2,0.9);
\end{tikzpicture}
}
\\
\;\hteqv\;
\boxed{
\shcr
\begin{tikzpicture}[menvone]
\draw[pttwo] (-.15,-0.6) rectangle ++(0.3,1.2);
\draw[thkln] (-1.35,0) -- (-.15,0)
node[near start,below] {$\scriptstyle \xca+1$ }
(.15,0) -- (1.35,0)
node [near end, below] {$\scriptstyle \xca+1$};
\end{tikzpicture}
\longrightarrow
\begin{tikzpicture}[menvone]
 \draw [color=white] (-0.15,-1.2) rectangle ++ (0.3,2.4);
\draw[ptzer] (-.15,-0.6) rectangle ++(0.3,1.2);
\draw[vthln] (-1,0) -- (-.15,0)
 node[near start, below] {$\scriptstyle \xca$}
(.15,0) -- (1,0)
 node[near end, below] {$\scriptstyle \xca$};
 \draw (-1,0.9) -- (1,0.9);
\end{tikzpicture}
}
\end{multline}
The first homotopy equivalence comes from \ex{eq:absorb}, the second follows from \ex{eq:projcn} and the last one follows from \eex{eq:xthrpr} and\rx{eq:cmppr}.
\end{proof}

\begin{theorem}
\label{thm:jwind}
The $(\xca+1)$-strand \cJWp\ can be presented as the following cone:
\begin{equation}
\label{eq:jwind}
\begin{tikzpicture}[menvone]
\draw[ptzer] (-.15,-0.6) rectangle ++(0.3,1.2);
\draw[thkln] (-1.35,0) -- (-.15,0)
node[near start,below] {$\scriptstyle \xca+1$ }
(.15,0) -- (1.35,0)
node [near end, below] {$\scriptstyle \xca+1$};
\end{tikzpicture}
\hteqv
\boxed{
\shcr^{2\xca+1}\shfr^{2}
\begin{tikzpicture}[menvone]
\draw[pttwo] (-.15,-0.6) rectangle ++(0.3,1.2);
\draw[thkln] (-1.35,0) -- (-.15,0)
node[near start,below] {$\scriptstyle \xca+1$ }
(.15,0) -- (1.35,0)
node [near end, below] {$\scriptstyle \xca+1$};
\end{tikzpicture}
\longrightarrow
\shcr^{\xca}
\begin{tikzpicture}[scale=0.5,baseline=-1.5]
\draw [xscale=-1]  (1, 0.9) to [out=180,in=0] (0,-0.4);
\draw [line width=6pt, color=white] (-1+0.15,0) -- (0,0);
\draw [thkln] (-1.75,0) -- (-1 - .15,0)
node [near start,below] {$\scriptstyle \xca$}
(-1+0.15,0) -- (0,0) -- (1-0.15,0) (1+.15,0) -- (1.75,0)
node [near end, below] {$\scriptstyle \xca$};
\draw [line width=6pt, color=white] (1, 0.9) to [out=180,in=0] (0,-0.4);
\draw (1, 0.9) to [out=180,in=0] (0,-0.4);
\draw (-1.75,0.9) -- (-1,0.9) (1,0.9) -- (1.75,0.9);
\draw [ptzer] (-1.15,-0.6) rectangle ++(0.3,1.2);
\draw [ptzer] (1 - 0.15,-0.6) rectangle ++(0.3,1.2);
\end{tikzpicture}
}
\end{equation}
\end{theorem}

\begin{lemma}
\label{lm:grwind}
There is a homotopy equivalence
\begin{equation}
\label{eq:grwind}
\begin{tikzpicture}[scale=0.5,baseline=-1.5]
\draw [xscale=-1]  (1, 0.9) to [out=180,in=0] (0,-0.4);
\draw [line width=6pt, color=white] (-1+0.15,0) -- (0,0);
\draw [thkln] (-1.75+0.3,0) -- (-1 +0.3 - .15,0)
node [near start,below] {$\scriptstyle \xca$}
(-1+0.15,0) -- (0,0) -- (1,0) (1.3,0) -- (2.15,0)
(2.05,0) -- (2.15,0)
node [near end, below] {$\scriptstyle \xca$};
\draw [line width=6pt, color=white] (1, 0.9) to [out=180,in=0] (0,-0.4);
\draw (1, 0.9) to [out=180,in=0] (0,-0.4);
\draw (-1.75+0.3,0.9) -- (-1,0.9) (1.3,0.9) -- (2.15,0.9); 
\draw [pttwo] (1,-1.2) rectangle ++(0.3,2.4);
\end{tikzpicture}
\hteqv
\shcr^{\xca}\shfr^2
\begin{tikzpicture}[menvone]
\draw[pttwo] (-.15,-0.6) rectangle ++(0.3,1.2);
\draw[thkln] (-1.35,0) -- (-.15,0)
node[near start,below] {$\scriptstyle \xca+1$ }
(.15,0) -- (1.35,0)
node [near end, below] {$\scriptstyle \xca+1$};
\end{tikzpicture}
\end{equation}
\end{lemma}
\begin{proof}
Consider the composition of the line winding around the $\xca$-cable with the left portion of the complex
$
\begin{tikzpicture}[scale=0.4,baseline=-1.5]
\draw [line width=6pt, color=white] (-1+0.15,0) -- (0,0);
\draw [line width=\cblth] (-1.75,0) -- (-1 - .15,0)
(-1+0.15,0) -- (0,0) -- (1-0.15,0) (1+.15,0) -- (1.75,0);
\draw (1, 0.9) to [out=180,in=90] (0.4,0.6) to [out=-90,in=180] (1-.15,0.3);
\draw [xscale=-1] (1, 0.9) to [out=180,in=90] (0.4,0.6) to [out=-90,in=180] (1-.15,0.3);
\draw (-1.75,0.9) -- (-1,0.9) (1,0.9) -- (1.75,0.9);
\draw [line width=\ljwp] (-1.15,-0.6) rectangle ++(0.3,1.2);
\draw [line width=\ljwp] (1 - 0.15,-0.6) rectangle ++(0.3,1.2);
\end{tikzpicture}
$
which generates the \tmcn\rx{eq:swgrpr}:
\begin{equation}
\label{eq:windjw}
\begin{tikzpicture}[scale=0.75,baseline=-1.5]
\draw [xscale=-1]  (1, 0.9) to [out=180,in=0] (0,-0.4);
\draw [line width=6pt, color=white] (-1+0.15,0) -- (0,0);
\draw [thkln] (-1.15,0) -- (-1 + .15,0)
node [near start,below] {$\scriptstyle \xca$}
(-1+0.15,0) -- (0,0) -- (1-0.15,0);
\draw [line width=6pt, color=white] (1, 0.9) to [out=180,in=0] (0,-0.4);
\draw (1, 0.9) to [out=180,in=0] (0,-0.4);
\draw (-1.15,0.9) -- (-1,0.9); 
\draw [line width=\ljwp] (1 - 0.15,-0.6) rectangle ++(0.3,1.2);
\begin{scope}[xshift=2cm]
\draw [thkln]
(-1+0.15,0) -- (-0.15+0.3,0) 
node [near end, below] {$\scriptstyle \xca-1$};
\draw [xscale=-1] (1, 0.9) to [out=180,in=90] (0.4,0.6) to [out=-90,in=180] (1-.15,0.3);
\end{scope}
\end{tikzpicture}
\;\hteqv\;
\begin{tikzpicture}[scale=0.75,baseline=-1.5]
\draw [xscale=-1]  (1, 0.9) to [out=180,in=0] (0,-0.4);
\draw [line width=4pt, color=white] (-1+0.15,0) -- (0,0);
\draw [thkln] (-1.75,0) -- (-1 - .15,0)
node [near start,below] {$\scriptstyle \xca$}
(-1+0.15,0) -- (0,0) -- (1+0.15,0);
\draw [line width=4pt, color=white] (-1+0.15,0.3) -- (0,0.3);
\draw (-1+0.15,0.3)  -- (1+0.15,0.3);
\draw [line width=6pt, color=white] (1, 0.9) to [out=180,in=0] (0,-0.4);
\draw (1, 0.9) to [out=180,in=0] (0,-0.4);
\draw (-1.75,0.9) -- (-1,0.9); 
\draw [ptzer] (-1.15,-0.6) rectangle ++(0.3,1.2);
\begin{scope}[xshift=2cm]
\draw [line width=\cblth]
(-1+0.15,0) -- (-0.15-0.2,0) 
node [near end, below] {$\scriptstyle \xca-1$};
\draw [xscale=-1] (1, 0.9) to [out=180,in=90] (0.4,0.6) to [out=-90,in=180] (1-.15,0.3);
\end{scope}
\end{tikzpicture}
\;\hteqv\;
\shcr \shfr^{2}
\begin{tikzpicture}[menvtwo]
\draw [thkln] (-0.75,0) -- (-0.15,0)
node [near start, below] {$\scriptstyle \xca$}
(0.15,0) -- (1,0);
\draw [lnovr] (0.15,0.3) to [out=0,in=180] (1,-0.4);
\draw (0.15,0.3) to [out=0,in=180] (1,-0.4);
\draw (1,-0.4) .. controls +(0:0.8) and +(0:1.2) .. (0.2,0.9) -- (-0.75,0.9);
\draw [lnovr] (1.1,0) -- (2,0);
\draw [thkln] (1,0) -- (2.2,0)
node [very near end, below] {$\scriptstyle \xca-1$};
\draw[ptzer] (-0.15,-0.6) rectangle (0.15,0.6);
\end{tikzpicture}
\;\hteqv\;
\shcr^{\xca}\shfr^{2}
\begin{tikzpicture}[menvtwo]
\draw [ptzer] (-0.15,-0.6) rectangle (0.15,0.6);
\draw [thkln] (-0.75,0) -- (-0.15,0)
node [near start, below] {$\scriptstyle \xca$} (0.15,0) -- (1,0)
node [near end,below] {$\scriptstyle \xca-1$};
\draw [xscale=-1,xshift=-1cm] (1, 0.9) to [out=180,in=90] (0.4,0.6) to [out=-90,in=180] (1-.15,0.3);
\draw (-0.75,0.9) -- (0,0.9);
\end{tikzpicture}
\end{equation}
Here the first equivalence is purely topological: the projector is moved left along the cable, the second equivalence uses \ex{eq:frsh} to remove two framing kinks on the single line and
 the third equivalence follows from \ex{eq:projtw}.
%
The equivalence\rx{eq:grwind} comes from applying equivalence\rx{eq:windjw} to every constituent complex in the \tmcn\rx{eq:swgrpr}.
\end{proof}
\begin{proof}[Proof of Theorem\rw{thm:jwind}]
Eq.\rx{eq:jwind} follows from a sequence of homotopy equivalences:
\begin{multline*}
\begin{tikzpicture}[scale=0.5,baseline=-1.5]
\draw [ptzer] (1.75,-0.6) rectangle ++(0.3,1.2);
\draw [thkln] (1.75-0.85-0.3,0) -- (1.75,0)
node [very near start, below] {$\scriptstyle \xca+1$}
(2.05,0) -- (2.9+0.3,0)
node [very near end, below] {$\scriptstyle \xca+1$};
\end{tikzpicture}
\;\hteqv\;
\shcr^{\xca}
\begin{tikzpicture}[scale=0.5,baseline=-1.5]
\draw [xscale=-1]  (1, 0.9) to [out=180,in=0] (0,-0.4);
\draw [line width=6pt, color=white] (-1+0.15,0) -- (0,0);
\draw [thkln] (-1.75+0.3,0) -- (-1 +0.3 - .15,0)
node [near start,below] {$\scriptstyle \xca$}
(-1+0.15,0) -- (0,0) -- (1,0) (1.3,0) -- (2.15,0)
(2.05,0) -- (2.15,0)
node [near end, below] {$\scriptstyle \xca$};
\draw [line width=6pt, color=white] (1, 0.9) to [out=180,in=0] (0,-0.4);
\draw (1, 0.9) to [out=180,in=0] (0,-0.4);
\draw (-1.75+0.3,0.9) -- (-1,0.9) (1.3,0.9) -- (2.15,0.9); 
\draw [ptzer] (1,-1.2) rectangle ++(0.3,2.4);
\end{tikzpicture}
\;\hteqv\;
\boxed{
\shcr^{\xca+1}
\begin{tikzpicture}[scale=0.5,baseline=-1.5]
\draw [xscale=-1]  (1, 0.9) to [out=180,in=0] (0,-0.4);
\draw [line width=6pt, color=white] (-1+0.15,0) -- (0,0);
\draw [thkln] (-1.75+0.3,0) -- (-1 +0.3 - .15,0)
node [near start,below] {$\scriptstyle \xca$}
(-1+0.15,0) -- (0,0) -- (1,0) (1.3,0) -- (2.15,0)
(2.05,0) -- (2.15,0)
node [near end, below] {$\scriptstyle \xca$};
\draw [line width=6pt, color=white] (1, 0.9) to [out=180,in=0] (0,-0.4);
\draw (1, 0.9) to [out=180,in=0] (0,-0.4);
\draw (-1.75+0.3,0.9) -- (-1,0.9) (1.3,0.9) -- (2.15,0.9); 
\draw [pttwo] (1,-1.2) rectangle ++(0.3,2.4);
\end{tikzpicture}
\longrightarrow
\shcr^{\xca}
\begin{tikzpicture}[scale=0.5,baseline=-1.5]
\draw [xscale=-1]  (1, 0.9) to [out=180,in=0] (0,-0.4);
\draw [line width=6pt, color=white] (-1+0.15,0) -- (0,0);
\draw [thkln] (-1.75+0.3,0) -- (-1 - .15+0.3,0)
node [near start,below] {$\scriptstyle \xca$}
(-1+0.15,0) -- (0,0) -- (1-0.15,0) (1+.15,0) -- (1.75,0)
node [near end, below] {$\scriptstyle \xca$};
\draw [line width=6pt, color=white] (1, 0.9) to [out=180,in=0] (0,-0.4);
\draw (1, 0.9) to [out=180,in=0] (0,-0.4);
\draw (-1.75+0.3,0.9) -- (-1,0.9) (1,0.9) -- (1.75,0.9);
\draw [line width=\ljwp] (1 - 0.15,-0.6) rectangle ++(0.3,1.2);
\end{tikzpicture}
}
\\
\;\hteqv\;
\boxed{
\shcr^{2\xca+1}\shfr^{2}
\begin{tikzpicture}[menvone]
\draw[pttwo] (-.15,-0.6) rectangle ++(0.3,1.2);
\draw[thkln] (-1.35,0) -- (-.15,0)
node[near start,below] {$\scriptstyle \xca+1$ }
(.15,0) -- (1.35,0)
node [near end, below] {$\scriptstyle \xca+1$};
\end{tikzpicture}
\longrightarrow
\shcr^{\xca}
\begin{tikzpicture}[scale=0.5,baseline=-1.5]
\draw [xscale=-1]  (1, 0.9) to [out=180,in=0] (0,-0.4);
\draw [line width=6pt, color=white] (-1+0.15,0) -- (0,0);
\draw [thkln] (-1.75,0) -- (-1 - .15,0)
node [near start,below] {$\scriptstyle \xca$}
(-1+0.15,0) -- (0,0) -- (1-0.15,0) (1+.15,0) -- (1.75,0)
node [near end, below] {$\scriptstyle \xca$};
\draw [line width=6pt, color=white] (1, 0.9) to [out=180,in=0] (0,-0.4);
\draw (1, 0.9) to [out=180,in=0] (0,-0.4);
\draw (-1.75,0.9) -- (-1,0.9) (1,0.9) -- (1.75,0.9);
\draw [ptzer] (-1.15,-0.6) rectangle ++(0.3,1.2);
\draw [ptzer] (1 - 0.15,-0.6) rectangle ++(0.3,1.2);
\end{tikzpicture}
}
\end{multline*}
Here the first equivalence follows from \ex{eq:projtw}, the second equivalence follows from \ex{eq:jwpar}, the third equivalence follows from \eex{eq:grwind} and\rx{eq:cmppr}.
\end{proof}

\begin{theorem}
\label{thm:projcr}
The $(\xca+1)$-strand \cJWp\ is homotopy equivalent to a cone
\begin{equation}
\label{eq:jwcrs}
\begin{tikzpicture}[menvone]
\draw[ptzer] (-.15,-0.6) rectangle ++(0.3,1.2);
\draw[thkln] (-1.35,0) -- (-.15,0)
node[near start,below] {$\scriptstyle \xca+1$ }
(.15,0) -- (1.35,0)
node [near end, below] {$\scriptstyle \xca+1$};
\end{tikzpicture}
\hteqv
\boxed{
\shcr^{2\xca}\shfr
\begin{tikzpicture}[menvone]
\draw[ptthr] (-.15,-0.6) rectangle ++(0.3,1.2);
\draw[thkln] (-1.35,0) -- (-.15,0)
node[near start,below] {$\scriptstyle \xca+1$ }
(.15,0) -- (1.35,0)
node [near end, below] {$\scriptstyle \xca+1$};
\end{tikzpicture}
\longrightarrow
\shcr^{\hlf}
\begin{tikzpicture}[menvone]
\draw [line width=6pt, color=white] (-1+0.15,0) -- (0,0);
\draw [thkln] (-1.75,0) -- (-1 - .15,0)
node [near start, below] {$\scriptstyle \xca$}
(-1+0.15,0) -- (0,0) -- (1-0.15,0) (1+.15,0) -- (1.75,0)
node [near end, below] {$\scriptstyle \xca$};
\draw [xscale=-1] (-1+0.15,0.3) -- (-0.5,0.3) to [out=0, in=180] (0.5,0.9) -- (1,0.9);
\draw [lnovr] (-0.5,0.3) to [out=0, in=180] (0.5,0.9);
\draw (-1+0.15,0.3) -- (-0.5,0.3) to [out=0, in=180] (0.5,0.9) -- (1,0.9);
\draw (-1.75,0.9) -- (-1,0.9) (1,0.9) -- (1.75,0.9);
\draw [ptzer] (-1.15,-0.6) rectangle ++(0.3,1.2);
\draw [ptzer] (1 - 0.15,-0.6) rectangle ++(0.3,1.2);
\end{tikzpicture}
}
\end{equation}
in which the complex
$
\begin{tikzpicture}[scale=0.4,baseline=-3]
\draw[ptthr] (-.15,-0.6) rectangle ++(0.3,1.2);
\draw[vthln] (-0.75,0) -- (-.15,0) (.15,0) -- (.75,0);
\end{tikzpicture}
$
has the following \tmcn\ presentation:
\begin{equation}
\label{eq:swnepr}
\begin{tikzpicture}[menvone]
\draw[ptthr] (-.15,-0.6) rectangle ++(0.3,1.2);
\draw[thkln] (-1.35,0) -- (-.15,0)
node[near start,below] {$\scriptstyle \xca+1$ }
(.15,0) -- (1.35,0)
node [near end, below] {$\scriptstyle \xca+1$};
\end{tikzpicture}
\;
\hteqv
\;
\boxed{
\cdots\longrightarrow
\shcr^i \bigoplus_{j=0}^i
\tprmltij \shfr^j
\begin{tikzpicture}[menvone]
\draw [line width=6pt, color=white] (-1+0.15,0) -- (0,0);
\draw [thkln] (-1.75,0) -- (-1 - .15,0)
node [near start, below] {$\scriptstyle \xca$}
(-1+0.15,0) -- (0,0) -- (1-0.15,0) (1+.15,0) -- (1.75,0)
node [near end, below] {$\scriptstyle \xca$};
\draw (1, 0.9) to [out=180,in=90] (0.4,0.6) to [out=-90,in=180] (1-.15,0.3);
\draw [xscale=-1] (1, 0.9) to [out=180,in=90] (0.4,0.6) to [out=-90,in=180] (1-.15,0.3);
\draw (-1.75,0.9) -- (-1,0.9) (1,0.9) -- (1.75,0.9);
\draw [ptzer] (-1.15,-0.6) rectangle ++(0.3,1.2);
\draw [ptzer] (1 - 0.15,-0.6) rectangle ++(0.3,1.2);
\end{tikzpicture}
\longrightarrow\cdots
}_{\;i=0}^{\infty},
\quad
\tprmltij =
\begin{cases}
\prmltv{i-1,j-1} & \text{if $i\geq 1$},
\\
1 & \text{if $i=0$}.
\end{cases}
\end{equation}
\end{theorem}

\begin{lemma}
There is a homotopy equivalence
\begin{equation}
\label{eq:wdcr}
\begin{tikzpicture}[scale=0.5,baseline=-1.5]
\draw [xscale=-1]  (1, 0.9) to [out=180,in=0] (0,-0.4);
\draw [line width=6pt, color=white] (-1+0.15,0) -- (0,0);
\draw [thkln] (-1.75,0) -- (-1 - .15,0)
node [near start,below] {$\scriptstyle \xca$}
(-1+0.15,0) -- (0,0) -- (1-0.15,0) (1+.15,0) -- (1.75,0)
node [near end, below] {$\scriptstyle \xca$};
\draw [line width=6pt, color=white] (1, 0.9) to [out=180,in=0] (0,-0.4);
\draw (1, 0.9) to [out=180,in=0] (0,-0.4);
\draw (-1.75,0.9) -- (-1,0.9) (1,0.9) -- (1.75,0.9);
\draw [ptzer] (-1.15,-0.6) rectangle ++(0.3,1.2);
\draw [ptzer] (1 - 0.15,-0.6) rectangle ++(0.3,1.2);
\end{tikzpicture}
\;\hteqv\;
\boxed{
\shcr^{\xca}\shfr
\begin{tikzpicture}[menvone]
\draw [line width=6pt, color=white] (-1+0.15,0) -- (0,0);
\draw [thkln] (-1.75,0) -- (-1 - .15,0)
node [near start, below] {$\scriptstyle \xca$}
(-1+0.15,0) -- (0,0) -- (1-0.15,0) (1+.15,0) -- (1.75,0)
node [near end, below] {$\scriptstyle \xca$};
\draw (1, 0.9) to [out=180,in=90] (0.4,0.6) to [out=-90,in=180] (1-.15,0.3);
\draw [xscale=-1] (1, 0.9) to [out=180,in=90] (0.4,0.6) to [out=-90,in=180] (1-.15,0.3);
\draw (-1.75,0.9) -- (-1,0.9) (1,0.9) -- (1.75,0.9);
\draw [ptzer] (-1.15,-0.6) rectangle ++(0.3,1.2);
\draw [ptzer] (1 - 0.15,-0.6) rectangle ++(0.3,1.2);
\end{tikzpicture}
\longrightarrow
\shcr^{-\xca+\hlf}
\begin{tikzpicture}[menvone]
\draw [line width=6pt, color=white] (-1+0.15,0) -- (0,0);
\draw [thkln] (-1.75,0) -- (-1 - .15,0)
node [near start, below] {$\scriptstyle \xca$}
(-1+0.15,0) -- (0,0) -- (1-0.15,0) (1+.15,0) -- (1.75,0)
node [near end, below] {$\scriptstyle \xca$};
\draw [xscale=-1] (-1+0.15,0.3) -- (-0.5,0.3) to [out=0, in=180] (0.5,0.9) -- (1,0.9);
\draw [lnovr] (-0.5,0.3) to [out=0, in=180] (0.5,0.9);
\draw (-1+0.15,0.3) -- (-0.5,0.3) to [out=0, in=180] (0.5,0.9) -- (1,0.9);
\draw (-1.75,0.9) -- (-1,0.9) (1,0.9) -- (1.75,0.9);
\draw [line width=\ljwp] (-1.15,-0.6) rectangle ++(0.3,1.2);
\draw [line width=\ljwp] (1 - 0.15,-0.6) rectangle ++(0.3,1.2);
\end{tikzpicture}
}
\end{equation}
\end{lemma}
\begin{proof}
The lemma is proved by applying \tKhbr\ formula\rx{eq:dkhbr} to one of the elementary crossings in the \lhs diagram:
\def\cht{0.4}
\[
\begin{tikzpicture}[scale=0.75,baseline=-1.5]
\draw [xscale=-1]  (1, 0.9) to [out=180,in=0] (0,-0.4);
\draw [lnovr] (-1+0.15,0) -- (0,0);
\draw [lnovr] (-1+0.15,\cht) -- (0,\cht);
\draw [thkln] (-1.75,0) -- (-1 - .15,0)
node [near start,below] {$\scriptstyle \xca$}
(-1+0.15,0) -- (0,0) -- (1-0.15,0) (1+.15,0) -- (1.75,0)
node [near end, below] {$\scriptstyle \xca$};
\draw (-1+0.15,\cht) -- (1-0.15,\cht);
\draw [line width=6pt, color=white] (1, 0.9) to [out=180,in=0] (0,-0.4);
\draw (1, 0.9) to [out=180,in=0] (0,-0.4);
\draw (-1.75,0.9) -- (-1,0.9) (1,0.9) -- (1.75,0.9);
\draw [ptzer] (-1.15,-0.6) rectangle ++(0.3,1.2);
\draw [ptzer] (1 - 0.15,-0.6) rectangle ++(0.3,1.2);
\end{tikzpicture}
\hteqv
\boxed{
\shcr^{\hlf}
\begin{tikzpicture}[scale=0.75,baseline=-1.5]
\draw (0,-0.4) to [out=-180,in=-90] (-0.3,-0.05) to [out=90,in=-180] (0,\cht);
\draw [lnovr] (-1+0.15,0) -- (0,0);
\draw [thkln] (-1.75,0) -- (-1 - .15,0)
node [near start,below] {$\scriptstyle \xca$}
(-1+0.15,0) -- (0,0) -- (1-0.15,0) (1+.15,0) -- (1.75,0)
node [near end, below] {$\scriptstyle \xca$};
\draw (0,\cht) -- (1-0.15,\cht);
\draw (-1+0.15,\cht) to [out=0,in=-90] (-0.5,0.6) to [out=90,in=0] (-1,0.9);
\draw [line width=6pt, color=white] (1, 0.9) to [out=180,in=0] (0,-0.4);
\draw (1, 0.9) to [out=180,in=0] (0,-0.4);
\draw (-1.75,0.9) -- (-1,0.9) (1,0.9) -- (1.75,0.9);
\draw [ptzer] (-1.15,-0.6) rectangle ++(0.3,1.2);
\draw [ptzer] (1 - 0.15,-0.6) rectangle ++(0.3,1.2);
\end{tikzpicture}
\longrightarrow
\shcr^{-\hlf}
\begin{tikzpicture}[scale=0.75,baseline=-1.5]
\draw [xscale=-1]  (1-0.15, \cht) to [out=180,in=0] (0,-0.4);
\draw (-1,0.9) to [out=0,in=180] (0,\cht);
\draw [lnovr] (-1+0.15,0) -- (0,0);
\draw [thkln] (-1.75,0) -- (-1 - .15,0)
node [near start,below] {$\scriptstyle \xca$}
(-1+0.15,0) -- (0,0) -- (1-0.15,0) (1+.15,0) -- (1.75,0)
node [near end, below] {$\scriptstyle \xca$};
\draw (0,\cht) -- (1-0.15,\cht);
\draw [line width=6pt, color=white] (1, 0.9) to [out=180,in=0] (0,-0.4);
\draw (1, 0.9) to [out=180,in=0] (0,-0.4);
\draw (-1.75,0.9) -- (-1,0.9) (1,0.9) -- (1.75,0.9);
\draw [ptzer] (-1.15,-0.6) rectangle ++(0.3,1.2);
\draw [ptzer] (1 - 0.15,-0.6) rectangle ++(0.3,1.2);
\end{tikzpicture}
}
\]
The diagrams in the \rhs cone are reduced to those of \ex{eq:wdcr} with the help of \eex{eq:frsh} and\rx{eq:projtw}.
\end{proof}
\begin{proof}[Proof of Theorem\rw{thm:projcr}]
A substitution of \ex{eq:wdcr} into \ex{eq:jwind} yields the cone presentation\rx{eq:jwcrs} with
\[
\begin{tikzpicture}[menvone]
\draw[ptthr] (-.15,-0.6) rectangle ++(0.3,1.2);
\draw[thkln] (-1.35,0) -- (-.15,0)
node[near start,below] {$\scriptstyle \xca+1$ }
(.15,0) -- (1.35,0)
node [near end, below] {$\scriptstyle \xca+1$};
\end{tikzpicture}
=
\boxed{
\shcr\shfr
\begin{tikzpicture}[menvone]
\draw[pttwo] (-.15,-0.6) rectangle ++(0.3,1.2);
\draw[thkln] (-1.35,0) -- (-.15,0)
node[near start,below] {$\scriptstyle \xca+1$ }
(.15,0) -- (1.35,0)
node [near end, below] {$\scriptstyle \xca+1$};
\end{tikzpicture}
\longrightarrow
\begin{tikzpicture}[menvone]
\draw[color=white] (-0.15,-1.2) rectangle (0.15,1.2);
\draw [line width=6pt, color=white] (-1+0.15,0) -- (0,0);
\draw [thkln] (-1.75,0) -- (-1 - .15,0)
node [near start, below] {$\scriptstyle \xca$}
(-1+0.15,0) -- (0,0) -- (1-0.15,0) (1+.15,0) -- (1.75,0)
node [near end, below] {$\scriptstyle \xca$};
\draw (1, 0.9) to [out=180,in=90] (0.4,0.6) to [out=-90,in=180] (1-.15,0.3);
\draw [xscale=-1] (1, 0.9) to [out=180,in=90] (0.4,0.6) to [out=-90,in=180] (1-.15,0.3);
\draw (-1.75,0.9) -- (-1,0.9) (1,0.9) -- (1.75,0.9);
\draw [ptzer] (-1.15,-0.6) rectangle ++(0.3,1.2);
\draw [ptzer] (1 - 0.15,-0.6) rectangle ++(0.3,1.2);
\end{tikzpicture}
}
\]
and \ex{eq:swnepr} follows from \ex{eq:swgrpr}.
\end{proof}

If we connect the endpoints of the upper single line in \ex{eq:jwcrs} and apply the framing relation\rx{eq:frsh} to the last diagram in that relation, then we come to the following corollary
\begin{corollary}
There is a homotopy equivalence
\begin{equation}
\label{eq:htloop}
\begin{tikzpicture}[menvone]
\draw [ptzer] (-0.15,-0.6) rectangle (.15,.6);
\draw [thkln] (-0.75,0) -- (-0.15,0)
node [near start,below] {$\scriptstyle \xca$}
(0.15,0) -- (.75,0)
node [near end,below] {$\scriptstyle \xca$};
\draw (-0.15,0.4) to [out=180,in=-90] (-0.6,0.8) to [out=90,in=180] (0,1.4)
to [out=0,in=90] (0.6,0.8) to [out=-90,in=0] (0.15,0.4);
\end{tikzpicture}
\;\hteqv\;
\boxed{
\shcr^{2\xca}\shfr
\begin{tikzpicture}[menvone]
\draw [color=white] (-0.15,-1) rectangle (0.15,1.6);
\draw [ptthr] (-0.15,-0.6) rectangle (.15,.6);
\draw [thkln] (-0.75,0) -- (-0.15,0)
node [near start,below] {$\scriptstyle \xca$}
(0.15,0) -- (.75,0)
node [near end,below] {$\scriptstyle \xca$};
\draw (-0.15,0.4) to [out=180,in=-90] (-0.6,0.8) to [out=90,in=180] (0,1.4)
to [out=0,in=90] (0.6,0.8) to [out=-90,in=0] (0.15,0.4);
\end{tikzpicture}
\longrightarrow
\shfr^{-1}
\begin{tikzpicture}[menvone]
\draw [ptzer] (-0.15,-0.6) rectangle (.15,.6);
\draw [thkln] (-0.75,0) -- (-0.15,0)
node [near start,below] {$\scriptstyle \xca$}
(0.15,0) -- (.75,0)
node [near end,below] {$\scriptstyle \xca$};
\end{tikzpicture}
}
\end{equation}
where
\begin{equation}
\label{eq:grlp}
\begin{tikzpicture}[menvone]
\draw [color=white] (-0.15,-1) rectangle (0.15,1.6);
\draw [ptthr] (-0.15,-0.6) rectangle (.15,.6);
\draw [thkln] (-0.75,0) -- (-0.15,0)
node [near start,below] {$\scriptstyle \xca$}
(0.15,0) -- (.75,0)
node [near end,below] {$\scriptstyle \xca$};
\draw (-0.15,0.4) to [out=180,in=-90] (-0.6,0.8) to [out=90,in=180] (0,1.4)
to [out=0,in=90] (0.6,0.8) to [out=-90,in=0] (0.15,0.4);
\end{tikzpicture}
\;
\hteqv
\;
\boxed{
\cdots\longrightarrow
\shcr^i \bigoplus_{j=0}^i
\tprmltij \shfr^j
\begin{tikzpicture}[menvone]
\draw [ptzer] (-0.15,-0.6) rectangle (.15,.6);
\draw [thkln] (-0.75,0) -- (-0.15,0)
node [near start,below] {$\scriptstyle \xca$}
(0.15,0) -- (.75,0)
node [near end,below] {$\scriptstyle \xca$};
\end{tikzpicture}
\longrightarrow\cdots
}_{\;i=0}^{\infty}.
\end{equation}
\end{corollary}

\section{The morphism $\mnfN$}
\label{sct:mrph}

\subsection{General setup}

For a link diagram $\xD$ we give a precise definition of a diagram (a complex) $\xDclN$. $\xDclN$ is constructed by first $\xca$-cabling all components of $\xD$ and then placing a \cJWp\ at every edge of $\xD$, an edge being a piece of $\xca$-cabled strand between two crossings.

The map $\mnfN$ of \ex{eq:spmaps} is a composition of many maps between \tKhoms\ of a sequence of diagrams related by \ltrf s, the first diagram in that sequence being $\xDclNo$ and the last being
$\xDclN$ (recall that maps go backwards).

We use three types of \ltrf s, which are based on the following \lrpl s:
\begin{gather}
\def\mxsh{1.5}
\def\mysh{1.2}
\begin{tikzpicture}[menvone]
\node at (-12,0) {(I):};
\node (l) at (-1.8,0) {};
\node (r) at (1.8,0) {};
\path[commutative diagrams/.cd, every arrow, every label]
(l) edge[commutative diagrams/squiggly] (r);
\begin{scope}[xshift=-6cm]
\draw [thkln] (-\mxsh+0.15,\mysh) to [out=0,in=180] (\mxsh-0.15,-\mysh);
\draw [lnovr] (-\mxsh+0.15,-\mysh) to [out=0,in=180] (\mxsh-0.15,\mysh);
\draw [thkln] (-\mxsh+0.15,-\mysh) to [out=0,in=180] (\mxsh-0.15,\mysh);
\draw [bcrc] (-0.75-0.6-\mxsh,-\mysh) -- (0.75+0.6+\mxsh,-\mysh);
\draw [bcrc] (-0.75-0.6-\mxsh,\mysh) -- (0.75+0.6+\mxsh,\mysh);
%
%
\draw [ptzert] (-0.15-\mxsh,-0.6-\mysh) rectangle ++(.3,1.2);
\draw [thkln] (-0.75-0.6-\mxsh,-\mysh) -- (-0.15-\mxsh,-\mysh)
node [very near start,below] {$\scriptstyle \xca+1$};
\begin{scope}[xscale=-1]
\draw [ptzert] (-0.15-\mxsh,-0.6-\mysh) rectangle ++(.3,1.2);
\draw [thkln] (-0.75-0.6-\mxsh,-\mysh) -- (-0.15-\mxsh,-\mysh)
node [very near start,below] {$\scriptstyle \xca+1$};
\end{scope}
\begin{scope}[yscale=-1]
\draw [ptzert] (-0.15-\mxsh,-0.6-\mysh) rectangle ++(.3,1.2);
\draw [thkln] (-0.75-0.6-\mxsh,-\mysh) -- (-0.15-\mxsh,-\mysh)
node [very near start,above] {$\scriptstyle \xca+1$};
\end{scope}
\begin{scope}[yscale=-1,xscale=-1]
\draw [ptzert] (-0.15-\mxsh,-0.6-\mysh) rectangle ++(.3,1.2);
\draw [thkln] (-0.75-0.6-\mxsh,-\mysh) -- (-0.15-\mxsh,-\mysh)
node [near start,above] {$\scriptstyle \xca+1$};
\end{scope}
\end{scope}
\begin{scope}[xshift=6cm]
\draw [thkln] (-\mxsh+0.15,\mysh) to [out=0,in=180]
node [sloped, near end,above] {$\scriptstyle \xca$} (\mxsh-0.15,-\mysh);
\draw [lnovr] (-\mxsh+0.15,-\mysh) to [out=0,in=180] (\mxsh-0.15,\mysh);
\draw [thkln] (-\mxsh+0.15,-\mysh) to [out=0,in=180]
node [sloped, near start,above] {$\scriptstyle \xca$} (\mxsh-0.15,\mysh);
\draw (0.15-\mxsh,\mysh+0.35) -- (-0.15+\mxsh,\mysh+0.35);
\draw (0.15-\mxsh,-\mysh-0.35) -- (-0.15+\mxsh,-\mysh-0.35);
\draw [bcrc] (-0.75-0.6-\mxsh,-\mysh) -- (0.75+0.6+\mxsh,-\mysh);
\draw [bcrc] (-0.75-0.6-\mxsh,\mysh) -- (0.75+0.6+\mxsh,\mysh);
%
%
\draw [ptzert] (-0.15-\mxsh,-0.6-\mysh) rectangle ++(.3,1.2);
\draw [thkln] (-0.75-0.6-\mxsh,-\mysh) -- (-0.15-\mxsh,-\mysh)
node [very near start,below] {$\scriptstyle \xca+1$};
\begin{scope}[xscale=-1]
\draw [ptzert] (-0.15-\mxsh,-0.6-\mysh) rectangle ++(.3,1.2);
\draw [thkln] (-0.75-0.6-\mxsh,-\mysh) -- (-0.15-\mxsh,-\mysh)
node [very near start,below] {$\scriptstyle \xca+1$};
\end{scope}
\begin{scope}[yscale=-1]
\draw [ptzert] (-0.15-\mxsh,-0.6-\mysh) rectangle ++(.3,1.2);
\draw [thkln] (-0.75-0.6-\mxsh,-\mysh) -- (-0.15-\mxsh,-\mysh)
node [very near start,above] {$\scriptstyle \xca+1$};
\end{scope}
\begin{scope}[yscale=-1,xscale=-1]
\draw [ptzert] (-0.15-\mxsh,-0.6-\mysh) rectangle ++(.3,1.2);
\draw [thkln] (-0.75-0.6-\mxsh,-\mysh) -- (-0.15-\mxsh,-\mysh)
node [near start,above] {$\scriptstyle \xca+1$};
\end{scope}
\end{scope}
\end{tikzpicture}
\\
\label{eq:thtr}
\\
\nonumber
\begin{tikzpicture}[menvtwo]
\node at (-4,0) {(II):};
\node (l) at (-0.8,0) {};
\node (r) at (0.8,0) {};
\path[commutative diagrams/.cd, every arrow, every label]
(l) edge[commutative diagrams/squiggly] (r);
\begin{scope}[xshift=-2cm]
\draw [ptzer] (-0.15,-0.6) rectangle (.15,.6);
\draw [thkln] (-0.75,0) -- (-0.15,0)
node [near start,below] {$\scriptstyle \xca$}
(0.15,0) -- (.75,0)
node [near end,below] {$\scriptstyle \xca$};
\draw [bcrc] (-0.75,0) -- (-0.15,0)  (0.15,0) -- (0.75,0);
\draw (-0.75,0.4) -- (-0.15,0.4) (0.15,0.4) -- (0.75,0.4);
\end{scope}
\begin{scope}[xshift=2cm]
\draw [ptzer] (-0.15,-0.6) rectangle (.15,.6);
\draw [thkln] (-0.75,0) -- (-0.15,0)
node [near start,below] {$\scriptstyle \xca$}
(0.15,0) -- (.75,0)
node [near end,below] {$\scriptstyle \xca$};
\draw [bcrc] (-0.75,0) -- (-0.15,0)  (0.15,0) -- (0.75,0);
\draw (-0.75,0.8) -- (0.75,0.8);
\end{scope}
\end{tikzpicture}
\qquad
\begin{tikzpicture}[menvtwo]
\node at (-4,0) {(II):};
\node (l) at (-0.8,0) {};
\node (r) at (0.8,0) {};
\path[commutative diagrams/.cd, every arrow, every label]
(l) edge[commutative diagrams/squiggly] (r);
\begin{scope}[xshift=-2cm]
\draw [ptzer] (-0.15,-0.6) rectangle (.15,.6);
\draw [thkln] (-0.75,0) -- (-0.15,0)
node [near start,below] {$\scriptstyle \xca$}
(0.15,0) -- (.75,0)
node [near end,below] {$\scriptstyle \xca$};
\draw [bcrc] (-0.75,0) -- (-0.15,0)  (0.15,0) -- (0.75,0);
\draw (-0.15,0.4) to [out=180,in=-90] (-0.6,0.8) to [out=90,in=180] (0,1.4)
to [out=0,in=90] (0.6,0.8) to [out=-90,in=0] (0.15,0.4);
\end{scope}
\begin{scope}[xshift=2cm]
\draw [ptzer] (-0.15,-0.6) rectangle (.15,.6);
\draw [thkln] (-0.75,0) -- (-0.15,0)
node [near start,below] {$\scriptstyle \xca$}
(0.15,0) -- (.75,0)
node [near end,below] {$\scriptstyle \xca$};
\draw [bcrc] (-0.75,0) -- (-0.15,0)  (0.15,0) -- (0.75,0);
\end{scope}
\end{tikzpicture}
\;.
\end{gather}
The thick gray lines in these pictures mark the \tBcr s of the diagram $\sBD$.

The transition from the diagram $\xDclNo$ to $\xDclN$ is performed in two stages. At the first stage we apply the first replacement of\rx{eq:thtr} to every crossing of $\xDclNo$. The result is the diagram $\xtDN$, which consists of two parts connected at \tJWp s. The first part is the $\xca$-cabled diagram $\xDclN$ and the second part consists of non-intersecting circles formed by single lines appearing in the final diagrams of replacements (I) of\rx{eq:thtr}. These single line circles go along the \tBcr s. We orient them clockwise and assume that in our pictures the clockwise orientation corresponds to the direction from the left to the right. The circles are attached to $\xDclN$ at the \tJWp s and those junctions have four possible forms:
\begin{equation}
\label{eq:frmprj}
\begin{tikzpicture}[menvtwo]
\draw [ptzer] (-0.15,-0.6) rectangle (.15,.6);
\draw [thkln] (-0.75,0) -- (-0.15,0)
node [near start,below] {$\scriptstyle \xca$}
(0.15,0) -- (.75,0)
node [near end,below] {$\scriptstyle \xca$};
\draw [bcrc] (-0.75,0) -- (-0.15,0)  (0.15,0) -- (0.75,0);
\draw (-0.75,0.4) -- (-0.15,0.4) (0.15,0.4) -- (0.75,0.4);
\end{tikzpicture},
\qquad
\begin{tikzpicture}[menvtwo]
\draw [ptzer] (-0.15,-0.6) rectangle (.15,.6);
\draw [thkln] (-0.75,0) -- (-0.15,0)
node [near start,above] {$\scriptstyle \xca$}
(0.15,0) -- (.75,0)
node [near end,above] {$\scriptstyle \xca$};
\draw [bcrc] (-0.75,0) -- (-0.15,0)  (0.15,0) -- (0.75,0);
\draw (-0.75,-0.4) -- (-0.15,-0.4) (0.15,-0.4) -- (0.75,-0.4);
\end{tikzpicture},
\qquad
\begin{tikzpicture}[menvtwo]
\draw [ptzer] (-0.15,-0.6) rectangle (.15,.6);
\draw [thkln] (-0.75,0) -- (-0.15,0)
node [near start,above] {$\scriptstyle \xca$}
(0.15,0) -- (.75,0)
node [near end,below] {$\scriptstyle \xca$};
\draw [bcrc] (-0.75,0) -- (-0.15,0)  (0.15,0) -- (0.75,0);
\draw (-0.75,-0.4) -- (-0.15,-0.4) (0.15,0.4) -- (0.75,0.4);
\end{tikzpicture},
\qquad
\begin{tikzpicture}[menvtwo]
\draw [thkln] (-0.75,0) -- (-0.15,0)
node [near start,below] {$\scriptstyle \xca$}
(0.15,0) -- (.75,0)
node [near end,above] {$\scriptstyle \xca$};
\draw [bcrc] (-0.75,0) -- (-0.15,0)  (0.15,0) -- (0.75,0);
\draw (-0.75,0.4) -- (-0.15,0.4) (0.15,-0.4) -- (0.75,-0.4);
\draw [ptzer] (-0.15,-0.6) rectangle (.15,.6);
\end{tikzpicture}.
\end{equation}

At the second stage of the transition from $\xDclNo$ to $\xDclN$ we remove the single circle lines
of $\xtDN$ one-by-one. In order to remove a particular circle we select an `initial' \tJWp\ on it and then detach the single lines from other projectors going clockwise. During this process, the single line between the initial and current \tJWp s are kept on the same side of the \tBcr s. If the current projector has the incoming and outgoing single lines on the opposite sides of the \tBcr\ (third and fourth type of\rx{eq:frmprj}) then, prior to detachment, we perform the following transformation
for the junction of the third type (and a similar transformation for the fourth type):
\begin{equation}
\label{eq:twflp}
\begin{tikzpicture}[menvtwo]
\begin{scope}[xshift=1.75cm]
\draw [thkln] (-1.15,0) -- (-0.15,0)
(0.15,0) -- (.75,0)
node [near end,below] {$\scriptstyle \xca$};
\draw [bcrc] (-4.25,0) -- (0.75,0);
\draw 
(0.15,0.4) -- (0.75,0.4);
\draw [ptzert] (-0.15,-0.6) rectangle (.15,.6);
\end{scope}
\begin{scope}[xshift=-1.75cm,xscale=-1]
\draw [thkln] (-1.15,0) -- (-0.15,0)
(0.15,0) -- (.75,0)
node [near end,above] {$\scriptstyle \xca$};
\draw 
(0.15,-0.4) -- (0.75,-0.4);
\draw [ptzert] (-0.15,-0.6) rectangle (.15,.6);
\end{scope}
\node (0,0) {$\cdots$};
\draw (-1.6,-0.4) -- (1.6,-0.4);
\end{tikzpicture}
\;\hteqv\;
\begin{tikzpicture}[menvtwo]
\begin{scope}[xshift=1.75cm]
\draw [bcrc] (-4.25,0) -- (0.75,0);
\draw [thkln] (-1.15,0) -- (-0.15,0)
(0.15,0) -- (.75,0)
node [near end,below] {$\scriptstyle \xca$};
\draw [lnovr]  (-1.05,0.4) to [out=0,in=180] (-0.15,-0.4);
\draw (-1.05,0.4) to [out=0,in=180] (-0.15,-0.4) (0.15,0.4) -- (0.75,0.4);
\draw [ptzert] (-0.15,-0.6) rectangle (.15,.6);
\end{scope}
\begin{scope}[xshift=-1.75cm,xscale=-1]
\draw [thkln] (-1.15,0) -- (-0.15,0)
(0.15,0) -- (.75,0)
node [near end,above] {$\scriptstyle \xca$};
\draw [lnovr]  (-1.05,0.4) to [out=0,in=180] (-0.15,-0.4);
\draw (-1.05,0.4) to [out=0,in=180] (-0.15,-0.4) (0.15,-0.4) -- (0.75,-0.4);
\draw [ptzert] (-0.15,-0.6) rectangle (.15,.6);
\end{scope}
\node (0,0) {$\cdots$};
\draw (-0.7,0.4) -- (0.7,0.4);
\end{tikzpicture}
\;\hteqv\;
\begin{tikzpicture}[menvtwo]
\begin{scope}[xshift=1.75cm]
\draw [bcrc] (-4.25,0) -- (0.75,0);
\draw [thkln] (-1.15,0) -- (-0.15,0)
(0.15,0) -- (.75,0)
node [near end,below] {$\scriptstyle \xca$};
\draw 
(0.15,0.4) -- (0.75,0.4);
\draw [ptzert] (-0.15,-0.6) rectangle (.15,.6);
\end{scope}
\begin{scope}[xshift=-1.75cm,xscale=-1]
\draw [thkln] (-1.15,0) -- (-0.15,0)
(0.15,0) -- (.75,0)
node [near end,above] {$\scriptstyle \xca$};
\draw 
(0.15,-0.4) -- (0.75,-0.4);
\draw [ptzert] (-0.15,-0.6) rectangle (.15,.6);
\end{scope}
\node (0,0) {$\cdots$};
\draw (-1.6,0.4) -- (1.6,0.4);
\end{tikzpicture},
\end{equation}
In these pictures the left projector is initial, the right projector is current, the first homotopy equivalence comes from the Reidemeister moves, while the second equivalence comes from \ex{eq:projtw}. Note that the single line between the initial and current projectors is kept always above the rest of the diagram.

After the single lines attached to the current projector are brought to the same side of the \tBcr, we detach the single line from that projector by the \lrpl\ (II) of\rx{eq:thtr}
and pass to the next projector on the single line.

The single line is kept above the rest of the diagram, so once it is detached from all projectors except the initial one, it can be contracted to a small loop attached to that initial projector with the help of Reidemeister moves. The final step is the removal of that loop by the replacement (III) of\rx{eq:thtr}. After all single line circles are removed, the diagram $\xtDN$ becomes $\xDclN$.

Our transition from $\xDclNo$ to $\xDclN$ is generally similar to that used by C.~Armond\cite{Arm11}, but the details are different. In particular, we do not replace $(\xca+1)$-cable crossings by projectors, but rather apply replacements (I) of\rx{eq:thtr} directly to the crossings.

\subsection{\Ltrf s generate isomorphisms at low \thdgr s}
\label{sct:hest}

We describe the \ltrf s related to replacements\rx{eq:thtr} and show that the corresponding maps\rx{eq:dgprmg} between shifted homologies are isomorphisms at low \thdgr s, thus proving Theorem\rw{thm:leviso}.

\subsubsection{\Ltrf\ \xltone}
Set
\begin{equation}
\label{eq:trone}
\def\mxsh{1.5}
\def\mysh{1.2}
\ytngsi =
\begin{tikzpicture}[menvone]
\draw [thkln] (-\mxsh+0.15,\mysh) to [out=0,in=180] (\mxsh-0.15,-\mysh);
\draw [lnovr] (-\mxsh+0.15,-\mysh) to [out=0,in=180] (\mxsh-0.15,\mysh);
\draw [thkln] (-\mxsh+0.15,-\mysh) to [out=0,in=180] (\mxsh-0.15,\mysh);
\draw [bcrc] (-0.75-0.6-\mxsh,-\mysh) -- (0.75+0.6+\mxsh,-\mysh);
\draw [bcrc] (-0.75-0.6-\mxsh,\mysh) -- (0.75+0.6+\mxsh,\mysh);
%
%
\draw [ptzert] (-0.15-\mxsh,-0.6-\mysh) rectangle ++(.3,1.2);
\draw [thkln] (-0.75-0.6-\mxsh,-\mysh) -- (-0.15-\mxsh,-\mysh)
node [very near start,below] {$\scriptstyle \xca+1$};
\begin{scope}[xscale=-1]
\draw [ptzert] (-0.15-\mxsh,-0.6-\mysh) rectangle ++(.3,1.2);
\draw [thkln] (-0.75-0.6-\mxsh,-\mysh) -- (-0.15-\mxsh,-\mysh)
node [very near start,below] {$\scriptstyle \xca+1$};
\end{scope}
\begin{scope}[yscale=-1]
\draw [ptzert] (-0.15-\mxsh,-0.6-\mysh) rectangle ++(.3,1.2);
\draw [thkln] (-0.75-0.6-\mxsh,-\mysh) -- (-0.15-\mxsh,-\mysh)
node [very near start,above] {$\scriptstyle \xca+1$};
\end{scope}
\begin{scope}[yscale=-1,xscale=-1]
\draw [ptzert] (-0.15-\mxsh,-0.6-\mysh) rectangle ++(.3,1.2);
\draw [thkln] (-0.75-0.6-\mxsh,-\mysh) -- (-0.15-\mxsh,-\mysh)
node [near start,above] {$\scriptstyle \xca+1$};
\end{scope}
\end{tikzpicture}\;,
\quad
\ytngsf=
\begin{tikzpicture}[menvone]
\draw [thkln] (-\mxsh+0.15,\mysh) to [out=0,in=180]
node [sloped, near end,above] {$\scriptstyle \xca$} (\mxsh-0.15,-\mysh);
\draw [lnovr] (-\mxsh+0.15,-\mysh) to [out=0,in=180] (\mxsh-0.15,\mysh);
\draw [thkln] (-\mxsh+0.15,-\mysh) to [out=0,in=180]
node [sloped, near start,above] {$\scriptstyle \xca$} (\mxsh-0.15,\mysh);
\draw (0.15-\mxsh,\mysh+0.35) -- (-0.15+\mxsh,\mysh+0.35);
\draw (0.15-\mxsh,-\mysh-0.35) -- (-0.15+\mxsh,-\mysh-0.35);
\draw [bcrc] (-0.75-0.6-\mxsh,-\mysh) -- (0.75+0.6+\mxsh,-\mysh);
\draw [bcrc] (-0.75-0.6-\mxsh,\mysh) -- (0.75+0.6+\mxsh,\mysh);
%
%
\draw [ptzert] (-0.15-\mxsh,-0.6-\mysh) rectangle ++(.3,1.2);
\draw [thkln] (-0.75-0.6-\mxsh,-\mysh) -- (-0.15-\mxsh,-\mysh)
node [very near start,below] {$\scriptstyle \xca+1$};
\begin{scope}[xscale=-1]
\draw [ptzert] (-0.15-\mxsh,-0.6-\mysh) rectangle ++(.3,1.2);
\draw [thkln] (-0.75-0.6-\mxsh,-\mysh) -- (-0.15-\mxsh,-\mysh)
node [very near start,below] {$\scriptstyle \xca+1$};
\end{scope}
\begin{scope}[yscale=-1]
\draw [ptzert] (-0.15-\mxsh,-0.6-\mysh) rectangle ++(.3,1.2);
\draw [thkln] (-0.75-0.6-\mxsh,-\mysh) -- (-0.15-\mxsh,-\mysh)
node [very near start,above] {$\scriptstyle \xca+1$};
\end{scope}
\begin{scope}[yscale=-1,xscale=-1]
\draw [ptzert] (-0.15-\mxsh,-0.6-\mysh) rectangle ++(.3,1.2);
\draw [thkln] (-0.75-0.6-\mxsh,-\mysh) -- (-0.15-\mxsh,-\mysh)
node [near start,above] {$\scriptstyle \xca+1$};
\end{scope}
\end{tikzpicture}\;,
\quad
\ytngsc= \shcr^{\xca+\hlf}
\begin{tikzpicture}[menvone]
\draw [thkln] (-\mxsh+0.15,\mysh) to [out=0,in=180]
node [sloped, very near end,below] {$\scriptstyle \xca$} (\mxsh-0.15,-\mysh);
\draw [lnovr] (-\mxsh+0.15,-\mysh) to [out=0,in=180] (\mxsh-0.15,\mysh);
\draw [thkln] (-\mxsh+0.15,-\mysh) to [out=0,in=180]
node [sloped, very near start,below] {$\scriptstyle \xca$} (\mxsh-0.15,\mysh);
\draw (0.15-\mxsh,\mysh-0.35) to [out=0,in=90] (-\mxsh*0.5,0) to [out=-90,in=0] (0.15-\mxsh,-\mysh+0.35);
\draw [xscale=-1] (0.15-\mxsh,\mysh-0.35) to [out=0,in=90] (-\mxsh*0.5,0) to [out=-90,in=0] (0.15-\mxsh,-\mysh+0.35);
\draw [bcrc] (-0.75-0.6-\mxsh,-\mysh) -- (0.75+0.6+\mxsh,-\mysh);
\draw [bcrc] (-0.75-0.6-\mxsh,\mysh) -- (0.75+0.6+\mxsh,\mysh);
%
%
\draw [ptzert] (-0.15-\mxsh,-0.6-\mysh) rectangle ++(.3,1.2);
\draw [thkln] (-0.75-0.6-\mxsh,-\mysh) -- (-0.15-\mxsh,-\mysh)
node [very near start,below] {$\scriptstyle \xca+1$};
\begin{scope}[xscale=-1]
\draw [ptzert] (-0.15-\mxsh,-0.6-\mysh) rectangle ++(.3,1.2);
\draw [thkln] (-0.75-0.6-\mxsh,-\mysh) -- (-0.15-\mxsh,-\mysh)
node [very near start,below] {$\scriptstyle \xca+1$};
\end{scope}
\begin{scope}[yscale=-1]
\draw [ptzert] (-0.15-\mxsh,-0.6-\mysh) rectangle ++(.3,1.2);
\draw [thkln] (-0.75-0.6-\mxsh,-\mysh) -- (-0.15-\mxsh,-\mysh)
node [very near start,above] {$\scriptstyle \xca+1$};
\end{scope}
\begin{scope}[yscale=-1,xscale=-1]
\draw [ptzert] (-0.15-\mxsh,-0.6-\mysh) rectangle ++(.3,1.2);
\draw [thkln] (-0.75-0.6-\mxsh,-\mysh) -- (-0.15-\mxsh,-\mysh)
node [near start,above] {$\scriptstyle \xca+1$};
\end{scope}
\end{tikzpicture}\;,
\end{equation}
while $\ytngkfp=\shcr^{-\xca+\hlf}\ytngkf$. Theorem\rw{lm:sss} provides the exact triangle relation\rx{eq:cnrel}.

\begin{proposition}
\label{prp:bndfs}
Let $\xDsi$ be the diagram constructed by performing \lrpl s I of\rx{eq:thtr} on some vertices of $\xDclNo$ and let $\xDsf$ be the diagram constructed by performing the \lrpl\ I on the `current' vertex in $\xDsi$. Then the \tdgpr\ map\rx{eq:dgprmg} is an isomorphism on $\tKHmvv{i}{\hem}$ for $i\leq 2\xca-1$.
\end{proposition}
\begin{proof}
Let $\xDsc$ be the diagram constructed by performing the \lrpl\ $\ytngsi \rightsquigarrow\ytngsc$ on the current vertex. We estimate the homological order of $\KHm(\xDsc)$ with the help of Theorem\rw{thm:smfr}: since $\yncrv{\xDsc} = \yncrv{\xDsi}-2\xca-1$, then $\KHmvv{i}{\hem}(\xDsc)=0$ for $i\leq -\shlf\yncrv{\xDsi} + 2\xca$ (we took into account the shift $\shcr^{\xca+\hlf}$ of $\ytngsc$ in \ex{eq:trone}) and the claim of the theorem follows from Proposition\rw{prp:gestdg}.
\end{proof}
%
%

\subsubsection{\Ltrf\ \xlttwo}
Set
\begin{equation}
\label{eq:lrdpr}
\ytngsi =
\begin{tikzpicture}[menvtwo]
\draw [bcrc] (-0.75,0) -- (0.75,0);
\draw [ptzert] (-0.15,-0.6) rectangle (.15,.6);
\draw [thkln] (-0.75,0) -- (-0.15,0)
node [near start,below] {$\scriptstyle \xca$}
(0.15,0) -- (.75,0)
node [near end,below] {$\scriptstyle \xca$};
\draw (-0.75,0.4) -- (-0.15,0.4) (0.15,0.4) -- (0.75,0.4);
\end{tikzpicture},
\qquad
\ytngsf =
\begin{tikzpicture}[menvtwo]
\draw [bcrc] (-0.75,0) -- (0.75,0);
\draw [ptzert] (-0.15,-0.6) rectangle (.15,.6);
\draw [thkln] (-0.75,0) -- (-0.15,0)
node [near start,below] {$\scriptstyle \xca$}
(0.15,0) -- (.75,0)
node [near end,below] {$\scriptstyle \xca$};
\draw (-0.75,0.8) -- (0.75,0.8);
\end{tikzpicture},
\qquad
\ytngsc =
\shcr
\begin{tikzpicture}[menvone]
\draw [bcrc] (-1.35,0) -- (1.35,0);
\draw[pttwo] (-.15,-0.6) rectangle ++(0.3,1.2);
\draw[thkln] (-1.35,0) -- (-.15,0)
node[near start,below] {$\scriptstyle \xca+1$ }
(.15,0) -- (1.35,0)
node [near end, below] {$\scriptstyle \xca+1$};
\end{tikzpicture}\;,
\end{equation}
while $\ytngkfp=\ytngkf$. The exact triangle relation\rx{eq:cnrel} is provided by Theorem\rw{thm:projpr}.
\begin{proposition}
\label{prp:bndsc}
Let $\xDsi$ be a diagram constructed by removing some single line circles from $\xtDN$ and by detaching the `current' single line circle from the projectors which lie between the initial one and the current one and let $\xDsf$ be the diagram constructed from $\xDsi$ by detaching the single line from the current projector.
Then the \tdgpr\ map\rx{eq:dgprmg} is an isomorphism on $\tKHmvv{i}{\hem}$ for $i\leq \xca-1$.
\end{proposition}
The proof uses the following
\begin{lemma}
\label{lm:cmbnd}
The tangle
\def\mxsh{2}
\def\mysh{1.5}
\begin{equation}
\label{eq:tngcr}
\ztau \; = \;
\begin{tikzpicture}[menvone]
\draw [thkln] (-\mxsh+0.15,\mysh) to [out=0,in=180] (\mxsh-0.15,-\mysh);
\draw [lnovr] (-\mxsh+0.15,-\mysh) to [out=0,in=180] (\mxsh-0.15,\mysh);
\draw [thkln] (-\mxsh+0.15,-\mysh) to [out=0,in=180] (\mxsh-0.15,\mysh);
\draw [lnovr] (-0.75-\mxsh,0) to [out=0,in=180] (0,-\mysh) to [out=0,in=180] (0.75+\mxsh,0);
\draw (-0.75-\mxsh,0) to [out=0,in=180] (0,-\mysh) to [out=0,in=180] (0.75+\mxsh,0);
\draw [ptzer] (-0.15-\mxsh,-0.6-\mysh) rectangle ++(.3,1.2);
\draw [thkln] (-0.75-\mxsh,-\mysh) -- (-0.15-\mxsh,-\mysh)
node [near start,below] {$\scriptstyle \xca$};
\begin{scope}[xscale=-1]
\draw [ptzer] (-0.15-\mxsh,-0.6-\mysh) rectangle ++(.3,1.2);
\draw [thkln] (-0.75-\mxsh,-\mysh) -- (-0.15-\mxsh,-\mysh)
node [near start,below] {$\scriptstyle \xca$};
\end{scope}
\begin{scope}[yscale=-1]
\draw [ptzer] (-0.15-\mxsh,-0.6-\mysh) rectangle ++(.3,1.2);
\draw [thkln] (-0.75-\mxsh,-\mysh) -- (-0.15-\mxsh,-\mysh)
node [near start,above] {$\scriptstyle \xca$};
\end{scope}
\begin{scope}[yscale=-1,xscale=-1]
\draw [ptzer] (-0.15-\mxsh,-0.6-\mysh) rectangle ++(.3,1.2);
\draw [thkln] (-0.75-\mxsh,-\mysh) -- (-0.15-\mxsh,-\mysh)
node [near start,above] {$\scriptstyle \xca$};
\end{scope}
\end{tikzpicture}
\end{equation}
has a homological bound
$\hmord{\ztau} \geq -\shlf\xca^2.$
\end{lemma}
\begin{remark}
\label{rmk:ignr}
This bound is better than the crude bound of Theorem\rw{thm:smfr}. In fact, it coincides with that bound, if we neglect the intersections between the single line and the $\xca$-cables.
\end{remark}
\begin{proof}[Proof of Lemma\rw{lm:cmbnd}]
Applying \ex{eq:colKhbr} to the $\xca$-cable crossing in $\ztau$ we get the presentation
\[
\def\mxsh{2}
\def\mysh{1.5}
\xKhv{\ztau} \hteqv
\shcr^{-\hlf\xca^2}\,
\Pcnv{
\bigoplus_{i=0}^{\xca}
\shcr^{i^2}{\xca \brace i}_{\shcr}\,
\ztau_i}
\qquad
\xKhv{\ztau_i}\;=\;
\begin{tikzpicture}[menvone]
\draw [thkln] (-\mxsh+0.15,\mysh-0.2) to [out=0,in=90] (-\mxsh+1.1,0)
 to [out=-90,in=0] (-\mxsh+0.15,-\mysh+0.2) ;
\begin{scope}[xscale=-1]
\draw [thkln] (-\mxsh+0.15,\mysh-0.2) to [out=0,in=90] (-\mxsh+1.1,0) to [out=-90,in=0] (-\mxsh+0.15,-\mysh+0.2);
\end{scope}
\node at (-\mxsh+1.5,0.4) {$\scriptstyle i$};
\node at (\mxsh-1.5,0.4) {$\scriptstyle i$};
\draw [thkln] (-\mxsh+0.15,\mysh+0.2) -- (\mxsh-0.15,\mysh+0.2) node [midway,above] {$\scriptstyle \xca-i$};
\draw [thkln] (-\mxsh+0.15,-\mysh-0.2) -- (\mxsh-0.15,-\mysh-0.2) node [midway, below] {$\scriptstyle \xca-i$};
\draw [lnovr] (-0.75-\mxsh,0) to [out=0,in=180] (0,-\mysh+0.7) to [out=0,in=180] (0.75+\mxsh,0);
\draw (-0.75-\mxsh,0) to [out=0,in=180] (0,-\mysh+0.7) to [out=0,in=180] (0.75+\mxsh,0);
\draw [ptzer] (-0.15-\mxsh,-0.6-\mysh) rectangle ++(.3,1.2);
\draw [thkln] (-0.75-\mxsh,-\mysh) -- (-0.15-\mxsh,-\mysh)
node [near start,below] {$\scriptstyle \xca$};
\begin{scope}[xscale=-1]
\draw [ptzer] (-0.15-\mxsh,-0.6-\mysh) rectangle ++(.3,1.2);
\draw [thkln] (-0.75-\mxsh,-\mysh) -- (-0.15-\mxsh,-\mysh)
node [near start,below] {$\scriptstyle \xca$};
\end{scope}
\begin{scope}[yscale=-1]
\draw [ptzer] (-0.15-\mxsh,-0.6-\mysh) rectangle ++(.3,1.2);
\draw [thkln] (-0.75-\mxsh,-\mysh) -- (-0.15-\mxsh,-\mysh)
node [near start,above] {$\scriptstyle \xca$};
\end{scope}
\begin{scope}[yscale=-1,xscale=-1]
\draw [ptzer] (-0.15-\mxsh,-0.6-\mysh) rectangle ++(.3,1.2);
\draw [thkln] (-0.75-\mxsh,-\mysh) -- (-0.15-\mxsh,-\mysh)
node [near start,above] {$\scriptstyle \xca$};
\end{scope}
\end{tikzpicture}
\]
The homological order of $\ztau_i$ can be estimated with the help of Theorem\rw{thm:smfr}: $\hmord{\xKhv{\ztau_i}} \geq - i$. Since the polynomial ${\xca \brace i}_{\shcr}$ has only non-negative powers of $\shcr$ and $i^2 - i \geq 0$ for all integer $i$, we come to the estimate of Lemma\rw{lm:cmbnd}.
\end{proof}

\begin{proof}[Proof of Proposition\rw{prp:bndsc}]
Let $\xDsc$ be the diagram constructed from $\xDsi$ by replacing the current projector ($\ytngsi$ of \ex{eq:lrdpr}) with the tangle complex $\ytngsc$ of \ex{eq:lrdpr}.
By Proposition\rw{prp:gestdg}, we have to prove the bound:
\[
\KHmvv{i}{\hem}(\xDsc)=0 \quad\text{for $i\leq-\shlf\yncrv{\xDsi}+\xca.$}
\]
estimate
\[\hmord{\xDsc}\geq-\hlf\yncrv{\xDsi} + \xca.\]
Since the complex $\ytngsc$ of \ex{eq:lrdpr} is a \tmcn\rx{eq:swgrpr} generated by an `elementary' tangle
\begin{equation}
\label{eq:eltng}
\ytngse \;=\;
\begin{tikzpicture}[menvone]
\draw [line width=6pt, color=white] (-1+0.15,0) -- (0,0);
\draw [line width=\cblth] (-1.75,0) -- (-1 - .15,0)
node [near start, below] {$\scriptstyle \xca$}
(-1+0.15,0) -- (0,0) -- (1-0.15,0) (1+.15,0) -- (1.75,0)
node [near end, below] {$\scriptstyle \xca$};
\draw (1, 0.9) to [out=180,in=90] (0.4,0.6) to [out=-90,in=180] (1-.15,0.3);
\draw [xscale=-1] (1, 0.9) to [out=180,in=90] (0.4,0.6) to [out=-90,in=180] (1-.15,0.3);
\draw (-1.75,0.9) -- (-1,0.9) (1,0.9) -- (1.75,0.9);
\draw [line width=\ljwp] (-1.15,-0.6) rectangle ++(0.3,1.2);
\draw [line width=\ljwp] (1 - 0.15,-0.6) rectangle ++(0.3,1.2);
\end{tikzpicture}
\end{equation}
then, according to Remark\rw{rmk:bndss},
it is sufficient to prove
\begin{equation}
\label{eq:scest}
\KHmvv{i}{\hem}(\xDse) = 0\quad\text{for $i\leq -\shlf\yncrv{\xDsi} + \xca-1$,}
\end{equation}
where $\xDse$ is the diagram constructed by replacing $\ytngsi$ in $\xDsi$ with $\ytngse$.

Consider a tangle within $\xDse$ which consists of the right half of $\ytngse$ and the cable crossing which follows the current projector and transform its complex with the help of two homotopy equivalences:
\begin{equation}
\label{eq:prsld}
\def\mxsh{2}
\def\mysh{1.5}
\begin{tikzpicture}[menvone]
\draw [thkln] (-\mxsh+0.15,\mysh) to [out=0,in=180] (\mxsh-0.15,-\mysh);
\draw [lnovr] (-\mxsh+0.15,-\mysh) to [out=0,in=180] (\mxsh-0.15,\mysh);
\draw [thkln] (-\mxsh+0.15,-\mysh) to [out=0,in=180] (\mxsh-0.15,\mysh);
\draw (\mxsh-0.15,\mysh+0.3) to [out=180,in=0]
(-\mxsh, \mysh+0.9) to [out=180,in=90] (-\mxsh-0.6,\mysh+0.6) to [out=-90,in=180] (-\mxsh-.15,\mysh+0.3);
\draw (\mxsh+0.15,\mysh+0.3) -- (\mxsh+0.75,\mysh+0.3);
\draw [dashed] (-\mxsh-0.75,-\mysh-0.3) -- (-\mxsh-0.15,-\mysh-0.3)
(-\mxsh+0.15,-\mysh-0.3) -- (\mxsh-0.15,-\mysh-0.3) (\mxsh+0.15,-\mysh-0.3) -- (\mxsh+0.75,-\mysh-0.3);
\draw [ptzert] (-0.15-\mxsh,-0.6-\mysh) rectangle ++(.3,1.2);
\draw [thkln] (-0.75-\mxsh,-\mysh) -- (-0.15-\mxsh,-\mysh)
node [near start,above] {$\scriptstyle \xca$};
\begin{scope}[xscale=-1]
\draw [ptzert] (-0.15-\mxsh,-0.6-\mysh) rectangle ++(.3,1.2);
\draw [thkln] (-0.75-\mxsh,-\mysh) -- (-0.15-\mxsh,-\mysh)
node [near start,above] {$\scriptstyle \xca$};
\end{scope}
\begin{scope}[yscale=-1]
\draw [ptzert] (-0.15-\mxsh,-0.6-\mysh) rectangle ++(.3,1.2);
\draw [thkln] (-0.75-\mxsh,-\mysh) -- (-0.15-\mxsh,-\mysh)
node [at start,below] {$\scriptstyle \xca-1\;\;\;$};
\end{scope}
\begin{scope}[yscale=-1,xscale=-1]
\draw [ptzert] (-0.15-\mxsh,-0.6-\mysh) rectangle ++(.3,1.2);
\draw [thkln] (-0.75-\mxsh,-\mysh) -- (-0.15-\mxsh,-\mysh)
node [near start,below] {$\scriptstyle \xca$};
\end{scope}
\end{tikzpicture}
\;\hteqv\;
\begin{tikzpicture}[menvone]
\draw (\mxsh+0.15,\mysh+0.3) -- (\mxsh+0.75,\mysh+0.3);
\draw [dashed] (-\mxsh-0.75,-\mysh-0.3) -- (-\mxsh-0.15,-\mysh-0.3)
(-\mxsh+0.15,-\mysh-0.3) -- (\mxsh-0.15,-\mysh-0.3) (\mxsh+0.15,-\mysh-0.3) -- (\mxsh+0.75,-\mysh-0.3);
%
\draw (\mxsh-0.15,\mysh+0.3) to [out=180,in=90]
(\mxsh-1.2,0) to [out=-90,in=180] (\mxsh-0.15,-\mysh+0.3);
\draw [thkln] (-\mxsh+0.15,\mysh) to [out=0,in=180] (\mxsh-0.15,-\mysh);
\draw [lnovr] (-\mxsh+0.15,-\mysh) to [out=0,in=180] (\mxsh-0.15,\mysh);
\draw [thkln] (-\mxsh+0.15,-\mysh) to [out=0,in=180] (\mxsh-0.15,\mysh);
\draw [ptzert] (-0.15-\mxsh,-0.6-\mysh) rectangle ++(.3,1.2);
\draw [thkln] (-0.75-\mxsh,-\mysh) -- (-0.15-\mxsh,-\mysh)
node [near start,above] {$\scriptstyle \xca$};
\begin{scope}[xscale=-1]
\draw [ptzert] (-0.15-\mxsh,-0.6-\mysh) rectangle ++(.3,1.2);
\draw [thkln] (-0.75-\mxsh,-\mysh) -- (-0.15-\mxsh,-\mysh)
node [near start,above] {$\scriptstyle \xca$};
\end{scope}
\begin{scope}[yscale=-1]
\draw [ptzert] (-0.15-\mxsh,-0.6-\mysh) rectangle ++(.3,1.2);
\draw [thkln] (-0.75-\mxsh,-\mysh) -- (-0.15-\mxsh,-\mysh)
node [at start,below] {$\scriptstyle \xca-1\;\;\;$};
\end{scope}
\begin{scope}[yscale=-1,xscale=-1]
\draw [ptzert] (-0.15-\mxsh,-0.6-\mysh) rectangle ++(.3,1.2);
\draw [thkln] (-0.75-\mxsh,-\mysh) -- (-0.15-\mxsh,-\mysh)
node [near start,below] {$\scriptstyle \xca$};
\end{scope}
\end{tikzpicture}
\;\hteqv\;
\shcr^{\hlf\xca}\,
\begin{tikzpicture}[menvone]
\draw (\mxsh+0.15,\mysh+0.3) -- (\mxsh+0.75,\mysh+0.3);
\draw [dashed] (-\mxsh-0.75,-\mysh-0.3) -- (-\mxsh-0.15,-\mysh-0.3)
(-\mxsh+0.15,-\mysh-0.3) -- (\mxsh-0.15,-\mysh-0.3) (\mxsh+0.15,-\mysh-0.3) -- (\mxsh+0.75,-\mysh-0.3);
\draw (\mxsh-0.15,\mysh-0.3) to [out=180,in=90]
(\mxsh-1.2,0) to [out=-90,in=180] (\mxsh-0.15,-\mysh+0.3);
\draw [thkln] (-\mxsh+0.15,\mysh) to [out=0,in=180] (\mxsh-0.15,-\mysh);
\draw [lnovr] (-\mxsh+0.15,-\mysh) to [out=0,in=180] (\mxsh-0.15,\mysh);
\draw [thkln] (-\mxsh+0.15,-\mysh) to [out=0,in=180] (\mxsh-0.15,\mysh);
\draw [ptzer] (-0.15-\mxsh,-0.6-\mysh) rectangle ++(.3,1.2);
\draw [thkln] (-0.75-\mxsh,-\mysh) -- (-0.15-\mxsh,-\mysh)
node [near start,above] {$\scriptstyle \xca$};
\begin{scope}[xscale=-1]
\draw [ptzer] (-0.15-\mxsh,-0.6-\mysh) rectangle ++(.3,1.2);
\draw [thkln] (-0.75-\mxsh,-\mysh) -- (-0.15-\mxsh,-\mysh)
node [near start,above] {$\scriptstyle \xca$};
\end{scope}
\begin{scope}[yscale=-1]
\draw [ptzer] (-0.15-\mxsh,-0.6-\mysh) rectangle ++(.3,1.2);
\draw [thkln] (-0.75-\mxsh,-\mysh) -- (-0.15-\mxsh,-\mysh)
node [at start,below] {$\scriptstyle \xca-1\;\;\;$};
\end{scope}
\begin{scope}[yscale=-1,xscale=-1]
\draw [ptzer] (-0.15-\mxsh,-0.6-\mysh) rectangle ++(.3,1.2);
\draw [thkln] (-0.75-\mxsh,-\mysh) -- (-0.15-\mxsh,-\mysh)
node [near start,below] {$\scriptstyle \xca$};
\end{scope}
\end{tikzpicture}
\end{equation}
The first equivalence comes from sliding the upper left projector down right along its $\xca$-cable, and the second equivalence comes from \ex{eq:projtw}. The dashed line indicates the possible presence of another single line which has not been removed yet, however, it plays no role in these calculations.

Let $\xDsep$ denote the diagram $\xDse$ in which the left tangle of \ex{eq:prsld} has been replaced by the right tangle, then
\begin{equation}
\label{eq:scsh}
\xDse \hteqv \shcr^{\hlf\xca}\xDsep.
\end{equation}
We would like to estimate $\hmord{\xDsep}$ with the help of Theorem\rw{thm:smfr}. In doing so we would have to take into account possible crossings coming from the stretch of the single line between the initial projector and the left projector of the tangle $\ytngse$ of \ex{eq:eltng} and $\xca$-cables participating in the crossings attached to the current single line circle. These new crossings are generated by the Reidemeister moves involved in the first homotopy equivalence of \ex{eq:twflp}: when a single line is flipped to the other side of the circle, it may come across the $\xca$-cable crossings, from which parts of this line originate through replacements I of\rx{eq:thtr} (see the picture\rx{eq:tngcr} of the tangle $\ztau$). However, Remark\rw{rmk:ignr} indicates that these crossings between the single line and the $\xca$-cables may be ignored when applying the estimate of Theorem\rw{thm:smfr}, so $\hmord{\xDsep} \geq -\hlf\yncrpv{\xDsep}$, where $\yncrpv{\xDsep}$ is the number of single line intersections within $\xDsep$, except those which we can ignore.

The cable intersection of the left tangle of \ex{eq:prsld} involves two $\xca$-cables, while the same intersection in the right tangle involves a $\xca$-cable and a $(\xca-1)$-cable, hence $\yncrpv{\xDsep} = \yncrv{\xDsi} - \xca$ and the inequality\rx{eq:scest} follows from \ex{eq:scsh}.
\end{proof}

\subsubsection{\Ltrf\ \xltthree}
Set
\begin{equation}
\label{eq:lrpth}
\ytngsi=
\begin{tikzpicture}[menvtwo]
\draw [bcrc] (-0.75,0) -- (0.75,0);
\draw [ptzert] (-0.15,-0.6) rectangle (.15,.6);
\draw [thkln] (-0.75,0) -- (-0.15,0)
node [near start,below] {$\scriptstyle \xca$}
(0.15,0) -- (.75,0)
node [near end,below] {$\scriptstyle \xca$};
\draw (-0.15,0.4) to [out=180,in=-90] (-0.6,0.8) to [out=90,in=180] (0,1.4)
to [out=0,in=90] (0.6,0.8) to [out=-90,in=0] (0.15,0.4);
\end{tikzpicture},
\qquad
\ytngsf =
\begin{tikzpicture}[menvtwo]
\draw [bcrc] (-0.75,0) -- (0.75,0);
\draw [ptzer] (-0.15,-0.6) rectangle (.15,.6);
\draw [thkln] (-0.75,0) -- (-0.15,0)
node [near start,below] {$\scriptstyle \xca$}
(0.15,0) -- (.75,0)
node [near end,below] {$\scriptstyle \xca$};
\end{tikzpicture},
\qquad
\ytngsc = \shcr^{2\xca}\shfr
\begin{tikzpicture}[menvtwo]
\draw [bcrc] (-0.75,0) -- (-0.15,0) (0.15,0) -- (0.75,0);
\draw [color=white] (-0.15,-1) rectangle (0.15,1.6);
\draw [ptthr] (-0.15,-0.6) rectangle (.15,.6);
\draw [thkln] (-0.75,0) -- (-0.15,0)
node [near start,below] {$\scriptstyle \xca$}
(0.15,0) -- (.75,0)
node [near end,below] {$\scriptstyle \xca$};
\draw (-0.15,0.4) to [out=180,in=-90] (-0.6,0.8) to [out=90,in=180] (0,1.4)
to [out=0,in=90] (0.6,0.8) to [out=-90,in=0] (0.15,0.4);
\end{tikzpicture},
\end{equation}
while $\ytngkfp = \shfr^{-1}\ytngkf$
and the cone relation\rx{eq:cnrel} is \ex{eq:htloop}.
\begin{proposition}
\label{prp:sckin}
Let $\xDsi$ be a diagram constructed by removing some single line circles from $\xtDN$ and by detaching the `current' single line circle from the all of its projectors, except the initial one, to which it is attached as in the picture\rx{eq:lrpth} of tangle $\ytngsi$. Let $\xDsf$ be the diagram $\xDsi$ from which this circle is completely removed.
Then the \tdgpr\ map\rx{eq:dgprmg} is an isomorphism on $\tKHmvv{i}{\hem}$ for $i\leq 2\xca-2$.
\end{proposition}
\begin{proof}
%
Since the complex $\ytngsc$ of \ex{eq:lrpth} is a \tmcn\rx{eq:grlp} generated by the elementary tangle
$\ytngsf$ of \ex{eq:lrpth}, then, according to Remark\rw{rmk:bndss}, the claim of this proposition
would follow from the bound
\[
\KHmvv{i}{\hem}(\xDsf)=0\quad\text{for $i\leq-\shlf\yncrv{\xDsi}-1 $.}
\]
The latter follows from Theorem\rw{thm:smfr} coupled with an obvious relation
$\yncrv{\xDsi}=\yncrv{\xDsf}$.
\end{proof}

\subsection{Proof of Theorem\rw{thm:leviso}}
Of all three types of \ltrf s considered in Propositions\rw{prp:bndfs},\rw{prp:bndsc} and\rw{prp:sckin}, it is the \ltrf\rx{eq:lrdpr} which yields the weakest estimate of the homological degrees at which the map\rx{eq:spmaps} is an isomorphism, and this is the estimate of Theorem\rw{thm:leviso}\qed

\section{Proof of Theorem\rw{thm:kqbnd}}

In order to compute the shifted homology $\tKHm(\xDclN)$, we apply the colored Khovanov bracket\rx{eq:colKhbr} to all crossings of $\xDclN$. As a result, this diagram turns into a \tmcn\ of flat diagrams of a special kind. Let $\svrt$ be the set of crossings of $\xD$.
A \emph{\tstt} of $\xDclN$ is a map $\spmp\colon \svrt\rightarrow \{0,1,\ldots,\xca\}$, $\xvrt\mapsto\sipvr$. It
determines a diagram $\xDs$ constructed by performing the following \ltrf s at each crossing $\xvrt\in\svrt$:
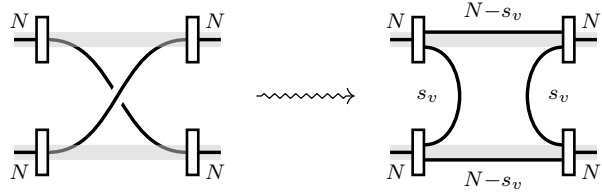
\begin{equation*}
\def\mxsh{2}
\def\mysh{1.5}
\begin{tikzpicture}[menvone]
\node (l) at (-1.6,0) {};
\node (r) at (1.6,0) {};
\path[commutative diagrams/.cd, every arrow, every label]
(l) edge[commutative diagrams/squiggly] (r);
\begin{scope}[xshift=-5cm]
\draw [thkln] (-\mxsh+0.15,\mysh) to [out=0,in=180] (\mxsh-0.15,-\mysh);
\draw [lnovr] (-\mxsh+0.15,-\mysh) to [out=0,in=180] (\mxsh-0.15,\mysh);
\draw [thkln] (-\mxsh+0.15,-\mysh) to [out=0,in=180] (\mxsh-0.15,\mysh);
%
%
\draw [bcrc] (-\mxsh-0.75,\mysh) -- (\mxsh+0.75,\mysh);
\draw [bcrc] (-\mxsh-0.75,-\mysh) -- (\mxsh+0.75,-\mysh);
\draw [ptzert] (-0.15-\mxsh,-0.6-\mysh) rectangle ++(.3,1.2);
\draw [thkln] (-0.75-\mxsh,-\mysh) -- (-0.15-\mxsh,-\mysh)
node [near start,below] {$\scriptstyle \xca$};
\begin{scope}[xscale=-1]
\draw [ptzert] (-0.15-\mxsh,-0.6-\mysh) rectangle ++(.3,1.2);
\draw [thkln] (-0.75-\mxsh,-\mysh) -- (-0.15-\mxsh,-\mysh)
node [near start,below] {$\scriptstyle \xca$};
\end{scope}
\begin{scope}[yscale=-1]
\draw [ptzert] (-0.15-\mxsh,-0.6-\mysh) rectangle ++(.3,1.2);
\draw [thkln] (-0.75-\mxsh,-\mysh) -- (-0.15-\mxsh,-\mysh)
node [near start,above] {$\scriptstyle \xca$};
\end{scope}
\begin{scope}[yscale=-1,xscale=-1]
\draw [ptzert] (-0.15-\mxsh,-0.6-\mysh) rectangle ++(.3,1.2);
\draw [thkln] (-0.75-\mxsh,-\mysh) -- (-0.15-\mxsh,-\mysh)
node [near start,above] {$\scriptstyle \xca$};
\end{scope}
\end{scope}
\begin{scope}[xshift=5cm]
\draw [bcrc] (-0.75-\mxsh,-\mysh) -- (0.75+\mxsh,-\mysh);
\draw [bcrc] (-0.75-\mxsh,\mysh) -- (0.75+\mxsh,\mysh);
\draw [thkln] (-\mxsh+0.15,\mysh-0.2) to [out=0,in=90] (-\mxsh+1.1,0)
 to [out=-90,in=0] (-\mxsh+0.15,-\mysh+0.2) ;
\begin{scope}[xscale=-1]
\draw [thkln] (-\mxsh+0.15,\mysh-0.2) to [out=0,in=90] (-\mxsh+1.1,0) to [out=-90,in=0] (-\mxsh+0.15,-\mysh+0.2);
\end{scope}
\node at (-\mxsh+0.25,0) {$\scriptstyle \spvr$};
\node at (\mxsh-0.25,0) {$\scriptstyle \spvr$};
\draw [thkln] (-\mxsh+0.15,\mysh+0.2) -- (\mxsh-0.15,\mysh+0.2) node [midway,above] {$\scriptstyle \xca-\spvr$};
\draw [thkln] (-\mxsh+0.15,-\mysh-0.2) -- (\mxsh-0.15,-\mysh-0.2) node [midway, below] {$\scriptstyle \xca-\spvr$};
%
%
\draw [ptzert] (-0.15-\mxsh,-0.6-\mysh) rectangle ++(.3,1.2);
\draw [thkln] (-0.75-\mxsh,-\mysh) -- (-0.15-\mxsh,-\mysh)
node [near start,below] {$\scriptstyle \xca$};
\begin{scope}[xscale=-1]
\draw [ptzert] (-0.15-\mxsh,-0.6-\mysh) rectangle ++(.3,1.2);
\draw [thkln] (-0.75-\mxsh,-\mysh) -- (-0.15-\mxsh,-\mysh)
node [near start,below] {$\scriptstyle \xca$};
\end{scope}
\begin{scope}[yscale=-1]
\draw [ptzert] (-0.15-\mxsh,-0.6-\mysh) rectangle ++(.3,1.2);
\draw [thkln] (-0.75-\mxsh,-\mysh) -- (-0.15-\mxsh,-\mysh)
node [near start,above] {$\scriptstyle \xca$};
\end{scope}
\begin{scope}[yscale=-1,xscale=-1]
\draw [ptzert] (-0.15-\mxsh,-0.6-\mysh) rectangle ++(.3,1.2);
\draw [thkln] (-0.75-\mxsh,-\mysh) -- (-0.15-\mxsh,-\mysh)
node [near start,above] {$\scriptstyle \xca$};
\end{scope}
\end{scope}
\end{tikzpicture}
\end{equation*}
The gray strips in $\xDs$ combine into \tBcr s in the background of this diagram.

The diagrams $\xDs$ for all \tstt s $\spmp$ generate a \tmcn\ presentation of $\xDclN$, hence $\tKHm(\xDclN)$ can be computed by spectral sequence, and its $E_1$ term is a sum of appropriately shifted homologies $\KHm(\xDs)$:
\[
\xEo = \shfr^{\xca\gvD}\bigoplus_{\substack{\spmp\\ k\geq 0 }}m_{\spmp,k}\,\shcr^{\xabms+k}\,  \KHm(\xDs),
\]
where $\xabms=\sum_{\xvrt\in\svrt}\spvr^2$.
Hence, a component $\xEoij$ of bi-degree $i,j$ (both are homological and have nothing to do with filtration!) has the form
\begin{equation}
\label{eq:eogr}
\xEoij = \bigoplus_{\substack{\spmp\\ k\geq 0 }}m_{\spmp,k} \KHmvv{i-\xabms-k}{j-\xca\gvD}(\xDs)
\end{equation}

As we already noted in Remark\rw{rmk:bndss}, further steps of spectral sequence may only reduce homology, hence $\xEoij=0$ implies $\tKHmvv{i}{j}(\xDclN)=0$. Moreover, all differentials have bi-grading (-1,1), hence  $\xEovv{i+1}{j-1}=\xEovv{i-1}{j+1}=0$ implies $\tKHmvv{i}{j}(\xDclN) = \xEoij$. These arguments imply that Theorem\rw{thm:kqbnd} follows from the proposition
\begin{proposition}
\label{prp:lm1}
$\xEoij=0$ if one of the following conditions is satisfied:
\begin{align}
\label{eq:bd1b}
i &<0,
\\
\label{eq:bd2b}
j&<- \shlf i - \shlf\ncrD - \sthlf\ncriD,
\\
\label{eq:bd3b}
j& < -i  - \ncriD - \xca\gvD.
\end{align}
Moreover, if $\xD$ is \tBadq, then
\begin{equation}
\label{eq:bd4b}
\xEovv{i}{-i} =
\begin{cases}
0,&\text{if $i\neq 0$,}
\\
\xalg,&\text{if $i=0$.}
\end{cases}
\end{equation}
\end{proposition}
Thus we proved Theorem\rw{thm:leviso}\qed

In view of \ex{eq:eogr}, Proposition\rw{prp:lm1} follows from the next one:
\begin{proposition}
\label{prp:lm2}
$\KHmvv{i}{j}(\xDs)=0$ if one of the following two conditions is satisfied:
\begin{align}
\label{eq:bd1d}
i &<0,
\\
\label{eq:bd2d}
j&<- \shlf \xabms - \shlf\ncrD - \sthlf\ncriD - \xca \gvD,
\\
\label{eq:bd3d}
j& < -\xabms  - \ncriD - \xca\gvD.
\end{align}
Furthermore, if a diagram $\xD$ is \tBadq, then $\KHmvv{i}{j}(\xDs)=0$ for
\begin{equation}
\label{eq:bd4d}
\KHmvv{0}{-\xabms-\xca\gvD}(\xDs) =
\begin{cases}
0,&\text{if $\xabms>0$},
\\
\IQ,&\text{if $\xabms=0$.}
\end{cases}
\end{equation}
\end{proposition}
\begin{proof}[Proof of Proposition\rw{prp:lm1}]
Conditions\rx{eq:bd1b}--\rxw{eq:bd3b} follow easily from the conditions\rx{eq:bd1d}--\rxw{eq:bd3d}. In order to prove \ex{eq:bd4b}, observe that according to \ex{eq:bd4b}, $\xEovv{i}{-i}$ is a sum of homologies $\KHmvv{i'}{j'}(\xDs)$ with $i' =i-\xabms-k$, $j'=-i-\xca\gvD$, hence
\[
j' = -\xabms -\xca\gvD - i' - k.
\]
Since $i'\geq 0$ by \ex{eq:bd1d} and $k\geq 0$ by \ex{eq:eogr}, then in view of the bound\rx{eq:bd3d} with $\ncriD=0$ we conclude that non-trivial contributing homology exists only for $i'=k=0$, so $j'= - \xabms-\xca\gvD$ and $i=\xabms$. Thus we proved that $\xEovv{i}{-i}$ is a sum of homologies $\KHmvv{0}{-\xabms-\xca\gvD}(\xDs)$ with $\xabms=i$, hence \ex{eq:bd4b} follows from \ex{eq:bd4d} and from the fact that the \tstt\ $\spmp$ with $\xabms=0$ is unique (it corresponds to \Bsplng\ all crossings in $\xDclN$) and its multiplicity in the presentation of $\xEovv{i}{-i}$ is one, because complete \Bsplng\ has multiplicity one in \ex{eq:colKhbr}.
\end{proof}

\begin{proof}[Proof of Proposition\rw{prp:lm2}]
First of all, we observe that the bound\rx{eq:bd1d} follows from the fact that $\xDs$ has no crossings, while the formulas\rx{eq:projcn},\rx{eq:grproj} for the \cJWp\ contain only non-negative shifts of \thdgr.

The proof of other bounds requires a simplification of the complex, whose homology yields Khovanov homology $\KHm(\xDs)$. We cut the diagram $\xDs$ into pieces (tangles), simplify their Khovanov complexes and then glue those complexes back together.

%
%
%
%

Consider a neighborhood of a \tBcr\ $c$ within a diagram $\xDs$ and cut in the middle all \tstrtl s which are attached to it. We are going to simplify the complex of the resulting colored tangle $\xtusc$ by inserting two extra \tJWp s in it and then purging all other (preexisting) projectors.

For $a\geq b$ let the box
$
\;
\begin{tikzpicture}[menvthree]
\draw [ptzer] (0.4,0.8) rectangle ++(-0.8,-1.6);
\draw [thkc] (-0.1,-0.2) arc (-90:90:0.2);
\draw [thkln] (-1.2,0) -- (-0.4,0)   node [near start,below] {$\scriptscriptstyle a$}
(0.4,0) -- (1.2,0)  node [near end,below] {$\scriptscriptstyle b$};
\end{tikzpicture}
\;
$
denote any \taTLt\ with the property
\[
\swdv{
\begin{tikzpicture}[menvone]
\draw [ptzer] (0.4,0.8) rectangle ++(-0.8,-1.6);
\draw [thkc] (-0.1,-0.2) arc (-90:90:0.2);
\draw [thkln] (-1,0) -- (-0.4,0) node [near start,below] {$\scriptstyle a$}
(0.4,0) -- (1,0) node [near end,below] {$\scriptstyle b$};
\end{tikzpicture}
} = b.
\]
In other words, a tangle
$
\;
\begin{tikzpicture}[menvthree]
\draw [ptzer] (0.4,0.8) rectangle ++(-0.8,-1.6);
\draw [thkc] (-0.1,-0.2) arc (-90:90:0.2);
\draw [thkln] (-1,0) -- (-0.4,0)  
(0.4,0) -- (1,0);  
\;
\end{tikzpicture}
$
contains no cups, but only caps and \txstrs s.
\begin{lemma}
The Khovanov categorification complex of the colored tangle diagram $\xtusc$ can be presented in the form
\begin{equation}
\label{eq:prcmp}
\xKhv{\xtusc} \hteqv
\boxed{
\cdots\longrightarrow\shcr^i
\bigoplus_{0\leq j\leq i}\shfr^j\left( \bigoplus_{\tau} m_{ij,\tau} \xKhv{\tau}\right)
\longrightarrow\cdots
}_{\;i=0}^{\;\infty}
\end{equation}
where the diagrams $\tau$ are of one of two types depicted in \fg{fg:twodg}:
\begin{figure}
\begin{equation*}
\begin{tikzpicture}[menvtwo]
\draw [bcrct] (-1.25,0) -- (1.25,0) (1.55,0) to [out=0,in=90] (4,-2) to [out=-90,in=0] (0,-5)
(-1.55,0) to [out=180,in=90] (-4,-2) to [out=-90,in=180] (0,-5) ;
\draw [thkln] (-1.25,0) -- (1.25,0) (1.55,0) to [out=0,in=90] (4,-2) to [out=-90,in=0] (0,-5)
(-1.55,0) to [out=180,in=90] (-4,-2) to [out=-90,in=180] (0,-5) ;
\node at (0,-5.5) {$\scriptstyle  \xca_1$};
\node at (0,0.5) {$\scriptstyle \xca_2$};
\draw [ptzert] (-1.55,-0.6) rectangle ++(0.3,1.2);
\draw [ptzert] (1.55,-0.6) rectangle ++(-0.3,1.2);
\draw [ptzer] (0.8,2) rectangle ++(-1.6,-0.8);
\draw [thkc] (0.2,0.1+1.6) arc (0:-180:0.2);
\draw [thkln]  (-0.5,2) -- (-0.5,2.6);
\draw [thkln]  (0.5,2) -- (0.5,2.6);
\node at (0.05,2.4) {$\scriptstyle \cdots$};
\draw [ptzer] (0.8,-2) rectangle ++(-1.6,0.8);
\draw [thkc] (0.2,-0.1-1.6) arc (0:180:0.2);
\draw [thkln]  (-0.5,-2) -- (-0.5,-2.6);
\draw [thkln]  (0.5,-2) -- (0.5,-2.6);
\node at (0.05,-2.4) {$\scriptstyle \cdots$};
\draw [thkln] (-1.25,0.35) to [out=0,in=-90] (-0.5,1.2);
\draw [thkln,xscale=-1] (-1.25,0.35) to [out=0,in=-90] (-0.5,1.2);
\draw [thkln,yscale=-1] (-1.25,0.35) to [out=0,in=-90] (-0.5,1.2);
\draw [thkln,xscale=-1,yscale=-1] (-1.25,0.35) to [out=0,in=-90] (-0.5,1.2);
\draw [decorate,decoration={brace,amplitude=4pt},xshift=0,yshift=-5pt]
(0.5cm+5pt,-2.6) -- (-0.5cm-5pt,-2.6)
node [black,midway,yshift=-0.4cm]
{\scriptsize to struts};
\draw [decorate,decoration={brace,amplitude=4pt},xshift=0,yshift=5pt]
(-0.5cm-5pt,2.6) -- (0.5cm+5pt,2.6)
node [black,midway,yshift=0.4cm]
{\scriptsize to struts};
\end{tikzpicture}
\qquad
\qquad
\begin{tikzpicture}[menvtwo]
\draw [bcrct] (-1.25,0) -- (1.25,0) (1.55,0) to [out=0,in=90] (4,-2) to [out=-90,in=0] (0,-5)
(-1.55,0) to [out=180,in=90] (-4,-2) to [out=-90,in=180] (0,-5) ;
\draw [thkln] (0,1.2) -- (0,-1.2) (1.55,0) to [out=0,in=90] (4,-2) to [out=-90,in=0] (0,-5)
(-1.55,0) to [out=180,in=90] (-4,-2) to [out=-90,in=180] (0,-5) ;
\node at (0,-5.5) {$\scriptstyle  \xca_1$};
\draw [ptzert] (-1.55,-0.6) rectangle ++(0.3,1.2);
\draw [ptzert] (1.55,-0.6) rectangle ++(-0.3,1.2);
\draw [ptzer] (0.8,2) rectangle ++(-1.6,-0.8);
\draw [thkc] (0.2,0.1+1.6) arc (0:-180:0.2);
\draw [thkln]  (-0.5,2) -- (-0.5,2.6);
\draw [thkln]  (0.5,2) -- (0.5,2.6);
\node at (0.05,2.4) {$\scriptstyle \cdots$};
\draw [ptzer] (0.8,-2) rectangle ++(-1.6,0.8);
\draw [thkc] (0.2,-0.1-1.6) arc (0:180:0.2);
\draw [thkln]  (-0.5,-2) -- (-0.5,-2.6);
\draw [thkln]  (0.5,-2) -- (0.5,-2.6);
\node at (0.05,-2.4) {$\scriptstyle \cdots$};
\draw [thkln] (-1.25,0.35) to [out=0,in=-90] (-0.5,1.2);
\draw [thkln,xscale=-1] (-1.25,0.35) to [out=0,in=-90] (-0.5,1.2);
\draw [thkln,yscale=-1] (-1.25,0.35) to [out=0,in=-90] (-0.5,1.2);
\draw [thkln,xscale=-1,yscale=-1] (-1.25,0.35) to [out=0,in=-90] (-0.5,1.2);
%
%
\draw [decorate,decoration={brace,amplitude=4pt},xshift=0,yshift=-5pt]
(0.5cm+5pt,-2.6) -- (-0.5cm-5pt,-2.6)
node [black,midway,yshift=-0.4cm]
{\scriptsize to struts};
\draw [decorate,decoration={brace,amplitude=4pt},xshift=0,yshift=5pt]
(-0.5cm-5pt,2.6) -- (0.5cm+5pt,2.6)
node [black,midway,yshift=0.4cm]
{\scriptsize to struts};
\end{tikzpicture}
\end{equation*}
\caption{A purged vicinity of a \tBcr}
\label{fg:twodg}
\end{figure}
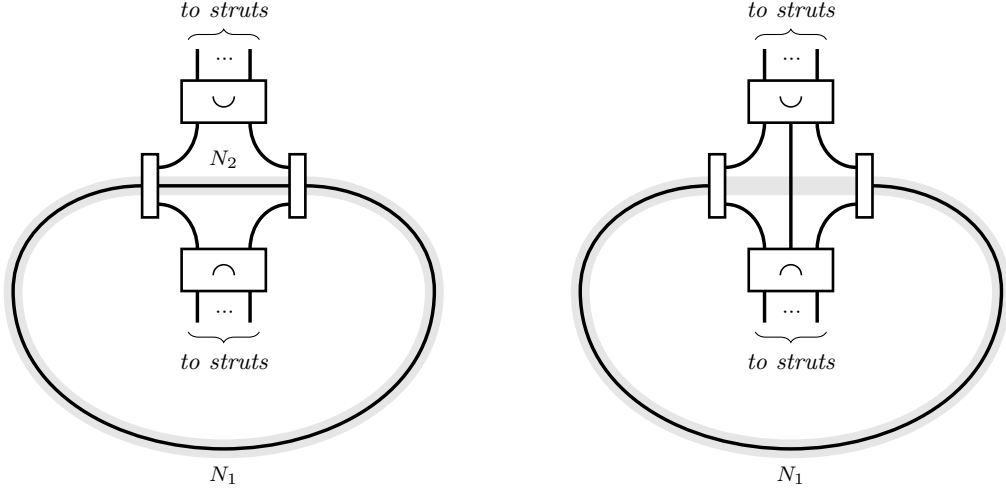
with $\xca_1,\xca_2\leq \xca$.
\end{lemma}

By gluing the complexes\rx{eq:prcmp} back together
we get a complex
\[
\xKhv{\xDs} \hteqv
\boxed{
\cdots\longrightarrow\shcr^i
\bigoplus_{0\leq j\leq i}\shfr^j\left( \bigoplus_{\xDcir} m_{ij,\tau} \xKhv{\xDcir}\right)
\longrightarrow\cdots
}_{\;i=0}^{\;\infty}
\]
such that $\KHm(\xDs) = \Hm(\xKhv{\xDs})$. The `circle diagrams' $\xDcir$ which result from gluing the diagrams of \fg{fg:twodg} at the \tstrtl\ cutting points and replacing projectors with complexes\rx{eq:projcn} and\rx{eq:grproj}, consist of multiple single line circles.
%
In view of the second formula of \ex{eq:dkhbr}, the lowest \tqdgr\ in the homology $\KHm(\xDs)$ may be bounded by the highest number of circles in those circle diagrams.

The circles in circle diagrams are of three types. The first type is \emph{\tjmp} circles: they contain at least one \tstrt\ line. The circles of the second and third type stay within the same \tBcr. A \emph{\txstr} circle goes along its \tBcr, passing straight through each \tJWp\ on its way. A \emph{\twndg} circle changes its direction at least twice, because it contains at least one cup and one cap of a constituent \taTLt\ coming from one of projectors.

Let us prove the inequalities\rx{eq:bd2d} and\rx{eq:bd3d} by finding upper bounds for the numbers $\njcr$, $\nscr$ and $\nwcr$ of \tjmp, \txstr\ and \twndg\ circles respectively in a circle diagram.

We begin with $\njcr$.
A \tjmpc\ must contain at least two \tstrt\ lines of an \tadq\ crossing or at least one \tstrt\ line of an inadequate crossing, so the number of \tjmpc s $\njcr$ has a bound:
\begin{equation}
\label{eq:sbnd}
\njcr \leq \sum_{\xvrt\in\svrta}\spvr + 2\sum_{\xvrt\in\svrti}\spvr,
\end{equation}
where $\svrta,\svrti\subset\svrt$ are the subsets of \tBadq\ and \tBiadq\ crossings.
The obvious inequalities
\[
\spvr\leq\shlf\spvr^2+\shlf,\qquad 2\spvr\leq\shlf\spvr^2 +2,\qquad
\spvr\leq\spvr^2,\qquad 2\spvr\leq\spvr^2 + 1,
\]
(the third inequality uses the fact that $\spvr$ is integer) indicate that the bound\rx{eq:sbnd} implies to other bounds:
\begin{equation}
\label{eq:mbnd}
\njcr \leq \shlf\xabms + \shlf\ncrD + \sthlf\ncriD,\qquad
\njcr \leq \xabms + \ncriD,
\end{equation}
which means that the first three terms in the \rhs of the inequality\rx{eq:bd2d} and the first two terms in the inequality\rx{eq:bd3d} bound the negative contribution of \tjmpc s to the \tqdgr\ of $\KHm(\xDs)$.

Next we prove the bound
\begin{equation}
\label{eq:wbd}
\nscrb + \nwcrb \leq \xca,
\end{equation}
where $\nscrb$ and $\nwcrb$ are the numbers of \txstr\ and \twndg\ circles within any given \tBcr\ $\crb$. It implies the bound $\nscr + \nwcr \leq \xca \gvD$ and combined with the bounds\rx{eq:mbnd} they imply the bounds\rx{eq:bd2d} and\rx{eq:bd3d}. In order to prove the bound\rx{eq:wbd}, we observe that
in the first diagram a \tstrghtc\ contains one strand from the $\xca_1$-cable and one strand from the $\xca_2$-cable, while a \twndgc\ contains at least two strands from one of these cables, hence there is a bound
\begin{equation}
\label{eq:bndNt}
\nscrb+\nwcrb\leq\shlf(\xca_1+\xca_2)\leq\xca.
\end{equation}
The second diagram is treated similarly, if we set $\xca_2=0$ in the previous argument.
Thus we proved the bounds\rx{eq:bd2d} and\rx{eq:bd3d}.

It remains to prove \ex{eq:bd4d}. Since this time $\xD$ is \tBadq, the second inequality of\rx{eq:mbnd} becomes
$ \njcr \leq \xabms$. Since we consider only homology of zeroth \thdgr, then according to \eex{eq:projcn} and\rx{eq:grproj}, we may replace \tJWp s with identity braids, so there is only one circle diagram $\xDcir$ contributing to  $\KHmvv{0}{-\xabms-\xca\gvD}(\xDs)$, and this circle diagram has no \twndg\ circles: $\nwcr=0$. Furthermore, $\nscr\leq\xca\gvD$, but if $\xabms\neq 0$, then there is at least one pair of \tstrtl s in $\xDs$, so $\nscr<\xca$ and $\njcr+\nscr<\xabms + \xca\gvD$, hence $\KHmvv{0}{-\xabms-\xca\gvD}(\xDs) = 0$. If $\xabms=0$, then $\xDs$ has no \tstrtl s and consists of disjoint $\xca$-cabled circles, so the relevant circle diagram $\xDcir$
consists
of $\xca\gvD$ single-line circles, and $\KHmvv{0}{-\xca\gvD}(\xDs)=\IQ$ follows from the second equation of\rx{eq:dkhbr}.
\end{proof}

\begin{proof}[Proof of Lemma\rw{eq:prcmp}]
We prove the lemma by `purging' \cJWp s appearing in the tangle $\xtusc$.
%
In order to bring the complex $\xKhv{\xtusc}$ to the form\rx{eq:prcmp} with diagrams $\tau$ depicted in \fg{fg:twodg}, we insert two extra \tJWp s side-by-side at any place on the cable which runs along the \tBcr. Then we go from the front one (relative to the clockwise orientation) to the back one in the clockwise direction, purging each preexisting projector that appears on our way. It is easy to prove by induction that after every projector purge we get a \tmcn\ presentation
\begin{equation*}
\xKhv{\xtusc} \hteqv
\boxed{
\cdots\longrightarrow\shcr^i
\bigoplus_{0\leq j\leq i}\shfr^j\left( \bigoplus_{\tau} m'_{ij,\tau} \xKhv{\tau}\right)
\longrightarrow\cdots
}_{\;i=0}^{\;\infty}
\end{equation*}
whose constituent diagrams $\tau$ have one of two possible forms between the front projector and the first unpurged  projector (which lies in the pictures to the left of the dashed line) depicted in \fg{fg:prgd}.
\begin{figure}[h]
\[
\begin{tikzpicture}[menvtwo]
\draw [bcrct] (-1.55-0.6-0.5,0) -- (-1.55,0) (-1.25,0) -- (1.25,0) (1.55,0) -- (4.05,0) (4.35,0) -- (4.95+0.5,0);
\draw [thkln] (-1.55-0.6-0.5,0) -- (-1.55,0) (-1.25,0) -- (1.25,0) (1.55,0) -- (4.05,0) (4.35,0) -- (4.95+0.5,0);
\draw [ptzert] (-1.55,-0.6) rectangle ++(0.3,1.2);
\draw [ptzert] (1.55,-0.6) rectangle ++(-0.3,1.2);
\draw [ptzert] (4.05,-0.6) rectangle ++(0.3,1.2);
\draw [ptzer] (0.8,2) rectangle ++(-1.6,-0.8);
\draw [thkc] (0.2,0.1+1.6) arc (0:-180:0.2);
\draw [thkln]  (-0.5,2) -- (-0.5,2.6);
\draw [thkln]  (0.5,2) -- (0.5,2.6);
\node at (0.05,2.4) {$\scriptstyle \cdots$};
\draw [ptzer] (0.8,-2) rectangle ++(-1.6,0.8);
\draw [thkc] (0.2,-0.1-1.6) arc (0:180:0.2);
\draw [thkln]  (-0.5,-2) -- (-0.5,-2.6);
\draw [thkln]  (0.5,-2) -- (0.5,-2.6);
\node at (0.05,-2.4) {$\scriptstyle \cdots$};
\draw [thkln] (-1.25,0.35) to [out=0,in=-90] (-0.5,1.2);
\draw [thkln,xscale=-1] (-1.25,0.35) to [out=0,in=-90] (-0.5,1.2);
\draw [thkln,yscale=-1] (-1.25,0.35) to [out=0,in=-90] (-0.5,1.2);
\draw [thkln,xscale=-1,yscale=-1] (-1.25,0.35) to [out=0,in=-90] (-0.5,1.2);
\draw [thkln,xscale=-1] (-1.55,0.35) to [out=180,in=-90] ++(-0.75,0.85) -- (-2.3,2.6);
\draw [thkln,xshift=5.6cm] (-1.55,0.35) to [out=180,in=-90] ++(-0.75,0.85) -- (-2.3,2.6);
\draw [thkln,xscale=-1,yscale=-1,xshift=-2.8cm] (-1.55,0.35) to [out=180,in=-90] ++(-0.75,0.85) -- (-2.3,2.6);
\draw [decorate,decoration={brace,amplitude=4pt},xshift=0,yshift=-5pt]
(5.1cm+5pt,-2.6) -- (-0.5cm-5pt,-2.6)
node [black,midway,yshift=-0.4cm]
{\scriptsize to struts};
\draw [decorate,decoration={brace,amplitude=4pt},xshift=0,yshift=5pt]
(-0.5cm-5pt,2.6) -- (3.3cm+5pt,2.6)
node [black,midway,yshift=0.4cm]
{\scriptsize to struts};
\draw [dashed] (2.8,2.6) -- (2.8,-2.6);
\end{tikzpicture}
\;,\qquad
\begin{tikzpicture}[menvtwo]
\draw [bcrct] (-1.55-0.6-0.5,0) -- (-1.55,0) (-1.25,0) -- (1.25,0) (1.55,0) -- (4.05,0) (4.35,0) -- (4.95+0.5,0);
\draw [thkln] (-1.55-0.6-0.5,0) -- (-1.55,0) (0,1.2) -- (0,-1.2) (1.55,0) -- (4.05,0) (4.35,0) -- (4.95+0.5,0);
\draw [ptzert] (-1.55,-0.6) rectangle ++(0.3,1.2);
\draw [ptzert] (1.55,-0.6) rectangle ++(-0.3,1.2);
\draw [ptzert] (4.05,-0.6) rectangle ++(0.3,1.2);
\draw [ptzer] (0.8,2) rectangle ++(-1.6,-0.8);
\draw [thkc] (0.2,0.1+1.6) arc (0:-180:0.2);
\draw [thkln]  (-0.5,2) -- (-0.5,2.6);
\draw [thkln]  (0.5,2) -- (0.5,2.6);
\node at (0.05,2.4) {$\scriptstyle \cdots$};
\draw [ptzer] (0.8,-2) rectangle ++(-1.6,0.8);
\draw [thkc] (0.2,-0.1-1.6) arc (0:180:0.2);
\draw [thkln]  (-0.5,-2) -- (-0.5,-2.6);
\draw [thkln]  (0.5,-2) -- (0.5,-2.6);
\node at (0.05,-2.4) {$\scriptstyle \cdots$};
\draw [thkln] (-1.25,0.35) to [out=0,in=-90] (-0.5,1.2);
\draw [thkln,xscale=-1] (-1.25,0.35) to [out=0,in=-90] (-0.5,1.2);
\draw [thkln,yscale=-1] (-1.25,0.35) to [out=0,in=-90] (-0.5,1.2);
\draw [thkln,xscale=-1,yscale=-1] (-1.25,0.35) to [out=0,in=-90] (-0.5,1.2);
\draw [thkln,xscale=-1] (-1.55,0.35) to [out=180,in=-90] ++(-0.75,0.85) -- (-2.3,2.6);
\draw [thkln,xshift=5.6cm] (-1.55,0.35) to [out=180,in=-90] ++(-0.75,0.85) -- (-2.3,2.6);
\draw [thkln,xscale=-1,yscale=-1,xshift=-2.8cm] (-1.55,0.35) to [out=180,in=-90] ++(-0.75,0.85) -- (-2.3,2.6);
\draw [decorate,decoration={brace,amplitude=4pt},xshift=0,yshift=-5pt]
(5.1cm+5pt,-2.6) -- (-0.5cm-5pt,-2.6)
node [black,midway,yshift=-0.4cm]
{\scriptsize to struts};
\draw [decorate,decoration={brace,amplitude=4pt},xshift=0,yshift=5pt]
(-0.5cm-5pt,2.6) -- (3.3cm+5pt,2.6)
node [black,midway,yshift=0.4cm]
{\scriptsize to struts};
\draw [dashed] (2.8,2.6) -- (2.8,-2.6);
\end{tikzpicture}
\]
\caption{Purging \tJWp s along a \tBcr}
\label{fg:prgd}
\end{figure}
In both diagrams the left projector on the grey strip is the front one, the middle projector is the first unpurged one and the right projector is the second unpurged one. It is not hard to see that if we purge the middleprojector, then we get similar diagrams with the third projector becoming the first unpurged one (the left diagram may become of either left or right type after the purge, while the right diagram remains of the same type). The \tqdgr\ shifts remain non-negative, because the purging does not produce any circles: it just makes explicit various line connections that were hidden inside the constituent \taTLt s of the purged projector.
\end{proof}

\section{Proof of invariance of the tail homology under \tBrdc}
\begin{proof}[Proof of Theorem\rw{thm:rddg}]
A removal of a \tBcr\ connected to the rest of the \tBdg\ by a single \tstrt\ corresponds to the first Reidemeister move, hence the invariance of the tail homology under this removal follows from the fact that tail homology of a \tBadq\ link is determined by shifted Khovanov homologies of its \uclrd\ diagrams and the latter are invariant under this type of first Reidemeister moves.

In view of Corollary\rw{cor:mpis}, the invariance of the tail homology under the removal of `extra' \tstrt s follows from the next lemma.
\end{proof}
\begin{lemma}
\label{lm:strm}
Suppose that two distinct \tBcr s of a link diagram $\xD$ are connected by multiple \tstrt s and the diagram $\xDp$ is constructed by removal of one of those \tstrt s. Then there exists a \tdgpr\ map $\tKHm(\xDclN')\xrightarrow{\xmg}\tKHm(\xDclN)$ which is an isomorphism on $\tKHmvv{i}{\hem}$ for $i\leq \xca-1$.
\end{lemma}
The proof of this lemma is similar to proofs of Section\rw{sct:mrph}:
we show that $\xDclN$ can be constructed from $\xDclN'$ with the help of a \ltrf\ and prove the homological smallness of the correction diagram.

We need a simple corollary of Theorem\rw{thm:colkhovbr}
\begin{corollary}
The \tKhbr\ formula\rx{eq:lmtcn} for the colored crossing can be recast in the form
\def\mxsh{2}
\def\mysh{1.2}
\begin{equation}
\label{eq:ccrc}
\begin{tikzpicture}[menvone]
\draw [thkln] (-\mxsh+0.15,\mysh) to [out=0,in=180] (\mxsh-0.15,-\mysh);
\draw [lnovr] (-\mxsh+0.15,-\mysh) to [out=0,in=180] (\mxsh-0.15,\mysh);
\draw [thkln] (-\mxsh+0.15,-\mysh) to [out=0,in=180] (\mxsh-0.15,\mysh);
%
%
\draw [ptzer] (-0.15-\mxsh,-0.6-\mysh) rectangle ++(.3,1.2);
\draw [thkln] (-0.75-\mxsh,-\mysh) -- (-0.15-\mxsh,-\mysh)
node [near start,below] {$\scriptstyle \xca$};
\begin{scope}[xscale=-1]
\draw [ptzer] (-0.15-\mxsh,-0.6-\mysh) rectangle ++(.3,1.2);
\draw [thkln] (-0.75-\mxsh,-\mysh) -- (-0.15-\mxsh,-\mysh)
node [near start,below] {$\scriptstyle \xca$};
\end{scope}
\begin{scope}[yscale=-1]
\draw [ptzer] (-0.15-\mxsh,-0.6-\mysh) rectangle ++(.3,1.2);
\draw [thkln] (-0.75-\mxsh,-\mysh) -- (-0.15-\mxsh,-\mysh)
node [near start,above] {$\scriptstyle \xca$};
\end{scope}
\begin{scope}[yscale=-1,xscale=-1]
\draw [ptzer] (-0.15-\mxsh,-0.6-\mysh) rectangle ++(.3,1.2);
\draw [thkln] (-0.75-\mxsh,-\mysh) -- (-0.15-\mxsh,-\mysh)
node [near start,above] {$\scriptstyle \xca$};
\end{scope}
\end{tikzpicture}
\;\hteqv\;
\shcr^{-\hlf\xca^2}\;
\boxed{
\begin{tikzpicture}[menvone]
\draw [thkln] (-\mxsh+0.15,\mysh) -- (-\mxsh*0.3,\mysh)
node [near start,above] {$\scriptstyle \xca$}
(\mxsh*0.3,\mysh) -- (\mxsh-0.15,\mysh)
node [near end,above] {$\scriptstyle \xca$};
\draw [thkln] (-\mxsh+0.15,-\mysh) -- (-\mxsh*0.3,-\mysh)
node [near start,below] {$\scriptstyle \xca$} (\mxsh*0.3,-\mysh) -- (\mxsh-0.15,-\mysh)
node [near end,below] {$\scriptstyle \xca$};
\draw[ptzer] (-\mxsh*0.3,-\mysh*1-0.6) rectangle ++(\mxsh*0.6,\mysh*2+1.2);
\node at (0,0) {$*$};
\end{tikzpicture}
\longrightarrow
\begin{tikzpicture}[menvone]
\draw [thkln] (-\mxsh+0.15,\mysh) --  (\mxsh-0.15,\mysh);
\draw [thkln] (-\mxsh+0.15,-\mysh) -- (\mxsh-0.15,-\mysh);
\draw [ptzer] (-0.15-\mxsh,-0.6-\mysh) rectangle ++(.3,1.2);
\draw [thkln] (-0.75-\mxsh,-\mysh) -- (-0.15-\mxsh,-\mysh)
node [near start,below] {$\scriptstyle \xca$};
\begin{scope}[xscale=-1]
\draw [ptzer] (-0.15-\mxsh,-0.6-\mysh) rectangle ++(.3,1.2);
\draw [thkln] (-0.75-\mxsh,-\mysh) -- (-0.15-\mxsh,-\mysh)
node [near start,below] {$\scriptstyle \xca$};
\end{scope}
\begin{scope}[yscale=-1]
\draw [ptzer] (-0.15-\mxsh,-0.6-\mysh) rectangle ++(.3,1.2);
\draw [thkln] (-0.75-\mxsh,-\mysh) -- (-0.15-\mxsh,-\mysh)
node [near start,above] {$\scriptstyle \xca$};
\end{scope}
\begin{scope}[yscale=-1,xscale=-1]
\draw [ptzer] (-0.15-\mxsh,-0.6-\mysh) rectangle ++(.3,1.2);
\draw [thkln] (-0.75-\mxsh,-\mysh) -- (-0.15-\mxsh,-\mysh)
node [near start,above] {$\scriptstyle \xca$};
\end{scope}
\end{tikzpicture}
}\;,
\end{equation}
where
\begin{equation}
\label{eq:stbx}
\begin{tikzpicture}[menvone]
\draw [thkln] (-\mxsh+0.15,\mysh) -- (-\mxsh*0.3,\mysh)
node [near start,above] {$\scriptstyle \xca$}
(\mxsh*0.3,\mysh) -- (\mxsh-0.15,\mysh)
node [near end,above] {$\scriptstyle \xca$};
\draw [thkln] (-\mxsh+0.15,-\mysh) -- (-\mxsh*0.3,-\mysh)
node [near start,below] {$\scriptstyle \xca$} (\mxsh*0.3,-\mysh) -- (\mxsh-0.15,-\mysh)
node [near end,below] {$\scriptstyle \xca$};
\draw[ptzer] (-\mxsh*0.3,-\mysh*1-0.6) rectangle ++(\mxsh*0.6,\mysh*2+1.2);
\node at (0,0) {$*$};
\end{tikzpicture}
\;\hteqv\;
\Pcnv{
\bigoplus_{\xki=1}^{\xca}
\shcr^{\xki^2}{a \brace \xki}_{\shcr}
\begin{tikzpicture}[menvone]
\draw [thkln] (-\mxsh+0.15,\mysh-0.2) to [out=0,in=90] (-\mxsh+1.1,0)
 to [out=-90,in=0] (-\mxsh+0.15,-\mysh+0.2) ;
\begin{scope}[xscale=-1]
\draw [thkln] (-\mxsh+0.15,\mysh-0.2) to [out=0,in=90] (-\mxsh+1.1,0) to [out=-90,in=0] (-\mxsh+0.15,-\mysh+0.2);
\end{scope}
\node at (-\mxsh+0.5,0) {$\scriptstyle \xki$};
\node at (\mxsh-0.5,0) {$\scriptstyle \xki$};
\draw [thkln] (-\mxsh+0.15,\mysh+0.2) -- (\mxsh-0.15,\mysh+0.2) node [midway,above] {$\scriptstyle \xca-\xki$};
\draw [thkln] (-\mxsh+0.15,-\mysh-0.2) -- (\mxsh-0.15,-\mysh-0.2) node [midway, below] {$\scriptstyle \xca-\xki$};
\draw [ptzert] (-0.15-\mxsh,-0.6-\mysh) rectangle ++(.3,1.2);
\draw [thkln] (-0.75-\mxsh,-\mysh) -- (-0.15-\mxsh,-\mysh)
node [near start,below] {$\scriptstyle \xca$};
\begin{scope}[xscale=-1]
\draw [ptzert] (-0.15-\mxsh,-0.6-\mysh) rectangle ++(.3,1.2);
\draw [thkln] (-0.75-\mxsh,-\mysh) -- (-0.15-\mxsh,-\mysh)
node [near start,below] {$\scriptstyle \xca$};
\end{scope}
\begin{scope}[yscale=-1]
\draw [ptzert] (-0.15-\mxsh,-0.6-\mysh) rectangle ++(.3,1.2);
\draw [thkln] (-0.75-\mxsh,-\mysh) -- (-0.15-\mxsh,-\mysh)
node [near start,above] {$\scriptstyle \xca$};
\end{scope}
\begin{scope}[yscale=-1,xscale=-1]
\draw [ptzert] (-0.15-\mxsh,-0.6-\mysh) rectangle ++(.3,1.2);
\draw [thkln] (-0.75-\mxsh,-\mysh) -- (-0.15-\mxsh,-\mysh)
node [near start,above] {$\scriptstyle \xca$};
\end{scope}
\end{tikzpicture}
}
\end{equation}
\end{corollary}
\begin{proof}
According to the second inequality of\rx{eq:twineq}, $\yki(0)=0$, hence the the right tangle of the cone\rx{eq:ccrc} is the one that appears at $i=0$ in the \tmcn\rx{eq:lmtcn}. According to \ex{eq:colKhbr}, this tangle has multiplicity one, so this is the only place where it may appear in that \tmcn, and it has a zero shift of \thdgr.
\end{proof}

\begin{proof}[Proof of Lemma\rw{lm:strm}]
Since $\xDp$ is constructed from $\xD$ by a removal of a single crossing (\tstrt\ of $\sBD$), we set
\def\mxsh{2}
\def\mysh{1.2}
\[
\ytngsi =
\begin{tikzpicture}[menvone]
\draw [thkln] (-\mxsh+0.15,\mysh) to [out=0,in=180] (\mxsh-0.15,-\mysh);
\draw [lnovr] (-\mxsh+0.15,-\mysh) to [out=0,in=180] (\mxsh-0.15,\mysh);
\draw [thkln] (-\mxsh+0.15,-\mysh) to [out=0,in=180] (\mxsh-0.15,\mysh);
%
%
\draw [ptzer] (-0.15-\mxsh,-0.6-\mysh) rectangle ++(.3,1.2);
\draw [thkln] (-0.75-\mxsh,-\mysh) -- (-0.15-\mxsh,-\mysh)
node [near start,below] {$\scriptstyle \xca$};
\begin{scope}[xscale=-1]
\draw [ptzer] (-0.15-\mxsh,-0.6-\mysh) rectangle ++(.3,1.2);
\draw [thkln] (-0.75-\mxsh,-\mysh) -- (-0.15-\mxsh,-\mysh)
node [near start,below] {$\scriptstyle \xca$};
\end{scope}
\begin{scope}[yscale=-1]
\draw [ptzer] (-0.15-\mxsh,-0.6-\mysh) rectangle ++(.3,1.2);
\draw [thkln] (-0.75-\mxsh,-\mysh) -- (-0.15-\mxsh,-\mysh)
node [near start,above] {$\scriptstyle \xca$};
\end{scope}
\begin{scope}[yscale=-1,xscale=-1]
\draw [ptzer] (-0.15-\mxsh,-0.6-\mysh) rectangle ++(.3,1.2);
\draw [thkln] (-0.75-\mxsh,-\mysh) -- (-0.15-\mxsh,-\mysh)
node [near start,above] {$\scriptstyle \xca$};
\end{scope}
\end{tikzpicture}\;,
\qquad
\ytngsf =
\begin{tikzpicture}[menvone]
\draw [thkln] (-\mxsh+0.15,\mysh) --  (\mxsh-0.15,\mysh);
\draw [thkln] (-\mxsh+0.15,-\mysh) -- (\mxsh-0.15,-\mysh);
\draw [ptzer] (-0.15-\mxsh,-0.6-\mysh) rectangle ++(.3,1.2);
\draw [thkln] (-0.75-\mxsh,-\mysh) -- (-0.15-\mxsh,-\mysh)
node [near start,below] {$\scriptstyle \xca$};
\begin{scope}[xscale=-1]
\draw [ptzer] (-0.15-\mxsh,-0.6-\mysh) rectangle ++(.3,1.2);
\draw [thkln] (-0.75-\mxsh,-\mysh) -- (-0.15-\mxsh,-\mysh)
node [near start,below] {$\scriptstyle \xca$};
\end{scope}
\begin{scope}[yscale=-1]
\draw [ptzer] (-0.15-\mxsh,-0.6-\mysh) rectangle ++(.3,1.2);
\draw [thkln] (-0.75-\mxsh,-\mysh) -- (-0.15-\mxsh,-\mysh)
node [near start,above] {$\scriptstyle \xca$};
\end{scope}
\begin{scope}[yscale=-1,xscale=-1]
\draw [ptzer] (-0.15-\mxsh,-0.6-\mysh) rectangle ++(.3,1.2);
\draw [thkln] (-0.75-\mxsh,-\mysh) -- (-0.15-\mxsh,-\mysh)
node [near start,above] {$\scriptstyle \xca$};
\end{scope}
\end{tikzpicture}\;,
\qquad
\ytngsc =
\shcr^{-\hlf\xca^2}
\begin{tikzpicture}[menvone]
\draw [thkln] (-\mxsh+0.15,\mysh) -- (-\mxsh*0.3,\mysh)
node [near start,above] {$\scriptstyle \xca$}
(\mxsh*0.3,\mysh) -- (\mxsh-0.15,\mysh)
node [near end,above] {$\scriptstyle \xca$};
\draw [thkln] (-\mxsh+0.15,-\mysh) -- (-\mxsh*0.3,-\mysh)
node [near start,below] {$\scriptstyle \xca$} (\mxsh*0.3,-\mysh) -- (\mxsh-0.15,-\mysh)
node [near end,below] {$\scriptstyle \xca$};
\draw[ptzer] (-\mxsh*0.3,-\mysh*1-0.6) rectangle ++(\mxsh*0.6,\mysh*2+1.2);
\node at (0,0) {$*$};
\end{tikzpicture}\;,
\]
so that $\xDsi = \xD$ and $\xDsf = \xDp$, while the relation\rx{eq:cnrel} comes from\rx{eq:ccrc} if we set $\xhshf=0$ and $\xqshf= -\hlf \xca^2$. Since $\yncrv{\ytngsi}=\xca^2$ , in view of Proposition\rw{prp:gestdg} it remains to establish the bound
$
\tKHmvv{i}{\hem}(\xDsc) = 0,
$
if $i\leq\xca$. Let $\xDsci$ be the diagram $\xDsc$ in which $\ytngsc$ is replaced by a constituent tangle $\ytngsci$ from the \rhs of \ex{eq:stbx}, in which one strand of a $\xki$-cable is separated from the others: 
\def\esh{2.5}
\[
\ytngsci =
\begin{tikzpicture}[menvone]
\draw [thkln] (-\mxsh+0.15,\mysh-0.2) to [out=0,in=90] (-\mxsh+1.1,0)
 to [out=-90,in=0] (-\mxsh+0.15,-\mysh+0.2) ;
\begin{scope}[xscale=-1]
\draw [thkln] (-\mxsh+0.15,\mysh-0.2) to [out=0,in=90] (-\mxsh+1.1,0) to [out=-90,in=0] (-\mxsh+0.15,-\mysh+0.2);
\end{scope}
\node at (-\mxsh+0.2,0) {$\scriptstyle \xki-1$};
\node at (\mxsh-0.5,0) {$\scriptstyle \xki$};
\draw [thkln] (-\mxsh+0.15,\mysh+0.2) -- (\mxsh-0.15,\mysh+0.2) node [midway,above] {$\scriptstyle \xca-\xki$};
\draw [thkln] (-\mxsh+0.15,-\mysh-0.2) -- (\mxsh-0.15,-\mysh-0.2) node [midway, below] {$\scriptstyle \xca-\xki$};
\draw [ptzert] (-0.15-\mxsh,-0.6-\mysh) rectangle ++(.3,1.2);
\draw [thkln] (-\esh-0.75-\mxsh,-\mysh) -- (-\esh-0.15-\mxsh,-\mysh)
node [near start,below] {$\scriptstyle \xca$};
\draw [thkln] (0.15-\esh-\mxsh,-\mysh-0.2) -- (-0.15-\mxsh,-\mysh-0.2)
node [midway, below] {$\scriptstyle \xca-1$};
\begin{scope}[xscale=-1]
\draw [ptzert] (-0.15-\mxsh,-0.6-\mysh) rectangle ++(.3,1.2);
\draw [thkln] (-0.75-\mxsh,-\mysh) -- (-0.15-\mxsh,-\mysh)
node [near start,below] {$\scriptstyle \xca$};
\end{scope}
\begin{scope}[yscale=-1]
\draw [ptzert] (-0.15-\mxsh,-0.6-\mysh) rectangle ++(.3,1.2);
\draw [thkln] (-\esh-0.75-\mxsh,-\mysh) -- (-\esh-0.15-\mxsh,-\mysh)
node [near start,above] {$\scriptstyle \xca$};
\draw [thkln] (0.15-\esh-\mxsh,-\mysh-0.2) -- (-0.15-\mxsh,-\mysh-0.2)
node [midway, above] {$\scriptstyle \xca-1$};
\end{scope}
\begin{scope}[yscale=-1,xscale=-1]
\draw [ptzert] (-0.15-\mxsh,-0.6-\mysh) rectangle ++(.3,1.2);
\draw [thkln] (-0.75-\mxsh,-\mysh) -- (-0.15-\mxsh,-\mysh)
node [near start,above] {$\scriptstyle \xca$};
\end{scope}
\draw [ptzert] (-\esh-0.15-\mxsh,-0.6-\mysh) rectangle ++(.3,1.2);
\draw [ptzert] (-\esh-0.15-\mxsh,0.6+\mysh) rectangle ++(.3,-1.2);
\draw (-\esh-\mxsh+0.15,\mysh-0.2) to [out=0,in=90] (-\esh-\mxsh+1.1,0) to [out=-90,in=0] (-\esh-\mxsh+0.15,-\mysh+0.2);
\end{tikzpicture}
\;\hteqv\;
\begin{tikzpicture}[menvone]
\draw [thkln] (-\mxsh+0.15,\mysh-0.2) to [out=0,in=90] (-\mxsh+1.1,0)
 to [out=-90,in=0] (-\mxsh+0.15,-\mysh+0.2) ;
\begin{scope}[xscale=-1]
\draw [thkln] (-\mxsh+0.15,\mysh-0.2) to [out=0,in=90] (-\mxsh+1.1,0) to [out=-90,in=0] (-\mxsh+0.15,-\mysh+0.2);
\end{scope}
\node at (-\mxsh+0.5,0) {$\scriptstyle \xki$};
\node at (\mxsh-0.5,0) {$\scriptstyle \xki$};
\draw [thkln] (-\mxsh+0.15,\mysh+0.2) -- (\mxsh-0.15,\mysh+0.2) node [midway,above] {$\scriptstyle \xca-\xki$};
\draw [thkln] (-\mxsh+0.15,-\mysh-0.2) -- (\mxsh-0.15,-\mysh-0.2) node [midway, below] {$\scriptstyle \xca-\xki$};
\draw [ptzert] (-0.15-\mxsh,-0.6-\mysh) rectangle ++(.3,1.2);
\draw [thkln] (-0.75-\mxsh,-\mysh) -- (-0.15-\mxsh,-\mysh)
node [near start,below] {$\scriptstyle \xca$};
\begin{scope}[xscale=-1]
\draw [ptzert] (-0.15-\mxsh,-0.6-\mysh) rectangle ++(.3,1.2);
\draw [thkln] (-0.75-\mxsh,-\mysh) -- (-0.15-\mxsh,-\mysh)
node [near start,below] {$\scriptstyle \xca$};
\end{scope}
\begin{scope}[yscale=-1]
\draw [ptzert] (-0.15-\mxsh,-0.6-\mysh) rectangle ++(.3,1.2);
\draw [thkln] (-0.75-\mxsh,-\mysh) -- (-0.15-\mxsh,-\mysh)
node [near start,above] {$\scriptstyle \xca$};
\end{scope}
\begin{scope}[yscale=-1,xscale=-1]
\draw [ptzert] (-0.15-\mxsh,-0.6-\mysh) rectangle ++(.3,1.2);
\draw [thkln] (-0.75-\mxsh,-\mysh) -- (-0.15-\mxsh,-\mysh)
node [near start,above] {$\scriptstyle \xca$};
\end{scope}
\end{tikzpicture}
\]
The \lumps\ presentation\rx{eq:stbx} of $\ytngsc$ allows us to use Remark\rw{rmk:bndss}: it is sufficient to establish a bound
\begin{equation}
\label{eq:bdck}
\tKHmvv{i}{\hem}(\xDsci) = 0,\quad\text{if $i\leq \xca-\xki $},
\end{equation}
because the tangle $\ytngsci$ has an extra \thdgr\ shift $\shcr^{\xki^2}$ in the \lumps\ formula\rx{eq:stbx} and $\xki^2-\xki\geq 0$.

Consider the portion of $\xDsci$ between the left $\xki$-cable of $\ytngsci$ and another crossing which connects the same \tBcr s:
\def\mxash{6}
\def\mxbsh{1}
\def\myash{2}
\begin{equation}
\label{eq:ordg}
\def\mxsh{3}
\def\mysh{1.6}
\begin{tikzpicture}[menvone]
\draw [bcrc] (-0.75-\mxsh-0.3-\mxash,-\mysh) -- (\mxash+1.6,-\mysh);
\draw [bcrc] (-0.75-\mxsh-0.3-\mxash,\mysh) -- (\mxash+1.6,\mysh);
\draw [thkln] (-\mxbsh+0.15,\mysh+\myash) to [out=0,in=180] (\mxbsh-0.15,\mysh);
\draw [lnovr,middle segment=0.3cm] (-\mxbsh+0.15,\mysh) to [out=0,in=180] (\mxbsh-0.15,\mysh+\myash);
\draw [thkln] (-\mxbsh+0.15,\mysh) to [out=0,in=180] (\mxbsh-0.15,\mysh+\myash);
\draw [ptzert] (-\mxbsh-0.15,\mysh-0.6) rectangle ++(.3,1.2);
\draw [ptzert] (\mxbsh-0.15,\mysh-0.6) rectangle ++(.3,1.2);
\begin{scope}[yshift=\myash cm]
\draw [ptzert] (-\mxbsh-0.15,\mysh-0.6) rectangle ++(.3,1.2);
\draw [ptzert] (\mxbsh-0.15,\mysh-0.6) rectangle ++(.3,1.2);
\draw [thkln] (-0.75-\mxbsh,\mysh) -- (-0.15-\mxbsh,\mysh) node [near start,above] {$\scriptstyle \xca$};
\draw [thkln] (0.75+\mxbsh,\mysh) -- (0.15+\mxbsh,\mysh) node [near start,above] {$\scriptstyle \xca$};
\end{scope}
\begin{scope}[xshift=-\mxash cm-0.3cm]
\draw [thkln] (-\mxsh+0.15,\mysh) to [out=0,in=180] (0,-\mysh);
\draw [lnovr,middle segment=0.3cm] (-\mxsh+0.15,-\mysh) to [out=0,in=180] (0,\mysh);
\draw [thkln] (-\mxsh+0.15,-\mysh) to [out=0,in=180] (0,\mysh);
\draw [ptzert] (-0.15-\mxsh,-0.6-\mysh) rectangle ++(.3,1.2);
\draw [ptzert] (0,-0.6-\mysh) rectangle ++(.3,1.2);
\draw [ptzert] (0,-0.6+\mysh) rectangle ++(.3,1.2);
\draw [thkln] (-0.75-\mxsh,-\mysh) -- (-0.15-\mxsh,-\mysh)
node [near start,below] {$\scriptstyle \xca$};
\begin{scope}[yscale=-1]
\draw [ptzert] (-0.15-\mxsh,-0.6-\mysh) rectangle ++(.3,1.2);
\draw [thkln] (-0.75-\mxsh,-\mysh) -- (-0.15-\mxsh,-\mysh)
node [near start,above] {$\scriptstyle \xca$};
\end{scope}
\end{scope}
\begin{scope}[xshift=\mxash cm]
\draw [ptzert] (0,-0.6-\mysh) rectangle ++(.3,1.2);
\draw [ptzert] (0,-0.6+\mysh) rectangle ++(.3,1.2);
\draw [thkln] (0.3,\mysh-0.2) to [out=0,in=90] (1,0) to [out=-90,in=0] (0.3,-\mysh+0.2);
\node at (1.5,0) {$\scriptstyle \xki$};
\draw [thkln] (0.3,\mysh+0.2) -- (1.6,\mysh+0.2) node [near end, above] {$\scriptstyle \xca - \xki $};
\draw [thkln] (0.3,-\mysh-0.2) -- (1.6,-\mysh-0.2) node [near end, below] {$\scriptstyle \xca - \xki $};
\end{scope}
\draw [thkln] (-\mxash,-\mysh) -- (-0.7,-\mysh) (0.7,-\mysh) -- (\mxash,-\mysh);
\node at (0,-\mysh) {$\cdots$};
\draw [thkln] (-\mxash,\mysh) -- (-0.5*\mxash-0.5*\mxbsh-0.7,\mysh)
(-0.5*\mxash-0.5*\mxbsh+0.7,\mysh) -- (-\mxbsh-0.15,\mysh);
\node at (-0.5*\mxash-0.5*\mxbsh,\mysh) {$\cdots$ };
\draw [thkln] (\mxash,\mysh) -- (0.5*\mxash+0.5*\mxbsh+0.7,\mysh)
(0.5*\mxash+0.5*\mxbsh-0.7,\mysh) -- (\mxbsh+0.15,\mysh);
\node at (0.5*\mxash+0.5*\mxbsh,\mysh) {$\cdots$ };
\end{tikzpicture}
\end{equation}
As usual, gray strips indicate \tBcr s of the \tBdg. We showed explicitly one of the crossings attached to a \tBcr s. Consider a modification of this diagram which results from a repeated application of
of Propositions\rw{prp:bndfs} and\rw{prp:bndsc} to these crossings:
\begin{equation}
\label{eq:mddg}
\def\mxsh{3}
\def\mysh{1.6}
\def\smshy{0.35}
\begin{tikzpicture}[menvone]
\draw [bcrc] (-0.75-\mxsh-0.3-\mxash,-\mysh) -- (\mxash+1.6,-\mysh);
\draw [bcrc] (-0.75-\mxsh-0.3-\mxash,\mysh) -- (\mxash+1.6,\mysh);
\begin{scope}[yshift=\smshy cm+0.3 cm]
\draw [thkln] (-\mxbsh+0.15,\mysh+\myash-0.2) to [out=0,in=180] (\mxbsh-0.15,\mysh);
\draw [lnovr,middle segment=0.3cm] (-\mxbsh+0.15,\mysh) to [out=0,in=180] (\mxbsh-0.15,\mysh+\myash-0.2);
\draw [thkln] (-\mxbsh+0.15,\mysh) to [out=0,in=180] (\mxbsh-0.15,\mysh+\myash-0.2);
\draw [thkln] (-\mxbsh+0.15,\mysh+\myash+0.2) -- (\mxbsh-0.15,\mysh+\myash+0.2)
node [midway,above] {$\scriptstyle \xki$};
\draw [ptzert] (-\mxbsh-0.15,\mysh-0.6) rectangle ++(.3,1.2);
\draw [ptzert] (\mxbsh-0.15,\mysh-0.6) rectangle ++(.3,1.2);
\begin{scope}[yshift=\myash cm]
\draw [ptzert] (-\mxbsh-0.15,\mysh-0.6) rectangle ++(.3,1.2);
\draw [ptzert] (\mxbsh-0.15,\mysh-0.6) rectangle ++(.3,1.2);
\draw [thkln] (-0.75-\mxbsh,\mysh) -- (-0.15-\mxbsh,\mysh) node [near start,above] {$\scriptstyle \xca$};
\draw [thkln] (0.75+\mxbsh,\mysh) -- (0.15+\mxbsh,\mysh) node [near start,above] {$\scriptstyle \xca$};
\end{scope}
\end{scope}
\begin{scope}[xshift=-\mxash cm-0.3cm]
\draw [thkln] (-\mxsh+0.15,\mysh) to [out=0,in=180] (0,-\mysh);
\draw [lnovr,middle segment=0.3cm] (-\mxsh+0.15,-\mysh) to [out=0,in=180] (0,\mysh);
\draw [thkln] (-\mxsh+0.15,-\mysh) to [out=0,in=180] (0,\mysh);
\draw [ptzert] (-0.15-\mxsh,-0.6-\mysh) rectangle ++(.3,1.2);
\draw [ptzert] (0,-0.6-\mysh) rectangle ++(.3,1.2);
\draw [ptzert] (0,-0.6+\mysh) rectangle ++(.3,1.2);
\draw [thkln] (-0.75-\mxsh,-\mysh) -- (-0.15-\mxsh,-\mysh)
node [near start,below] {$\scriptstyle \xca$};
\begin{scope}[yscale=-1]
\draw [ptzert] (-0.15-\mxsh,-0.6-\mysh) rectangle ++(.3,1.2);
\draw [thkln] (-0.75-\mxsh,-\mysh) -- (-0.15-\mxsh,-\mysh)
node [near start,above] {$\scriptstyle \xca$};
\end{scope}
\end{scope}
\begin{scope}[xshift=\mxash cm]
\draw [ptzert] (0,-0.4-\mysh-\smshy) rectangle ++(.3,0.8);
\draw [ptzert] (0,-0.4+\mysh+\smshy) rectangle ++(.3,0.8);
\draw [thkln] (0.3,\mysh-\smshy) to [out=0,in=90] (1,0) to [out=-90,in=0] (0.3,-\mysh+\smshy);
\node at (1.5,0) {$\scriptstyle \xki$};
\draw [thkln] (0.3,\mysh+\smshy) -- (1.6,\mysh+\smshy) node [near end, above] {$\scriptstyle \xca - \xki $};
\draw [thkln] (0.3,-\mysh-\smshy) -- (1.6,-\mysh-\smshy) node [near end, below] {$\scriptstyle \xca - \xki $};
\end{scope}
\begin{scope}[yshift=-\smshy cm]
\draw [thkln] (-\mxash,-\mysh) -- (-0.7,-\mysh) (0.7,-\mysh) -- (\mxash,-\mysh);
\node at (0,-\mysh) {$\cdots$};
\end{scope}
\draw [thkln] (-\mxash,-\mysh+\smshy) -- (0.3+\mxash,-\mysh+\smshy);;
\begin{scope}[yshift=\smshy cm]
\draw [thkln] (-\mxash,\mysh) -- (-0.5*\mxash-0.5*\mxbsh-0.7,\mysh)
(-0.5*\mxash-0.5*\mxbsh+0.7,\mysh) -- (-\mxbsh-0.15,\mysh);
\node at (-0.5*\mxash-0.5*\mxbsh,\mysh) {$\cdots$ };
\draw [thkln] (\mxash,\mysh) -- (0.5*\mxash+0.5*\mxbsh+0.7,\mysh)
(0.5*\mxash+0.5*\mxbsh-0.7,\mysh) -- (\mxbsh+0.15,\mysh);
\node at (0.5*\mxash+0.5*\mxbsh,\mysh) {$\cdots$ };
\end{scope}
\draw [thkln] (-\mxash,\mysh-\smshy) -- (\mxash+0.3,\mysh-\smshy);
\end{tikzpicture}
\end{equation}
Let $\xDscip$ be the diagram constructed from $\xDsci$ by replacing the subdiagram\rx{eq:ordg} with\rx{eq:mddg}. According to Propositions\rw{prp:bndfs} and\rw{prp:bndsc}
there is a map of shifted Khovanov homologies $\tKHm(\xDscip)\rightarrow \tKHm(\xDsci)$, which is an isomorphism on $\tKHmvv{i}{\hem}(\xDsci)$ for $i\leq \xca-\xki$. We are going to show that
\begin{equation}
\label{eq:bdckb}
\tKHmvv{i}{\hem}(\xDscip) = 0\quad \text{if $i\leq \xca-\xki $},
\end{equation}
hence this will imply the bound\rx{eq:bdck}.

Consider a sequence of homotopy equivalences:
\begin{equation}
\label{eq:xthm}
\begin{tikzpicture}[menvone]
\begin{scope}[xshift=-\mxsh cm]
\draw [thkln] (-\mxsh+0.15,\mysh) to [out=0,in=180] (0,-\mysh);
\draw [lnovr] (-\mxsh+0.15,-\mysh) to [out=0,in=180] (0,\mysh);
\draw [thkln] (-\mxsh+0.15,-\mysh) to [out=0,in=180] (0,\mysh);
\draw [ptzer] (-0.15-\mxsh,-0.6-\mysh) rectangle ++(.3,1.2);
\draw [thkln] (-0.75-\mxsh,-\mysh) -- (-0.15-\mxsh,-\mysh)
node [near start,below] {$\scriptstyle \xca$};
\begin{scope}[yscale=-1]
\draw [ptzer] (-0.15-\mxsh,-0.6-\mysh) rectangle ++(.3,1.2);
\draw [thkln] (-0.75-\mxsh,-\mysh) -- (-0.15-\mxsh,-\mysh)
node [near start,above] {$\scriptstyle \xca$};
\end{scope}
\end{scope}
\draw [thkln] (-\mxsh+0.15,\mysh-0.2) to [out=0,in=90] (-\mxsh+1.1,0)
 to [out=-90,in=0] (-\mxsh+0.15,-\mysh+0.2);
\node at (-\mxsh+1.5,0) {$\scriptstyle \xki$};
\draw [thkln] (-\mxsh+0.15,\mysh+0.2) -- (0.3,\mysh+0.2) node [midway,above] {$\scriptstyle \xca-\xki$};
\draw [thkln] (-\mxsh+0.15,-\mysh-0.2) -- (0.3,-\mysh-0.2) node [midway, below] {$\scriptstyle \xca-\xki$};
\draw [ptzert] (-0.15-\mxsh,-0.6-\mysh) rectangle ++(.3,1.2);
\begin{scope}[yscale=-1]
\draw [ptzert] (-0.15-\mxsh,-0.6-\mysh) rectangle ++(.3,1.2);
\end{scope}
\end{tikzpicture}
\;\hteqv\;
\begin{tikzpicture}[menvone]
\begin{scope}[xshift=-\mxsh cm]
\draw [thkln] (-\mxsh+0.15,\mysh-0.2) to [out=0,in=180] (0.15+0.75,-\mysh-0.2);
\draw [thkln,yshift=0.4cm] (-\mxsh+0.15,\mysh-0.2) to [out=0,in=180] (0.15+0.75,-\mysh-0.2+0.6);
\draw [lnovrtw,yshift=-0.4cm] (-\mxsh+0.15,-\mysh+0.2) to [out=0,in=180] (0.15+0.75,\mysh+0.2-0.6);
\draw [thkln,yshift=-0.4cm] (-\mxsh+0.15,-\mysh+0.2) to [out=0,in=180] (0.15+0.75,\mysh+0.2-0.6);
\draw [lnovrtw] (-\mxsh+0.15,-\mysh+0.2) to [out=0,in=180] (0.15+0.75,\mysh+0.2);
\draw [thkln] (-\mxsh+0.15,-\mysh+0.2) to [out=0,in=180] (0.15+0.75,\mysh+0.2);
\draw [ptzer] (-0.15-\mxsh,-0.6-\mysh) rectangle ++(.3,1.2);
\draw [thkln] (-0.75-\mxsh,-\mysh) -- (-0.15-\mxsh,-\mysh)
node [near start,below] {$\scriptstyle \xca$};
\begin{scope}[yscale=-1]
\draw [ptzer] (-0.15-\mxsh,-0.6-\mysh) rectangle ++(.3,1.2);
\draw [thkln] (-0.75-\mxsh,-\mysh) -- (-0.15-\mxsh,-\mysh)
node [near start,above] {$\scriptstyle \xca$};
\end{scope}
\end{scope}
\draw [thkln,xshift=0.75cm] (-\mxsh+0.15,\mysh-0.2-0.6) to [out=0,in=90] (-\mxsh+0.5,0)
 to [out=-90,in=0] (-\mxsh+0.15,-\mysh+0.2+0.6) ;
\node at (-\mxsh+1.65,0) {$\scriptstyle \xki$};
\draw [thkln] (-\mxsh+0.15+0.75,\mysh+0.2) -- (0.3,\mysh+0.2) node [midway,above] {$\scriptstyle \xca-\xki$};
\draw [thkln] (-\mxsh+0.15+0.75,-\mysh-0.2) -- (0.3,-\mysh-0.2) node [midway, below] {$\scriptstyle \xca-\xki$};
%
\begin{scope}[yscale=-1]
\end{scope}
\end{tikzpicture}
\;
\hteqv
\;
\shcr^{\xca \xki - \shlf \xki^2}\shfr^\xki
\begin{tikzpicture}[menvone]
\draw [thkln] (-\mxsh+0.15,\mysh-0.2) to [out=0,in=90] (-\mxsh+0.8,0)
 to [out=-90,in=0] (-\mxsh+0.15,-\mysh+0.2) ;
\node at (-\mxsh+0.4,0) {$\scriptstyle \xki$};
\draw [thkln] (-\mxsh+0.15,\mysh+0.2) to [out=0,in=180]
(1.2,-\mysh);
\draw [thkln] (1.2,-\mysh) -- (2,-\mysh)
node [near start,below] {$\scriptstyle \xca-\xki$};
\draw [lnovr] (-\mxsh+0.15,-\mysh-0.2) to [out=0,in=180] (1.2,\mysh);
\draw [thkln] (-\mxsh+0.15,-\mysh-0.2) to [out=0,in=180]
(1.2,\mysh);
\draw [thkln] (1.2,\mysh) -- (2,\mysh)
node [near start,above] {$\scriptstyle \xca-\xki$};
%
\draw [ptzert] (-0.15-\mxsh,-0.6-\mysh) rectangle ++(.3,1.2);
\draw [thkln] (-0.75-\mxsh,-\mysh) -- (-0.15-\mxsh,-\mysh)
node [near start,below] {$\scriptstyle \xca$};
\begin{scope}[yscale=-1]
\draw [ptzert] (-0.15-\mxsh,-0.6-\mysh) rectangle ++(.3,1.2);
\draw [thkln] (-0.75-\mxsh,-\mysh) -- (-0.15-\mxsh,-\mysh)
node [near start,above] {$\scriptstyle \xca$};
\end{scope}
\end{tikzpicture}
\end{equation}
Here the first homotopy equivalence comes from sliding $\xki$-cable projectors to the left along $\xca$-cables and then contracting double projectors into single ones, while the second homotopy comes from \eex{eq:cfrsh} and\rx{eq:projtw}. Let $\xDscipp$ be the diagram constructed from $\xDscip$ by replacing the left tangle of \ex{eq:xthm} with the right tangle. Since $\yncrv{\xDscipp} = \yncrv{\xDscip} + (\xca-\xki)^2 - \xca^2$, homotopy equivalence\rx{eq:xthm} implies the isomorphism of shifted Khovanov homologies
\[
\tKHm(\xDscip) = \shcr^{\xki(2\xca-\xki)}
\shfr^\xki\,\tKHm(\xDscipp)
\]
and the bound\rx{eq:bdckb} follows from Theorem\rw{thm:eaest}.
\end{proof}

\appendix
\section{A single crossing of colored strands approximates a \tJWp}

Let $\xD$ be a diagram which may include both single and cabled lines as well as \tJWp s. Suppose that $\xD$ has a crossing of two $\xca$-cables with a projector on each. Let $\xDp$ be a diagram, in which the crossing is replaced with a \tJWp:
\def\mxsh{2}
\def\mysh{1.2}
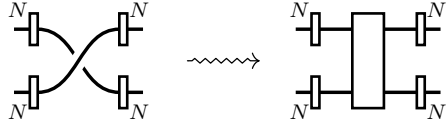
\[
\begin{tikzpicture}[menvthree]
\node (l) at (-1.8,0) {};
\node (r) at (1.8,0) {};
\path[commutative diagrams/.cd, every arrow, every label]
(l) edge[commutative diagrams/squiggly] (r);
\begin{scope}[xshift=-5.5cm]
\def\mxsh{1.7}
\draw [thkln] (-\mxsh+0.15,\mysh) to [out=0,in=180] (\mxsh-0.15,-\mysh);
\draw [lnovr] (-\mxsh+0.15,-\mysh) to [out=0,in=180] (\mxsh-0.15,\mysh);
\draw [thkln] (-\mxsh+0.15,-\mysh) to [out=0,in=180] (\mxsh-0.15,\mysh);
%
%
\draw [ptzer] (-0.15-\mxsh,-0.6-\mysh) rectangle ++(.3,1.2);
\draw [thkln] (-0.75-\mxsh,-\mysh) -- (-0.15-\mxsh,-\mysh)
node [near start,below] {$\scriptstyle \xca$};
\begin{scope}[xscale=-1]
\draw [ptzer] (-0.15-\mxsh,-0.6-\mysh) rectangle ++(.3,1.2);
\draw [thkln] (-0.75-\mxsh,-\mysh) -- (-0.15-\mxsh,-\mysh)
node [near start,below] {$\scriptstyle \xca$};
\end{scope}
\begin{scope}[yscale=-1]
\draw [ptzer] (-0.15-\mxsh,-0.6-\mysh) rectangle ++(.3,1.2);
\draw [thkln] (-0.75-\mxsh,-\mysh) -- (-0.15-\mxsh,-\mysh)
node [near start,above] {$\scriptstyle \xca$};
\end{scope}
\begin{scope}[yscale=-1,xscale=-1]
\draw [ptzer] (-0.15-\mxsh,-0.6-\mysh) rectangle ++(.3,1.2);
\draw [thkln] (-0.75-\mxsh,-\mysh) -- (-0.15-\mxsh,-\mysh)
node [near start,above] {$\scriptstyle \xca$};
\end{scope}
\end{scope}
\begin{scope}[xshift=5.5cm]
\draw [thkln] (-\mxsh+0.15,\mysh) -- (-\mxsh*0.3,\mysh) (\mxsh*0.3,\mysh) -- (\mxsh-0.15,\mysh);
\draw [thkln] (-\mxsh+0.15,-\mysh) -- (-\mxsh*0.3,-\mysh) (\mxsh*0.3,-\mysh) -- (\mxsh-0.15,-\mysh);
%
%
\draw [ptzer] (-0.15-\mxsh,-0.6-\mysh) rectangle ++(.3,1.2);
\draw [thkln] (-0.75-\mxsh,-\mysh) -- (-0.15-\mxsh,-\mysh)
node [near start,below] {$\scriptstyle \xca$};
\begin{scope}[xscale=-1]
\draw [ptzer] (-0.15-\mxsh,-0.6-\mysh) rectangle ++(.3,1.2);
\draw [thkln] (-0.75-\mxsh,-\mysh) -- (-0.15-\mxsh,-\mysh)
node [near start,below] {$\scriptstyle \xca$};
\end{scope}
\begin{scope}[yscale=-1]
\draw [ptzer] (-0.15-\mxsh,-0.6-\mysh) rectangle ++(.3,1.2);
\draw [thkln] (-0.75-\mxsh,-\mysh) -- (-0.15-\mxsh,-\mysh)
node [near start,above] {$\scriptstyle \xca$};
\end{scope}
\begin{scope}[yscale=-1,xscale=-1]
\draw [ptzer] (-0.15-\mxsh,-0.6-\mysh) rectangle ++(.3,1.2);
\draw [thkln] (-0.75-\mxsh,-\mysh) -- (-0.15-\mxsh,-\mysh)
node [near start,above] {$\scriptstyle \xca$};
\end{scope}
\draw[ptzer] (-\mxsh*0.3,-\mysh*1-0.6) rectangle ++(\mxsh*0.6,\mysh*2+1.2);
\end{scope}
\end{tikzpicture}
\]

\begin{theorem}
\label{thm:crpr}
There exists a map $\tKHm(\xD)\xrightarrow{\xmg}\tKHm(\xDp)$ which is an isomorphism on $\tKHmvv{i}{\hem}$ for $i\leq 2\xca-2$.
\end{theorem}
\begin{proof}
Consider three tangles
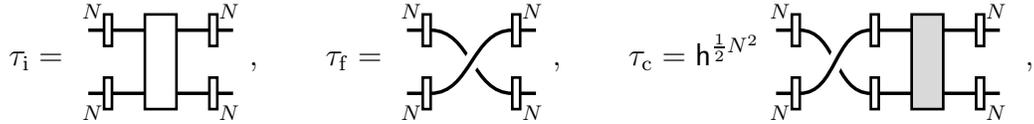
\[
\ytngsi =
\begin{tikzpicture}[menvthree]
\draw [thkln] (-\mxsh+0.15,\mysh) -- (-\mxsh*0.3,\mysh) (\mxsh*0.3,\mysh) -- (\mxsh-0.15,\mysh);
\draw [thkln] (-\mxsh+0.15,-\mysh) -- (-\mxsh*0.3,-\mysh) (\mxsh*0.3,-\mysh) -- (\mxsh-0.15,-\mysh);
\draw [ptzer] (-0.15-\mxsh,-0.6-\mysh) rectangle ++(.3,1.2);
\draw [thkln] (-0.75-\mxsh,-\mysh) -- (-0.15-\mxsh,-\mysh)
node [near start,below] {$\scriptstyle \xca$};
\begin{scope}[xscale=-1]
\draw [ptzer] (-0.15-\mxsh,-0.6-\mysh) rectangle ++(.3,1.2);
\draw [thkln] (-0.75-\mxsh,-\mysh) -- (-0.15-\mxsh,-\mysh)
node [near start,below] {$\scriptstyle \xca$};
\end{scope}
\begin{scope}[yscale=-1]
\draw [ptzer] (-0.15-\mxsh,-0.6-\mysh) rectangle ++(.3,1.2);
\draw [thkln] (-0.75-\mxsh,-\mysh) -- (-0.15-\mxsh,-\mysh)
node [near start,above] {$\scriptstyle \xca$};
\end{scope}
\begin{scope}[yscale=-1,xscale=-1]
\draw [ptzer] (-0.15-\mxsh,-0.6-\mysh) rectangle ++(.3,1.2);
\draw [thkln] (-0.75-\mxsh,-\mysh) -- (-0.15-\mxsh,-\mysh)
node [near start,above] {$\scriptstyle \xca$};
\end{scope}
\draw[ptzer] (-\mxsh*0.3,-\mysh*1-0.6) rectangle ++(\mxsh*0.6,\mysh*2+1.2);
\end{tikzpicture}
,\qquad
\def\mxsh{1.7}
\ytngsf =
\begin{tikzpicture}[menvthree]
\draw [thkln] (-\mxsh+0.15,\mysh) to [out=0,in=180] (\mxsh-0.15,-\mysh);
\draw [lnovr] (-\mxsh+0.15,-\mysh) to [out=0,in=180] (\mxsh-0.15,\mysh);
\draw [thkln] (-\mxsh+0.15,-\mysh) to [out=0,in=180] (\mxsh-0.15,\mysh);
%
%
\draw [ptzer] (-0.15-\mxsh,-0.6-\mysh) rectangle ++(.3,1.2);
\draw [thkln] (-0.75-\mxsh,-\mysh) -- (-0.15-\mxsh,-\mysh)
node [near start,below] {$\scriptstyle \xca$};
\begin{scope}[xscale=-1]
\draw [ptzer] (-0.15-\mxsh,-0.6-\mysh) rectangle ++(.3,1.2);
\draw [thkln] (-0.75-\mxsh,-\mysh) -- (-0.15-\mxsh,-\mysh)
node [near start,below] {$\scriptstyle \xca$};
\end{scope}
\begin{scope}[yscale=-1]
\draw [ptzer] (-0.15-\mxsh,-0.6-\mysh) rectangle ++(.3,1.2);
\draw [thkln] (-0.75-\mxsh,-\mysh) -- (-0.15-\mxsh,-\mysh)
node [near start,above] {$\scriptstyle \xca$};
\end{scope}
\begin{scope}[yscale=-1,xscale=-1]
\draw [ptzer] (-0.15-\mxsh,-0.6-\mysh) rectangle ++(.3,1.2);
\draw [thkln] (-0.75-\mxsh,-\mysh) -- (-0.15-\mxsh,-\mysh)
node [near start,above] {$\scriptstyle \xca$};
\end{scope}
\end{tikzpicture}
,\qquad
\def\mxsh{2}
\ytngsc= \shcr^{\hlf\xca^2}
\begin{tikzpicture}[menvthree]
\draw [thkln] (-\mxsh+0.15,\mysh) -- (-\mxsh*0.3,\mysh) (\mxsh*0.3,\mysh) -- (\mxsh-0.15,\mysh);
\draw [thkln] (-\mxsh+0.15,-\mysh) -- (-\mxsh*0.3,-\mysh) (\mxsh*0.3,-\mysh) -- (\mxsh-0.15,-\mysh);
\def\xmlf{2.5}
\draw [thkln] (-0.15-\mxsh,-\mysh) to [out=180, in=0] (-\xmlf*\mxsh+0.15,\mysh);
\draw [lnovr] (-0.15-\mxsh,\mysh) to [out=180, in=0] (-\xmlf*\mxsh+0.15,-\mysh);
\draw [thkln] (-0.15-\mxsh,\mysh) to [out=180, in=0] (-\xmlf*\mxsh+0.15,-\mysh);
\draw [ptzer] (-0.15-\xmlf*\mxsh,-0.6-\mysh) rectangle ++(.3,1.2);
\draw [thkln] (-0.75-\xmlf*\mxsh,-\mysh) -- (-0.15-\xmlf*\mxsh,-\mysh)
node [near start,below] {$\scriptstyle \xca$};
\begin{scope}[yscale=-1]
\draw [ptzer] (-0.15-\xmlf*\mxsh,-0.6-\mysh) rectangle ++(.3,1.2);
\draw [thkln] (-0.75-\xmlf*\mxsh,-\mysh) -- (-0.15-\xmlf*\mxsh,-\mysh)
node [near start,above] {$\scriptstyle \xca$};
\end{scope}
\draw [ptzer] (-0.15-\mxsh,-0.6-\mysh) rectangle ++(.3,1.2);
\begin{scope}[xscale=-1]
\draw [ptzer] (-0.15-\mxsh,-0.6-\mysh) rectangle ++(.3,1.2);
\draw [thkln] (-0.75-\mxsh,-\mysh) -- (-0.15-\mxsh,-\mysh)
node [near start,below] {$\scriptstyle \xca$};
\end{scope}
\begin{scope}[yscale=-1]
\draw [ptzer] (-0.15-\mxsh,-0.6-\mysh) rectangle ++(.3,1.2);
\end{scope}
\begin{scope}[yscale=-1,xscale=-1]
\draw [ptzer] (-0.15-\mxsh,-0.6-\mysh) rectangle ++(.3,1.2);
\draw [thkln] (-0.75-\mxsh,-\mysh) -- (-0.15-\mxsh,-\mysh)
node [near start,above] {$\scriptstyle \xca$};
\end{scope}
\draw[ptone] (-\mxsh*0.3,-\mysh*1-0.6) rectangle ++(\mxsh*0.6,\mysh*2+1.2);
\end{tikzpicture}\;,
\]
where the complex
$
\begin{tikzpicture}[scale=0.4,baseline=-2.5]
\draw[line width=\ljwp,fill=gray!30] (-0.15,-0.6) rectangle ++(0.3,1.2);
\draw [line width=\cblth] (-.65,0) -- (-0.15,0) node [near start, above] {$\scriptstyle a$}
(0.15,0) -- (.65,0);
\end{tikzpicture}
$
is defined by \ex{eq:grproj}, and set $\xKhv{\ytngsfp}=\shcr^{\hlf\xca^2}\xKhv{\ytngsf}$. These tangles have a relation\rx{eq:cnrel} which comes from a sequence of homotopy equivalences:
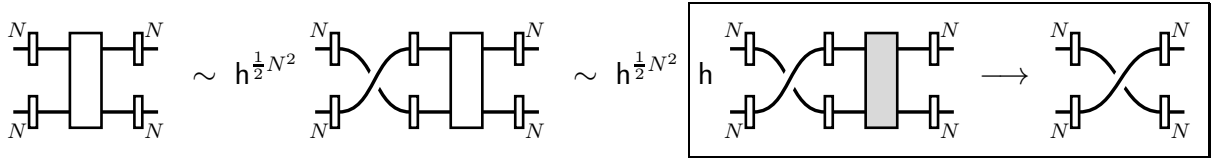
\[
\begin{tikzpicture}[menvthree]
\draw [thkln] (-\mxsh+0.15,\mysh) -- (-\mxsh*0.3,\mysh) (\mxsh*0.3,\mysh) -- (\mxsh-0.15,\mysh);
\draw [thkln] (-\mxsh+0.15,-\mysh) -- (-\mxsh*0.3,-\mysh) (\mxsh*0.3,-\mysh) -- (\mxsh-0.15,-\mysh);
\draw [ptzer] (-0.15-\mxsh,-0.6-\mysh) rectangle ++(.3,1.2);
\draw [thkln] (-0.75-\mxsh,-\mysh) -- (-0.15-\mxsh,-\mysh)
node [near start,below] {$\scriptstyle \xca$};
\begin{scope}[xscale=-1]
\draw [ptzer] (-0.15-\mxsh,-0.6-\mysh) rectangle ++(.3,1.2);
\draw [thkln] (-0.75-\mxsh,-\mysh) -- (-0.15-\mxsh,-\mysh)
node [near start,below] {$\scriptstyle \xca$};
\end{scope}
\begin{scope}[yscale=-1]
\draw [ptzer] (-0.15-\mxsh,-0.6-\mysh) rectangle ++(.3,1.2);
\draw [thkln] (-0.75-\mxsh,-\mysh) -- (-0.15-\mxsh,-\mysh)
node [near start,above] {$\scriptstyle \xca$};
\end{scope}
\begin{scope}[yscale=-1,xscale=-1]
\draw [ptzer] (-0.15-\mxsh,-0.6-\mysh) rectangle ++(.3,1.2);
\draw [thkln] (-0.75-\mxsh,-\mysh) -- (-0.15-\mxsh,-\mysh)
node [near start,above] {$\scriptstyle \xca$};
\end{scope}
\draw[ptzer] (-\mxsh*0.3,-\mysh*1-0.6) rectangle ++(\mxsh*0.6,\mysh*2+1.2);
\end{tikzpicture}
\;\hteqv\;
\shcr^{\hlf\xca^2}
\begin{tikzpicture}[menvthree]
\draw [thkln] (-\mxsh+0.15,\mysh) -- (-\mxsh*0.3,\mysh) (\mxsh*0.3,\mysh) -- (\mxsh-0.15,\mysh);
\draw [thkln] (-\mxsh+0.15,-\mysh) -- (-\mxsh*0.3,-\mysh) (\mxsh*0.3,-\mysh) -- (\mxsh-0.15,-\mysh);
\def\xmlf{2.5}
\draw [thkln] (-0.15-\mxsh,-\mysh) to [out=180, in=0] (-\xmlf*\mxsh+0.15,\mysh);
\draw [lnovr] (-0.15-\mxsh,\mysh) to [out=180, in=0] (-\xmlf*\mxsh+0.15,-\mysh);
\draw [thkln] (-0.15-\mxsh,\mysh) to [out=180, in=0] (-\xmlf*\mxsh+0.15,-\mysh);
\draw [ptzer] (-0.15-\xmlf*\mxsh,-0.6-\mysh) rectangle ++(.3,1.2);
\draw [thkln] (-0.75-\xmlf*\mxsh,-\mysh) -- (-0.15-\xmlf*\mxsh,-\mysh)
node [near start,below] {$\scriptstyle \xca$};
\begin{scope}[yscale=-1]
\draw [ptzer] (-0.15-\xmlf*\mxsh,-0.6-\mysh) rectangle ++(.3,1.2);
\draw [thkln] (-0.75-\xmlf*\mxsh,-\mysh) -- (-0.15-\xmlf*\mxsh,-\mysh)
node [near start,above] {$\scriptstyle \xca$};
\end{scope}
\draw [ptzer] (-0.15-\mxsh,-0.6-\mysh) rectangle ++(.3,1.2);
\begin{scope}[xscale=-1]
\draw [ptzer] (-0.15-\mxsh,-0.6-\mysh) rectangle ++(.3,1.2);
\draw [thkln] (-0.75-\mxsh,-\mysh) -- (-0.15-\mxsh,-\mysh)
node [near start,below] {$\scriptstyle \xca$};
\end{scope}
\begin{scope}[yscale=-1]
\draw [ptzer] (-0.15-\mxsh,-0.6-\mysh) rectangle ++(.3,1.2);
\end{scope}
\begin{scope}[yscale=-1,xscale=-1]
\draw [ptzer] (-0.15-\mxsh,-0.6-\mysh) rectangle ++(.3,1.2);
\draw [thkln] (-0.75-\mxsh,-\mysh) -- (-0.15-\mxsh,-\mysh)
node [near start,above] {$\scriptstyle \xca$};
\end{scope}
\draw[ptzer] (-\mxsh*0.3,-\mysh*1-0.6) rectangle ++(\mxsh*0.6,\mysh*2+1.2);
\end{tikzpicture}
\;\hteqv\;
\shcr^{\hlf\xca^2}\;
\boxed{
\shcr
\begin{tikzpicture}[menvthree]
\draw [thkln] (-\mxsh+0.15,\mysh) -- (-\mxsh*0.3,\mysh) (\mxsh*0.3,\mysh) -- (\mxsh-0.15,\mysh);
\draw [thkln] (-\mxsh+0.15,-\mysh) -- (-\mxsh*0.3,-\mysh) (\mxsh*0.3,-\mysh) -- (\mxsh-0.15,-\mysh);
\def\xmlf{2.5}
\draw [thkln] (-0.15-\mxsh,-\mysh) to [out=180, in=0] (-\xmlf*\mxsh+0.15,\mysh);
\draw [lnovr] (-0.15-\mxsh,\mysh) to [out=180, in=0] (-\xmlf*\mxsh+0.15,-\mysh);
\draw [thkln] (-0.15-\mxsh,\mysh) to [out=180, in=0] (-\xmlf*\mxsh+0.15,-\mysh);
\draw [ptzer] (-0.15-\xmlf*\mxsh,-0.6-\mysh) rectangle ++(.3,1.2);
\draw [thkln] (-0.75-\xmlf*\mxsh,-\mysh) -- (-0.15-\xmlf*\mxsh,-\mysh)
node [near start,below] {$\scriptstyle \xca$};
\begin{scope}[yscale=-1]
\draw [ptzer] (-0.15-\xmlf*\mxsh,-0.6-\mysh) rectangle ++(.3,1.2);
\draw [thkln] (-0.75-\xmlf*\mxsh,-\mysh) -- (-0.15-\xmlf*\mxsh,-\mysh)
node [near start,above] {$\scriptstyle \xca$};
\end{scope}
\draw [ptzer] (-0.15-\mxsh,-0.6-\mysh) rectangle ++(.3,1.2);
\begin{scope}[xscale=-1]
\draw [ptzer] (-0.15-\mxsh,-0.6-\mysh) rectangle ++(.3,1.2);
\draw [thkln] (-0.75-\mxsh,-\mysh) -- (-0.15-\mxsh,-\mysh)
node [near start,below] {$\scriptstyle \xca$};
\end{scope}
\begin{scope}[yscale=-1]
\draw [ptzer] (-0.15-\mxsh,-0.6-\mysh) rectangle ++(.3,1.2);
\end{scope}
\begin{scope}[yscale=-1,xscale=-1]
\draw [ptzer] (-0.15-\mxsh,-0.6-\mysh) rectangle ++(.3,1.2);
\draw [thkln] (-0.75-\mxsh,-\mysh) -- (-0.15-\mxsh,-\mysh)
node [near start,above] {$\scriptstyle \xca$};
\end{scope}
\draw[ptone] (-\mxsh*0.3,-\mysh*1-0.6) rectangle ++(\mxsh*0.6,\mysh*2+1.2);
\end{tikzpicture}
\longrightarrow
\begin{tikzpicture}[menvthree]
\def\mxsh{1.7}
\draw [thkln] (-\mxsh+0.15,\mysh) to [out=0,in=180] (\mxsh-0.15,-\mysh);
\draw [lnovr] (-\mxsh+0.15,-\mysh) to [out=0,in=180] (\mxsh-0.15,\mysh);
\draw [thkln] (-\mxsh+0.15,-\mysh) to [out=0,in=180] (\mxsh-0.15,\mysh);
%
%
\draw [ptzer] (-0.15-\mxsh,-0.6-\mysh) rectangle ++(.3,1.2);
\draw [thkln] (-0.75-\mxsh,-\mysh) -- (-0.15-\mxsh,-\mysh)
node [near start,below] {$\scriptstyle \xca$};
\begin{scope}[xscale=-1]
\draw [ptzer] (-0.15-\mxsh,-0.6-\mysh) rectangle ++(.3,1.2);
\draw [thkln] (-0.75-\mxsh,-\mysh) -- (-0.15-\mxsh,-\mysh)
node [near start,below] {$\scriptstyle \xca$};
\end{scope}
\begin{scope}[yscale=-1]
\draw [ptzer] (-0.15-\mxsh,-0.6-\mysh) rectangle ++(.3,1.2);
\draw [thkln] (-0.75-\mxsh,-\mysh) -- (-0.15-\mxsh,-\mysh)
node [near start,above] {$\scriptstyle \xca$};
\end{scope}
\begin{scope}[yscale=-1,xscale=-1]
\draw [ptzer] (-0.15-\mxsh,-0.6-\mysh) rectangle ++(.3,1.2);
\draw [thkln] (-0.75-\mxsh,-\mysh) -- (-0.15-\mxsh,-\mysh)
node [near start,above] {$\scriptstyle \xca$};
\end{scope}
\end{tikzpicture}
}
\]
Here the first homotopy equivalence comes from \ex{eq:projtw}, while the second equivalence comes from \ex{eq:projcn}.

In order to put a bound on the homological order of $\xDsc$, we purge the gray box in $\ytngsc$, that is, we contract all constituent \taTLt s, whose cup or cap is connected directly to \tJWp s sitting on $\xca$-cables. After the purge, the complex of $\ytngsc$ takes the form
\begin{eqnarray*}
\ytngsc &\hteqv&
\boxed{\cdots\longrightarrow \bigoplus_{\substack{0\leq\zmj\leq \zmi \\ \xki\geq 1} }
\prmltijk\,\shcr^\zmi \shfr^\zmj
\begin{tikzpicture}[menvthree]
\begin{scope}[xshift=-\mxsh cm]
\draw [thkln] (-\mxsh+0.15,\mysh) to [out=0,in=180] (0,-\mysh);
\draw [lnovr] (-\mxsh+0.15,-\mysh) to [out=0,in=180] (0,\mysh);
\draw [thkln] (-\mxsh+0.15,-\mysh) to [out=0,in=180] (0,\mysh);
\draw [ptzer] (-0.15-\mxsh,-0.6-\mysh) rectangle ++(.3,1.2);
\draw [thkln] (-0.75-\mxsh,-\mysh) -- (-0.15-\mxsh,-\mysh)
node [near start,below] {$\scriptstyle \xca$};
\begin{scope}[yscale=-1]
\draw [ptzer] (-0.15-\mxsh,-0.6-\mysh) rectangle ++(.3,1.2);
\draw [thkln] (-0.75-\mxsh,-\mysh) -- (-0.15-\mxsh,-\mysh)
node [near start,above] {$\scriptstyle \xca$};
\end{scope}
\end{scope}
\draw [thkln] (-\mxsh+0.15,\mysh-0.2) to [out=0,in=90] (-\mxsh+1.1,0)
 to [out=-90,in=0] (-\mxsh+0.15,-\mysh+0.2) ;
\begin{scope}[xscale=-1]
\draw [thkln] (-\mxsh+0.15,\mysh-0.2) to [out=0,in=90] (-\mxsh+1.1,0) to [out=-90,in=0] (-\mxsh+0.15,-\mysh+0.2);
\end{scope}
\node at (-\mxsh+0.5,0) {$\scriptstyle \xki$};
\node at (\mxsh-0.5,0) {$\scriptstyle \xki$};
\draw [thkln] (-\mxsh+0.15,\mysh+0.2) -- (\mxsh-0.15,\mysh+0.2) node [midway,above] {$\scriptstyle \xca-\xki$};
\draw [thkln] (-\mxsh+0.15,-\mysh-0.2) -- (\mxsh-0.15,-\mysh-0.2) node [midway, below] {$\scriptstyle \xca-\xki$};
\draw [ptzert] (-0.15-\mxsh,-0.6-\mysh) rectangle ++(.3,1.2);
\begin{scope}[xscale=-1]
\draw [ptzert] (-0.15-\mxsh,-0.6-\mysh) rectangle ++(.3,1.2);
\draw [thkln] (-0.75-\mxsh,-\mysh) -- (-0.15-\mxsh,-\mysh)
node [near start,below] {$\scriptstyle \xca$};
\end{scope}
\begin{scope}[yscale=-1]
\draw [ptzert] (-0.15-\mxsh,-0.6-\mysh) rectangle ++(.3,1.2);
\end{scope}
\begin{scope}[yscale=-1,xscale=-1]
\draw [ptzert] (-0.15-\mxsh,-0.6-\mysh) rectangle ++(.3,1.2);
\draw [thkln] (-0.75-\mxsh,-\mysh) -- (-0.15-\mxsh,-\mysh)
node [near start,above] {$\scriptstyle \xca$};
\end{scope}
\end{tikzpicture}
\longrightarrow\cdots
}_{\;\xki=1}^{\;\infty}
\\
&\hteqv&
\def\mxsh{3}
\boxed{
\cdots\longrightarrow \bigoplus_{\substack{0\leq\zmj\leq \zmi \\ \xki\geq 1} }
\prmltijk\,\shcr^{\zmi+\xki(\xca-\xki)+\hlf\xki^2} \shfr^{\zmj+\xki}
\begin{tikzpicture}[menvthree]
\draw [thkln] (-\mxsh+0.15,\mysh-0.2) to [out=0,in=90] (-\mxsh+0.8,0)
 to [out=-90,in=0] (-\mxsh+0.15,-\mysh+0.2) ;
\begin{scope}[xscale=-1]
\draw [thkln] (-\mxsh+0.15,\mysh-0.2) to [out=0,in=90] (-\mxsh+0.8,0) to [out=-90,in=0] (-\mxsh+0.15,-\mysh+0.2);
\end{scope}
\node at (-\mxsh+0.4,0) {$\scriptstyle \xki$};
\node at (\mxsh-0.4,0) {$\scriptstyle \xki$};
\draw [thkln] (-\mxsh+0.15,\mysh+0.2) to [out=0,in=180]
node [sloped, near start,above] {$\scriptstyle \xca-\xki$} (\mxsh-0.15,-\mysh-0.2);
\draw [lnovr] (-\mxsh+0.15,-\mysh-0.2) to [out=0,in=180] (\mxsh-0.15,\mysh+0.2);
\draw [thkln] (-\mxsh+0.15,-\mysh-0.2) to [out=0,in=180]
node [sloped,near start,below] {$\scriptstyle \xca-\xki$}
(\mxsh-0.15,\mysh+0.2);
%
\draw [ptzert] (-0.15-\mxsh,-0.6-\mysh) rectangle ++(.3,1.2);
\draw [thkln] (-0.75-\mxsh,-\mysh) -- (-0.15-\mxsh,-\mysh)
node [near start,below] {$\scriptstyle \xca$};
\begin{scope}[xscale=-1]
\draw [ptzert] (-0.15-\mxsh,-0.6-\mysh) rectangle ++(.3,1.2);
\draw [thkln] (-0.75-\mxsh,-\mysh) -- (-0.15-\mxsh,-\mysh)
node [near start,below] {$\scriptstyle \xca$};
\end{scope}
\begin{scope}[yscale=-1]
\draw [ptzert] (-0.15-\mxsh,-0.6-\mysh) rectangle ++(.3,1.2);
\draw [thkln] (-0.75-\mxsh,-\mysh) -- (-0.15-\mxsh,-\mysh)
node [near start,above] {$\scriptstyle \xca$};
\end{scope}
\begin{scope}[yscale=-1,xscale=-1]
\draw [ptzert] (-0.15-\mxsh,-0.6-\mysh) rectangle ++(.3,1.2);
\draw [thkln] (-0.75-\mxsh,-\mysh) -- (-0.15-\mxsh,-\mysh)
node [near start,above] {$\scriptstyle \xca$};
\end{scope}
\end{tikzpicture}
\longrightarrow\cdots
}_{\;\xki=1}^{\;\infty}
\end{eqnarray*}
We used homotopy equivalence\rx{eq:xthm}.
Note that there are no tangles with $k=0$, because the complex\rx{eq:grproj} does not contain identity braids.

Let $\xDsci$ be the diagram $\xDsc$ in which the complex $\ytngsc$ is replaced by the tangle diagram
\[
\def\mxsh{3}
\ytngsci \;=\;
\begin{tikzpicture}[menvthree]
\draw [thkln] (-\mxsh+0.15,\mysh-0.2) to [out=0,in=90] (-\mxsh+0.8,0)
 to [out=-90,in=0] (-\mxsh+0.15,-\mysh+0.2) ;
\begin{scope}[xscale=-1]
\draw [thkln] (-\mxsh+0.15,\mysh-0.2) to [out=0,in=90] (-\mxsh+0.8,0) to [out=-90,in=0] (-\mxsh+0.15,-\mysh+0.2);
\end{scope}
\node at (-\mxsh+0.4,0) {$\scriptstyle \xki$};
\node at (\mxsh-0.4,0) {$\scriptstyle \xki$};
\draw [thkln] (-\mxsh+0.15,\mysh+0.2) to [out=0,in=180]
node [sloped, near start,above] {$\scriptstyle \xca-\xki$} (\mxsh-0.15,-\mysh-0.2);
\draw [lnovr] (-\mxsh+0.15,-\mysh-0.2) to [out=0,in=180] (\mxsh-0.15,\mysh+0.2);
\draw [thkln] (-\mxsh+0.15,-\mysh-0.2) to [out=0,in=180]
node [sloped,near start,below] {$\scriptstyle \xca-\xki$}
(\mxsh-0.15,\mysh+0.2);
%
\draw [ptzert] (-0.15-\mxsh,-0.6-\mysh) rectangle ++(.3,1.2);
\draw [thkln] (-0.75-\mxsh,-\mysh) -- (-0.15-\mxsh,-\mysh)
node [near start,below] {$\scriptstyle \xca$};
\begin{scope}[xscale=-1]
\draw [ptzert] (-0.15-\mxsh,-0.6-\mysh) rectangle ++(.3,1.2);
\draw [thkln] (-0.75-\mxsh,-\mysh) -- (-0.15-\mxsh,-\mysh)
node [near start,below] {$\scriptstyle \xca$};
\end{scope}
\begin{scope}[yscale=-1]
\draw [ptzert] (-0.15-\mxsh,-0.6-\mysh) rectangle ++(.3,1.2);
\draw [thkln] (-0.75-\mxsh,-\mysh) -- (-0.15-\mxsh,-\mysh)
node [near start,above] {$\scriptstyle \xca$};
\end{scope}
\begin{scope}[yscale=-1,xscale=-1]
\draw [ptzert] (-0.15-\mxsh,-0.6-\mysh) rectangle ++(.3,1.2);
\draw [thkln] (-0.75-\mxsh,-\mysh) -- (-0.15-\mxsh,-\mysh)
node [near start,above] {$\scriptstyle \xca$};
\end{scope}
\end{tikzpicture}
\]
According to Theorem\rw{thm:smfr}, $\KHmvv{i}{\hem}(\xDsci)=0$ for $i\leq -\shlf (\xca-\xki)^2 -\shlf\yncrv{\xDsi}-1$, so, by Remark\rw{rmk:bndss}, $\KHmvv{i}{\hem}(\xDsc)=0$ for $i\leq 
-\shlf\yncrv{\xDsi}+2\xca-2$ (here we used inequality $2\xca\xki - \xki^2\geq 2\xca-1$ for $\xki\geq 1$). Now the claim of Theorem\rw{thm:crpr} follows from \ex{eq:ineqo}.
\end{proof}

\end{document}

Let $\xDsci$ be the diagram $\xDsc$ in which the \tmcn\rx{eq:stbx} is replaced by its $i$-th constituent diagram.
\begin{proposition}
\label{prp:smit}
There is a bound $\KHmvv{i}{\hem}(\xDsci)=0$ for $i\geq$.
\end{proposition}

\begin{proof}[Proof of Proposition\rw{prp:smdsc}]

\end{proof}
\begin{proof}[Proof of Proposition\rw{prp:smit}]
Consider a homotopy equivalence transformation of a constituent diagram of the \tmcn\rx{eq:stbx} composed with a crossing tangle on the left:
\[
\begin{tikzpicture}[menvone]
\begin{scope}[xshift=-\mxsh cm]
\draw [thkln] (-\mxsh+0.15,\mysh) to [out=0,in=180] (0,-\mysh);
\draw [lnovr] (-\mxsh+0.15,-\mysh) to [out=0,in=180] (0,\mysh);
\draw [thkln] (-\mxsh+0.15,-\mysh) to [out=0,in=180] (0,\mysh);
\draw [ptzer] (-0.15-\mxsh,-0.6-\mysh) rectangle ++(.3,1.2);
\draw [thkln] (-0.75-\mxsh,-\mysh) -- (-0.15-\mxsh,-\mysh)
node [near start,below] {$\scriptstyle \xca$};
\begin{scope}[yscale=-1]
\draw [ptzer] (-0.15-\mxsh,-0.6-\mysh) rectangle ++(.3,1.2);
\draw [thkln] (-0.75-\mxsh,-\mysh) -- (-0.15-\mxsh,-\mysh)
node [near start,above] {$\scriptstyle \xca$};
\end{scope}
\end{scope}
\draw [thkln] (-\mxsh+0.15,\mysh-0.2) to [out=0,in=90] (-\mxsh+1.1,0)
 to [out=-90,in=0] (-\mxsh+0.15,-\mysh+0.2) ;
\begin{scope}[xscale=-1]
\draw [thkln] (-\mxsh+0.15,\mysh-0.2) to [out=0,in=90] (-\mxsh+1.1,0) to [out=-90,in=0] (-\mxsh+0.15,-\mysh+0.2);
\end{scope}
\node at (-\mxsh+0.5,0) {$\scriptstyle i$};
\node at (\mxsh-0.5,0) {$\scriptstyle i$};
\draw [thkln] (-\mxsh+0.15,\mysh+0.2) -- (\mxsh-0.15,\mysh+0.2) node [midway,above] {$\scriptstyle \xca-i$};
\draw [thkln] (-\mxsh+0.15,-\mysh-0.2) -- (\mxsh-0.15,-\mysh-0.2) node [midway, below] {$\scriptstyle \xca-i$};
\draw [ptzert] (-0.15-\mxsh,-0.6-\mysh) rectangle ++(.3,1.2);
\begin{scope}[xscale=-1]
\draw [ptzert] (-0.15-\mxsh,-0.6-\mysh) rectangle ++(.3,1.2);
\draw [thkln] (-0.75-\mxsh,-\mysh) -- (-0.15-\mxsh,-\mysh)
node [near start,below] {$\scriptstyle \xca$};
\end{scope}
\begin{scope}[yscale=-1]
\draw [ptzert] (-0.15-\mxsh,-0.6-\mysh) rectangle ++(.3,1.2);
\end{scope}
\begin{scope}[yscale=-1,xscale=-1]
\draw [ptzert] (-0.15-\mxsh,-0.6-\mysh) rectangle ++(.3,1.2);
\draw [thkln] (-0.75-\mxsh,-\mysh) -- (-0.15-\mxsh,-\mysh)
node [near start,above] {$\scriptstyle \xca$};
\end{scope}
\end{tikzpicture}
\;\hteqv\;
\begin{tikzpicture}[menvone]
\begin{scope}[xshift=-\mxsh cm]
\draw [thkln] (-\mxsh+0.15,\mysh-0.2) to [out=0,in=180] (0.15+0.75,-\mysh-0.2);
\draw [thkln,yshift=0.4cm] (-\mxsh+0.15,\mysh-0.2) to [out=0,in=180] (0.15+0.75,-\mysh-0.2+0.6);
\draw [lnovrtw,yshift=-0.4cm] (-\mxsh+0.15,-\mysh+0.2) to [out=0,in=180] (0.15+0.75,\mysh+0.2-0.6);
\draw [thkln,yshift=-0.4cm] (-\mxsh+0.15,-\mysh+0.2) to [out=0,in=180] (0.15+0.75,\mysh+0.2-0.6);
\draw [lnovrtw] (-\mxsh+0.15,-\mysh+0.2) to [out=0,in=180] (0.15+0.75,\mysh+0.2);
\draw [thkln] (-\mxsh+0.15,-\mysh+0.2) to [out=0,in=180] (0.15+0.75,\mysh+0.2);
\draw [ptzer] (-0.15-\mxsh,-0.6-\mysh) rectangle ++(.3,1.2);
\draw [thkln] (-0.75-\mxsh,-\mysh) -- (-0.15-\mxsh,-\mysh)
node [near start,below] {$\scriptstyle \xca$};
\begin{scope}[yscale=-1]
\draw [ptzer] (-0.15-\mxsh,-0.6-\mysh) rectangle ++(.3,1.2);
\draw [thkln] (-0.75-\mxsh,-\mysh) -- (-0.15-\mxsh,-\mysh)
node [near start,above] {$\scriptstyle \xca$};
\end{scope}
\end{scope}
\draw [thkln,xshift=0.75cm] (-\mxsh+0.15,\mysh-0.2-0.6) to [out=0,in=90] (-\mxsh+0.5,0)
 to [out=-90,in=0] (-\mxsh+0.15,-\mysh+0.2+0.6) ;
\begin{scope}[xscale=-1]
\draw [thkln] (-\mxsh+0.15,\mysh-0.2) to [out=0,in=90] (-\mxsh+1.1,0) to [out=-90,in=0] (-\mxsh+0.15,-\mysh+0.2);
\end{scope}
\node at (-\mxsh+1.65,0) {$\scriptstyle i$};
\node at (\mxsh-0.65,0) {$\scriptstyle i$};
\draw [thkln] (-\mxsh+0.15+0.75,\mysh+0.2) -- (\mxsh-0.15,\mysh+0.2) node [midway,above] {$\scriptstyle \xca-i$};
\draw [thkln] (-\mxsh+0.15+0.75,-\mysh-0.2) -- (\mxsh-0.15,-\mysh-0.2) node [midway, below] {$\scriptstyle \xca-i$};
%
\begin{scope}[xscale=-1]
\draw [ptzert] (-0.15-\mxsh,-0.6-\mysh) rectangle ++(.3,1.2);
\draw [thkln] (-0.75-\mxsh,-\mysh) -- (-0.15-\mxsh,-\mysh)
node [near start,below] {$\scriptstyle \xca$};
\end{scope}
\begin{scope}[yscale=-1]
\end{scope}
\begin{scope}[yscale=-1,xscale=-1]
\draw [ptzert] (-0.15-\mxsh,-0.6-\mysh) rectangle ++(.3,1.2);
\draw [thkln] (-0.75-\mxsh,-\mysh) -- (-0.15-\mxsh,-\mysh)
node [near start,above] {$\scriptstyle \xca$};
\end{scope}
\end{tikzpicture}
\;
\hteqv
\;
\shcr^{\xca i - \shlf i^2}\shfr^i
\begin{tikzpicture}[menvone]
\draw [thkln] (-\mxsh+0.15,\mysh-0.2) to [out=0,in=90] (-\mxsh+0.8,0)
 to [out=-90,in=0] (-\mxsh+0.15,-\mysh+0.2) ;
\begin{scope}[xscale=-1]
\draw [thkln] (-\mxsh+0.15,\mysh-0.2) to [out=0,in=90] (-\mxsh+0.8,0) to [out=-90,in=0] (-\mxsh+0.15,-\mysh+0.2);
\end{scope}
\node at (-\mxsh+0.4,0) {$\scriptstyle i$};
\node at (\mxsh-0.4,0) {$\scriptstyle i$};
\draw [thkln] (-\mxsh+0.15,\mysh+0.2) to [out=0,in=180]
node [sloped, near start,above] {$\scriptstyle \xca-i$} (\mxsh-0.15,-\mysh-0.2);
\draw [lnovr] (-\mxsh+0.15,-\mysh-0.2) to [out=0,in=180] (\mxsh-0.15,\mysh+0.2);
\draw [thkln] (-\mxsh+0.15,-\mysh-0.2) to [out=0,in=180]
node [sloped,near start,below] {$\scriptstyle \xca-i$}
(\mxsh-0.15,\mysh+0.2);
%
\draw [ptzert] (-0.15-\mxsh,-0.6-\mysh) rectangle ++(.3,1.2);
\draw [thkln] (-0.75-\mxsh,-\mysh) -- (-0.15-\mxsh,-\mysh)
node [near start,below] {$\scriptstyle \xca$};
\begin{scope}[xscale=-1]
\draw [ptzert] (-0.15-\mxsh,-0.6-\mysh) rectangle ++(.3,1.2);
\draw [thkln] (-0.75-\mxsh,-\mysh) -- (-0.15-\mxsh,-\mysh)
node [near start,below] {$\scriptstyle \xca$};
\end{scope}
\begin{scope}[yscale=-1]
\draw [ptzert] (-0.15-\mxsh,-0.6-\mysh) rectangle ++(.3,1.2);
\draw [thkln] (-0.75-\mxsh,-\mysh) -- (-0.15-\mxsh,-\mysh)
node [near start,above] {$\scriptstyle \xca$};
\end{scope}
\begin{scope}[yscale=-1,xscale=-1]
\draw [ptzert] (-0.15-\mxsh,-0.6-\mysh) rectangle ++(.3,1.2);
\draw [thkln] (-0.75-\mxsh,-\mysh) -- (-0.15-\mxsh,-\mysh)
node [near start,above] {$\scriptstyle \xca$};
\end{scope}
\end{tikzpicture}
\]
Here the first homotopy equivalence comes from sliding the middle projectors to the left along $\xca$-cables and then contracting double projectors into single ones, while the second homotopy comes from \eex{eq:cfrsh} and\rx{eq:projtw}.
\end{proof}

\end{document}

In order to prove the remaining bound\rx{eq:bd4a} in the case of \tBadq\ diagram $\xD$, we observe that  apart from inequality it corresponds to the `equal' part of the bound\rx{eq:bd3a}, so relation $\xca \gvD+j =-i -\xabms$ holds only if there are equalities in the second inequality of\rx{eq:mbnd} and in the inequality \ex{eq:bndNt}:
\begin{equation}
\label{eq:threq}
\njcr = \xabms,\qquad \nscrb+\nwcrb=\xca,\qquad \xca_1+\xca_2=2\xca.
\end{equation}
Now consider the case of $\xabms>0$ and the case of $\xabms=0$. If $\xabms>0$, then the diagram $\xDs$ must have at least one \tstrtl. If that line is attached to a \tBcr\ $c$, then it is easy to see that in that circle either $\xca_1<\xca$ or $\xca_2<\xca$, and the third equation of\rx{eq:threq} can't hold. If $\xabms=0$, then $\spvr=0$ for all crossings $\xvrt$ of $\xD$, that is, $\xDs$ contains no \tstrtl s. This means that the diagram $\xDs$ consists of $\xca$-colored circles, each carrying one projector. Let us convert them into circle diagrams by using \eex{eq:projcn} and\rx{eq:grproj} for the projectors. A \twndg\ circle in a \tBcr\ $c$ contains at least two strands of the corresponding $\xca$-cable, so in this case $\nscrb + 2\nwcrb\leq\xca$, hence the second relation of\rx{eq:threq} may hold only if $\nwcrb=0$.
%
%
%
However, since $\xabms=0$, the condition $i+\xabms\neq 0$ of the bound\rx{eq:bd4a} implies $i>0$. A non-trivial \thdgr\ means that the circle diagrams contributing to it must have at least one \taTLt\ $\gamma$ from \ex{eq:grproj} from one of projectors. If this happens on a \tBadq\ circle $c$, then since $\gamma$ has at least one cap and one cup, there must be at least one \twndg\ circle, which contradicts $\nwcrb=0.$\qed

However, the constituent \tTLt s of \tJWp s on a \tBcr\ $c$ produce too many cups and caps, thus potentially creating too many \twndg\ circles in circle diagrams and violating the bound\rx{eq:wbd}.

In order to resolve this difficulty, we will contract the circle diagrams with too many cups and caps within the complex $\xKhv{\xDs}$. Consider a neighborhood of a particular \tBcr\ $c$. We are going to cut temporarily the \tstrtl s connecting it to the rest of the diagram, thus creating a colored tangle diagram $\xtusc$, and then perform some contractions within its Khovanov complex $\xKhv{\xtusc}$.

If the \tstt\ $\spmp$ was such that $c$ had no \tstrtl s attached to it, then we could simply use the projector property of \tJWp s and contract them all to a single projector. In the presence of \tstrtl s we, roughly speaking, contract only parts of the projectors which are connected directly to each other.

We use the following `purging' procedure. We select the initial projector on $c$, the preceding projector being called `final', and purge all other projectors on $c$ one-by-one going clockwise. Purging means that we consider the complex\rx{eq:projcn},\rx{eq:grproj} for the current projector, thus presenting $\xKhv{\xtusc}$ as a \tmcn\ of complexes, each coming from a particular \tTLt\ $\gamma$ in \ex{eq:grproj} (and the identity braid of \ex{eq:projcn}) , and contract all complexes, whose graph $\gamma$ has a cap connected directly to the initial projector or a cup connected directly to the next projector.
%
%
%
%
We claim that after we perform this procedure on all projectors, except the initial and final ones, we get a \tmcn
\begin{equation}
\label{eq:prcmp}
\xKhv{\xtusc} \hteqv
\boxed{
\cdots\longrightarrow\shcr^i
\bigoplus_{0\leq j\leq i}\shfr^j\left( \bigoplus_{\tau} m_{ij,\tau} \xKhv{\tau}\right)
\longrightarrow\cdots
}_{\;i=0}^{\;\infty}
\end{equation}
where the diagrams $\tau$ are of one of two types:
\begin{equation}
\label{eq:twpdgs}
\begin{tikzpicture}[menvtwo]
\draw [bcrct] (-1.25,0) -- (1.25,0) (1.55,0) to [out=0,in=90] (4,-2) to [out=-90,in=0] (0,-5)
(-1.55,0) to [out=180,in=90] (-4,-2) to [out=-90,in=180] (0,-5) ;
\draw [thkln] (-1.25,0) -- (1.25,0) (1.55,0) to [out=0,in=90] (4,-2) to [out=-90,in=0] (0,-5)
(-1.55,0) to [out=180,in=90] (-4,-2) to [out=-90,in=180] (0,-5) ;
\node at (0,-5.5) {$\scriptstyle  \xca_1$};
\node at (0,0.5) {$\scriptstyle \xca_2$};
\draw [ptzert] (-1.55,-0.6) rectangle ++(0.3,1.2);
\draw [ptzert] (1.55,-0.6) rectangle ++(-0.3,1.2);
\draw [ptzer] (0.8,2) rectangle ++(-1.6,-0.8);
\draw [thkc] (0.2,0.1+1.6) arc (0:-180:0.2);
\draw [thkln]  (-0.5,2) -- (-0.5,2.6);
\draw [thkln]  (0.5,2) -- (0.5,2.6);
\node at (0.05,2.4) {$\scriptstyle \cdots$};
\draw [ptzer] (0.8,-2) rectangle ++(-1.6,0.8);
\draw [thkc] (0.2,-0.1-1.6) arc (0:180:0.2);
\draw [thkln]  (-0.5,-2) -- (-0.5,-2.6);
\draw [thkln]  (0.5,-2) -- (0.5,-2.6);
\node at (0.05,-2.4) {$\scriptstyle \cdots$};
\draw [thkln] (-1.25,0.35) to [out=0,in=-90] (-0.5,1.2);
\draw [thkln,xscale=-1] (-1.25,0.35) to [out=0,in=-90] (-0.5,1.2);
\draw [thkln,yscale=-1] (-1.25,0.35) to [out=0,in=-90] (-0.5,1.2);
\draw [thkln,xscale=-1,yscale=-1] (-1.25,0.35) to [out=0,in=-90] (-0.5,1.2);
\draw [thkln] (-1.55,0.35) to [out=180,in=-90] ++(-0.75,0.85) -- (-2.3,2.6);
\draw [thkln,xscale=-1] (-1.55,0.35) to [out=180,in=-90] ++(-0.75,0.85) -- (-2.3,2.6);
\draw [decorate,decoration={brace,amplitude=4pt},xshift=0,yshift=-5pt]
(0.5cm+5pt,-2.6) -- (-0.5cm-5pt,-2.6)
node [black,midway,yshift=-0.4cm]
{\scriptsize to struts};
\draw [decorate,decoration={brace,amplitude=4pt},xshift=0,yshift=5pt]
(-2.3cm-5pt,2.6) -- (2.3cm+5pt,2.6)
node [black,midway,yshift=0.4cm]
{\scriptsize to struts};
\end{tikzpicture}
\qquad
\qquad
\begin{tikzpicture}[menvtwo]
\draw [bcrct] (-1.25,0) -- (1.25,0) (1.55,0) to [out=0,in=90] (4,-2) to [out=-90,in=0] (0,-5)
(-1.55,0) to [out=180,in=90] (-4,-2) to [out=-90,in=180] (0,-5) ;
\draw [thkln] (0,1.2) -- (0,-1.2) (1.55,0) to [out=0,in=90] (4,-2) to [out=-90,in=0] (0,-5)
(-1.55,0) to [out=180,in=90] (-4,-2) to [out=-90,in=180] (0,-5) ;
\node at (0,-5.5) {$\scriptstyle  \xca_1$};
\draw [ptzert] (-1.55,-0.6) rectangle ++(0.3,1.2);
\draw [ptzert] (1.55,-0.6) rectangle ++(-0.3,1.2);
\draw [ptzer] (0.8,2) rectangle ++(-1.6,-0.8);
\draw [thkc] (0.2,0.1+1.6) arc (0:-180:0.2);
\draw [thkln]  (-0.5,2) -- (-0.5,2.6);
\draw [thkln]  (0.5,2) -- (0.5,2.6);
\node at (0.05,2.4) {$\scriptstyle \cdots$};
\draw [ptzer] (0.8,-2) rectangle ++(-1.6,0.8);
\draw [thkc] (0.2,-0.1-1.6) arc (0:180:0.2);
\draw [thkln]  (-0.5,-2) -- (-0.5,-2.6);
\draw [thkln]  (0.5,-2) -- (0.5,-2.6);
\node at (0.05,-2.4) {$\scriptstyle \cdots$};
\draw [thkln] (-1.25,0.35) to [out=0,in=-90] (-0.5,1.2);
\draw [thkln,xscale=-1] (-1.25,0.35) to [out=0,in=-90] (-0.5,1.2);
\draw [thkln,yscale=-1] (-1.25,0.35) to [out=0,in=-90] (-0.5,1.2);
\draw [thkln,xscale=-1,yscale=-1] (-1.25,0.35) to [out=0,in=-90] (-0.5,1.2);
\draw [thkln] (-1.55,0.35) to [out=180,in=-90] ++(-0.75,0.85) -- (-2.3,2.6);
\draw [thkln,xscale=-1] (-1.55,0.35) to [out=180,in=-90] ++(-0.75,0.85) -- (-2.3,2.6);
%
%
\draw [decorate,decoration={brace,amplitude=4pt},xshift=0,yshift=-5pt]
(0.5cm+5pt,-2.6) -- (-0.5cm-5pt,-2.6)
node [black,midway,yshift=-0.4cm]
{\scriptsize to struts};
\draw [decorate,decoration={brace,amplitude=4pt},xshift=0,yshift=5pt]
(-2.3cm-5pt,2.6) -- (2.3cm+5pt,2.6)
node [black,midway,yshift=0.4cm]
{\scriptsize to struts};
\end{tikzpicture}
\end{equation}
with $\xca_1,\xca_2\leq \xca$.
A box with an arc denotes a `\ntwd' \tTLt: for $a \geq b$,
\[
\swdv{
\begin{tikzpicture}[menvone]
\draw [ptzer] (0.4,0.8) rectangle ++(-0.8,-1.6);
\draw [thkc] (-0.1,-0.2) arc (-90:90:0.2);
\draw [thkln] (-1,0) -- (-0.4,0) node [near start,below] {$\scriptstyle a$}
(0.4,0) -- (1,0) node [near end,below] {$\scriptstyle b$};
\end{tikzpicture}
} = b.
\]
In other words, a tangle
$
\begin{tikzpicture}[menvthree,rotate=90]
\draw [ptzer] (0.4,0.8) rectangle ++(-0.8,-1.6);
\draw [thkc] (-0.1,-0.2) arc (-90:90:0.2);
\draw [thkln] (-1,0) -- (-0.4,0)  
(0.4,0) -- (1,0);  
\end{tikzpicture}
$
contains no cups, but only caps and \txstrs s.

Establishing a bound $\njrc\leq\xca$ within diagrams\rx{eq:twpdgs} is easy.  In the first diagram a \tstrghtc\ contains one strand from the $\xca_1$-cable and one strand from the $\xca_2$-cable, while a \twndgc\ contains at least two strands from one of these cables, hence the total number of \trlxc s
has a bound $\njrc\leq\hlf(\xca_1+\xca_2)\leq\xca$. The second diagram contains only \twndgc s, each of them contains at least two strands of the $\xca_1$-cable, hence there $\njrc\leq\hlf\xca_1\leq\xca$.

It remains to show that the purging process ends up with the complex\rx{eq:prcmp} generated by diagrams\rx{eq:twpdgs} with non-negative shifts of \tqdgr. We prove inductively that after we purged projectors between the initial and the current one, we get a \tmcn\
\begin{equation*}
\xKhv{\xDscr} \hteqv
\boxed{
\cdots\longrightarrow\shcr^i
\bigoplus_{0\leq j\leq i}\shfr^j\left( \bigoplus_{\tau} m'_{ij,\tau} \xKhv{\tau}\right)
\longrightarrow\cdots
}_{\;i=0}^{\;\infty}
\end{equation*}
whose constituent diagrams $\tau$ have one of two possible forms between the initial and current projector (which lies to the left of the dashed line):
\[
\begin{tikzpicture}[menvtwo]
\draw [bcrct] (-1.55-0.6-0.5,0) -- (-1.55,0) (-1.25,0) -- (1.25,0) (1.55,0) -- (4.05,0) (4.35,0) -- (4.95+0.5,0);
\draw [thkln] (-1.55-0.6-0.5,0) -- (-1.55,0) (-1.25,0) -- (1.25,0) (1.55,0) -- (4.05,0) (4.35,0) -- (4.95+0.5,0);
\draw [ptzert] (-1.55,-0.6) rectangle ++(0.3,1.2);
\draw [ptzert] (1.55,-0.6) rectangle ++(-0.3,1.2);
\draw [ptzert] (4.05,-0.6) rectangle ++(0.3,1.2);
\draw [ptzer] (0.8,2) rectangle ++(-1.6,-0.8);
\draw [thkc] (0.2,0.1+1.6) arc (0:-180:0.2);
\draw [thkln]  (-0.5,2) -- (-0.5,2.6);
\draw [thkln]  (0.5,2) -- (0.5,2.6);
\node at (0.05,2.4) {$\scriptstyle \cdots$};
\draw [ptzer] (0.8,-2) rectangle ++(-1.6,0.8);
\draw [thkc] (0.2,-0.1-1.6) arc (0:180:0.2);
\draw [thkln]  (-0.5,-2) -- (-0.5,-2.6);
\draw [thkln]  (0.5,-2) -- (0.5,-2.6);
\node at (0.05,-2.4) {$\scriptstyle \cdots$};
\draw [thkln] (-1.25,0.35) to [out=0,in=-90] (-0.5,1.2);
\draw [thkln,xscale=-1] (-1.25,0.35) to [out=0,in=-90] (-0.5,1.2);
\draw [thkln,yscale=-1] (-1.25,0.35) to [out=0,in=-90] (-0.5,1.2);
\draw [thkln,xscale=-1,yscale=-1] (-1.25,0.35) to [out=0,in=-90] (-0.5,1.2);
\draw [thkln] (-1.55,0.35) to [out=180,in=-90] ++(-0.75,0.85) -- (-2.3,2.6);
\draw [thkln,xscale=-1] (-1.55,0.35) to [out=180,in=-90] ++(-0.75,0.85) -- (-2.3,2.6);
\draw [thkln,xshift=5.6cm] (-1.55,0.35) to [out=180,in=-90] ++(-0.75,0.85) -- (-2.3,2.6);
\draw [thkln,xscale=-1,yscale=-1,xshift=-2.8cm] (-1.55,0.35) to [out=180,in=-90] ++(-0.75,0.85) -- (-2.3,2.6);
\draw [decorate,decoration={brace,amplitude=4pt},xshift=0,yshift=-5pt]
(5.1cm+5pt,-2.6) -- (-0.5cm-5pt,-2.6)
node [black,midway,yshift=-0.4cm]
{\scriptsize to struts};
\draw [decorate,decoration={brace,amplitude=4pt},xshift=0,yshift=5pt]
(-2.3cm-5pt,2.6) -- (3.3cm+5pt,2.6)
node [black,midway,yshift=0.4cm]
{\scriptsize to struts};
\draw [dashed] (2.8,2.6) -- (2.8,-2.6);
\end{tikzpicture}
\;,\qquad
\begin{tikzpicture}[menvtwo]
\draw [bcrct] (-1.55-0.6-0.5,0) -- (-1.55,0) (-1.25,0) -- (1.25,0) (1.55,0) -- (4.05,0) (4.35,0) -- (4.95+0.5,0);
\draw [thkln] (-1.55-0.6-0.5,0) -- (-1.55,0) (0,1.2) -- (0,-1.2) (1.55,0) -- (4.05,0) (4.35,0) -- (4.95+0.5,0);
\draw [ptzert] (-1.55,-0.6) rectangle ++(0.3,1.2);
\draw [ptzert] (1.55,-0.6) rectangle ++(-0.3,1.2);
\draw [ptzert] (4.05,-0.6) rectangle ++(0.3,1.2);
\draw [ptzer] (0.8,2) rectangle ++(-1.6,-0.8);
\draw [thkc] (0.2,0.1+1.6) arc (0:-180:0.2);
\draw [thkln]  (-0.5,2) -- (-0.5,2.6);
\draw [thkln]  (0.5,2) -- (0.5,2.6);
\node at (0.05,2.4) {$\scriptstyle \cdots$};
\draw [ptzer] (0.8,-2) rectangle ++(-1.6,0.8);
\draw [thkc] (0.2,-0.1-1.6) arc (0:180:0.2);
\draw [thkln]  (-0.5,-2) -- (-0.5,-2.6);
\draw [thkln]  (0.5,-2) -- (0.5,-2.6);
\node at (0.05,-2.4) {$\scriptstyle \cdots$};
\draw [thkln] (-1.25,0.35) to [out=0,in=-90] (-0.5,1.2);
\draw [thkln,xscale=-1] (-1.25,0.35) to [out=0,in=-90] (-0.5,1.2);
\draw [thkln,yscale=-1] (-1.25,0.35) to [out=0,in=-90] (-0.5,1.2);
\draw [thkln,xscale=-1,yscale=-1] (-1.25,0.35) to [out=0,in=-90] (-0.5,1.2);
\draw [thkln] (-1.55,0.35) to [out=180,in=-90] ++(-0.75,0.85) -- (-2.3,2.6);
\draw [thkln,xscale=-1] (-1.55,0.35) to [out=180,in=-90] ++(-0.75,0.85) -- (-2.3,2.6);
\draw [thkln,xshift=5.6cm] (-1.55,0.35) to [out=180,in=-90] ++(-0.75,0.85) -- (-2.3,2.6);
\draw [thkln,xscale=-1,yscale=-1,xshift=-2.8cm] (-1.55,0.35) to [out=180,in=-90] ++(-0.75,0.85) -- (-2.3,2.6);
\draw [decorate,decoration={brace,amplitude=4pt},xshift=0,yshift=-5pt]
(5.1cm+5pt,-2.6) -- (-0.5cm-5pt,-2.6)
node [black,midway,yshift=-0.4cm]
{\scriptsize to struts};
\draw [decorate,decoration={brace,amplitude=4pt},xshift=0,yshift=5pt]
(-2.3cm-5pt,2.6) -- (3.3cm+5pt,2.6)
node [black,midway,yshift=0.4cm]
{\scriptsize to struts};
\draw [dashed] (2.8,2.6) -- (2.8,-2.6);
\end{tikzpicture}
\]
In both diagrams the left projector on the grey strip is the initial one, the middle projector is the current one and the right projector is the next after the current. It is not hard to see that if we purge the current projector, then we get similar diagrams with the next projector becoming the current one (the first type diagram may become of either type after the purge, while the second type diagram remains of the same type). The \tqdgr\ shifts remain non-negative, because the purging does not produce any circles: it just makes explicit various line connections that were hidden inside the constituent \tTLt s of the current projector.

Our goal is to remove the single line circles one-by-one. Consider a particular circle. We designate one of its projectors as initial and then, going clockwise, we detach the circle from other projectors step-by-step. First, if the current projector is of the third or fourth type of\rx{eq:frmprj}, then we replace it with the first or second projector with the help of the homotopy equivalences
\begin{equation}
\label{eq:smsdr}
\begin{tikzpicture}[menvtwo]
\draw [ptzer] (-0.15,-0.6) rectangle (.15,.6);
\draw [thkln] (-0.75,0) -- (-0.15,0)
node [near start,above] {$\scriptstyle \xca$}
(0.15,0) -- (.75,0)
node [near end,below] {$\scriptstyle \xca$};
\draw (-0.75,-0.4) -- (-0.15,-0.4) (0.15,0.4) -- (0.75,0.4);
\end{tikzpicture}
\;\hteqv\;
\shcr^{\hlf\xca}
\begin{tikzpicture}[menvtwo]
\draw [thkln] (-1.25,0) -- (-0.15,0)
node [very near start,above] {$\scriptstyle \xca$}
(0.15,0) -- (.75,0)
node [near end,below] {$\scriptstyle \xca$};
\draw [lnovr]  (-1.25,-0.4) to [out=0,in=180] (-0.15,0.4);
\draw (-1.25,-0.4) to [out=0,in=180] (-0.15,0.4) (0.15,0.4) -- (0.75,0.4);
\draw [ptzer] (-0.15,-0.6) rectangle (.15,.6);
\end{tikzpicture},
\qquad
\begin{tikzpicture}[menvtwo]
\draw [thkln] (-0.75,0) -- (-0.15,0)
node [near start,below] {$\scriptstyle \xca$}
(0.15,0) -- (.75,0)
node [near end,above] {$\scriptstyle \xca$};
\draw (-0.75,0.4) -- (-0.15,0.4) (0.15,-0.4) -- (0.75,-0.4);
\draw [ptzer] (-0.15,-0.6) rectangle (.15,.6);
\end{tikzpicture}
\;\hteqv\;
\shcr^{-\hlf\xca}
\begin{tikzpicture}[menvtwo]
\draw [thkln] (-1.25,0) -- (-0.15,0)
node [very near start,below] {$\scriptstyle \xca$}
(0.15,0) -- (.75,0)
node [near end,above] {$\scriptstyle \xca$};
\draw [lnovr]  (-1.25,0.4) to [out=0,in=180] (-0.15,-0.4);
\draw (-1.25,0.4) to [out=0,in=180] (-0.15,-0.4) (0.15,-0.4) -- (0.75,-0.4);
\draw [ptzer] (-0.15,-0.6) rectangle (.15,.6);
\end{tikzpicture}
\end{equation}
(note that we keep the level of the outgoing single lines and keep single lines above the $\xca$-cables).

Now suppose that our projector junction is of the first type (the second type is treated similarly). We perform a \ltrf
\begin{equation}
\label{eq:lrdpr_o}
\ytngsi =
\begin{tikzpicture}[menvtwo]
\draw [ptzer] (-0.15,-0.6) rectangle (.15,.6);
\draw [thkln] (-0.75,0) -- (-0.15,0)
node [near start,below] {$\scriptstyle \xca$}
(0.15,0) -- (.75,0)
node [near end,below] {$\scriptstyle \xca$};
\draw (-0.75,0.4) -- (-0.15,0.4) (0.15,0.4) -- (0.75,0.4);
\end{tikzpicture},
\qquad
\ytngsf =
\begin{tikzpicture}[menvtwo]
\draw [ptzer] (-0.15,-0.6) rectangle (.15,.6);
\draw [thkln] (-0.75,0) -- (-0.15,0)
node [near start,below] {$\scriptstyle \xca$}
(0.15,0) -- (.75,0)
node [near end,below] {$\scriptstyle \xca$};
\draw (-0.75,0.8) -- (0.75,0.8);
\end{tikzpicture},
\qquad
\ytngsc =
\shcr
\begin{tikzpicture}[menvtwo]
\draw[ptfour] (-.15,-0.6) rectangle ++(0.3,1.2);
\draw[thkln] (-1.35,0) -- (-.15,0)
node[near start,below] {$\scriptstyle \xca+1$ }
(.15,0) -- (1.35,0)
node [near end, below] {$\scriptstyle \xca+1$};
\end{tikzpicture}
\end{equation}
and the cone relation\rx{eq:cnrel} is \ex{eq:jwpar}.

After we performed the \ltrf\rx{eq:lrdpr} on all projectors of a given single line, except the initial one, the circle remains attached to $\xDclN$ only at the initial projector. We perform the replacement\rx{eq:smsdr} on it, if necessary. Suppose that as a result the projector is of the first type of\rx{eq:frmprj}. The single line circle passes above the $\xca$-cable in all their intersections, hence the circle can be contracted towards the remaining junction, and the junctions has the form
$
\begin{tikzpicture}[menvthree]
\draw [ptzer] (-0.15,-0.6) rectangle (.15,.6);
\draw [thkln] (-0.75,0) -- (-0.15,0)
node [near start,below] {$\scriptstyle \xca$}
(0.15,0) -- (.75,0)
node [near end,below] {$\scriptstyle \xca$};
\draw (-0.15,0.4) to [out=180,in=-90] (-0.6,0.8) to [out=90,in=180] (0,1.4)
to [out=0,in=90] (0.6,0.8) to [out=-90,in=0] (0.15,0.4);
\end{tikzpicture}
$
Now we perform a \ltrf
\[
\ytngsi=
\begin{tikzpicture}[menvtwo]
\draw [ptzer] (-0.15,-0.6) rectangle (.15,.6);
\draw [thkln] (-0.75,0) -- (-0.15,0)
node [near start,below] {$\scriptstyle \xca$}
(0.15,0) -- (.75,0)
node [near end,below] {$\scriptstyle \xca$};
\draw (-0.15,0.4) to [out=180,in=-90] (-0.6,0.8) to [out=90,in=180] (0,1.4)
to [out=0,in=90] (0.6,0.8) to [out=-90,in=0] (0.15,0.4);
\end{tikzpicture},
\qquad
\ytngsf = \shfr^{-1}
\begin{tikzpicture}[menvtwo]
\draw [ptzer] (-0.15,-0.6) rectangle (.15,.6);
\draw [thkln] (-0.75,0) -- (-0.15,0)
node [near start,below] {$\scriptstyle \xca$}
(0.15,0) -- (.75,0)
node [near end,below] {$\scriptstyle \xca$};
\end{tikzpicture},
\qquad
\ytngsc =
\shcr^{2\xca}\shfr
\begin{tikzpicture}[menvtwo]
\draw [color=white] (-0.15,-1) rectangle (0.15,1.6);
\draw [ptthr] (-0.15,-0.6) rectangle (.15,.6);
\draw [thkln] (-0.75,0) -- (-0.15,0)
node [near start,below] {$\scriptstyle \xca$}
(0.15,0) -- (.75,0)
node [near end,below] {$\scriptstyle \xca$};
\draw (-0.15,0.4) to [out=180,in=-90] (-0.6,0.8) to [out=90,in=180] (0,1.4)
to [out=0,in=90] (0.6,0.8) to [out=-90,in=0] (0.15,0.4);
\end{tikzpicture}
\]
and the cone relation\rx{eq:cnrel} is \ex{eq:htloop}.

Performing these steps on each single line circle of $\xtDNo$ we transform it into $\xDclN$, and the composition of elementary maps yields the map
\begin{equation}
\label{eq:spmapsp}
\KHm(\xDclN)\xrightarrow{\;\mntfN\;}
\end{equation}

 Composing all maps, which lead from $\xDclN$ to $\xDclNo$, we get a \tdgpr\ map
\begin{equation}
\label{eq:spmapsp_o}
\KHm(\xDclN)\xrightarrow{\;\mntfN\;} \shcr^{(\xca+\hlf)\ncrD}\, \shfr^{\gvD}\,\KHm(\xDclNo).
\tKHm(\xDclN)\xrightarrow{\;\mntfN\;} \tKHm(\xDclNo).
\end{equation}
Note that the degree shifts of \ex{eq:smsdr} cancel each other, since a single line circle has equal number of the junctions of the third and fourth type of\rx{eq:frmprj}.
After the degree shifts\rx{eq:shdgt}, the map\rx{eq:spmapsp} becomes the map\rx{eq:spmaps}.

\section{Estimates of \thdgr\ bounds on correction diagrams}

Let $\xD$ be a diagram of a link which may involve \tJWp s. Let $\yncrD$ denote the number of single line crossings in $\xD$. If $\xD$ involves cables, then we consider them as multiple single lines for the purpose of defining $\yncrD$: a crossing of an $a$-cable and a $b$-cable contributes $ab$ to $\yncrD$.

Our estimates of \thdgr\ bounds are based on the following:
\begin{proposition}
\label{prp:dgrest}
$\KHmvv{i}{\hem}(\xD)=0$ for $i\leq-\hlf \yncrD-1$, where $\yncrD$ is the number of single line crossings in $\xD$.
\end{proposition}
\begin{proof}
The claim follows easily from the defining relations\rx{eq:dkhbr} of the \tKhbr\ and from the fact that a \cJWp\ has a presentation as a complex with only non-negative homological degrees.
\end{proof}

\subsection{\Ltrf s of the first type}

\def\xpsh{3}

\begin{proof}[Proof of Proposition\rw{prp:bndfs}]

The diagram $\xtDN$ is the result of applying \Arpl s to all crossings of $\xDclNo$. A composition of maps $\xmg$ corresponding to \ltrf s\rx{eq:frstst} produces a map of \tbdgr\ zero
\[
\shcr^{\hlf\yncrv{\xtDN}}
\shfr^{\gvD(\xca+1)}
\KHm(\xtDN) \xrightarrow{\;\xmgN\;} \tKHm(\xDclNo)
\]

We perform a finite induction over the number of \Arpl s on the diagram $\xDclNo$.
Suppose that we performed \Arpl s on first $\incra$ crossings, thus obtaining the diagram $\xDsi$, and consider the \ltrf\rx{eq:frstst} on the $(\incra+1)$-st vertex. By the assumption of induction, we have a 
map
\[
\shcr^{\hlf\yncrv{\xDsi}} \shfr^{\gvD(\xca+1)}
\KHm(\xDsi) \xrightarrow{\mnfi} \tKHm(\xDclNo),
\]
where
\[
\yncrv{\xDsi} = \clN^2\ncrD + (2\xca+1)(\ncrD-\incra)
\]
is the number of single line crossings in $\xDsi$. The map $\mnfi$ has zero \tbdgr\ and it is an isomorphism on $\tKHmvb{i}(\xDclNo)$ for $i\leq 2\xca-1$.

In diagram $\xDsf$ the $(\incra+1)$-st vertex is \Arpld, while is the diagram $\xDsc$ this vertex is \Brpld. Hence the diagrams $\xDsf$ and $\xDsc$ have equal number of single line crossings:
\[
\yncrv{\xDsf} = \yncrv{\xDsc} = \yncrv{\xDsi} - (2\xca + 1).
\]
By Proposition\rw{prp:dgrest},
$
\KHmvb{i}(\xDsc) = 0
$
for $i\leq-\hlf\yncrv{\xDsc}+\xca-\hlf$, where $\xca+\hlf$ accounts for the degree shift of the correction tangle in\rx{eq:frstst}. Hence by Proposition\rw{prp:shest} the map $\KHmvb{i}(\xDsf)\xrightarrow{\xmg}\KHmvb{i}(\xDsi)$ is an isomorphism for
$i \leq-\hlf\yncrv{\xDsc}+\xca-\thlf$ and the composition map
\[
\shcr^{\hlf\yncrv{\xDsf}+\xca+\hlf}
\shfr^{\gvD(\xca+1)}\KHmvb{i}(\xDsf) \xrightarrow{\mnff}\tKHmvb{i}(\xDclNo)
\]
is an isomorphism for $i\leq 2\xca-1$.
\end{proof}

\subsection{\Ltrf s of the second type}

The resulting diagram is $\xDsi$.
 and, by the assumption of the induction we have a map

 of $\xDclNo$. The resulting

Let $\xDNoa$ denote a diagram constructed from $\xDclNo$ by performing \Brpl s on crossings $\xlbb$, $1\leq \incrb \leq \incra$ (in particular,  $\xtDNo= \xDpNof$), and let $\xDpNoa$ denote a diagram constructed by performing \Brpl s on crossings $\xlbb$, $1\leq \incrb \leq \incra$ and \Arpl\ on the crossing $\xlbao$. We define their shifted \tKhom:
\[\tKHm(\xDNoa) = \shcr^{abcd} \KHm(\xDNoa).
\]

According to Lemma\rw{lm:sss}, the complexes of the picture\rx{eq:abrpl} form an exact triangle, hence there is a canonical map
\[
\shcr^{-(\xca+\hlf)}\KHm(\xDNoa)\xrightarrow{\;\xmgNa\;}
\KHm(\xDNoamo)
\]
which has zero degree with respect to all gradings. A composition of maps $\xmgNa$ for $\incra=0,\ldots,\ncrD-1$ is a zero-degree map
\[
\shcr^{-\ncrD(\xca+\hlf)}\KHm(\xDclNo) \xrightarrow{\;\xmgN\;} \KHm(\xtDNo).
\]

%

and let $\xDpNoa$ denote a diagram constructed by performing \Brpl s on crossings $\xlbb$, $1\leq \incrb \leq \incra$ and \Arpl\ on the crossing $\xlba$.

\subsection{Second stage}
Consider the structure of the diagram $\xtDNo$. It consists of two parts. The first part is $\xca$-cabled diagram $\xDclN$. The second part consists of non-intersecting circles formed by single lines, which appeared when \Brpl s were applied to all vertices of $\xDclNo$. These circles are the same as the circles of the diagram $\sBD$ which is constructed by \Bsplng\ all crossings of the original diagram $\xD$.
We orient single line circles clockwise and on our pictures the clockwise orientation goes from the left to the right.

Both parts of $\xtDNo$ are joined by \cJWp s, the junctions having four possible forms:
\begin{equation}
\label{eq:frmprj1}
\begin{tikzpicture}[menvtwo]
\draw [ptzer] (-0.15,-0.6) rectangle (.15,.6);
\draw [thkln] (-0.75,0) -- (-0.15,0)
node [near start,above] {$\scriptstyle \xca$}
(0.15,0) -- (.75,0)
node [near end,above] {$\scriptstyle \xca$};
\draw (-0.75,-0.4) -- (-0.15,-0.4) (0.15,-0.4) -- (0.75,-0.4);
\end{tikzpicture},
\qquad
\begin{tikzpicture}[menvtwo]
\draw [ptzer] (-0.15,-0.6) rectangle (.15,.6);
\draw [thkln] (-0.75,0) -- (-0.15,0)
node [near start,below] {$\scriptstyle \xca$}
(0.15,0) -- (.75,0)
node [near end,below] {$\scriptstyle \xca$};
\draw (-0.75,0.4) -- (-0.15,0.4) (0.15,0.4) -- (0.75,0.4);
\end{tikzpicture},
\qquad
\begin{tikzpicture}[menvtwo]
\draw [ptzer] (-0.15,-0.6) rectangle (.15,.6);
\draw [thkln] (-0.75,0) -- (-0.15,0)
node [near start,above] {$\scriptstyle \xca$}
(0.15,0) -- (.75,0)
node [near end,below] {$\scriptstyle \xca$};
\draw (-0.75,-0.4) -- (-0.15,-0.4) (0.15,0.4) -- (0.75,0.4);
\end{tikzpicture},
\qquad
\begin{tikzpicture}[menvtwo]
\draw [thkln] (-0.75,0) -- (-0.15,0)
node [near start,below] {$\scriptstyle \xca$}
(0.15,0) -- (.75,0)
node [near end,above] {$\scriptstyle \xca$};
\draw (-0.75,0.4) -- (-0.15,0.4) (0.15,-0.4) -- (0.75,-0.4);
\draw [ptzer] (-0.15,-0.6) rectangle (.15,.6);
\end{tikzpicture}.
\end{equation}

Our strategy is to remove the circles one by one by though the following procedure applied to each circle: we choose the zeroth projector on the circle and then go from it clockwise, detaching the single line from the projectors. After we detach the single line circle from all projectors except the initial one, we

 First, if the projector is of the third or fourth type of\rx{eq:frmprj}, then we turn it into the first or second form by the homotopy equivalences
\begin{equation}
\label{eq:smsdr1}
\begin{tikzpicture}[menvtwo]
\draw [ptzer] (-0.15,-0.6) rectangle (.15,.6);
\draw [thkln] (-0.75,0) -- (-0.15,0)
node [near start,above] {$\scriptstyle \xca$}
(0.15,0) -- (.75,0)
node [near end,below] {$\scriptstyle \xca$};
\draw (-0.75,-0.4) -- (-0.15,-0.4) (0.15,0.4) -- (0.75,0.4);
\end{tikzpicture}
\;\hteqv\;
\shcr^{\hlf\xca}
\begin{tikzpicture}[menvtwo]
\draw [thkln] (-1.25,0) -- (-0.15,0)
node [very near start,above] {$\scriptstyle \xca$}
(0.15,0) -- (.75,0)
node [near end,below] {$\scriptstyle \xca$};
\draw [lnovr]  (-1.25,-0.4) to [out=0,in=180] (-0.15,0.4);
\draw (-1.25,-0.4) to [out=0,in=180] (-0.15,0.4) (0.15,0.4) -- (0.75,0.4);
\draw [ptzer] (-0.15,-0.6) rectangle (.15,.6);
\end{tikzpicture},
\qquad
\begin{tikzpicture}[menvtwo]
\draw [thkln] (-0.75,0) -- (-0.15,0)
node [near start,below] {$\scriptstyle \xca$}
(0.15,0) -- (.75,0)
node [near end,above] {$\scriptstyle \xca$};
\draw (-0.75,0.4) -- (-0.15,0.4) (0.15,-0.4) -- (0.75,-0.4);
\draw [ptzer] (-0.15,-0.6) rectangle (.15,.6);
\end{tikzpicture}
\;\hteqv\;
\shcr^{-\hlf\xca}
\begin{tikzpicture}[menvtwo]
\draw [thkln] (-1.25,0) -- (-0.15,0)
node [very near start,below] {$\scriptstyle \xca$}
(0.15,0) -- (.75,0)
node [near end,above] {$\scriptstyle \xca$};
\draw [lnovr]  (-1.25,0.4) to [out=0,in=180] (-0.15,-0.4);
\draw (-1.25,0.4) to [out=0,in=180] (-0.15,-0.4) (0.15,-0.4) -- (0.75,-0.4);
\draw [ptzer] (-0.15,-0.6) rectangle (.15,.6);
\end{tikzpicture}
\end{equation}
(note that we keep the level of the outgoing single lines and keep single lines above the $\xca$-cables).
Second, we replace the projectors of the first two types of\rx{eq:frmprj} by the diagrams
\[
\begin{tikzpicture}[menvtwo]
\draw [ptzer] (-0.15,-0.6) rectangle (.15,.6);
\draw [thkln] (-0.75,0) -- (-0.15,0)
node [near start,above] {$\scriptstyle \xca$}
(0.15,0) -- (.75,0)
node [near end,above] {$\scriptstyle \xca$};
\draw (-0.75,-0.8) -- (0.75,-0.8);
\end{tikzpicture},
\qquad
\begin{tikzpicture}[menvtwo]
\draw [ptzer] (-0.15,-0.6) rectangle (.15,.6);
\draw [thkln] (-0.75,0) -- (-0.15,0)
node [near start,below] {$\scriptstyle \xca$}
(0.15,0) -- (.75,0)
node [near end,below] {$\scriptstyle \xca$};
\draw (-0.75,0.8) -- (0.75,0.8);
\end{tikzpicture}
\]
After we finish, the single line circle is attached only at the zeroth projector. The final step is to `suck' it into the projector.
%
%

We label the single line circles of $\xtDNo$ as $\ylbb$, $1\leq \crcb \leq \gvD$. On each circle $\ylbb$ we choose the zeroth \tJWp\ and then label other projectors clockwise as $\zlbbg$, $0 \leq \prjg \leq\pncrb-1 $, where $\pncrb$ is the total number projectors on the circle $\ylbb$ ($\pncrb$ is also equal to the number of crossings ... ). $\xtDNob$ denotes the diagram $\xtDNo$ in which the single line circles $\ylbbp$ with $\crcbp < \crcb$ are removed. Further, $\xtDNobg$ denotes the diagram $\xtDNob$ in which the circle $\ylbb$ is detached from the projectors $\zlbbgp$ with $1\leq\prjgp<\prjg$.
The homotopy equivalence\rx{eq:jwpar} determines a canonical map
\begin{equation}
\label{eq:mpgo}
\KHm(\xtDNobgo)\xrightarrow{\;\xtmgbg\;}
\KHm(\xtDNobg).
\end{equation}

In the diagram $\xtDNobf$ the cirlce $\ylbb$ is detached only at the zeroth projector, and whenever it crosses the cable (due to replacements\rx{eq:smsdr}) it passes above it. Hence the circle can be contracted towards the zeroth projector. Let $\xhDNob$ denote the resulting diagram, the vicinity of the zeroth projector there has the form:
\[
\begin{tikzpicture}[menvtwo]
\draw [ptzer] (-0.15,-0.6) rectangle (.15,.6);
\draw [thkln] (-0.75,0) -- (-0.15,0)
node [near start,below] {$\scriptstyle \xca$}
(0.15,0) -- (.75,0)
node [near end,below] {$\scriptstyle \xca$};
\draw (-0.15,0.4) to [out=180,in=-90] (-0.6,0.8) to [out=90,in=180] (0,1.4)
to [out=0,in=90] (0.6,0.8) to [out=-90,in=0] (0.15,0.4);
\end{tikzpicture}
\]
(the single line loop may also appear below, but this distinction is not important).
Thus there is an isomorphism $\KHm(\xtDNobf)\cong\KHm(\xhDNob)$ and the composition of all maps\rx{eq:mpgo} for $1\leq\prjg\leq \pncrb-1$ produces a map
\[
\KHm(\xhDNob)\xrightarrow{}\KHm(\xtDNob).
\]
The homotopy equivalence\rx{eq:htloop} determines a canonical map
\[
\KHm(\xtDNobo) \xrightarrow{} \KHm(\xhDNob).
\]

Let $\xhDNob=$ denote the diagram in which the circle $\ylbb$ is attached


\subsection{The tail structure}
\subsubsection{Tail conjecture}

We want to study the structure of the \thead\ of the colored Jones polynomial which is the dominant part of $\pJqaL$ at $q\rightarrow 0$. In other words, we want to describe the coefficients at the highest negative powers of $q$ in $\pJqaL$. If $\xLb$ is the mirror image of $\xL$, then $\pJvv{\xca}{\xLb}(q) = \pJvv{\xca}{\xL}(\qi)$, hence all results about the behavior of $\pJqaL$ at $q\rightarrow 0$ can be applied at $q\rightarrow \infty$ by replacing $\xL$ with $\xLb$.


It has been observed experimentally that at the large values of the color $\xca$ the \thead\ of $\pJqaL$ stabilizes: if $\xca\p - \xca$ is a positive even integer, then $\pJqvv{\xca\p}{\xL}$ and $\pJqaL$ share the same first $\xca$ coefficients.
\begin{conjecture}
\label{conj.tail}
For every link $\xL$ and for $\xnu\in\{0,1\}$ there exists a  series $\zpolnqnL\in\ZZ[[q]]$  and an integer-valued function $\xfsLN$ such that
\[
\pJqaL = q^{-\xfsLN} (\zpolnqnL + O(q^{\xca})),
\]
if $\xnu = \xca \mod 2$.
\end{conjecture}


The strongest result so far regarding this conjecture has been obtained by C.~Armond.
\begin{theorem}[C.~Armond]
If $\xL$ is an alternating link, then Conjecture\rw{conj.tail}  holds true and the \thead\ does not depend on the parity of $\xca$: $\zpolnqvv{0}{\xL}=\zpolnqvv{1}{\xL}$.
\end{theorem}

\subsubsection{Multi-layered structure}
A closer look at the experimental data suggests that the tail of $\pJqaL$ can be better described by the double Laurent series of the form
\[
q^{-\xfsLN} \sum_{m,n} c_{m,n}\, q^{m + \xn\xca}.
\]
 We are going to prove the following:
%
%
%
\begin{theorem}
\label{thm.main}
For a framed link $\xL$ with the minimal crossing number $\ncrL$
there exists a non-negative integers $\gvL$, $\geL$ and $\glL$
and a two-variable
Laurent series
\[
\ypolqtL = \xt^{-\hlf\gvL}\smxnzi\xt^\xn\, \ypolqnL,
\]
where $\ypolqnL$ is a Laurent series with degree bounds
\begin{equation}
\label{eq.bnd}
\degq \ypolqnL\geq  -\shlf\geL \xn(\xn+1) - \shlf\glL\left(2 n + 1\right)^2 ,
\end{equation}
such that for any integer $\xM\geq 1$
\begin{equation}
\label{eq.mrel}
\pJqaL = \Nsgn q^{-\frth\ncrL\xca^2 - \gvL\xca }
\left(
\sum_{\xn = 0}^{\xM-1} q^{\xn\xca}\ypolqnL
+ O\left( q^{\xca(\xM-\hlf)}\right)
\right)
\end{equation}
if $\xca \geq 16\, \ncrtD^2\xM^3$.
The series $\ypolqtL$ is a topological invariant of the framed link $\xL$.
\end{theorem}

This theorem is an easy corollary of the following
\begin{theorem}
\label{thm.ma}
A link diagram $\xD$ determines a two-variable Laurent series
\[
\ypolqtD = t^{-\hlf\nBD} \smxnzi \xt^\xn\,\ypolqnD
\]
where $\ypolqnD$ are Laurent series with degree bounds
\[
\degq \ypolqnL\geq -\shlf\,\geD \xn(\xn+1) - \shlf\,\glD\left( 2n + 1\right)^2,
\]
such that for any integer $\xM\geq 1$
\[
\pJqaD =\Nsgn q^{-\frth\ncrD\xca^2 - \hlf\gvD \xca }
\left(
\sum_{\xn = 0}^{\xM-1} q^{\xn\xca}\ypolqnL
+ O\left( q^{\xca(\xM-\hlf)}\right)
\right)
\]
if $\xca \geq 16\, \ncrtD^2\xM^3$.
\end{theorem}

\begin{proof}[Proof of Theorem]
By choosing the diagram $\xD$ of Theorem\rw{thm.ma} to be minimal and setting
$\ypolqtL = \ypolqtD$, $\gvL=\gvD$, $\geL=\geD$ and $\glL=\glD$ we are led to the bounds and relations.

It is easy to see from the relation\rx{eq.mrel} that the colored Jones polynomial $\pJqaL$ determines the series $\ypolqtL$, hence that series is a topological invariant of $\xL$.
\end{proof}

\subsubsection{Adequate links}

The main drawback of Theorem\rw{thm.main} is that for many links the series $\ypolqtL$
is identically zero: $\ypolqtL\equiv 0$. Examples of this phenomenon are identified with the help of the following simple corollary of Theorem\rw{thm.ma}:
%
\begin{corollary}
\label{cor.equiv}
If for a framed link $\xL$ there exist numbers $\bdA$ and $\xca_0$ such that for $\xca\geq \xca_0$ the colored Jones polynomial $\pJqaL$ has a degree bound
\begin{equation}
\label{eq.dcond1}
\degq \pJqaL \geq  -\sfrth(\ncrD-1)\xca^2 + \bdA \xca,
\end{equation}
then $\ypolqnL\equiv 0$.
\end{corollary}
\begin{proof}
It is easy to see from \ex{eq.mrel} that the Jones polynomials $\pJqaL$  at $\xca\geq \xca_0$ determine the series $\ypolqnL$ and that the condition\rx{eq.dcond1} implies $\ypolqtL\equiv 0$.
\end{proof}

As an example of this Corollary, consider torus knots $\yTmn$, $1< m < n$, their framing being determined by the standard diagram of the twisted $m$-cable of the unknot. Then
\[
\degq \pJqvv{\xca}{\yTmn} \geq -\sfrth n \xca^2 + \shlf n\xca.
\]
Since $\ncrv{\yTmn} = (m-1)n$, then according to Corollary\rw{cor.equiv}, $\ypolqtv{\yTmn}\equiv 0$ if $m\geq 3$. In fact, we conjecture, that this happens to all \tnBadq\ links:
\begin{conjecture}
If a framed link $\xL$ is \tnBadq, then $\ypolqtL\equiv 0$.
\end{conjecture}
The relation $\ypolqtL=0$ means that for the link $\xL$, the series $\ypolqtL$ computed according to the prescription of , misses the actual \thead\ of $\pJqaL$ and this \thead\ has to be computed differently.

For \tBadq\ links the method of captures the actual \thead. Indeed, a theorem by O.~Dasbach and X.-S.~Lin\cx{asdf} implies that if $\xL$ is \tBadq, then $\degq\pJqaL=-\frth\ncrL\xca^2 - \gvL\xca$ while the coefficient at the lowest power of $q$ is $1$. Hence Theorem\rw{thm.main} has the following corollary:
\begin{corollary}
If a framed link $\xL$ is \tBadq, then $\ypolqtL\not\equiv 0$ and $\ypolqvv{0}{\xL} = 1 + O(q)$.
\end{corollary}

\subsubsection{The \thead s of torus links $\yTmmn$}
The bound\rx{eq.bnd} for a \tBadq\ link is stronger:
\[
\degq \ypolqnL\geq  -\shlf\geL \xn(\xn+1).
\]
The computation of the \thead s of torus links indicates  that this bound is sharp. Indeed, according to~\cite{Mort95}, the colored Jones polynomial of $\yTmmn$ is

\subsubsection{The \thead\ and Habiro series}

Still it does not permit us to perform a direct substitution of $\xt = q^\xca$ in the whole double series $\ypolqtL$.

\end{document}

Let us recall the definitions of a pair of categories and a pair of 2-categories related to the Jones polynomial and its categorification. The first category in each pair is purely topological while the second category is associated to $\xSo$ by the 3-dimensional Chern-Simons-Witten \TQFT\ and by the 4-dimensional \TQFT\ corresponding to Khovanov's categorification.

An object of the topological tangle category $\cTng$ is an oriented 2-disc $\xDon$ with $n$ marked points in its interior
placed on a special oriented `equatorial' diameter in the disc. The marked points are framed, that is, there is a choice of a tangent vector at each point. We assume that the framing vectors are tangent to the \eqdiam. The set of morphisms $\Hom_{\cTng}(\xDom,\xDon)$ consists of \ttnglmn s, that is, we glue the discs $\xDom$ and $\xDon$ together along their circle boundaries so that the orientations and end-points of \eqdiam s match, and consider a 3-ball $\xBmpn$, whose boundary is the resulting 2-sphere. Then an \ttnglmn\ is an embedding of framed (unoriented) segments and cirlces into $\xBmpn$ such that the end-points of segments map to the $m+n$ marked points on its boundary. The morphism composition rule is the obvious composition of tangles.

The category $\cTng$ can be promoted to the 2-category $\ctTng$ if  the morphism sets $\Hom(\xDom,\xDon)$ are defined as categories, the morphisms between two tangles being cobordisms.

Let $q$ be a commutative variable and let $\QQqqi$ be the algebra of Laurent polynomials of $q$.
The \tTL\ (\taTL) category $\cTL$ is an additive category over the algebra $\QQqqi$ of Laurent polynomials of $q$. Its objects are the same as those of $\cTng$, but the set of morphisms
$\aTLmn = \Hom_{\cTL}(\xDom,\xDon)$ is a module over $\QQqqi$ generated freely by \taTLt s, which are  tangles consisting of segments embedded (crossinglessly) into the equatorial disc of $\xBmpn$. The composition of \taTLt s may produce disjoint circles, but each disjoint circle is replaced by the factor $\mqpqi$. The sum of all modules $\aTL = \bigoplus_{m,n\geq 0}\aTLmn$
is called the \tTLa.

In order to work with \tJWp s we have to introduce the algebra of formal Laurent series $\QQqqip$. The corresponding \tTLc\ is denoted as $\cTLp$.

The \tKbr\
$
\xKbrd\colon \cTng\rightarrow\cTL
$
is the functor between categories defined by the relation
\[
\xKbrBv{\xcrsp}
\;\;=\;\;
\qvh\;\;
\xKbrBv{\xpver}
\;\;+\;\;
\qvmh\;\;
\xKbrBv{\xphor}
\]
and by the rule that disjoint circles are converted into the factors $\mqpqi$. In particular, the \tKbr\ maps a framed link into its \tJpol.

Khovanov categorified the Jones polynomial as well as the \tKbr\ for tangles. We find it convenient to use Bar-Natan's canopoly version of this categorification with slight modifications. Thus we consider the 2-category $\ctTL$, whose objects are again the discs $\xDon$. For two objects $\xDom$ and $\xDon$ consider the $\ZZZtt$-graded additive category $\tHom_{\ctTL}(\xDom,\xDon)$ generated by objects
$\xKhl$ corresponding to \TLttnglmn s $\xlam$ and by their translations $\xKhl\hgrshklm$ with respect to the $\ZZZtt$-grading. Note that $k$ and $l$ may take half-integer values, but integer and half-integer degrees `do not mix', so we can safely pretend that the grading is $\ZZZtt$.

\subsection{The center of the 2-category $\ctTL$}

Our main goal is the study of the center of the 2-category $\ctTL$. Recall that an element $\elcf$ of the center $\Zcatv{\caC}$ of a category $\caC$ is a choice of an endomorphism $\elcfv{A}$ for every object $A$ of $\caC$ such that for any morphism $\gAB\in\Hom_{\caC}(A,B)$ there is an equality $\elcfv{B} \gAB = \gAB \elcfv{A}$.
 The elements of the center can be composed, so $\Zcatv{\caC}$ is a monoid and  the center $\Zcatv{\ctaC}$ of a 2-category $\ctaC$ is a monoidal category.

The centers of the \TQFT\ categories $\cTL$ and $\ctTL$ are especially interesting, since they are related to the module and, respectively, category associated within these \TQFT s to the 2-torus $\xTt$. To illustrate this point, consider first the centers of the topological categories $\cTng$ and $\ctTng$.

Let $\xIdn$ denote the \TLttnglnn\ corresponding to the $n$-strand identity braid. For a link $\xL$ in $\xStSo$, let $\xIdnL$ denote the \ttnglnn\ constructed by wrapping the link $\xL$ as a meridian (that is, as a band) around $\xIdn$. It is easy to see that a collection of endomorphisms $\xcztL$ defined as
$\xczt_{\xDon} = \xIdnL$ is an element in the center of both $\cTL$ and $\ctTL$, because the $\xL$-band can be slid along any tangle around which it is wrapped. From the \TQFT\ perspective, since the boundary of $\xStSo$ is $\xTt$, the solid torus $\xStSo$ containing a link determines an element (an object) in the module (category) associated with $\xTt$.

Notably, $\xczt$ are not the only elements (objects) in the center of topological categories. Let
$\gbrmn$ be the \ttnglnn\ corresponding to a \trbr\ with $m$ full \clckw\ rotations of $n$ strands having zero framing. Then a collection of endomorphisms

\section{Multi-cones}

Let $\caCt$ be an additive category generated freely by a finite set of objects, that is, the objects of $\caCt$ are finite sums of generators  (we have categories $\tHom_{\ctTL}(\xDom,\xDon)$ in mind). Let $\caC = \cKomm{\caCt} $ be the homotopy category of complexes bounded from above: an object of $\caC$ is a complex
$\cmA = (\cdots \rightarrow A_i \rightarrow A_{i+1} \rightarrow\cdots\rightarrow A_k )$ and morphisms are chain maps up to homotopy.


A \JWp\ $\jwpn$ is a special idempotent element of the $n$-strand
\TLa\ $\cTLn$, whose defining property is
the annihilation of cap and cup tangles.
The coefficients in its expression in terms of \TLb\ tangles are
rational (rather than polynomial) functions of $q$. This suggests
that the categorification $\ctjwn$ of $\jwpn$ in the
universal tangle category $\dTLn$ constructed by D.~Bar-Natan\cx{BN1}
should be presented by a semi-infinite \chcpl. In fact, there are
two mutually dual categorifications: the complex $\ctjwpn$ which is bound from
above and the complex $\ctjwmn$ which is bound from below. We will
consider only $\ctjwpn$ in detail, since the story of $\ctjwmn$ is totally
similar.

%

The construction of $\ctjwmn$ by
B.~Cooper and S.~Krushkal\cx{CK} is based upon the
Frenkel-Khovanov formula for $\jwpn$ and requires the invention of morphisms
between constituent \TL\ tangles as well as non-trivial `thickening'
of the complex. An alternative `representation-theoretic'
approach to the categorification of the \JWp\ is developed by Igor Frenkel,
Catharina Stroppel, and Joshua Sussan\cx{FSS}.

Our approach is rather straightforward: the
categorified projector $\ctjwpn$ is a direct limit of
appropriately shifted
categorification complexes of \cbr s
(\ie braid analogs of torus links) with high \clckw\ twist (the
other projector $\ctjwmn$ comes from high \cclckw\ twists).
The limit  $\ctjwpn$
can be presented as a cone:
\xlee{eq:int1}
\ctjwpn\hteqv
\CnBv{\Ohp\big(2m(n-1)\big)\longrightarrow\cbrmns},
\xeee
where
$\gbrmn$ is a \cbr\ with $m$ full \clckw\ rotations of $n$ strands,
$\symcats{-}$ is the
categorification complex with a special grading shift, and
$\Ohp(k)$ denotes a \chcpl\ which ends at the homological degree
$-k$. Theorem\rw{th:cnpr} imposes even stronger restrictions on
the complex $\Ohp\big(2m(n-1)\big)$ in \ex{eq:int1}.


The advantage of our approach is that one can use \cbr s with high
twist as approximations to $\ctjwpn$ in a computation of \Kh\ of a
spin network which involves \JWp s:
if a spin network $\xnu$ is constructed by connecting $\jwpn$ to an \ttngnn\
$\xtau$ such that $\symcat{\xtau}\hteqv\Ohp(k)$, while a spin network $\xnum$ is constructed
by replacing $\jwpn$ in $\xnu$ with $\gbrmn$, then the homology of
$\symcat{\xnu}$ coincides with the shifted homology of
$\symcat{\xnum}$ in all homological degrees $i$
such that $i> -k - 2m(n-1)$. Thus one may say that
there is a stable limit
\xlee{eq:stlimsn}
\symcat{\xnu} =
\lim_{m\rightarrow+\infty}\symcats{\xnum}.
\xeee
We will define homological limits more precisely in
subsection\rw{sss.homcal}.

The practical
importance of the relation between $\symcat{\xnu}$
and $\symcat{\xnum}$ stems from the fact that $\xnum$ is an
ordinary link and its homology
can be computed with the help of
existing efficient computer programs even for high values
of $m$.
%
%

The simplest example of a spin network  is the unknot
`colored' by the $(n+1)$-dimensional representation of
$\mathrm{SU}(2)$ with the help of the projector $\jwpn$. Its
\Kh\ is approximated by the homology
of torus links $\mathrm{T}_{n,-mn}$ which appear as cyclic closures of
$\gbrmn$. The \Kh\ of torus links has been studied by Marko Stosic\cx{St}, who
observed that it stabilizes at lower degrees as $m$ grows. This is
a particular case of the `stable limit'\rx{eq:stlimsn}.

In Section\rw{s:notres}
we explain all notations and conventions
which are used in the paper. In particular, in
subsection\rw{sss:trgr} we define a non-traditional grading of
\Kh, which is convenient for our computations.
Then we formulate our results.


In Section\rw{s:elhomcal} we review basic facts about homological
`calculus' required to work with limits of sequences of
complexes in a homotopy category. In Section\rw{s:cbr} we construct
a sequence of categorification complexes of
\cbr s related by special \chmp s. This sequence yields $\ctjwpn$
as its direct limit. In Section\rw{s:prfs} we use homological
calculus of Section\rw{s:elhomcal} in order to prove that
$\ctjwpn$ is a categorification of the \JWp.

\def\bbS{ \mathbb{S} }
\def\So{ \bbS^1 }
\def\St{ \bbS^2 }
\def\Sot{ \So\times\St }

\subsection*{Acknowledgements}

This paper is a spinoff of a joint project with Mikhail
Khovanov\cx{KhRS} which is
dedicated to the study of categorification complexes of \cbr s and their
relation to the categorification of the Witten-Reshetikhin-Turaev
invariant of links in $\Sot$. I am deeply indebted to Mikhail for
numerous discussions and suggestions.

I would like to thank Slava Krushkal for sharing the results of
his ongoing research. I am also indebted to organizers of the M.S.R.I.
workshop `Homology Theories of Knots and Links' which stimulated
me to write this paper.

This work is supported by the NSF grant DMS-0808974.

\end{document}

\section{Notations and results}
\label{s:notres}
\subsection{Notations}
\label{ss:not}
\subsubsection{Tangles and \TLa}

All tangles in this paper are framed and we assume the blackboard
framing in pictures. We use the symbol 
$\xygraph{
!{0;/r1.5pc/:}
[u(0.5)]
!{\xcapv@(0)}
[u(0.45)r(0.23)]
*{\symfr\;\scriptstyle{k}}
[u(1.5)]
}
$to indicate an
addition of $k$ framing twists to a tangle strand:
\xlee{ae1.1b}
\xygraph{
!{0;/r1.5pc/:}
[u(0.5)]
!{\hover}
!{\hcap}
[u(0.5)l(0.25)]
}
\;\; = \;\;
\xygraph{
!{0;/r1.5pc/:}
[u(0.5)]
!{\xcapv@(0)}
[u(0.45)r(0.23)]
*{\symfr\;\scriptstyle{1}}
[u(1.5)]
}
\xeee

A tangle is called \emph{\plnr} if it can be presented by a diagram
without crossings. A \plnr\ tangle is called \emph{connected} or
\emph{\TLb} (\TLba) if
it does not contain disjoint circles. Let $\rTNG$ denote the set of all
framed tangles,
$\rTNGmn$ -- the set of \ttngmn s and $\rTNGn$ -- the set of
\ttngnn s.
We adopt similar notations for the set  $\rTL$ of \TLba-tangles.

We use the symbol $\tcmp$ to denote the composition of tangles:
$\xtauo\tcmp\xtaut$. The same symbol is used to denote the
multiplication in \TLa\ and the composition bifunctor in the
category $\dTL$.

A \TLa\ $\cTL$ over the ring of Laurent polynomials $\Zqqi$\footnote{It is clear from our normalization of the
Kauffman  bracket relation\rx{ae1.2} that we should rather use the
ring $\Zqqhi$. However, in all expressions in this paper the
half-integer power of $q$ appears only as a common factor, so the terms with integer
and half-integer powers of $q$ do not mix. Hence
we refer to $\Zqqi$, while keeping in mind that $\qh$ may
appear as a common factor is some expressions.}
is a quiver ring. The vertices $v_n$ of the quiver are indexed by
non-negative integers $n$ and each pair of vertices $v_m$, $v_n$,
such that $m-n$ is even, is connected
by an edge $e_{mn}$. To a vertex $v_n$ we associate a ring $\cTLnn$ (also denoted as
$\cTLn$)
and to an edge $e_{mn}$ we associate a
$\cTLn\otimes\cTLm^{\mathrm{op}}$-module $\cTLmn$. As a module,
$\cTLmn$ is generated freely by elements $\clam$ corresponding to \TL\ \ttngmn s $\xlam$, while ring
and module structures come from the composition of tangles modulo
the relation
\xlee{ae1.1}
\Bsymalg{\lcir}
 = -(\qpqi),
\xeee
which is needed to remove disjoint circles that may appear in the composition
of \TLb\ tangles.

The map $\rTNG\xrightarrow{\symalg{-}}\cTL$
associates an element $\ctau$ to a tangle $\xtau$ with the help of
\ex{ae1.1} and the Kauffman bracket relation
\xlee{ae1.2}
\Bsymalg{\xcrsp}
\;\;=\;\;
\qvh\;\;
\Bsymalg{\xpver}
\;\;+\;\;
\qvmh\;\;
\Bsymalg{\xphor}.
\xeee
This relation removes crossings and disjoint circles from the
diagram of $\xtau$, hence
\xlee{ae1.2a0}
\ctau = \sltln \xcalt\, \clam,\qquad
\xcalt = \sum_{i\in\ZZ}\xcalit\,q^i
\xeee
with only finitely many coefficients $\xcalit$ being non-zero.
%

If two tangles differ only by the framing of their strands, then
the corresponding algebra elements differ by the $q$
power factor coming from the following relation associated with
the first Reidemeister move:
\xlee{ae1.2a}
\Bsymalg{\xvfro\hspace*{-0.2cm}}
\;\; = \;\;
-q^{\frac{3}{2}}\;\;
\Bsymalg{\;\xvert\hspace*{-0.5cm}}
\xeee

A \ttngzz\ $\xL$ is a framed link, so $\symalg{\xL}$
is
the framing dependent Jones polynomial defined by the
Kauffman bracket.



We use the notations $\QcTL$ and $\cTLpinf$ for \TLa s defined over
the field $\Qq$ of rational functions of $q$ and over the field
$\Zsqqi$ of Laurent power series.
A sequence of injective homomorphisms
$\Zqqi\hookrightarrow\Qq\hookrightarrow\Zsqqi$, the latter one
generated by the expansion in powers of $q$,
produce a sequence of injective homomorphisms of the corresponding
\TLa s.

%

%

\subsubsection{The \JWp}

Let $\gcupni\in\rTLvv{n-2}{n}$ and $\gcapni\in\rTLvv{n}{n-2}$,
$1\leq i\leq n-1$, denote the following \TL\ tangles:
\ylee{ae1.3}
\gcupni=\xygraph{
!{0;/r1.5pc/:}
[r(0.25)u(0.5)]
!{\xcapv@(0)}
[u(0.5)r(1)]
*{\cdots}
[r(01)u(0.5)]
!{\xcapv@(0)}
[r(0.5)u(1)]
!{\vcap-}
[r(1.5)]
!{\xcapv@(0)}
[u(0.5)r(1)]
*{\cdots}
[r(01)u(0.5)]
!{\xcapv@(0)}
[u(1.5)l(3.5)]
*{\scriptstyle{i}}
[r(1)]
*{\scriptstyle{i+1}}
[l(3.5)]
*{\scriptstyle{1}}
[r(6)]
*{\scriptstyle{n}}
}
,
\quad\quad
\gcapni=
\xygraph{
!{0;/r1.5pc/:}
[r(0.25)u(0.5)]
!{\xcapv@(0)}
[u(0.5)r(1)]
*{\cdots}
[r(01)u(0.5)]
!{\xcapv@(0)}
[r(0.5)]
!{\vcap}
[r(1.5)u(1)]
!{\xcapv@(0)}
[u(0.5)r(1)]
*{\cdots}
[r(01)u(0.5)]
!{\xcapv@(0)}
[d(0.5)l(3.5)]
*{\scriptstyle{i}}
[r(1)]
*{\scriptstyle{i+1}}
[l(3.5)]
*{\scriptstyle{1}}
[r(6)]
*{\scriptstyle{n}}
}
\yeee
Their compositions $\xUni = \gcupni\tcmp \gcapni$ are standard
generators of the \TLa\ $\cTLn$.

The \JWp\ $\jwpn\in\QcTLn$ is the unique non-trivial idempotent element satisfying the
condition
\xlee{ae1.4}
\acapni\;\tcmp\jwpn =0,\qquad 1\leq i\leq n-1.
\xeee
The \JWp\ also satisfies the relation
\xlee{ae1.4a}
\jwpn\tcmp\;\acupni =0,\qquad 1\leq i\leq n-1.
\xeee

We denote the idempotent element of
$\cTLpinf_n$  corresponding to $\jwpn$ as $\jwpnp$.

\subsubsection{Basic notions of homological algebra}
Let $\xChA$ be a category of \chcpls\ associated with an additive
category $\xctA$. An object of $\xChA$ is a \chcpl\
\ylee{ae1.ch1}
\xbA  = (\cdots \rightarrow
\xAi\xrightarrow{\xdi}\xAio\rightarrow\cdots),
\yeee
and a morphism
between two chain complexes is a \chmp\ defined as a \mmp
\xlee{ae1.10d}
\vcenter{\xymatrix{
\xbA \ar[d]^-{\xbf} &&
\cdots\ar[r]^-{\xdimo} & \xAi \ar[r]^-{\xdi} \ar[d]^-{\yfi} &
\xAio
\ar[r]^-{\xdio} \ar[d]^{\yfio} & \cdots
\\
\xbB &&
\cdots\ar[r]^-{\xdpimo} & \xBi \ar[r]^{\xdpi} & \xBio
\ar[r]^-{\xdpio} & \cdots
}
}
\xeee
which commutes with the chain differential: $\xdpi\,\yfi = \yfio\,\xdi$ for all $i$.
The cone of a \chmp\ $\xbA\xrightarrow{\xbf}\xbB$ is a complex
\ylee{ae1.10b1}
\Cnbf
=
\lrbc{
\vcenter{
\xymatrix@C=1.5cm@R=0.5cm{
\cdots \ar[dr] \ar[r] & \xAi
\ar@{}[d] |{\oplus} \ar[r]^-{-\xdi} \ar[dr]^{-\xfi} &
\xAio
\ar@{}[d] |{\oplus}
\ar[r] \ar[dr]& \cdots
\\
\cdots \ar[r] & \xBimo\ar[r]_{\xdpimo} & \xBi \ar[r] & \cdots
}
}
}
\yeee
in which the object $\xAio\oplus\xBi$ has the homological degree
$i$.
There are two special \chmp s
$\xbB\xrightarrow{\idlbf}\Cnbf$ and
$\Cnbf[1]\xrightarrow{\chdlbf}\xbA$ associated to
the cone:
\ylee{ae1.10b2}
\xymatrix{
\xbB \ar[d]^-{\idlbf}&&
\cdots \ar[r] &
\xBi \ar[r] \ar[d]^-{0\oplus \xId}
&
\xBio \ar[r] \ar[d]^-{0\oplus \xId}
&
\cdots
\\
\Cnbf \ar[d]^-{\chdlbf} &&
\cdots \ar[r] &
\xAio \oplus \xBi \ar[r] \ar[d]^-{\xId\oplus 0} &
\xAit\oplus \xBio \ar[r] \ar[d]^-{\xId\oplus 0} &
\cdots
\\
\xbA[-1] &&
\cdots \ar[r]
&
\xAio \ar[r]
&
\xAit \ar[r]
&
\cdots
}
\yeee
These complexes and \chmp s form a \dstt:
\xlee{ae1.ch2}
\xymatrix{
\xbA\ar[r]^-{\xbf} &
\xbB \ar[r]^-{\idlbf} &
\Cnbf \ar[r]^-{\chdlbf} &
\xbA[-1]
}.
\xeee

The homotopy category of complexes $\xKhA$ has the same objects as
$\xChA$ and the morphisms are the morphisms of $\xChA$ modulo
homotopies.
%
%
We denote
homotopy equivalence by the sign $\hteqv$.
The notion of a cone extends to $\xKhA$ and there
are additional relations in that category: $\Cnv{\idlbf} \hteqv \xbA[-1]$ and
$\Cnv{\chdlbf} \hteqv \xbB[-1]$, so all vertices of a \dstt\ have
equal properties.

\subsubsection{A triply graded categorification of the Jones
polynomial}
\label{sss:trgr}
In his famous paper\cx{Kh1}, M.~Khovanov
introduced a categorification of
the Jones polynomial of links. To a diagram $\xL$ of a
link he associates a complex of graded modules
\xlee{ae1.5}
\dL = \lrbc{ \cdots \rightarrow \dLi \rightarrow \dLio\rightarrow\cdots}
\xeee
so that
if two diagrams represent the same link then the corresponding
complexes are homotopy equivalent, and the graded Euler
characteristic of $\dL$ is equal to the Jones polynomial of $\xL$.

Thus, overall, the complex\rx{ae1.5} has two gradings: the first one
is
the grading related to powers of $q$ and the second one is the
homological grading of the complex itself, the corresponding
degree being equal to $i$.
In this paper we adopt a slightly different convention which is
convenient for working with framed links and tangles. It is
inspired by matrix factorization categorification\cx{KR1} and its
advantage is that it is no longer necessary to assign orientation to
link strands in order to obtain the grading of the categorification
complex\rx{ae1.5} which would make it invariant under the second
Reidemeister move.

To a framed link
diagram $\xL$ we associate a $\ZZ \oplus\ZZ\oplus\ZZ_2$-graded complex\rx{ae1.5} with
degrees $\dgo$, $\dgt$ and $\dgh$.
The first two gradings are of the same nature as in\cx{Kh1} and, in
particular, $\dgo\dLi=i$. The third grading is an inner grading of
chain modules defined modulo 2 and of homological
nature, that is, the homological parity of an element of $\dL$,
which affects various sign factors, is the sum of $\dgo$ and
$\dgh$. Both homological degrees are either integer or
half-integer simultaneously, so the homological parity is integer
and takes values in $\ZZ_2$. The $q$-degree $\dgt$ may also take
half-integer values.


Let $\tgrshv{l}{m}{n}$ denote the shift of three degrees by $l$,
$m$ and $n$ units respectively\footnote{
Our degree shift is defined in such a way that if an object $M$
has a homogeneous degree $n$, then the shifted object $M[1]$ has a
homogeneous degree $n+1$.
}. We use abbreviated notations
$$
\tgrsshv{m}{l} = \tgrshv{m}{l}{0},\qquad
\qshv{m} = \tgrshv{m}{0}{0}
$$
as well as the following `power' notation:
$$
\tgrshv{m}{l}{n}^k = \tgrshv{km}{kl}{kn}.
$$

With new grading conventions, the categorification
formulas of\cx{Kh1} take the following form:
the module associated with an unknot is still $\ZZ[x]/(x^2)$ but with
a different degree assignment:
\begin{eqnarray}
\label{ae1.6}
&
\Bsymcat{\lcir}=\ZZ[x]/(x^2)\,
\tgrshv{-1}{0}{1},
\\
&\dgt 1 = 0, \quad \dgt x = 2,
\quad\dgo 1 = \dgo x = \dgh 1=\dgh x =0,
\end{eqnarray}
and the categorification complex of a crossing is the same as
in\cx{Kh1} but with a different degree shift:
\xlee{ae1.7}
\Bsymcat{\xcrsp}
\;\;=\;\;
\Bigg(\;\;
\Bsymcat{\xpver}
\;\tgrshv{\vthf}{-\vthf}{\vthf}
\xrightarrow{\;\;\;\;\xmrf\;\;\;\;}
\Bsymcat{\xphor}
\;\tgrshv{-\vthf}{\vthf}{-\vthf}
\vspace*{18pt}
\;\;
\Bigg),
\xeee
where $f$ is either a multiplication or a comultiplication of the
ring $\ZZ[x]/(x^2)$ depending on how the arcs in the \rhs are
closed into circles.
The resulting categorification complex\rx{ae1.5} is invariant
up to homotopy under the second and third Reidemeister moves, but
it acquires a degree shift under the first Reidemeister move:
\xlee{ae1.8}
\Bsymcat{\xvfro\hspace*{-0.2cm}}
\;\; = \;\;
\Bsymcat{\;\xvert\hspace*{-0.5cm}}
\tgrshv{\vthh}{-\vthf}{-\vthf}.
\xeee
It is easy to see that the whole categorification complex\rx{ae1.5} has a
homogeneous degree $\dgh$.

\subsubsection{A universal categorification of the \TLa}
D.~Bar-Natan\cx{BN1} described the universal category $\dTL$, whose
Grothendieck \Kzg\ is $\cTL$ considered as a $\Zqqi$-module.
We will use this category with obvious adjustments required by the new
grading conventions.

Let $\dTLt$ be an additive category whose objects are in
one-to-one correspondence with \TLb\ tangles, morphisms being generated
by tangle cobordisms (see\cx{BN1} for details). The universal category
$\dTL$ is the homotopy category of bounded complexes associated with
$\dTLt$. In other words, an object of $\dTL$ is a complex
\xlee{ae1.8a}
\xbC =
\lrbc{\cdots\rightarrow\xCi\rightarrow\xCipo\rightarrow\cdots},\qquad
\xCi =
\bigoplus_{j,\xmu}
\oltln \cjilam\,\dlam \tgrshv{j}{0}{\mu},
\xeee
where
non-negative integers $\cjilam$ are multiplicities; since the
complex is bounded, they are non-zero for only finitely many
values of $i$.

%

A categorification map $\rTNG\xrightarrow{\mpcat}\dTL$ turns a framed
tangle diagram $\xtau$ into a complex $\dtau$ according to the
rules\rx{ae1.6} and\rx{ae1.7}, the morphism $\xmrf$ in the
complex\rx{ae1.7} being the saddle cobordism. A composition of
tangles becomes a composition bi-functor
$\dTL\times\dTL\rightarrow\dTL$ if we apply
the categorified version of the rule\rx{ae1.1} in order to remove
disjoint circles:
\xlee{ae1.01}
\Bsymcat{\lcir}= \cnot \tgrshv{1}{0}{1} + \cnot\tgrshv{-1}{0}{1},
\xeee
where $\xnot$ is the empty \TL\ \ttngzz.

A complex $\dtau$ associated to a tangle $\xtau$ is defined only up to
homotopy. We use a notation $\spcc{\dtau}$ for a particular complex
with special properties which represents $\dtau$.
%

Overall,
we have the following commutative diagram:
%
\begin{equation}
\xymatrix@C=1.5cm@R=0.3cm{
& {}\dTL \ar[dd]^{\Kz}
\\
{}\rTNG \ar[ur]^{\mpcat} \ar[dr]^{\mpalg}
\\
& {}\cTL
}
\end{equation}
where the map $\Kz$ turns the complex\rx{ae1.8a} into the
sum\rx{ae1.2a0}:
\xlee{ae1.9a}
\Kz(\xbC)
=
\sltln\sum_{j}
\xcalj \,q^j\,\clam,
\qquad
\xcalj = \sum_{i,\xmu}
(-1)^{i+\xmu}\, \cjilam.
\xeee
Since the complex is bounded, the sum in the expression for
$\xcalj$ is finite.


In addition to $\dTL$ we consider
a category $\dTLp$ of complexes
bounded from above, that is, the multiplicity coefficients in the
sum\rx{ae1.8a} are zero if $i$ is greater than certain value.
Define the $q^+$ order of a \qcmd\ $\xCi$:
$\yordq{\xCi} = \xinfv{j\colon \exists\mu\colon\cjilam\neq
0}$.
A complex $\xbC$ in $\dTLp$ is \emph{\qpb} if $\lmii\yordq{\xCmi} =
+\infty$.
%
For a \qpb\ complex,
the sum in the expression\rx{ae1.9a} for $\xcalj$ is
finite, hence the element $\Kz(\xbC)$ is well defined.


\subsection{Results}
\label{ss:res}

\subsubsection{Infinite \cbr\ as a \JWp\ in a \TLa}
A braid with $n$ strands is a particular example of a \ttngnn.
A \emph{\cbr} is a braid
that can be drawn on a cylinder $\So\times[0,1]$
without intersections. In fact, all \cbr s have the form
$\btcyln^m$, $m\in\ZZ$, where $\btcyln$ is the elementary
clockwise winding \cbr:
\xlee{ae1.10p}
\btcyln \;\;=\;\;
\xygraph{
!{0;/r1.5pc/:}
!{\vover}
[u(1.5)l(0.5)]
!{\xbendr[0.5]}
[u(1.25)l(1.25)]
!{\xbendd[-0.5]}
[u(1.25)l(1)]
!{\xcapv[0.25]@(0)}
[r(1.75)]
!{\xcaph@(0)}
[u(1)]
!{\vover[-1]}
[r(1)]
!{\xbendr[-0.5]}
[u(0.75)l(0.75)]
!{\xbendd[0.5]}
[u(0.5)l(0.5)]
!{\xcapv[0.25]@(0)}
[u(1)l(1.5)]
*{\cdots}
[u(1.5)l(1.5)]
*{\cdots}
[u(0.75)l(1.5)]
*{\scriptstyle{1}}
[r(2.5)]
*{\scriptstyle{n-1}}
[r(1.25)]
*{\scriptstyle{n}}
[d(3)l(3)]
*{\scriptstyle{1}}
[r(1)]
*{\scriptstyle{2}}
[r(2.75)]
*{\scriptstyle{n}}
[u(2.2)l(0.35)]
}
\xeee
%
%
%
%
%
%
We introduce a special notation for the \cbr\ which corresponds to
$m$ full rotations of $n$ strands:
\ylee{eq:brdf}
\gbrmn = \btcyl^{mn}.
\yeee

Let $\Opqm$ denote any element of $\cTLpinf$ of the form
$\sltln\sum_{j\geq m} \xcalj\,q^j\,\clam$.
We define a \emph{\qord} of an element $\yal\in\cTLpinf$ as
$\yordq{\yal} = \xsupv{m\colon \yal = \Opqm}$.

\begin{definition}
\label{df:qlm}
A sequence of elements
\wlee{ae1.2a1}
\yal_1,\yal_2,\ldots\in\cTLpinf%
\weee
has a limit
$\lim_{k\rightarrow \infty} \yal_k = \ybet$,
if $\lmii\yordq{\ybet-\yal_k} = +\infty$.
\end{definition}

The following theorem may be known, so we do not claim
credit for it. It appears here as a by-product
and it is an easy corollary of \ex{ae2.m4}.
\begin{theorem}
\label{th:alg}
The \TL\ element corresponding to the infinite \cbr\ equals the
\JWp:
\xlee{ae1.9}
\lim_{m\rightarrow+\infty} q^{\vthf mn(n-1)}\abrmn = \jwpnp,
\xeee
where $\jwpnp\in \cTLpinf_n$ corresponds to the \JWp\
$\jwpn\in\QcTLn$.
\end{theorem}
In fact, a more general statement is also true:
\xlee{ae1.10}
\lim_{m\rightarrow+\infty} q^{\vthf m(n-1)}\symalg{\btcyln^m} =
\jwpnp,
\xeee
but its proof is more technical and we omit it here.

\subsubsection{A bit of homological calculus}
\label{sss.homcal}

Let $\xKhA$ denote the homotopy category of complexes associated
with an additive category $\xctA$ (we have in mind a particular case of $\xKhA = \dTLp$).

A \chcpl\ is considered `homologically small' if it ends at a low
(that is, high negative)
homological degree.
Let
$\Ohpm$ denote a complex which ends at $(-m)$-th homological
degree: $\Ohpm = (\cdots \xAv{-m-1} \rightarrow\xAv{-m})$. We define
a homological order of a complex $\xbA$ as $\yordh{\xbA} =
\xsupv{m\colon\xbA\hteqv\Ohpm}$.

Two complexes connected by a \chmp: $\xbA\xrightarrow{\xbf}\xbB$
are considered `homologically close' if $\Cnbf$ is homologically
small.

A \emph{\chsq} is a sequence of complexes connected by \chmp s:
\ylee{ae1.10d1}
\scA = (\xbAz\xrightarrow{\xbfz} \xbAo
\xrightarrow{\xbfo}\cdots).
\yeee
\begin{definition}
\label{df:cauchy}
A \chsq\ $\scA$ is \emph{\Cch} if $\lmii \yordh{\Cnbfi} = \pinft$.
\end{definition}
\begin{definition}
\label{df:sqlm}
A \chsq\ has a limit
\footnote{This definition differs
from the standard categorical definition of a direct limit, however
Theorem\rw{pr:spmp} indicates that our definition implies the standard one. We expect that
both definitions are equivalent.}
: $\dlm\scA = \xbA$, where $\xbA$ is a \chcpl, if
there exist \chmp s $\xbAi\xrightarrow{\xbtfi}\xbA$ such that
they form commutative triangles
\xlee{ae1.10e1}
\cmtr{\xbfi}{\xbtfio}{\xbtfi}{\xbAi}{\xbAio}{\xbA}
\xeee
%
and $\lmii\yordhr{\Cnv{\xbtfi}} =\pinft$.
\end{definition}

In Section\rw{s:elhomcal} we prove the following homology versions
of standard theorems about limits
(Propositions\rw{pr:chlm},\rw{pr:lmch} and\rw{pr:lmun}):
\begin{theorem}
\label{th:lmt}
A \chsq\ $\scA$ has a limit if and only if it is \Cch.
\end{theorem}
\begin{theorem}
\label{th:lmt2}
The limit of a \chsq\ is unique up to homotopy equivalence.
\end{theorem}

\subsubsection{Infinite \cbr\ as a \JWp\ in the universal category}

For a tangle diagram $\xtau$ let
$\dtaus$ denote the categorification complex $\dtau$
with a degree shift proportional to the number $\crsv{\xtau}$ of crossings
in the diagram $\xtau$:
\xlee{ae1.10b1}
\dtaus = \dtau\tgrshv{\vthf}{-\vthf}{\vthf}^{\crsv{\xtau}}.
\xeee
%

In subsection\rw{ss:brchsq} we define a
\chsq\ of categorification complexes of \cbr s connected by special
\chmp s
%
%
%
%
\begin{multline}
\label{ae1.10c}
\xctBn =
\Big(
\cidbrn
\xratv{\mrfz}
\cbrons \xratv{\mrfo}
\cdots
\\
\cdots
\xrightarrow{\mrfmmo}
\cbrmns \xraov{\mrfm}
\cbrmons \xrightarrow{\mrfmo}\cdots\Big).
\end{multline}
We prove that
$\yordhr{\Cnv{\mrfm}}\geq2m(n-1) + 1$,
so
$\xctBn$ is a \Csq\ and by  Theorem\rw{th:lmt} it has a unique limit:
$\dlm\xctBn =\ctjwpn \in\dTLnp $.

\begin{theorem}
\label{th:enum}
The limiting complex $\ctjwpn$ has the following properties:
\begin{enumerate}

\item A composition of $\ctjwpn$ with cap- and \uptg s is contractible:
\ylee{eq:auc}
\ccapni \;\tcmp\ctjwpn \hteqv \ctjwpn\tcmp\; \ccupni\hteqv 0.
\yeee
\item The complex $\ctjwpn$ is idempotent with respect to tangle composition:
$\ctjwpn \tcmp\ctjwpn
\hteqv \ctjwpn$.

\end{enumerate}
\end{theorem}

We provide a glimpse into the structure of $\ctjwpn$.
A complex $\xbC$  in $\cTLn$ is called \emph{\odct} if
$\gidbrn$ never appears in \qcmds\ $\xCi$.
A complex $\xbC$ in $\cTLn$ is called \emph{\otbl} if the multiplicities
$\cjilam$ of \ex{ae1.8a} satisfy the property
\xlee{ae2.m1}
\cjmilam=0\qquad\text{if $i<0$, or $j<i$, or $j>2i$.}
\xeee
%


Let $\cbrmns\xraov{\xbtfm}\ctjwpn$ be \chmp s associated
with the limit $\dlm\xctBn = \ctjwpn$ in accordance with
Definition\rw{df:sqlm}.
\begin{theorem}
\label{th:cnpr}
There exist  \odct\ \otbl\ complexes $\wbCmn$ such that
$$\Cnv{\xbtfm}\hteqv\wbCmn\spshmnm\tgrsshv{1}{-1}.$$
\end{theorem}
\noindent
In other words, there exists a distinguished triangle
%
\ylee{ae2.m2}
\wbCmn\spshmnm\qsho \xratv{\chdlbtfm} \cbrmns \xrahv{\xbtfm} \ctjwpn
\xrightarrow{\;\;\;\;\;\;\;} \wbCmn\spshmnm\tgrsshv{1}{-1}
\yeee
so there is a presentation
%
\xlee{ae2.m4}
\ctjwpn \hteqv \CnBv{ \wbCmn\spshmnm\qsho
\xrahv{\chdlbtfm} \cbrmns},
\xeee
%
where the complex $\wbCmn$ is
\odct\ and \otbl.

At $m=0$ the formula\rx{ae2.m4} becomes
\xlee{ae2.m5}
\ctjwpn \hteqv \CnBv{ \wbCzn\qsho
\xrahv{\chdlbtfz} \cidbrn},
\xeee
where the complex $\wbCzn$ is \odct\ and \otbl.

%
%

Since $\wbCzn$ is \otbl, the complex $\Cnchdlbtfz$ is
also \otbl\ and consequently \qpb. Hence $\Kctjwpn$ is well-defined. Also
$\Kctjwpn\neq 0$, because it contains $\cidbrn$ with coefficient 1.
Theorem\rw{th:enum} indicates that
$\Kctjwpn$ satisfies
defining properties of the \JWp, hence by uniqueness it is the \JWp:
\begin{corollary}
The complex $\ctjwpn$ categorifies the \JWp\ in
$\cTLpinf$:
\xlee{eq:catKz}
\Kctjwpn = \jwpn.
\xeee
\end{corollary}




\section{Elementary homological calculus}
\label{s:elhomcal}

\subsection{Limits in the category of complexes}

Consider a category $\xChA$ of \chcpls\ associated with
an additive category $\xctA$.
An $i$-th \emph{\trnc} $\xtrniv{\xbA}$ of a \chcpl\ $\xbA$ is
the \chcpl\
$\xAmi\xrightarrow{\xdmi}\xAmio\rightarrow\cdots$.
An
\trnci\ of a \chmp\ $\xbf$ is defined similarly.

Define an \emph{\isor}
$\ysiobf$
of a chain map $\xbA\xrightarrow{\xbf}\xbB$  as
the largest number $i$ for which a truncated \chmp\ $\xtrniv{\xbf}$ is an
isomorphism of truncated complexes.

\begin{remark}
\label{rm:cnord}
Consider a \dstt\rx{ae1.ch2}.
If $\xbA\hteqv\Ohpm$, then $\ysiov{\idlbf}\geq m-1$.
\end{remark}

\begin{definition}
A \chsq\
$\scA = (\xbAo\xrightarrow{\xbfo}\xbAt\xrightarrow{\xbft}\cdots)$ is
\emph{\stblz} if
$\lim_{i\rightarrow\infty} \ysiobfi=\pinft$.
\end{definition}

\begin{definition}
A \chsq\ $\scA$ has a \tchlm\ $\chlm\scA=\xbA$ if there exist
\chmp s $\xbAi\xrightarrow{\xbtfi}\xbA$ such that
$\xbtfi = \xbtfio\,\xbfi$ and $\lmii \ysiorv{\xbtfi} = \pinft$.
\end{definition}

The following two theorems are easy to prove:
\begin{theorem}
A \chsq\ has a \tchlm\ if and only if it is \stblz.
If a \tchlm\ exists then it is unique.
\end{theorem}

%

\begin{theorem}
\label{pr:fnctchlm}
Suppose that $\chlm\scA = \xbA$.
Then for a complex $\xbB$ and \chmp s
$\xbAi\xrightarrow{\ybgi}\xbB$ such that $\ybgi = \ybgio\xbfi$,
 there exists a unique \chmp\ $\xbA\xrightarrow{\ybg}\xbB$
such that $\ybgi = \ybg\,\xbtfi$.
\end{theorem}

\begin{definition}
\label{df:chlmmp}
A sequence of \chmp s $\xbA\xrightarrow{\xbfz,\xbfo,\cdots}\xbB$
has a \tchlm\ $\lmii\xbfi = \xbf$ if for any $N$ there exists $N\p$
such that $\xtrnNv{\xbfi} = \xtrnNv{\xbf}$ for any $i\geq N\p$.
\end{definition}

\subsection{Limits in the homotopy category}

Definitions\rw{df:cauchy} and\rw{df:sqlm} extend the notion of a
\stblz\ \chsq\ and its limit to the homotopy category $\xKhA$:
obviously, a \stblz\ \chsq\ is \Cch, while $\chlm\scA=\xbA$ implies
$\dlm\scA=\xbA$.

\begin{proposition}
\label{pr:chlm}
A \Csq\ has a limit.
\end{proposition}
\proof
Consider a \Csq\ $\scA$. We construct a special \chcpl\
$\xbAs$ such that $\dlm\scA=\xbAs$ in accordance with
Definition\rw{df:sqlm}. Roughly speaking, we take $\xbAz$ and
attach to it the cones $\Cnbfi$ represented by homologically small
complexes, one by one. The result is a sequence $\scAs=\xbApz,\xbApo,\ldots$ of \stblz\
complexes $\xbApi$ such that $\xbApi\hteqv\xbAi$, and
$\xbAs=\chlm\scAs$ is their \tchlm.

Here is a detailed explanation.
By Definition\rw{df:cauchy}, there exist complexes $\ybCi$ such
that
\xlee{ae1.10a1}
\Cnv{\xbfi} \hteqv \ybCi[1] = \Ohpmi,\qquad\lmii m_i=+\infty.
\xeee
The complexes $\xbAi$, $\xbAio$ and $\ybCi$ form exact triangles:
\ylee{ae1.10a2}
\xymatrix{
\ybCi \ar[r]^-{\chdlbfi} & \xbAi \ar[r]^-{\xbfi} & \xbAio \ar[r] &
\ybCi[-1]
}
\yeee
and $\xbAio \hteqv \Cnv{\chdlbfi}$. We define recursively a new sequence
of complexes $\scAp = (\xbApz \xrightarrow{\idlbgz} \xbApo\xrightarrow{\idlbgo}\cdots)$ by
the relations $\xbApz = \xbAz$, $\xbApi\hteqv\xbAi$ and
$\xbApio = \Cnv{\ybgi}$, where the \chmp\
$\ybCi\xrightarrow{\ybgi}\xbApi$ is homotopy equivalent to the
\chmp\ $\chdlbfi$. In other words,
\xlee{ae1.10g1}
\xbApio = \Cnv{\ybCi\xrightarrow{\ybgi}
\underbrace{
\Cnv{\ybCimo\xrightarrow{\ybgimo}\cdots\xrightarrow{\ybgt}
\underbrace{
\Cnv{\ybCo\xrightarrow{\ybgo}
\underbrace{
\Cnv{\ybCz\xrightarrow{\chdlbfz}\xbAz
}}_{\xbApo}\;
}
}_{\xbApt}\;
}
}_{\xbApi}\;
}
\xeee

According to Remark\rw{rm:cnord}, $\ysiov{\idlbgi}\geq m_i$,
hence the sequence
$\scAp$ is \stblz, so there exists a chain limit
$\chlm\scAp = \xbAs$ and consequently there is a limit
$\dlm\scA=\xbAs$.\qed

Simply saying, the complex $\xbAs$ is an infinite \mtcn\ extension
of the complex\rx{ae1.10g1}:
\xlee{ae1.10g2}
\xbAs =
\cdots\xrightarrow{\ybgh}\Cnv{\ybCt
\xrightarrow{\ybgt}
\Cnv{\ybCo\xrightarrow{\ybgo}
\Cnv{\ybCz\xrightarrow{\chdlbfz}\xbAz
}
}
}.
\xeee

For our applications it is important to express $\Cnbtfz$ in terms
of complexes $\ybCi$. This can be done by rearranging
the infinite \mtcn\rx{ae1.10g2} with the help of associativity of
cone formation, which exists even within the category $\xChA$:
%
%
\xlee{ae1.10h1}
\xbAs = \Cnv{\ybtC\xrightarrow{\ybtg}\xbAz},\qquad
\ybtC =\cdots\xrightarrow{\ybht}
\Cnv{\ybCt[1]\xrightarrow{\ybho}\Cnv{\ybCo[1]\xrightarrow{\ybhz}\ybCz}},
\xeee
so that $\xbtfz \hteqv \idlv{\ybtg}$, and $\Cnbtfz\hteqv\ybtC[-1]$
is expressed in terms of complexes $\ybCi$
arranged into an infinite \mtcn\ $\ybtC$. Here is a more formal
statement.
\begin{theorem}
\label{th:rshfl}
For a \Csq\ $\scA$ there exists another \Csq\
$\yctC = (\ybCz
\xrightarrow{\ybhpz} \ybtCo\xrightarrow{\ybhpo}\cdots)$
and \chmp s
$\ybCi[1]\xrightarrow{\ybhi} \ybtCi$ such that
$\Cnv{\ybhi}=\ybtCio$, $\ybhpi = \idlv{\ybhi}$ and for the
limiting complex $\ybtC=\chlm\yctC$ there exists a \chmp\
$\ybtC\xrightarrow{\ybtg}\xbAz$ such that $\xbAs = \Cnv{\ybtg}$,
$\;\xbtfz \hteqv \idlv{\ybtg}$ and consequently $\Cnbtfz\hteqv\ybtC[-1]$.
\end{theorem}

\proof
%
%
Let us recall the associativity of cones in a general setting.
For a \chmp\ $\xbA\xrightarrow{\xbf}\xbB$,
a \chmp\ $\xbC\xrightarrow{\ybg} \Cnbf$
is a sum: $\ybg = \ybgA \oplus\ybgB$
\ylee{ae1.10f1}
\xymatrix{
 & \xbA\ar[d]^-{\xbf}
\\
\xbC \ar[ur]|{[1]}^-{\ybgA} \ar[r]_-{\ybgB}
& \xbB
}
\yeee
where $\xbC\xrightarrow{\ybgA}\xbA[-1]$ is a \chmp\ and
$\xbC\xrightarrow{\ybgB}\xbB$ is a \mmp. Now it is obvious
that
\xlee{ae1.10f2}
\Cnv{\xbC\xrightarrow{\ybg}\Cnv{\xbA\xrightarrow{\xbf}\xbB}}
=\Cnv{\Cnv{\xbC[1]\xrightarrow{\ybgA}\xbA}\xrightarrow{\ybgB\oplus\xbf}\xbB
}.
\xeee

We apply the associativity relation\rx{ae1.10f2} to
\mtcn s\rx{ae1.10g1} consecutively for $i=1,2,\ldots$ in order to rearrange
them, so that $\xbApi = \Cnv{\ybtCi\xrightarrow{\ybtgi}\xbAz}$,
while the complexes $\ybtCi$ and \chmp s $\ybtgi$ are defined
recursively: $\ybtCz=\ybCz$, $\ybtgz = \chdlbfz$, $\ybtCio = \Cnv{\ybhi}$,
while the \chmp s $\ybCi[1]\xrightarrow{\ybhi} \ybtCi$ and
$\ybtCio\xrightarrow{\ybtgio}\xbAz$
are defined by applying the associativity
relation\rx{ae1.10f2} to the  double cone on the second line of
the following equation:
\begin{equation}
\label{ae1.10f3}
\begin{split}
\xbApio & = \Cnv{\ybCi\xrightarrow{\ybgi}\xbApi}
\\
& = \Cnv{\ybCi\xrightarrow{\ybgi}\Cnv{\ybtCi\xrightarrow{\ybtgi}\xbAz} }
\\
& = \Cnv{\Cnv{\ybCi[1]\xrightarrow{\ybhi} \ybtCi
} \xrightarrow{\ybtgio}
\xbAz
}
\\
& = \Cnv{\ybtCio\xrightarrow{\ybtgio}\xbAz}.
\end{split}
\end{equation}
Distinguished triangles
\ylee{ae1.10f3}
\xymatrix{
\ybCi[1]\ar[r]^-{\ybhi}
&
\ybtCi \ar[r]^-{\idlbhi}
&
\ybtCio \ar[r]
&
\ybCi
}
\yeee
determine  \chmp s $\ybhpi=\idlbhi$ of the \chsq\  $\yctC = (\ybtCz
\xrightarrow{\ybhpz} \ybtCo\xrightarrow{\ybhpo}\cdots)$. By
Remark\rw{rm:cnord} it has a \tchlm: $\chlm\yctC =
\ybtC$, which is an infinite \mtcn:
\ylee{ae1.10f4}
\ybtC =\cdots\xrightarrow{\ybht}
\Cnv{\ybCt[1]\xrightarrow{\ybho}\Cnv{\ybCo[1]\xrightarrow{\ybhz}\ybCz}}.
\yeee
The \chmp s $\ybtCi\xrightarrow{\ybhpi}\ybtCio$ satisfy a relation
$\ybtgi = \ybtgio\,\ybhpi$, so by Theorem\rw{pr:fnctchlm} there
exists a unique \chmp\ $\ybtC \xrightarrow{\ybtg}\xbAz$ such that
$\ybtgi=\ybtg\,\ybhtpi$.
It is easy to show that $\xbAs = \Cnv{\ybtC\xrightarrow{\ybtg}\xbAz}$,
and $\xbtfz = \idlv{\ybtg}$, hence $\Cnbtfz\hteqv \ybtC$.\qed

It is easy to prove the analog of Theorem\rw{pr:fnctchlm}:
\begin{theorem}
\label{pr:spmp}
For a complex $\xbB$ and \chmp s
$\xbAi\xrightarrow{\ybgi}\xbB$ such that $\ybgi \hteqv \ybgio\xbfi$,
 there exists a unique (up to homotopy) \chmp\ $\xbAs\xrightarrow{\ybg}\xbB$
which forms commutative triangles
\xlee{ae1.10f1}
\cmtr{\xbtfi}{\ybg}{\ybgi}{\xbAi}{\xbAio}{\xbB}
\xeee
%
\end{theorem}

In order to complete the proof
of Theorems\rw{th:lmt} and\rw{th:lmt2},
we need two simple propositions. The first one establishes a
triangle inequality for homological orders of cones.
\begin{proposition}
If three \chmp s form a commutative triangle
\xlee{ae1.10c1}
\xymatrix{
\xbA \ar[r]_-{\xbfAB} \ar@/^1pc/[rr]^-{\xbfAC} &
\xbB \ar[r]_-{\xbfBC} &
\xbC
},\qquad \xbfAC\hteqv \xbfBC\xbfAB.
\xeee
then the homological orders of their cones satisfy the inequalities
\begin{align}
\label{eq:in1}
\yordch{\xbfAB}\geq \min\big(\yordch{\xbfAC},\yordch{\xbfBC}-1\big),
\\
\label{eq:in2}
\yordch{\xbfBC}\geq \min\big(\yordch{\xbfAB}+1,\yordch{\xbfAC}\big).
\end{align}
\end{proposition}
\proof
If \chmp s form a commutative triangle\rx{ae1.10c1}, then their cones form a \dstt
\ylee{ae1.10c2}
\Cnv{\xbfAB}\xrightarrow{\ybgo} \Cnv{\xbfAC} \xrightarrow{\ybgt} \Cnv{\xbfBC}
\xrightarrow{\ybgh}\Cnv{\xbfAB}[1],
\yeee
so the first inequality follows from the relation
$\Cnv{\xbfAB}\hteqv\Cnv{\ybgt}[1]$ and the second inequality
follows from the relation $\Cnv{\xbfBC} \hteqv \Cnv{\ybgo}$.\qed

The second proposition says that if a complex is homologically
infinitely small then it is contractible.
\begin{proposition}
\label{pr:ismc}
If $\yordh{\xbA} = +\infty$ then $\xbA$ is contractible.
\end{proposition}
\proof Since $\yordh{\xbA} = +\infty$, there exist complexes
$\xbAi\hteqv\xbA$, such that $\xbAi=\Ohpmi$ and $\lmii m_i =
+\infty$. Consider a sequence of \chmp s establishing homotopy
equivalence between the complexes:
\ylee{ae1.10c3}
\xymatrix@C=0.5cm{
\xbA \ar@<0.6ex>[r]^-{\xbfz}
&
\xbAo
\ar@<0.3ex>[l]^-{\ybgz}
\ar@<0.6ex>[r]^-{\xbfo}
&
\xbAt
\ar@<0.3ex>[l]^-{\ybgo}
\ar@<0.6ex>[r]
&
\cdots
\ar@<0.3ex>[l]
\ar@<0.6ex>[r]
&
\xbAi
\ar@<0.3ex>[l]
\ar@<0.6ex>[r]^-{\xbfi}
&
\xbAio
\ar@<0.3ex>[l]^-{\ybgi}
\ar@<0.6ex>[r]
&
\cdots
\ar@<0.3ex>[l]
}
,\qquad
\yIdAi-\ybgi\xbfi =
\atcmv{\xbdi}{\xbhi},
\yeee
where $\yIdAi$ is the identity
\chmp\ of $\xbAi$, while $\xbAi[1]\xrightarrow{\xbhi}\xbAi$ is a homotopy \chmp\ (it does
not commute with the chain differential $\xbdi$ in the complex $\xbAi$).

Consider the compositions $\xbhfi = \xbfi\cdots\xbfo\xbfz$,
$\xbhgi = \ybgz\ybgo\cdots\ybgi$ and $\xbhhi =
\xbhgimo\,\xbhi\,\xbhfimo$. It is easy to see
that $\xbhgimo\,\xbhfimo - \xbhgi\,\xbhfi = \atcmv{\xbd}{\xbhhi}$,
hence $\yIdA - \xbhgi\,\xbhfi = \atcmv{\xbd}{\xbchi}$, where
$\xbchi = \xbhhz + \xbhho +\cdots + \xbhhi$.
There is a limit (\cf Definition\rw{df:chlmmp}) $\lmii\xbchi = \xbch$, while
$\lmii\xbhgi\,\xbhfi = 0$, hence $\yIdA = \atcmv{\xbd}{\xbch}$
which means that the complex $\xbA$ is contractible.
\qed

\begin{proposition}
\label{pr:lmch}
If a \chsq\ $\scA$ has a limit, then it is \Cch.
\end{proposition}
\proof
The inequality\rx{eq:in1} applied to the commutative
triangle\rx{ae1.10e1} says that
$$\yordch{\xbfi}\geq
\min\lrbc{\zordch{\xbtfi},\zordch{\xbtfio}-1},$$
hence the limit $\lmii\zordch{\xbtfi} = +\infty$ implies the \Cch\
property of $\scA$.

\begin{proposition}
\label{pr:lmun}
If a \chsq\ $\scA$ has a limit then it is unique.
\end{proposition}
\proof
If $\scA$ has a limit then by Proposition\rw{pr:lmch} it is \Cch.
Hence it has a special limit $\xbAs$ described in the
proof of Proposition\rw{pr:chlm}. If $\scA$ has another limit
$\xbAp$ with \chmp s $\xbAi\xrightarrow{\xbtfpi}\xbAp$ then by
Theorem\rw{pr:spmp} there is a \chmp\
$\xbAs\xrightarrow{\ybg}\xbAp$ with commutative
triangles\rx{ae1.10f1}. The inequality\rx{eq:in2} says
\ylee{ae1.10c3}
\yordch{\ybg}
\geq \min\big(\zordch{\xbtfi}+1,\zordch{\ybgi}\big).
\yeee
Since both cones in the \rhs become homologically infinitely small
at $i\rightarrow +\infty$, the cone $\Cnv{\ybg}$ is also
homologically infinitely small. Then Proposition\rw{pr:ismc} says
that $\Cnv{\ybg}$ is contractible and as a result
$\xbAp\hteqv\xbAs$.\qed

We end this section with a theorem which follows easily from
Definition\rw{df:sqlm}.
\begin{theorem}
\label{th:zl}
If a \chsq\ $\scA$ satisfies the property $\lmii\yordh{\xbAi} =
+\infty$ then its limit is contractible: $\dlm\scA = 0$.
\end{theorem}

\section{A \chsq\ of
categorification complexes of  \cbr s} 
\label{s:cbr}
\subsection{A special categorification complex of a \ngbr}

Let $\xsgi$ denote an elementary negative $n$-strand braid:
\ylee{ae2.1}
\xsgi=\xygraph{
!{0;/r1.5pc/:}
[r(0.25)u(0.5)]
!{\xcapv@(0)}
[u(0.5)r(1)]
*{\cdots}
[r(01)u(0.5)]
!{\xcapv@(0)}
[r(0.5)u(1)]
!{\vcross}
[r(1.5)u(1)]
!{\xcapv@(0)}
[u(0.5)r(1)]
*{\cdots}
[r(01)u(0.5)]
!{\xcapv@(0)}
[d(0.5)l(3.5)]
*{\scriptstyle{i}}
[r(1)]
*{\scriptstyle{i+1}}
[l(3.5)]
*{\scriptstyle{1}}
[r(6)]
*{\scriptstyle{n}}
}
\yeee

\begin{theorem}
\label{th:prop}
If an $n$-strand braid $\brb$ can be presented as a product of elementary
negative braids: $\brb = \xsgiv{k}\cdots\xsgiv{2}\xsgiv{1}$, then
its categorification complex has a special presentation $\cbrbs$:
\xlee{aea2.1}
\cbrbas =
\Big(\ldots\rightarrow\xCmt\rightarrow\xCmo\rightarrow\cidbrn\Big)
\xeee
such that the complex
\xlee{aea2.2}
\xbC = (\ldots\rightarrow\xCmt\rightarrow\xCmo)\tgrsshv{-1}{1}
\xeee
is
\odct\ and \otbl.



\end{theorem}

More abstractly, the theorem says that there exists a
\odct\ and \otbl\ complex $\xbC$ and a \chmp\
$\xbC\rightarrow\cidbrn$ such that $\cbrba \hteqv
\CnBv{\xbC\qsho\rightarrow\cidbrn}$.

\begin{remark}
\label{rm:spbl}
Theorem\rw{th:prop} implies that the special complex
$\cbrbas$ is  \otbl.
\end{remark}

\pr{Theorem}{th:prop}
Let $\xlam$ be a \TL\ \ttngnn. Fix $i$ such that $1\leq i\leq n-1$. If the
composition $\gcapni\tcmp\xlam$ does not contain a disjoint circle,
then, in accordance with \ex{ae1.7},
we define the special categorification complex of $\xsgi\tcmp\xlam$ as
%
\xlee{ae2.3}
\symcatps{\xsgi\tcmp\xlam}  =
\Big(\symbcat{\xUni\tcmp\xlam}\tgrshv{1}{-1}{1}
\rightarrow \dlam \Big)
\xeee
%
If $\gcapni\tcmp\xlam$ contains a disjoint circle, then $\xlam$ must
have the form $\gcupni\tcmp\xlamp$. Hence
$\xsgi\tcmp\xlam=\xsgi\tcmp\gcupni\tcmp\xlamp$. The tangle $\xsgi\tcmp\gcupni$
is the same as $\gcupni$ with a positive framing twist, so
according to \ex{ae1.8},
$\bsymcat{\xsgi\tcmp\gcupni} = \ccupni \tgrshv{\vthh}{-\vthf}{-\vthf}$.
Hence in this case we define the special categorification complex
of $\xsgi\tcmp\xlam$ simply as shifted $\dlam$:
\xlee{ae2.4}
\symcatps{\xsgi\tcmp\xlam} = 
\dlam
\tgrshv{2}{-1}{0}.
\xeee
%

Now we define a recursive algorithm for constructing the complex
$\cbrbas$. For $\brb = \gidbrn$ we define $\cbrbas = \cidbrn$. Let
$\brb = \xsgiv{k}\tcmp\cdots\tcmp\xsgiv{1}$ and
suppose that we have defined its special complex $\cbrbas$. We
define the special categorification complex of a
braid $\brbp=\xsgikpo\tcmp\brb$ by applying the rules\rx{ae2.3}
and\rx{ae2.4} to all constituent tangles $\xlam$ in the complex
$\cbrbs$ (see the formula\rx{ae1.8a}).

We prove the properties of $\cbrbas$ by induction over $k$.
If $k=0$ then $\brb = \gidbrn$ and the properties of $\cbrbas$ are
obvious.

Suppose that the special categorification complex
$\cbrbas$ of a braid $\brb = \xsgiv{k}\tcmp\cdots\tcmp\xsgiv{1}$
has the form\rx{aea2.1} and its tail\rx{aea2.2} is \odct\ and
\otbl.
Consider
a longer braid
$\brbp=\xsgikpo\tcmp\brb$. The object $\cidbrn$ may appear in $\cbrbpas$ if
and only if $\xlam=\gidbrn$ and the extra crossing $\xsgikpo$ is
\nsplcd\ in \ex{ae2.3}, hence
$\cbrbpas$ has the form\rx{aea2.1} and its tail\rx{aea2.2} is
\odct.

If the negative crossing $\xsgikpo$ is composed with the head
$\cidbrn$ of the complex\rx{aea2.1}, then the formula\rx{ae2.3}
applies and the tangle $\xUv{n}{i_k+1}$ appearing in the tail of
$\cbrbpas$ satisfies the property\rx{ae2.m1}.


If the crossing $\xsgikpo$ is composed with a
\TL\ tangle $\xlam$ from the $(-i)$-th \qcmd\ $\xCmi$ (see \ex{ae1.8a})
in the tail of the complex $\cbrbas$ with the $q$-degree shift $j$
satisfying the inequality $i-1 \leq j-1 \leq 2(i-1)$, then the
shifted objects in the \rhs of \eex{ae2.3} and\rx{ae2.4}
also satisfy this inequality.\qed

The picture\rx{ae1.10p} presents a \cbr\ as a product of negative
crossings, hence
\begin{corollary}
\label{cr:otbl}
A \cbr\ $\gbrmn$ has a special \otbl\ categorification complex
$\cbrmnps$. In particular, for $m=1$
\xlee{ae2.5}
\cbronps
 = \CnBv{\xbCon\qsho\rightarrow\cidbrn},
\xeee
where the complex
$\xbCon$ is \odct\ and \otbl.
\end{corollary}

\subsection{Special morphisms between \cbr\ complexes}
\label{ss:brchsq}

Relation\rx{ae2.5}  indicates that there is a \dstt\
%
%
%
\ylee{ae2.6}
\xbCon\qsho \longrightarrow
\cidbrn \xratv{\mrfo}
\cbrons \longrightarrow
\xbCon\tgrsshomo
\yeee
and
\xlee{ae2.6a}
\Cnv{\mrfo} \hteqv \xbCon\tgrsshomo.
\xeee
%
%
Composing both sides of the morphism $\mrfo$ with
the \cbr\ complex
$\cbrmns$,
we get a morphism
%
\ylee{ae2.7}
\cbrmns \xratv{\mrfm}\cbrmons
\yeee
such that
\xlee{ae2.8}
\Cnv{\mrfm} \hteqv \Cnv{\mrfo}\tcmp\cbrmns.
\xeee
\begin{theorem}
\label{th:2.1}
The cone\rx{ae2.8} can be presented by a shifted complex
\ylee{ae2.9}
\Cnv{\mrfm} \hteqv \xbCmn \tgrsshnontm\tgrsshomo,
\yeee
%
such that $\xbCmn$ is \odct\ and \otbl.
\end{theorem}

The proof is based on a simple geometric lemma:
\begin{lemma}
\label{l:1}
For $n\geq 2$, the following two compositions of framed tangles are isotopic:
%
\xlee{ae2.b}
\gcapni \tcmp\;\gbron = \gbronmt\;\tcmp\gcapnit
\xeee
where $\gcapnit$ is the tangle $\gcapni$ with double framing twist
on the cap:
\ylee{ae2.10}
\gcapnik=
\xygraph{
!{0;/r1.5pc/:}
[r(0.25)u(0.5)]
!{\xcapv@(0)}
[u(0.5)r(1)]
*{\cdots}
[r(01)u(0.5)]
!{\xcapv@(0)}
[r(0.5)]
!{\vcap}
[r(1.5)u(1)]
!{\xcapv@(0)}
[u(0.5)r(1)]
*{\cdots}
[r(01)u(0.5)]
!{\xcapv@(0)}
[d(0.5)l(3.5)]
*{\scriptstyle{i}}
[r(1)]
*{\scriptstyle{i+1}}
[l(3.5)]
*{\scriptstyle{1}}
[r(6)]
*{\scriptstyle{n}}
[l(3)u(1)]
*{\symfr}
[u(0.5)]
*{\scriptstyle{k}}
}
\yeee
\end{lemma}
\proof
This lemma is geometrically obvious: a cap on a pair of adjacent strands slides down
through the \cbr\ to the
bottom.\qed

An immediate corollary of \ex{ae2.b} and of the framing change
formula\rx{ae1.8} is the following relation:
%
\xlee{ae2.11}
\bsymcats{\gcapni \tcmp\gbrmn} \hteqv \bsymcats{\gbrmnmt\;\tcmp\gcapni}
\tgrsshnontm.
\xeee


In order to prove Theorem\rw{th:2.1}, we need three simple
propositions.
For a positive integer $d\leq \frac{n}{2}$,
let $\stI=(i_1,\ldots,i_d)$ be a sequence of positive integer
numbers such that $i_k<n-2k+2$ for all $k\in\{1,\ldots,d\}$.
A \emph{\aptg} $\gcapnI$
is a $(n,n-2d)$-tangle which
can be presented as a product of $d$ tangles of the form
$\gcapv{m}{i}{0.75}$:
%
\ylee{aes2.1a}
\gcapnI =
\gcapv{n-2d+2}{i_d}{2}\tcmp\cdots\tcmp
\gcapv{n-2}{i_2}{1.5}
\tcmp
\gcapv{n}{i_1}{0.75}.
\yeee
A \emph{\uptg} $\gcupnI$ is defined similarly:
%
\ylee{aes2.2a}
\gcupnI =
\gcupv{n}{i_1}{-0.75}
\tcmp
\gcupv{n-2}{i_2}{-1.25}
\tcmp
\cdots
\tcmp
\gcupv{n-2d+2}{i_d}{-2.25}
.
\yeee
The first proposition is obvious:
\begin{proposition}
\label{pr:3}
Every \TL\ \ttngnn\ $\xlam$ has a presentation
\xlee{aes2.3a}
\xlam = \gcupnIp\tcmp
\gcapnI,\qquad
\nI=\nIp.
\xeee
\end{proposition}
The number $\cpdlam=\nI=\nIp$ is determined by the tangle $\xlam$
and we call it  the \apdg\ (or \updg) of $\xlam$.

The second proposition is also obvious:
\begin{proposition}
\label{pr:1}
If at least one of two complexes $\xbCo$ and $\xbCt$ in $\dTLn$ is
\odct\ then their composition $\xbCo\tcmp\xbCt$ is \odct.
\end{proposition}
Note that even if
both complexes are \otbl,  then their composition is not necessarily
\otbl. Indeed, in contrast to the homological degree,
the \qdgr\ is not additive with respect to the composition of
tangles:
%
if  the composition of two \TL\ tangles contains a disjoint
circle then the \qdgr\ shifts of the rule\rx{ae1.01}
violate additivity. However, if the upper tangle in the composition has no caps or the
lower tangle has no cups then no circles are created and the
\otblc\ is maintained:
\begin{proposition}
\label{pr:2}
If a complex $\xbC$ in $\dTLv{n-2\cpdlam}$ is \otbl, then the complexes
$\ccupnI\tcmp\xbC$
and
$\xbC\tcmp\ccapnI$
are also \otbl.
\end{proposition}


\pr{Theorem}{th:2.1}
In order to construct the \odct\ and \otbl\ complex $\xbCmn$, we
use the presentation
\xlee{ae2.12}
\Cnv{\mrfm} \hteqv
\xbCon\tcmp\cbrmns\tgrsshomo,
\xeee
which follows from \eex{ae2.8} and\rx{ae2.6a}.
We construct $\xbCmn$ by
simplifying the complexes
$\bsymcats{\xlam\tcmp\gbrmn}$
for \TL\ \ttngnn s $\xlam$
appearing in the \qcmds\ of $\xbCon$, with the help of the
relation\rx{ae2.11}, thus creating necessary degree shifts, and then
using Corollary\rw{cr:otbl} which says that emerging \cbr s have \otbl\ categorification
complexes.




Let
$\dlam\tgrsshjmi$
be an object appearing in the $(-i)$-th \qcmd\ of
$\xbCon$ with a non-zero multiplicity (we made its homological degree explicit by
including $-i$ in the shift).
%
We apply \ex{ae2.11} consequently to every cap $\gcapnk$ appearing
in the \aptg\ $\gcapnI$ in the presentation\rx{aes2.3a} of
$\xlam$:
\begin{multline}
\label{ae2.13}
\dlam\tgrsshjmi\tcmp\cbrmns
\\
\hteqv
\Bigg(
\ccupnIp\tcmp
\bsymcatps{\gbrv{m}{n-2\cpdlam}{2.5}}\tcmp\ccapnI
\tgrsshv{\alm}{-\blm}^{2m} \tgrsshjmi\Bigg)
\tgrsshnontm,
\end{multline}
where
\xlee{ae2.14}
\alm = \sum_{k=1}^{\cpdlam-1}(n-2k)
,\qquad
\blm = \sum_{k=1}^{\cpdlam-1}(n-2k-1).
\xeee

The object $\dlam$ comes from the \odct\ complex $\xbCon$,
hence $\cpdlam>0$ and
the complex
in big brackets in the \rhs of \ex{ae2.13} is \odct\ in view of
Proposition\rw{pr:1}. Proposition\rw{pr:2} implies that the complex
$\ccupnIp\tcmp
\bsymcatps{\gbrv{m}{n-2\cpdlam}{2.5}}\tcmp\ccapnI$
is also \otbl. Since $\dlam$
comes from the \otbl\ complex $\xbCon$, the numbers $i$ and $j$
satisfy inequalities $i\geq 0$ and $i\leq j\leq 2i$. It is easy to check that the numbers
$\alm$ and $\blm$ of \ex{ae2.14} satisfy the same inequalities:
$\blm\geq 0$, $\blm\leq \alm \leq 2\blm$, hence the complex in big
brackets in the \rhs of \ex{ae2.13} is also \otbl. The
complex
$\xbCon\tcmp\cbrmns$
in the \rhs of \ex{ae2.12} is
composed of complexes\rx{ae2.13}, hence Theorem\rw{th:2.1} is proved.
\qed

\section{A categorified \JWp}
\label{s:prfs}

Consider the \chsq\rx{ae1.10c}. Theorem\rw{th:2.1} implies that
$\yordhr{\Cnv{\mrfm}} \geq 2m(n-1)+1$,
hence $\xctBn$ is
\Cch\ and it has a unique limit $\dlm\xctBn =\ctjwpn \in\dTLnp$.

Now we prove Theorems\rw{th:enum} and Theorem\rw{th:cnpr} which
describe the properties of $\ctjwpn$.

\pr{Theorem}{th:cnpr}
Consider the \chsq\rx{ae1.10c} truncated from below:
\ylee{eq:np1}
\xctBmn =
\Big(
\cbrmns \xraov{\mrfm}
\cbrmons \xrightarrow{\mrfmo}\cdots\Big)\longrightarrow\ctjwpn.
\yeee
According to Theorem\rw{th:rshfl}, the limit $\ctjwpn$ can be
presented as a cone\rx{ae2.m4}, where $\wbCmnp = \wbCmn\spshmnm$
and $\wbCmn$ is an infinite
\mtcn:
\begin{multline}
\nonumber
\wbCmn =\cdots\rightarrow\Cnv{\xbCvn{m+k}\tgrsshv{2kn}{-2k(n-1)+1}
\rightarrow
\cdots
\\
\cdots
\rightarrow
\Cnv{
\xbCvn{m+1}\tgrsshv{2n}{-2n+3}\rightarrow\xbCmn}
}
\end{multline}
with \odct\ and \otbl\ complexes $\xbCmn$ introduced in
Theorem\rw{th:2.1}. Hence the complex $\wbCmn$ itself is \odct\ and
\otbl.\qed

%

\pr{part 1 of Theorem}{th:enum}
The tangle composition with $\ccapni$ is a `continuous'
functor, that is, it can be applied to both sides of
$\dlm\xctBn = \ctjwpn$, hence $ \dlm\;
\ccapni\tcmp\xctBn = \ccapni\tcmp\ctjwpn$. According to \ex{ae2.11},
\begin{equation}
\nonumber
\begin{split}
\ccapni\tcmp\xctBn & = \Big(\ccapni\tcmp\cidbrn\rightarrow \cdots\rightarrow
\ccapni\tcmp\cbrmns\rightarrow\cdots\Big)
\\
& = \Big( \ccapni\rightarrow\cdots\rightarrow
\cbrmnmts\tcmp\ccapni
\spshmnm
\rightarrow\cdots
\Big).
\end{split}
\end{equation}
Since
\ylee{ae3.3}
\yordhb{\cbrmnmts\tcmp\ccapni
\spshmnm}
= 2m(n-1)\xrightarrow[m\rightarrow +\infty]{}+ \infty,
\yeee
according
to Theorem\rw{th:zl}, $\dlm\ccapni\tcmp\xctBn=0$, hence
$\ccapni\tcmp\ctjwpn$ is contractible.\qed

\begin{remark}
The contractibility of $\ctjwpn\tcmp\ccupni$ is proved similarly.
\end{remark}

\begin{corollary}
\label{cr:odct}
If $\xbC$ is a \odct\ complex in $\dTLnp$, then $\xbC\tcmp\ctjwpn$ is
contractible.
\end{corollary}

\pr{part 2 of Theorem}{th:enum}
According to
\ex{ae2.m5},
\begin{multline}
\nonumber
\ctjwpn\tcmp\ctjwpn \hteqv \CnBv{\wbCzn\qsho\longrightarrow\cidbrn}\tcmp\ctjwpn
\\
\hteqv\CnBv{\wbCzn\tcmp\ctjwpn\qsho\longrightarrow\cidbrn\tcmp\ctjwpn}
\hteqv\ctjwpn,
\end{multline}
where we used the fact that $\wbCzn$ is \odct\ and Corollary\rw{cr:odct} in order to establish the last
homotopy equivalence.\qed

\pr{Theorem}{th:alg}
The complexes $\ctjwpn$, $\wbCmn$ and $\cbrmns$ in \ex{ae2.m4} are
\otbl, hence they are \qpb\ and their $\Kz$ images are
well-defined. Applying $\Kz$ to this equation and taking into account \ex{eq:catKz} and
the definition\rx{ae1.10b1}, we find
\ylee{ae3.5}
\jwpn = q^{\vthf mn(n-1)}\abrmn - q^{2mn+1}\Kz(\wbCmn).
\yeee
The complex $\wbCmn$ is \otbl, so $\yordq{\Kz(\wbCmn)}\geq 0$ and
by Definition\rw{df:qlm} there is a limit\rw{ae1.9}.\qed

\section{The other projector}
A dual of an \ttngmn\ $\xtau$ is the \ttngnm\ tangle
$\xtaud$ which is its mirror image. Duality extends to an
isomorphism $\cTL \xrightarrow{\dsym} \cTLop$ combined with the
isomorphism of the ground ring $\Zqqi\xrightarrow{\dsym}\Zqqi$, such that
$\dsymv{q} = q^{-1}$. Furthermore, duality establishes an
isomorphism $\cTLpinf\xrightarrow{\dsym}\cTLminfop$, where
$\cTLminf$ is the analog of $\cTLpinf$ constructed over the ring
$\Zsqiq$ of Laurent series in $q^{-1}$.

Since the relations\rx{ae1.4} and\rx{ae1.4a} are dual to each other,
while the idempotency condition $\jwpn\tcmp\jwpn=\jwpn$ is duality
invariant, the uniqueness of the \JWp\ implies that it is duality
invariant: $\dsymv{\jwpn} = \jwpn$. Hence the corresponding
idempotents $\jwpnp\in\cTLpinf$ and $\jwpnm\in\cTLminf$ are also
dual to each other: $\jwpnm = \dsymv{(\jwpnp)}$. Taking the dual
of \ex{ae1.9} we find that $\jwpnm$ is the limit of \cbr s with
high positive (that is, \cclckw) twist:
\xlee{ae1.9b}
\lim_{m\rightarrow+\infty} q^{-\vthf mn(n-1)}\aobrmn = \jwpnm,
\xeee
because $\dsymv{\Big(\gbrmn\Big)} = \gobrmn$.


Duality extends to a contravariant
equivalence functor $\dTL\xrightarrow{\dsym}\dTLop$, where
$\dTLop$ is the same category as $\dTL$, except that the
composition of tangles is performed in reversed order. The functor
$\dsym$ also switches the signs of all three gradings of $\dTL$.
Applying the duality functor to the construction of $\ctjwpn$ we
find that there exists a \chsq
\begin{multline}
\label{ae1.10f}
\xctBnd =
\Big(
\cidbrn
\rxratv{\dmrfz}
\cobrons \rxratv{\dmrfo}
\cdots
\\
\cdots
\xleftarrow{\dmrfmmo}
\cobrmns \rxraov{\dmrfm}
\cobrmons \xleftarrow{\dmrfmo}\cdots\Big),
\end{multline}
where $-\xspsh$ denotes the grading shift which is opposite
to\rx{ae1.10b1}.
The system\rx{ae1.10f} is
dual to the system\rx{ae1.10c}
and it has an inductive limit $\ilm \xctBnd =\ctjwmn $, which satisfies
projector properties
\ylee{ae1.10f1}
\ctjwmn\tcmp \ctjwmn \hteqv \ctjwmn,\qquad
\ccapni \tcmp\ctjwmn \hteqv \ctjwmn\tcmp \ccupni\hteqv 0
\yeee
and has a presentation
\ylee{ae1.10f2}
\ctjwpn \hteqv \CnBv{ \dsymv{\wbCmn}\ospshmnm\qshmo
\xrahv{\dsymv{\chdlbtfm}} \cobrmns},
\yeee
%
where the complex $\wbCmn$ is
\odct\ and \otbl. In particular, at $m=0$ we get the dual of
presentation\rx{ae2.m5}:
\ylee{ae1.10f3}
\ctjwmn \hteqv \CnBv{ \dsymv{\wbCzn}\qshmo
\xrahv{\dsymv{\chdlbtfz}} \cidbrn},
\yeee
where the complex $\wbCzn$ is \odct\ and \otbl.






\begin{bibdiv}
\begin{biblist}

\bib{Arm11}{article}
{
author={Armond, Cody}
title={The head and tail conjecture for alternating knots}
eprint={arXiv:1112.3995}
}

\bib{BN05}{article}
{
author={Bar-Natan, Dror}
title={Khovanov's homology for tangles and cobordisms}
journal={Geometry and Topology}
volume={9}
year={2005}
pages={1443-1499}
eprint={math.GT/0410495}
}

\bib{CK10}{article}
{
author={Cooper, Benjamin}
author={Krushkal, Vyacheslav}
title={Categorification of the Jones-Wenzl Projectors}
journal={Quantum Topology}
volume={3}
year={2012}
pages={139-180}
eprint={arXiv:1005.5117}
}

\bib{AD11}{article}
{
author={Armond, Cody}
author={Dasbach, Oliver}
title={Rogers-Ramanujan type identities and the head and tail of the colored Jones polynomial}
eprint={arXiv:1106.3948}
}

\bib{AD12}{article}
{
author={Armond, Cody}
author={Dasbach, Oliver}
year={2012}
status={in preparation}
}


\bib{FSS11}{article}
{
author={Frenkel, Igor}
author={Stroppel, Catharina}
author={Sussan, Joshua}
title={Categorifying fractional Euler characteristics, Jones-Wenzl projector and $3j$-symbols}
eprint={arXiv:1007.4680}
}

\bib{GL11}{article}
{
author={Garoufalidis, Stavros}
author={Le, Thang T. Q.}
title={Nahm sums, stability and the colored Jones polynomial}
eprint={arXiv:1112.3905}
}

\bib{GS11}{article}
{
author={Gukov, Sergei}
author={\Stosic, Marko}
title={Homological algebra of knots and BPS states}
eprint={arXiv:1112.0030}
}

\bib{Kh99}{article}
{
author={Khovanov, Mikhail}
title={A categorification of the Jones polynomial}
journal={Duke Mathematical Journal}
year={1999}
volume={101}
pages={359--426}
eprint={arXiv:math/9908171}
}

\bib{Ro11}{article}
{
author={Rozansky, Lev}
title={An infinite torus braid yields a categorified Jones-Wenzl projector}
eprint={arXiv:1005.3266}
}

\bib{Wi11}{article}
{
author={Witten, Edward}
title={Fivebranes and knots}
eprint={arXiv:1101.3216}
}
\end{biblist}
\end{bibdiv}

\begin{bibdiv}
\begin{biblist}

\bib{Mort95}{article}
{
author={Morton, H.R.}
title={The coloured Jones function and Alexander polynomial for torus knots}
journal={Math. Proc. Camb. Phil. Soc.}
volume={117}
year={1995}
pages={129-135}
}



\bib{BN1}{article}
{
author={Bar-Natan, Dror}
title={Khovanov's homology for tangles and cobordisms}
journal={Geometry and Topology}
volume={9}
year={2005}
pages={1443-1499}
eprint={arXiv:math.GT/0410495}
}

\bib{CK}{misc}
{
author={Cooper, Ben}
author={Krushkal, Slava}
title={Categorification of the Jones-Wenzl projectors}
note={in preparation}
}

\bib{FSS}{misc}
{
author={Frenkel, Igor}
author={Stroppel, Catharina}
author={Sussan, Joshua}
note={in preparation}
}



\bib{Kh1}{article}
{
author={Khovanov, Mikhail}
title={A categorification of the Jones polynomial}
journal={Duke Journal of Mathematics}
volume={101}
year={2000}
pages={359-426}
eprint={arXiv:math.QA/9908171}
}

\bib{KR1}{article}
{
author={Khovanov, Mikhail}
author={Rozansky, Lev}
title={Matrix factorizations and link homology}
journal={Fundamenta Mathematicae}
volume={199}
year={2008}
pages={1-91}
eprint={arXiv:math.QA/0401268}
}

\bib{KhRS}{misc}
{
author={Khovanov, Mikhail}
author={Rozansky, Lev}
note={in preparation}
}

\bib{St}{article}
{
author={Stosic, Marko}
title={Homological thickness and stability of torus knots}
journal={Algebraic and Geometric Topology}
volume={7}
year={2007}
pages={261-284}
eprint={arXiv:math.GT/0511532}
}

\end{biblist}
\end{bibdiv}

\begin{bibdiv}
\begin{biblist}


\bib{BN1}{article}
{
author={Bar-Natan, Dror}
title={Khovanov's homology for tangles and cobordisms}
journal={Geometry and Topology}
volume={9}
year={2005}
pages={1443-1499}
eprint={arXiv:math.GT/0410495}
}

\bib{CK}{misc}
{
author={Cooper, Ben}
author={Krushkal, Slava}
title={Categorification of the Jones-Wenzl projectors}
note={in preparation}
}

\bib{FSS}{misc}
{
author={Frenkel, Igor}
author={Stroppel, Catharina}
author={Sussan, Joshua}
note={in preparation}
}



\bib{Kh1}{article}
{
author={Khovanov, Mikhail}
title={A categorification of the Jones polynomial}
journal={Duke Journal of Mathematics}
volume={101}
year={2000}
pages={359-426}
eprint={arXiv:math.QA/9908171}
}

\bib{KR1}{article}
{
author={Khovanov, Mikhail}
author={Rozansky, Lev}
title={Matrix factorizations and link homology}
journal={Fundamenta Mathematicae}
volume={199}
year={2008}
pages={1-91}
eprint={arXiv:math.QA/0401268}
}

\bib{KhRS}{misc}
{
author={Khovanov, Mikhail}
author={Rozansky, Lev}
note={in preparation}
}

\bib{St}{article}
{
author={Stosic, Marko}
title={Homological thickness and stability of torus knots}
journal={Algebraic and Geometric Topology}
volume={7}
year={2007}
pages={261-284}
eprint={arXiv:math.GT/0511532}
}

\end{biblist}
\end{bibdiv}
\end{document}